%% file: Dissertation.tex
\documentclass[11pt, a4paper]{book} %% [11pt, a4paper, leqno]
\usepackage[sort]{cite}       % nicer citations
\usepackage{amssymb}          % various ams stuff
\usepackage{latexsym}         % dito
\usepackage[centertags]{amsmath}          % dito
\usepackage{amscd}            % nasty commutative ams diagrams
\usepackage[latin1]{inputenc} % smart input of special characters
\usepackage{eucal}            % nicer caligraphic fonts
\usepackage{bbm}              % nicer bbm fonts
\usepackage{theorem}          % nettere Theoreme
\usepackage{stmaryrd}         % noch ein paar mehr Symbole (Klammern)
\usepackage{mathrsfs}         % nette kaligraphische Schrift
\usepackage{paralist}         % kleinere Listen COOL!
\usepackage[all,dvips]{xy}
\usepackage{fancyhdr}
\usepackage{tabularx}
\renewcommand{\mathbb}[1]{\mathbbm{#1}} % use nicer bbm fonts

\newcounter{comment}

\theoremstyle{break}
\theoremheaderfont{\normalfont\bfseries}
\theorembodyfont{\itshape}
\newtheorem{lemma}{Lemma}[section]
\newtheorem{proposition}[lemma]{Proposition}
\newtheorem{theorem}[lemma]{Theorem}
\newtheorem{corollary}[lemma]{Corollary}

\theorembodyfont{\rmfamily}

\newtheorem{definition}[lemma]{Definition}
\newtheorem{example}[lemma]{Example}
\newtheorem{remark}[lemma]{Remark}
\theoremstyle{plain}

\newenvironment{proof}[1][{}]{\par\noindent\textsc{Proof{#1}: }}{\hspace*{\fill}$\blacksquare$\smallskip\noindent\par}

\newcommand{\cc}[1]      {\overline{{#1}}}              
\newcommand{\id}         {\operatorname{\mathsf{id}}}   
\newcommand{\image}      {\operatorname{{\mathrm{im}}}} 
 
\newcommand{\Lie}        {\operatorname{\mathscr{L}\!}}    
\newcommand{\supp}       {\operatorname{\mathrm{supp}}}

\newcommand{\tr}         {\operatorname{\mathsf{tr}}}    
\newcommand{\ad}         {\operatorname{\mathrm{ad}}}    
\newcommand{\Ad}         {\operatorname{\mathrm{Ad}}}    
\newcommand{\Conj}       {\operatorname{\mathrm{Conj}}}    
\newcommand{\Hom}        {\operatorname{\mathsf{Hom}}}   
\newcommand{\End}        {\operatorname{\mathsf{End}}}

\newcommand{\Aut}        {\operatorname{\mathsf{Aut}}} 
 
\newcommand{\cl}         {\mathrm{cl}}                    
 
\newcommand{\norm}[1]    {\left\|{#1}\right\|}            
 
\newcommand{\I}          {\mathrm{i}}
\newcommand{\E}          {\mathrm{e}}
\newcommand{\D}          {\operatorname{\mathrm{d}}}
\newcommand{\lie}[1]     {\mathfrak{#1}}

\newcommand{\ins}        {\operatorname{\mathrm{i}}}
\newcommand{\insa}       {\operatorname{\mathrm{i}_{\scriptscriptstyle\mathrm{a}}}}
\newcommand{\inss}       {\operatorname{\mathrm{i}_{\scriptscriptstyle\mathrm{s}}}}

\newcommand{\starp}     {\mathbin{\star'}}
\newcommand{\bulletp}   {\mathbin{\bullet'}}

\newcommand{\HH}            {\mathrm{HH}}
\newcommand{\Def}           {\mathrm{Def}}

\newcommand{\HC}            {\mathrm{HC}}
\newcommand{\HCtype}        {\mathrm{HC}_{\scriptstyle{\mathrm{type}}}}
\newcommand{\HCcont}        {\mathrm{HC}_{\scriptstyle{\mathrm{cont}}}}

\newcommand{\HCdiff}        {\mathrm{HC}_{\scriptstyle{\mathrm{diff}}}}

\newcommand{\HHtype}        {\mathrm{HH}_{\scriptstyle{\mathrm{type}}}}
\newcommand{\HHcont}        {\mathrm{HH}_{\scriptstyle{\mathrm{cont}}}}

\newcommand{\HHdiff}        {\mathrm{HH}_{\scriptstyle{\mathrm{diff}}}}

\newcommand{\Diffop}        {\operatorname{\mathrm{DiffOp}}}
\newcommand{\Diffopver}     {\operatorname{\mathrm{DiffOp}_{\mathrm{ver}}}}
\newcommand{\Diffopverg}[1] {\operatorname{\mathrm{DiffOp}_{\mathrm{ver}}^{\mathit{#1}}}}
\newcommand{\Union}         {\operatorname*{\mbox{$\bigcup$}}}

\newcommand{\bullett}   {\mathbin{\tilde{\bullet}}} 
\newcommand{\Loc}       {\operatorname{\mathrm{Loc}}}      % lokale Abbildungen bzw. Operatoren
\newcommand{\Mat}       {\mathrm{Mat}}                     % Matrizen z.B \Mat(N,\mathbb{R})
\newcommand{\acts}      {\mathbin{\rhd}}                   % "left action" bei Gruppenwirkungen
\newcommand{\pp}        {\mathsf{p}}                       % Fusspunktprojektion des HFB oder der surj. Sub. P                 
\newcommand{\pr}        {\mathsf{pr}}                      % Projektionsabbildung bei kartesischen Produkten
\newcommand{\pe}        {\mathsf{q}}                       % Fusspunktprojektion des VB E über P
\newcommand{\Homs}      {\operatorname{\it{Hom}}}          % Garben der Garbenhomomorphismen
\newcommand{\Homab}     {\operatorname{\it{Hom}_{\mathrm{ab}}}}     % Garben der Garbenhomomorphismen bei Abelschen Gruppen
\newcommand{\circi}     {\mathbin{\circ_i}}                % insertion after the i-th position
\newcommand{\circp}[1]  {\mathbin{\circ_{#1}}}             % insertion after variable position
\newcommand{\HCtypep}   {\mathrm{HC}_{\scriptstyle{\mathrm{type'}}}}
\newcommand{\opp}       {\mathrm{opp}}
\newcommand{\Ext}       {\mathrm{Ext}}
\newcommand{\rr}        {\mathsf{r}}                       % right action "r"
\newcommand{\rR}        {\mathsf{R}}                       % right action "R"
\newcommand{\lL}        {\mathsf{L}}                       % left action "L"
\newcommand{\LL}        {\mathsf{L}}                       % left multiplication used for Diffops
\newcommand{\Dach}[1]   {\mbox{$\bigwedge$}^{\!\!#1}}      % antisymm. Potenzen
\newcommand{\rep}       {\pi}                              % Darstellung einer Lie-Gruppe
\newcommand{\lieg}      {\mathfrak{g}}                     % Lie-Algebra "g"

\newcommand{\ddto}      {\left.\frac{\D}{\D t}\right|_{t=0}} %      ... an der Stelle t=0
\newcommand{\Fl}        {\mathrm{Fl}}                      % Fluß von Vektorfeldern
\newcommand{\hor}       {\mathrm{hor}}                     % "hor" für horizontale Formen
\newcommand{\hlift}     {\mathsf{h}}                       % "h" für horizontale Lifts
\newcommand{\Gau}       {\mathrm{Gau}}                     % Eichtransformationen des HFB
\newcommand{\gau}       {\mathfrak{gau}}                   % infinitesimale Eichtransformationen
\newcommand{\catopen}[1]{{{#1}\textrm{-}\mathbf{set}}}     % Kategorie der offenen nichtleeren Teilmengen von #1

\makeatletter
\def\cleardoublepage{\clearpage\if@twoside \ifodd\c@page\else
    \hbox{}
    \vspace*{\fill}

    \vspace{\fill}                                              %
    \thispagestyle{empty}                                       %
    \newpage                                                    %
    \if@twocolumn\hbox{}\newpage\fi\fi\fi}                      %
\makeatother  

\begin{document}

\pagestyle{empty}
\input{title}
%\cleardoublepage
%\input{joke}

\frontmatter
\pagestyle{fancy}
%\fancyhead[CE]{\slshape \nouppercase{\leftmark}}  
%\fancyhead[CO]{\slshape \nouppercase{\rightmark}} 

\tableofcontents
\mainmatter
%\addcontentsline{toc}{chapter}{Introduction}

\include{hauptteil}
%\include{intro}
%\include{gauge}   
%\include{defmodules}
%%\include{actions}
%%\include{projective}
%\include{differential}
%\include{sheaf}
%\include{barkoszul}
%\include{sursub}
%\include{symbol}
%%\include{fibration}
%%\include{commutants}
%\include{vector}
%%\include{techniques}
%%\include{HCdifferential}

\begin{appendix}
\renewcommand{\chaptername}{Appendix}
\include{anhang}

%\include{bundles}
%\include{homological}
%\include{multidiffop}
%\include{sheaves}
%\include{computations}
\end{appendix}
%\begin{appendix}

%\end{appendix}

\backmatter

\addcontentsline{toc}{chapter}{Bibliography}
%\bibliography{lit,litweiss,dqbook,dqarticle,dqpreprints,dqproceeding,dqprocentry,dqthesis,docdip,misc,preprints,script} 
%\bibliography{lit}
%\bibliographystyle{ewde}

%\addcontentsline{toc}{chapter}{Publication}

\include{paper}

%\addcontentsline{toc}{chapter}{Acknowledgement}
\include{danke}

\end{document}

%% file: title.tex
\begin{titlepage}
\thispagestyle{empty}
\vspace*{1.5cm}
\begin{center}
    \Huge Deformation Quantization\\ of \\Principal Fibre Bundles \\
    and \\ Classical Gauge Theories
\end{center}
\vspace{1.5cm}
\begin{center}
    \Large \textsc{Inaugural-Dissertation}\\
    \large zur \\ Erlangung des Doktorgrades\\ der\\
    \textsc{Fakult{\"{a}}t
      f{\"{u}}r Mathematik und Physik}\\ der\\
    \textsc{Albert-Ludwigs-Universit{\"{a}}t Freiburg}\\ [1.0cm]
    \large vorgelegt von\\[1.0cm]
    \Large \textbf{Stefan Wei\ss}\\ [1.0cm]
    \large aus\\ [1.0cm]
    St{\"{u}}hlingen\\
    \vspace{2cm} Oktober 2009
\end{center}

\newpage

\large

\vspace*{\fill}

\begin{tabular}{ll}
    Dekan: & Prof. Dr. Kay K{\"{o}}nigsmann \\
    Betreuer der Arbeit:& Apl. Prof. Dr. Stefan Waldmann\\
    Referent: & Apl. Prof. Dr. Stefan Waldmann\\
    Korreferent: & Prof. Dr. Hartmann R{\"{o}}mer
    \\
    \\
    \multicolumn{2}{l} {Tag der m{\"{u}}ndlichen Pr{\"{u}}fung: 14.12.2009}
\end{tabular}

\bigskip

\newpage

\vspace*{5cm}

\begin{center}
    \textit{To my family}
\end{center}

\end{titlepage}

% Local Variables: 
% TeX-master: "Dissertation"
% End: 

%% file: hauptteil.tex
\chapter*{Introduction}
\label{cha:introduction}

\addcontentsline{toc}{chapter}{Introduction}

\fancyhead[CE]{\slshape \nouppercase{Introduction}}
\fancyhead[CO]{\slshape \nouppercase{Introduction}}

As typical for a thesis in mathematical physics the aspiration of this
work is twofold.

From the physical point of view it is motivated by the aim to find a
global and geometric formulation of classical gauge theories on
arbitrary noncommutative space-times where the noncommutative
structure is given by a star product. In principle, such a formulation
is done by an appropriate adaption of all occurring geometric
structures to the noncommutative case. As a first step towards this
ambitious goal, this work presents and investigates a notion of a
deformation quantization of principal fibre bundles and related
geometric structures that play a role in physical
applications.

The investigation of the underlying mathematical problem, namely the
deformation of right modules with respect to a given deformed algebra
is of independent mathematical interest. From this point of view the
work provides a general algebraic approach to such deformation
problems. For certain geometric examples the existence and the
uniqueness up to equivalence of the investigated structures is shown
by an explicit computation of Hochschild cohomologies.

\section*{Motivation}
\label{sec:physical-motivation}

Quantum mechanics and general relativity are indisputably two of the
most important theories in physics and build the basis of the present
understanding of the fundamental laws of nature. Both concepts are
well-established and impressively confirmed by experiments. However,
it still is an unsolved problem to unify the quantum mechanical
description of nature at microscopic scales and the sophisticated
model of space-time geometry and its interaction with matter provided
by general relativity. The attempt to combine the two theories leads
to conceptual problems and reveals the necessity for new models of the
space-time. The promising concept of so-called noncommutative
space-times first arose in the different context of quantum
electrodynamics where it was supposed to be a way to handle the
occurring UV-divergences. In a letter to Peierls from 1930,
Heisenberg already suggested to introduce uncertainties in the space
coordinates and regretted not to be able to provide a reasonable
mathematical framework \cite{heisenberg:1930}. In 1947, Snyder was the
first who substantiated the idea and formulated a first notion of a
quantized space-time by substituting the space-time coordinates by
hermitian operators \cite{snyder:1947a,snyder:1947b}. However, with
the great achievements of other regularization methods and the concept
of renormalization the noncommutative space-times passed out of mind
for some years.

Meanwhile, the mathematical field of noncommutative geometry was
established in the works of Connes \cite{connes:1994a}, Madore
\cite{madore:1999a}, Gracia-Bond{\'{\i}}a et
al. \cite{gracia.bondia:2000a}, and many others. The developed notions
in this new area in particular clarified the mathematical definition
of noncommutative spaces and provided the appropriate mathematical
framework for all later investigations. Based on the observation that
many structures in differential geometry are determined by their
algebraic properties, noncommutative geometry extends the
correspondences between geometric data and commutative algebras to the
noncommutative case. In particular, any smooth manifold is already
determined by the commutative algebra of smooth functions on
it. Taking the algebraic structure as the fundamental one and dropping
the notion of a manifold as a set of points, `any' algebra can be
interpreted as the function algebra on a corresponding noncommutative
manifold.

That this new concept of noncommutative manifolds is physically
meaningful became clear with the works
\cite{doplicher.fredenhagen.roberts:1994,
  doplicher.fredenhagen.roberts:1995a} of Doplicher, Fredenhagen and
Roberts. Using the uncertainty relations of quantum mechanics and the
Einstein equations of general relativity they pointed out in a
gedankenexperiment that it is not possible to determine the space-time
coordinates of a point-like particle with arbitrary accuracy. Instead,
they deduced uncertainty relations for the space-time coordinates
where the Planck length $l_P=1.6\times 10^{-35} m$ has been perceived
as the scale where the new effects of noncommutativity should occur.
Consequently, the notion of a space-time consisting of single points
has no operational physical meaning. In the proposed new framework the
space-time coordinates have been replaced by noncommuting
operators. In \cite{doplicher.fredenhagen.roberts:1995a} this also led
to a first approach towards quantum field theories on noncommutative
space-times.

Starting with a noncommutative algebra describing a noncommutative
manifold it is one of the main purposes of noncommutative geometry to
find the algebraically motivated counterparts of any geometric
structure that is defined on an ordinary manifold. In particular,
vector bundles and principal bundles are of special interest since
they provide the appropriate framework to study classical gauge
theories. Vector bundles over a noncommutative manifold are simply
defined as finitely generated and projective modules over the
considered algebra. This definition based on the work of Swan
\cite{swan:1962a} is the natural generalization of the fact that
vector bundles are determined by their sections which are such modules
with respect to the functions on the base manifold.

The situation for principal fibre bundles is not that easy since the
geometric structures have no obvious algebraic properties determining
them. For trivial bundles and matrix Lie groups as structure groups
Dubois-Violette, Kerner and Madore
\cite{dubois-violette.kerner.madore:1990} presented first approaches
towards new geometric models of gauge theory. More general but
abstract notions of so-called quantum principal bundles were presented
by Brzezi{\'n}ski and Majid \cite{Brzezinski.majid:1993} and further
investigated by Hajac \cite{hajac:1996, hajac:1997}. Basically, the
structure group there is replaced by a Hopf algebra \cite{majid:1995a}
co-acting on the algebra which is related to the total space. This
definition of quantum principal bundles has later been realized to be
nothing but so-called Hopf-Galois extensions. A related definition
where the total space of the bundle is also replaced by an algebra was
given by Durdevi{\'c} in
\cite{durdevic:1992a,durdevic:1994,durdevic:1995c} who also studied
notions of corresponding gauge theories
\cite{durdevic:1995b}. Structures in the context of associated vector
bundles were also investigated from this point of view, for instance
in \cite{brzezinski:1996,coquereaux.garcia.trinchero:2000}. The
underlying geometrical and algebraic background of all these
approaches is summarized in the expository articles
\cite{baum.hajac.matthes.szymanski:2007} and \cite{masson:2008}. The
influence of noncommutative geometry gave rise to a huge amount of new
approaches in various realms of theoretical physics as it is outlined
in the articles \cite{douglas.nekrasov:2001} of Douglas and Nekrasov
and \cite{szabo:2003,szabo:2006} of Szabo.

However, due to the very abstract formulation of these concepts it is
often not clear how the occurring structures and features have to be
interpreted and how one can construct physically relevant models. For
the investigation of noncommutative space-times and their impact on
corresponding classical gauge field theories in the aspired geometric
setting, the methods of deformation quantization seem to be
appropriate. The basic concepts of deformation quantization initially
investigated by Bayen and coworkers in \cite{bayen.et.al:1978a} in
order to study the quantization of phase spaces and the relationship
between classical and quantum mechanics can be applied in the
following way. In the most general case of the pursued approach the
noncommutative space-time is assumed to be given by an associative
deformation of the pointwise product of functions on the considered
classical space-time manifold, for instance by a star product with
respect to a Poisson structure. In general, the new multiplication is
just defined for the formal power series of functions with respect to
some formal parameter. In the spirit of Gerstenhaber's deformation
theory \cite{gerstenhaber:1963a,gerstenhaber:1964a} which is the basic
framework for deformation quantization one should then try to deform
all relevant algebraic structures of a principal fibre bundle that are
related to the initially deformed multiplication of functions on the
base manifold such that the algebraic properties are
preserved. Although the convergence of the formal series is a
nontrivial problem this approach, that already gave rise to the notion
of deformation quantization of vector bundles, has some important
advantages. Since the classical objects are mostly just replaced by
formal power series their physical interpretation is clear. If the
basic deformed structures are established the found situation may
indeed lead to new perceptions that result in further generalizations.
Moreover, the framework of deformation quantization allows an
understanding of the classical limit. In the new context the formal
parameter encodes the influence of the noncommutative aspects of
space-time and can be interpreted as the Planck area if the formal
series converge. For macroscopic phenomena where this scale is
negligibly small, the parameter can be set to zero and one obtains the
ordinary gauge theory. These two aspects, the reformulation of gauge
theories in stages and the existing notion of a classical limit, are
the main reasons for this approach.

Although the notion of a star product is also clear for complicated
geometries as typically occurring in general relativity, most of the
attempts to formulate physical theories known so far use the
controversial assumption of a flat Minkowski space equipped with the
Weyl-Moyal product or other star products respecting the relevant
commutation relations of the coordinates. Within this setting the work
\cite{seiberg.witten:1999a} of Seiberg and Witten on string theory
initiated the investigations of Yang-Mills and other gauge theories on
noncommutative flat space-times which arise as low-energy limits of
string theories. Neglecting their geometric background all physically
relevant structures can be traced back to functions on the
space-time. In the approaches that make use of this simplification all
ordinary products are replaced by the considered star product. In the
works \cite{jurco.schupp.wess:2000a, jurco.schraml.schupp.wess:2000a,
  madore.schraml.schupp.wess:2000a,
  wess:2001a,jurco.moeller.schraml.schupp.wess:2001a, wess2005:} Wess
and coworkers pointed out that the Lie algebra defining the gauge
fields and the infinitesimal gauge transformations of such theories
has to be extended to its universal enveloping algebra. The initial
considerations in \cite{seiberg.witten:1999a} further motivated the
so-called Seiberg-Witten maps which express the noncommutative gauge
and matter fields together with their infinitesimal gauge
transformations by formal series of their commutative counterparts in
the deformation parameter. For quantum field theoretic models on a
flat space-time the works of Grosse and Wohlgenannt
\cite{Grosse.Wulkenhaar:2005a, Grosse.Wulkenhaar:2005b} contributed
results concerning the renormalizability of quantum field theories on
noncommutative spaces. An overview over this aspect was given by
Rivasseau \cite{Rivasseau:2007}. Based on these results there already
exist concrete formulations of the standard model on noncommutative
space-times and proposals how to measure the occurring parameters that
encode the noncommutativity in experiments \cite{alboteanu:2007,
  alboteanu.ohl.rueckl:2006, ohl.reuter:2004}.

Although the mentioned models have been investigated to an enormous
extent, there are various issues that remain unclear and thus require
further investigations and improvements. First of all, the physically
adequate notion of a noncommutative space-time is still not
well-understood. The use of star products immediately implies that the
space-time carries a Poisson structure and there do not exist any
physical arguments that distinguish a reasonable choice for it or even
guarantee its existence. There is no experimental evidence for a
globally defined noncommutative structure, in particular not for a
constant Poisson structure as used so far. In contrast to this it is
rather natural to assume that the noncommutative features of
space-time are a local effect. Waldmann and coworkers have discussed
these aspects in \cite{bahns.waldmann:2007a,
  heller.neumaier.waldmann:2007a, heller.neumaier.waldmann:2007b,
  waldmann:2008a} and presented an approach to locally noncommutative
space-times.

Second, and this is the issue basically motivating the present work,
the transparent and clear geometric formulation of classical gauge
field theories in terms of principal fibre bundles and associated
vector bundles has not been sufficiently clarified and investigated in
the context of deformation quantization. The need for a geometrical
generalization is supported by the fact that the physical effects at
the Planck scale which initially gave rise to the consideration of
noncommutative space-times can not be described in the flat Minkowski
space and require the full geometric framework of general
relativity. Moreover, many of the simplifications used so far are too
strong and not sustainable. The naive understanding of fields as
functions on the space-time is just a local description of the
geometric objects like vector-valued differential forms or sections of
vector bundles and depends on the used local charts. This is
especially important for nontrivial bundles which should be taken into
account out of the following obvious reasons. The relevant space-times
have non-vanishing curvature and due to possible singularities and
black holes they may carry nontrivial topologies admitting nontrivial
principal and vector bundles over them. Moreover, one usually demands
that the fields have a certain behaviour to decrease at infinity in
order to obtain a finite amount of energy. Using a one-point
compactification the fields that vanish at infinity can be defined as
fields on a sphere which also allows nontrivial bundles.

From this point of view it becomes clear that the procedures used in
the most models, namely to replace all products of functions by star
products has to be improved. The fields are no functions and thus it
is simply not possible to multiply them with a star product. If this
is done for the functions obtained in the local expressions of the
fields the new products heavily depend on the choice of the local
charts and have no global meaning. In the geometric formulation fields
are multiplied using fibre metrics which in the noncommutative case
have to be adapted properly. Besides establishing the deformation
quantization of vector bundles, Bursztyn and Waldmann provided a
notion of deformed hermitian fibre metrics and discussed some aspects
of noncommutative field theories in this setting \cite{bursztyn:2001a,
  bursztyn:2002a, bursztyn.waldmann:2000b, waldmann:2001b,
  waldmann:2002b, waldmann:2005a, waldmann:2009a}.

This well-understood framework for vector bundles only describes
matter fields. Gauge potentials and gauge transformations are not
taken into account. In order to do this one has to find the deformed
versions of principal fibre bundles. The appropriate notion of a
deformation quantization of principal fibre bundles was first given in
the authors diploma thesis \cite{weiss:2006a}. Using an adapted
Fedosov construction \cite{fedosov:1994a,fedosov:1996a} it was shown
that under the assumption of a symplectic base manifold the functions
on the total space can be deformed into an invariant right module with
respect to a Fedosov star product on the base.

\section*{Results of the work}
\label{sec:purpose}

As a continuation of the above developments it is the main concern of
this work to investigate the notion of deformation quantization of
principal fibre bundles with arbitrary Poisson manifolds as
basis. From the algebraic point of view the considered definition
encodes a deformation problem of right modules which also can be
discussed in a more general geometric framework. Due to this
observation a notion of deformation quantization of surjective
submersions is also taken into account.

\emph{It is the main result of this work to show that the functions on
  the total space of any surjective submersion can be deformed into a
  right module with respect to an arbitrary differential star product
  on the base manifold. These deformations defined as formal series of
  bidifferential operators with the pointwise module structure in the
  lowest order of the formal parameter, are further shown to be unique
  up to equivalence. The same result concerning the existence and the
  classification of deformation quantizations is further shown for the
  special case of principal fibre bundles where the deformations and
  their equivalence transformations are additionally required to be
  invariant under the action of the structure Lie group of the
  bundle.}

The attempt to construct deformations of algebraic structures and
correspondingly defined equivalence transformations between them order
by order in the formal parameter is typically opposed to obstructions
which are encoded in certain Hochschild cohomology groups. For star
products which are associative deformations of the algebra of smooth
functions the theorem of Hochschild, Kostant, and Rosenberg
\cite{hochschild.kostant.rosenberg:1962a} states that the relevant
Hochschild cohomology of the considered algebra is isomorphic to the
Gerstenhaber algebra of antisymmetric multivector fields. Due to this
well-known result further discussed in \cite{cahen.gutt.dewilde:1980a,
  cattaneo:2005a, nadaud:1999a, pflaum:1998a} the proof for the
existence and the classification of star products given by Kontsevich
\cite{kontsevich:2003a, kontsevich:1997a, kontsevich:1999a} has to
make use of a more sophisticated approach. Due to these facts and
other typical examples of algebraic structures where the mentioned
obstructions in general are nontrivial, the following facts are
rather surprising.

\emph{The general investigation of module deformations in this work in
  particular substantiates the obstructions against the desired
  orderwise constructions proving the above stated result concerning
  the classification and the existence of deformation quantizations of
  surjective submersions and principal fibre bundles.
  In both cases they are given by the first and the second cohomology
  group of the Hochschild complex of the algebra of functions on the
  base manifold with values in the corresponding bimodule of
  differential operators of the functions on the total space which are
  required to be invariant operators in the case of principal fibre
  bundles.
  An explicit and nontrivial computation points out that the first
  and all higher cohomology groups are trivial.}

The computation of the cohomology based on explicit local homotopies not only
provides an elegant proof of the above stated main assertion but, in
principle, also provides an explicit prescription how to construct the
relevant structures.

The mathematical framework which is necessary to describe and discuss
the arising deformation problems is provided in full generality. The
work presents a new definition for the deformation of right modules
with respect to a deformed algebra. This discussion in the sense of
Gerstenhaber's deformation theory \cite{gerstenhaber:1964a} extends the
purely algebraic notions in the works of Donald and Flanigan
\cite{donald.flanigan:1974a} and Yau \cite{yau:2005}. In the course of
the considerations it turns out that the notion of deformed module
structures consisting of bidifferential operators is a subtle
point. The aspired local formulas are obtained if the deformed module
structure is a formal series of one-cochains of the mentioned
Hochschild complex. The developed tools and techniques to perform the
crucial computation of the cohomology are also formulated in the most general
way and may possibly serve for further investigations of similar
problems. With respect to the commutant, this means the algebra of
endomorphisms of the functions on the total space respecting their
deformed module structure, the vanishing first cohomologies yield the
following result.

\emph{Depending on the geometric choices of a principal connection and
  an always existing invariant and torsion-free covariant derivative
  on the total space of the principal fibre bundle which respects the
  vertical bundle and which is related to a torsion-free covariant
  derivative on the base manifold, the commutant of a deformation
  quantization of a principal fibre bundle within the differential
  operators is isomorphic to the formal power series of vertical
  differential operators.  In consequence, this induces an invariant
  associative deformation of the algebra of vertical differential
  operators and an invariant deformation of the classical left module
  structure of the functions on the total space with respect to the
  commutant. For a fixed star product on the base manifold, all
  deformations in the resulting bimodule structure are unique up to
  invariant equivalence.  Moreover, the two occurring algebras turn out
  to be mutual commutants.}

For surjective submersions one obtains an analogous assertion where
the structures are not invariant. It is a basic feature of any
principal fibre bundle that a representation of its structure group on
a finite dimensional vector space induces an associated vector
bundle. Concerning this point the above result establishes a
connection to the works \cite{bursztyn.waldmann:2000b, waldmann:2001b,
  waldmann:2002b, waldmann:2005a, waldmann:2009a} of Bursztyn and
Waldmann. The well-known isomorphism of invariant vector-valued
functions on the principal fibre bundle and the sections of the
associated vector bundle lead to the following observation.

\emph{Every deformation quantization of a principal fibre bundle
  induces a deformation quantization of any associated vector
  bundle. Moreover, one always obtains a surjective algebra morphism
  from the deformed algebra of vertical differential operators to the
  commutant of the deformed associated vector bundle.}

Motivated by the aspired applications, the work presents a
generalization of the above statements. The notions concerning the
deformation with respect to the functions on the total space can also
be defined for the sections of arbitrary vector bundles over the total
space which, in the case of principal fibre bundles, are assumed to be
equivariant vector bundles. The above results then induce the
following generalization.

\emph{The existence and uniqueness up to equivalence can be proved in
  the same way for the right module of sections of any (equivariant)
  vector bundle over the total space.
  The assertion in particular holds for horizontal differential forms
  on a principal fibre bundle. Depending on the above geometric
  choices the invariant module deformation then induces a
  corresponding deformation of the module structure of the forms on
  the base manifold with values in any associated vector bundle.}

Although these results can be traced back to the above fundamental
notion of deformation quantization it is convenient in the main text
to investigate the sections and regard the functions as a special case
having further properties. Some of the mathematical results of this
work, basically those for the module of functions, have already been
published in the article
\cite{bordemann.neumaier.waldmann.weiss:2007a:pre} which will appear
in the \emph{Journal f{\"{u}}r die reine und angewandte Mathematik
  (Crelle's Journal)}.

\section*{Applications and outlook}
\label{sec:applications-and-outlook}

Due to the fact that surjective submersions and principal fibre
bundles are basic geometric structures that are omnipresent in
differential geometry, the above results have various applications in
mathematical physics.

\begin{itemize}
\item For the aspired geometric formulation of gauge theories on
    noncommutative space-times the above results provide a first step
    towards a better understanding of the algebraic structures which
    are crucial for the theory. All the occurring module structures
    with respect to the function algebra on the space-time allow an
    adapted deformation for any given star product. The matter fields
    of a field theory, given as the module of sections on an
    associated vector bundle can now be deformed in the gauge
    theoretic setting. The Lie algebra of infinitesimal gauge
    transformations which is required to respect the module structure
    then has to be extended to the commutant within the differential
    operators. This is the geometric and global analogue of the local
    observations made by Jur{\v{c}}o, Schraml, Schupp and Wess in
    \cite{jurco.schraml.schupp.wess:2000a}. It turns out in the
    geometric formulation that the gauge potentials, given by the
    principal connections are an affine vector space and have no
    direct module structure like the matter fields. This is only the
    case for the underlying vector space whose module structure can be
    deformed. The physical interpretation of these observations is
    still unclear. The role of gauge potentials and force fields in
    the noncommutative setting has to be investigated in more
    detail. In the symplectic case the adapted Fedosov construction of
    the deformed module structure for the matter fields reveals a
    functorial dependence of the necessary choice of a gauge
    potential, confer \cite{weiss:2006a}. This dependence shows a nice
    behaviour under local gauge transformations and gives rise to the
    conjecture that the Seiberg Witten maps for the gauge potentials
    in the global geometric framework could possibly be seen as the
    map assigning a deformed module structure to any gauge
    potential. In order to investigate this in full generality it
    would be desirable to find a more explicit construction of the
    module structures in the slightly more realistic case of
    space-times with a Poisson structure. For the formulation of a
    complete gauge theory one further has to clarify the notion of
    deformed fibre metrics as done for vector bundles in
    \cite{bursztyn.waldmann:2000b}. Then the newly defined actions
    should be invariant under the noncommutative analogues of gauge
    transformations.
\item The structure of a principal fibre bundle also occurs in the
    context of phase space reduction where the found results provide a
    deeper understanding of the quantization of the Marsden-Weinstein
    reduction \cite[Sect.~4.3]{abraham.marsden:1985a}. If a Lie group
    acts on a phase space by a Hamiltonian group action with an
    equivariant momentum map it also acts on the surface determined by
    the momentum level zero. In the case of a free and proper action
    this surface is the total space of a principal bundle over the
    quotient space given by the orbits. This base manifold inherits a
    Poisson structure of the initial phase space and is called the
    reduced phase space.  In the framework of deformation quantization
    one is then usually interested in star products of the initial
    phase space inducing star products on the reduced one. This
    problem has already been discussed using the BRST formalism
    \cite{bordemann:2000a, bordemann.herbig.waldmann:2000a,
      herbig:2006a} or more direct constructions of the reduced star
    product, confer \cite{bordemann.brischle.emmrich.waldmann:1996a,
      bordemann.brischle.emmrich.waldmann:1996b, gloessner:1998a:pre}.
    Another approach is provided by the fact that the functions on the
    momentum level surface are a bimodule with respect to the two
    algebras of functions on the two phase spaces. The left and right
    module structures are induced by the pullbacks of the embedding
    map of the surface and the bundle projection. If it is possible to
    deform the left module structure with respect to a given star
    product such that the corresponding commutant is isomorphic to the
    formal series of functions on the reduced phase space, the latter
    automatically inherits a reduced star product. In
    \cite{bordemann:2004a:pre, bordemann:2005a} this approach is
    presented together with first results in the symplectic case. For
    general Poisson manifolds there have been found sufficient
    conditions for a successful reduction \cite{cattaneo:2004a,
      cattaneo.felder:2004a, cattaneo.felder:2007a}. The results of
    this work now yield a contribution to this problem from the
    opposite direction. For every star product on a reduced phase
    space it is possible to construct a right module structure for the
    functions on momentum level surface. Together with the commutant
    one even obtains a deformed bimodule structure which is unique up
    to equivalence.
\item As discussed in the article
    \cite{bordemann.neumaier.waldmann.weiss:2007a:pre} the results
    also motivate further investigations with respect to Morita
    theory.
\end{itemize}

\section*{Outline of the work}
\label{sec:outline}

The present work is organized as follows.

\begin{itemize}
\item Chapter~\ref{cha:gauge-theory} provides a short introduction to
    the concepts of deformation quantization and it is explained in
    which sense a star product represents and describes the
    noncommutative aspects of space-time at the Planck scale. A
    discussion of the algebraic structures occurring in the geometric
    formulation of a classical gauge theory will then lead to the
    basic definitions.
\item A general and detailed investigation of deformations of right
    module structures with respect to deformed algebras is the content
    of Chapter~\ref{cha:deformation-algebras-modules}. The new
    algebraic notion is presented after a detailed discussion of the
    Hochschild complexes of algebras and modules. The cohomology
    groups of these complexes are then identified to encode the
    obstructions for the orderwise construction of such deformations
    and the corresponding equivalence transformations. These results
    can moreover be obtained for structures of particular types having
    further specific properties that should be respected by
    corresponding deformations. It is further shown that trivial
    cohomology groups not only imply the existence and uniqueness up
    to equivalence of deformed module structures but also allow a
    computation of the corresponding commutants which give rise to
    deformed bimodules. For modules which are invariant under group
    actions a separate discussion reveals further general
    results. Using an explicit homotopy it is finally shown that the
    obstructions for projective modules always vanish.
\item Chapter~\ref{cha:differential-G-invariant-structures} starts
    with a summary of the well-known facts concerning the algebraic
    notion of multidifferential operators. This is necessary to
    define the notion of differential algebra and module structures
    which then are seen to be examples of particular types.
\item In Chapter~\ref{cha:sheaves} it is pointed out how the general
    concepts of presheaves and sheaves encoding the relations between
    global and local data, can be used to compute (invariant)
    differential Hochschild cohomologies.
    % The main results there are the important
    % Propositions~\ref{proposition:sheaves-and-cohomology} and
    % \ref{proposition:sheaves-and-cohomology-G-inv}.
\item It is well-known from homological algebra that Hochschild
    cohomologies are certain values of derived functors and that this
    fact can be used to compute the cohomologies by doing it for
    certain other complexes. Motivated by this purely algebraic
    situation it is shown in Chapter~\ref{cha:bar-koszul} that such a
    result can also be obtained for a more specific situation. The
    differential Hochschild cohomologies of the algebra of smooth
    functions on a convex subset of $\mathbb{R}^n$ with values in
    certain bimodules are isomorphic to cohomologies involving the
    topological versions of the bar and Koszul resolutions. This
    result does not follow from abstract arguments. All the
    isomorphisms have to be established by explicitly given maps whose
    construction is a nontrivial issue.
\item The central result of this work is proved in
    Chapter~\ref{cha:deformation-on-bundles}. The Hochschild
    cohomologies which are crucial for the deformation quantization of
    surjective submersions and principal fibre bundles are computed
    explicitly and shown to be trivial. Based on the results of the
    previous chapters the computation makes use of a locally given
    homotopy. Besides proving the existence and the uniqueness up to
    equivalence the obtained result also presents a way to construct
    the considered structures explicitly.
\item The investigation of the commutant within the differential
    operators is the basic concern of
    Chapter~\ref{cha:symbol-calculus}. The obtained results are based
    on adapted versions of the symbol calculus of differential
    operators and the existence of particular covariant
    derivatives. This is moreover used to show that there always exist
    deformations which respect a further algebraic property.
\item In Chapter~\ref{cha:associated-vector-bundles} it is shown that
    each deformation quantization of a principal fibre bundle induces
    a deformation quantization of any associated vector bundle.
\item The two Appendices~\ref{cha:bundle-geometry} and
    \ref{cha:homological-algebra} provide the necessary basics of
    differential geometry and homological algebra which are used in
    the main text. Appendix~\ref{cha:technical-proofs} contains a
    technical but simple proof which has been omitted in
    Chapter~\ref{cha:bar-koszul}. 
\end{itemize}

\chapter{Deformation quantization and classical gauge theories}
\label{cha:gauge-theory}

\fancyhead[CE]{\slshape \nouppercase{\leftmark}}
\fancyhead[CO]{\slshape \nouppercase{\rightmark}}

As mentioned in the introduction, it is the main motivation of this
work to investigate how the methods of deformation quantization can be
used for a geometrically meaningful discussion and formulation of
classical gauge theories on noncommutative space-times. First of all,
it has to be clarified in which sense a manifold equipped with a star
product can be interpreted as a noncommutative space-time. This will
become clear after a short review of the basic notions and results of
deformation quantization. Based on the obtained interpretation of a
star product algebra, the geometric formulation of classical gauge
theories in terms of principal fibre bundles will then be
investigated. For the implementation of the new noncommutative
structure into the theory it turns out that the crucial right module
structures have to be adapted. This will finally lead to the
definition of a deformation quantization of principal fibre bundles.

\section{Deformation quantization and noncommutative space-times}
\label{sec:noncommutative-space-times}

In order to formulate the definition of a star product and all other
relevant structures occurring in the framework of deformation
quantization, one needs the concept of formal power series, confer
\cite[Chap.~IV, $\mathsection$9]{lang:1997a}.

\begin{definition}[Formal power series]
    Let $X$ be a module over a ring $\mathsf{R}$ and $\lambda$ a
    letter. Then the \emph{set $X[[\lambda]]$ of formal power series
      in the formal parameter $\lambda$ with coefficients in $X$} is
    given by the cartesian product $X[[\lambda]]= \prod_{r=0}^\infty
    X_r$ with $X_r=X$ for all $r\in \mathbb{N}_0$. In this context,
    all sequences $(x_r)_{r \in \mathbb{N}_0}= (x_0,x_1,\dots)$ in $X$
    are denoted as formal power series $x = \sum_{r=0}^\infty
    \lambda^r x_r$.
\end{definition}

In the obvious way $X[[\lambda]]$ is an $\mathsf{R}$-module again. If
$X,Y,Z$ are $\mathsf{R}$-modules any operation $\cdot: X\times
Y\longrightarrow Z$ induces a corresponding one $\cdot:
X[[\lambda]]\times Y[[\lambda]] \longrightarrow Z[[\lambda]]$ by the
formal Cauchy product
\begin{equation}
    \label{eq:formal-Cauchy-product}
    x\cdot y= \left( \sum_{r=0}^\infty \lambda^r x_r\right) \cdot
    \left(\sum_{l=0}^\infty \lambda^l y_l \right) =
    \sum_{r=0}^\infty\lambda^r \sum_{s=0}^r x_s\cdot y_{r-s}. 
\end{equation}
Following this purely algebraic prescription one can extend rings,
algebras and other algebraic structures to the formal power
series. Analogously, one defines the application of formal power
series of maps itself and one finds the isomorphism
$\Hom_{\mathsf{R}}(X,Y)[[\lambda]]\cong
\Hom_{\mathsf{R}[[\lambda]]}(X[[\lambda]],Y[[\lambda]])$ for the
$\mathsf{R}$-linear maps between two modules $X,Y$. For associative
and unital algebras $\mathcal{A}$ with unit $1$ it is an important
fact that every series $a=\sum_{r=0}^\infty\lambda^r a_r\in
\mathcal{A}[[\lambda]]$ starting with an invertible element $a_0$ is
invertible. Writing $a= a_0+\lambda b$ with $b\in
\mathcal{A}[[\lambda]]$ the element
\begin{equation}
    \label{eq:inverse-element}
    a^{-1}= a_0^{-1}\sum_{s=0}^\infty(\lambda b a_0^{-1})^s \in
    \mathcal{A}[[\lambda]] 
\end{equation}
satisfies $aa^{-1}=1= a^{-1}a$, confer
\cite[Lemma~2.2.2]{neumaier:1998a}. Further details and discussions of
formal power series can for example be found in
\cite[App.~A]{neumaier:2001a} and \cite[Bem.~4.2.36,
Sect.~6.2.1]{waldmann:2007a}.  Within this framework the definition of
a star product can be formulated as follows, confer
\cite{bayen.et.al:1978a}.

\begin{definition}[Formal star product]
    Let $M$ be a smooth manifold and $C^\infty(M)$ denote the smooth
    complex-valued functions. Then a formal \emph{star product} is an
    associative $\mathbb{C}[[\lambda]]$-bilinear product $\star$ for
    $C^\infty(M)[[\lambda]]$ of the form
    \begin{equation}
        \label{eq:star-product}
        a\star b = \sum_{r=0}^\infty\lambda^r C_r(a,b)
    \end{equation}
    for all $a,b\in C^\infty(M)[[\lambda]]$ with $\mathbb{C}$-bilinear
    maps $C_r: C^\infty(M)\times C^\infty (M)\longrightarrow
    C^\infty(M)$ such that
    \begin{enumerate}
    \item $C_0(a,b)= a\cdot b$ is the pointwise product of functions and
    \item $1\star a= a= a\star 1$ which means that
        $C_r(1,a)=0=C_r(a,1)$ for all $r\ge 1$.
    \end{enumerate}
    Two star products $\star$ and $\tilde{\star}$ on $M$ are said to be
    \emph{equivalent} if there exists a formal series
    $S=\id_{C^\infty(M)}+ \sum_{r=1}^\infty\lambda^r S_r$ of linear maps
    $S_r: C^\infty(M)\longrightarrow C^\infty(M)$ with 
    \begin{equation}
        \label{eq:star-product-equivalent}
        S(a\star b)= Sa \mathbin{\tilde{\star}} Sb \quad \textrm{and}
        \quad 
        S(1)=1. 
    \end{equation}
    A star product $\star$ is called \emph{differential} if the $C_r$
    are bidifferential operators. Then, the equivalence
    transformations $S$ have to consist of differential operators.
\end{definition}

The associativity of a star product shows that the manifold $M$ is
necessarily a Poisson manifold with a Poisson bracket
$\{\cdot,\cdot\}$ defined by
\begin{equation}
    \label{eq:first-order-commutator-Poisson}
    C_1(a,b)- C_1(b,a)= \I \{a,b\},
\end{equation}
where $\I\in \mathbb{C}$ is the imaginary unit.

The existence of star products on symplectic manifolds was
independently proved by DeWilde and Lecomte
\cite{dewilde.lecomte:1983b}, and Fedosov \cite{fedosov:1994a}. A
further contribution was given by Omori, Maeda and Yoshioka
\cite{omori.maeda.yoshioka:1991a}. With his famous formality theorem
Kontsevich showed that there exist star product for any Poisson
manifold such that \eqref{eq:first-order-commutator-Poisson} holds for
the given Poisson bracket. Moreover, this approach clarified the
classification of star products with respect to the given notion of
equivalence. Previously, this was done for the symplectic case by
Bertelson, Cahen and Gutt \cite{bertelson.cahen.gutt:1997a}, Gutt and
Rawnsley \cite{gutt.rawnsley:1999a}, as well as by Nest and Tsygan
\cite{nest.tsygan:1995a}. The topic was further pursued by Deligne
\cite{deligne:1995a} and Neumaier
\cite{neumaier:2001a,neumaier:2002a}.

The notion of star products first arose in the context of deformation
quantization as established by Bayen, Flato, Fr{\o}nsdal,
Lichnerowicz, and Sternheimer \cite{bayen.et.al:1978a}. It is a
remarkable physical phenomenon that the fundamental theory of quantum
mechanics has a classical limit. Due to this fact we only have a
physical intuition for the observables in classical mechanics, given
by the real-valued functions in the Poisson algebra
$(C^\infty(M),\{\cdot,\cdot\})$ of functions on a Poisson manifold
$M$. In contrast to this classical situation the observables in
quantum mechanics are defined as the self-adjoint elements in the
noncommutative algebra of operators on a Hilbert space. The reason for
this generally accepted axiomatic point of view can be seen in
Heisenbergs uncertainty relations $\Delta Q^i \Delta P^j \ge
\frac{\hbar}{2}\delta^{ij}$ for position and momentum which are
satisfied if the corresponding observables $Q^i$ and $P^j$ in quantum
mechanics satisfy the commutation relation $[Q^i,P^j]=\I \hbar
\delta^{ij}$. The resulting quantization problem, namely to construct
the quantum theory for any arbitrarily given classical system and to
establish a relation between the observables is not trivial. In
general, there exists no canonical way of quantization that in
particular provides the aspired relation between the Poisson bracket
$\{\cdot,\cdot\}$ for the classical observables and the commutator
$\frac{1}{\I\hbar}[\cdot,\cdot]$ for the assigned
operators. Deformation quantization in its original intention is a
geometrically motivated approach to this problem. A star product
deforming a Poisson bracket is seen as the noncommutative algebraic
structure of the observables in quantum mechanics expressed in terms
of the classical functions. Due to the fact that
\begin{equation}
    \label{eq:star-commutator}
    [a,b]_\star= a\star b-b\star a= \I \lambda \{a,b\}+ O(\lambda^2)
\end{equation}
it is clear that the deformation parameter can be interpreted as
Plancks constant $\hbar=6.6\times 10^{-34} \textit{Js}$ and that
deformation quantization yields a notion of the classical limit. A
detailed discussion of all this and an introduction to the topic can
in particular be found in the book \cite{waldmann:2007a} of Waldmann.

According to the considerations and results of Doplicher, Fredenhagen,
and Roberts in \cite{doplicher.fredenhagen.roberts:1994,
  doplicher.fredenhagen.roberts:1995a} the space time coordinates
$x^\mu$, $\mu=0,\dots, 3$, of an event in the usual Minkowski
space-time $M= \mathbb{R}^4$ are subjected to the uncertainty
relations
\begin{equation}
    \label{eq:uncertainties-space-time}
    \Delta x^0 \sum_{i=1}^3 \Delta x^i \ge l_P^2 
    \quad \textrm{and} \quad
    \sum_{j<k=1}^3 \Delta x^j \Delta x^k \ge l_P^2.
\end{equation}
The occurring Planck length $l_P=\sqrt{\frac{\gamma \hbar}{c^3}}= 1.6
\times 10^{-35} m$ is determined by the gravitational constant
$\gamma$, the Planck constant $\hbar$ and the velocity $c$ of light
and can be seen as the scale where the predictions of classical
mechanics and general relativity lead to contradictions and where a
new and more fundamental model of space-time is necessary. As shown in
the mentioned works the uncertainty relations can be satisfied in a
specific framework where the coordinate functions $x^\mu$ are replaced
by certain generators $q^\mu$ of an abstract algebra satisfying
commutation relations of the form
\begin{equation}
    \label{eq:commutation-relations-generators-space-time}
    [q^\mu,q^\nu]= \I l_P^2 \Theta^{\mu\nu}
\end{equation}
with a constant Poisson tensor $\Theta^{\mu\nu}$. Again, the problem
is to construct the fundamental theory of quantum space-time out of
its well-known and obviously existing classical limit. Completely
analogous to the above quantization problem one has to replace the
commutative algebra of functions on the space-time manifold $M$ by a
noncommutative one. With this analogy it is obvious that the general
framework of deformation quantization can also be used to describe the
physical effects of such a noncommutative space-time in terms of the
well-known one. The only additional problem is that there is no
physical evidence for a globally defined Poisson structure, confer the
discussions in \cite{bahns.waldmann:2007a,
  heller.neumaier.waldmann:2007b, waldmann:2008a}.

Nevertheless, in order to investigate the impact of a noncommutative
space-time structure on field theories it is still adequate to work
with noncommutative space-times that are given by a smooth manifold
$M$ together with a differential star product. The basic algebraic
structure for further investigations is thus the star product algebra
$(C^\infty(M)[[\lambda]], \star)$.

\section{The geometric formulation of a classical gauge theory}
\label{sec:geometry-gauge-theory}

It has become evident in the last century that classical gauge
theories have a simple and clear geometrical formulation which shall
be outlined in the following. Details, concrete applications and
further discussions of the well-known topic can be found in
\cite{bleecker:2005, daniel.viallet:1980, naber:1997, naber:2000,
  nakahara:1990a, schottenloher:1995a, weiss:2006a}. It is assumed
that the reader is familiar with the necessary concepts of
differential geometry which can also be found in the above references
or in the most books on differential geometry, for instance
\cite{kobayashi.nomizu:1963a, kobayashi.nomizu:1969a, lee:2003a,
  michor:2008}.  Appendix~\ref{cha:bundle-geometry} contains a summary
of the most important facts together with an explanation of the used
notation.

For the description of the kinematics of a classical gauge theory one
needs the following crucial data.
\begin{enumerate}
\item A principal fibre bundle $\pp:P\longrightarrow M$ with structure
    Lie group $G$ and principal right action $\rr$.
\item A representation $\rep: G\longrightarrow \Aut(V)$ of $G$ on a
    finite dimensional vector space $V$ from the left.
\end{enumerate}
These fundamental structures have the following interpretation. The
base manifold $M$ plays the role of the space-time. Thus it is
equipped with a lorentzian metric $h$, which in particular is a
non-degenerate symmetric tensor field $h\in \Gamma^\infty(M,S^2TM)$.
The total space $P$ is sometimes referred to as the \emph{space of
  generalized phase factors}. Considering the fibres one can say that
at any point $p\in M$ there is an additional degree of freedom which
is strongly related with the group $G$. The name comes from classical
electrodynamics where the elements of the relevant group $G= U(1)=
\{\E^{i\Theta}\:|\: \Theta\in (0,2\pi]\}$ are determined by the phase
$\Theta$.  The structure group $G$ is the \emph{internal symmetry
  group} of the theory and controls the gauge fields by the
representations on the vector space $V$ and the Lie algebra $\lieg$ of
$G$. The meaning of this will become clear with the following
notions.

In the global geometric context the $G$-invariant vector-valued
functions
\begin{equation}
    \label{eq:matter-fields}
    f\in C^\infty(P,V)^G
\end{equation} 
satisfying $f \circ \rr_g= \rep_{g^{-1}} \circ f$ for all $g\in G$ are
interpreted as the \emph{matter} or \emph{particle fields} of the
theory. The principal connections in the form of their connection
one-forms
\begin{equation}
    \label{eq:gauge-fields}
    \omega\in \mathcal{C}\subseteq (\Gamma^\infty(P,T^*P) \otimes
    \lieg)^G ,
\end{equation}
which are $\lieg$-valued one-forms that are $G$-invariant with respect
to the adjoint representation, $\rr_g^*\omega= \Ad_{g^{-1}}\circ
\omega$, and reproduce the generators of the fundamental vector
fields, confer Definition~\ref{definition:connection-one-form}, are
seen as the \emph{gauge potentials}. The induced covariant exterior
derivative $\D_\omega$ of such a connection then encodes the
description of the changing rate of any field. Using this derivative
one obtains the corresponding curvature forms
\begin{equation}
    \label{eq:force-fields}
    \Omega= \D_\omega \omega \in (\Gamma^\infty_\hor(P,\Dach{2}T^*P)
    \otimes \lieg)^G ,
\end{equation}
which are horizontal and invariant $\lieg$-valued two-forms. They are
seen as the \emph{force fields} or \emph{field strength tensors}. An
immediate consequence of the above definition are the \emph{structure
  equation} and the \emph{Bianchi identity},
\begin{eqnarray}
    \label{eq:structure-equation}
    \Omega&=& \D_\omega \omega= \D \omega + \frac{1}{2}
    [\omega,\omega]_\wedge,\\
    \label{eq:Bianchi-global}
    \D_\omega \Omega &=& 0,
\end{eqnarray}
where $[\cdot,\cdot]_\wedge$ is the bracket for $\lieg$-valued forms
defined in \eqref{eq:bracket-g-valued-forms}.
The covariant exterior derivative further satisfies
\begin{eqnarray}
    \label{eq:covariant-exterior-derivative-matter-fields}
    \D_\omega f&=& \D f+ \rep'(\omega)f,\\
    \label{eq:covariant-exterior-derivative-inv-hor-forms}
    \D_\omega \alpha &= & \D\alpha + [\omega,\alpha]_\wedge
\end{eqnarray}
for all $f\in C^\infty(P,V)^G$ and $\alpha\in
(\Gamma^\infty(P,\Dach{\bullet}T^*P)\otimes \lieg)^G$.
% In the considered setting the kinematic variables given by the matter,
% gauge and force fields and the covariant exterior derivatives thereof
% allow a unified description as invariant vector-valued forms
% \begin{equation}
%     \label{eq:kinematic-variables-vector-valued-forms}
%     \alpha \in (\Gamma^\infty(P,\Dach{\bullet}T^*P)\otimes V)^G
% \end{equation}
% for some appropriate vector space $V$ and a representation $\rep:
% G\longrightarrow \Aut(V)$.\\
The \emph{gauge group} or \emph{group of gauge transformations} is
given by the group of $G$-invariant $G$-valued functions
\begin{equation}
    \label{eq:gauge-group}
    C^\infty(P,G)^G
\end{equation}
with respect to the action of the Lie group $G$ on itself via
conjugation. An element $H\in C^\infty(P,G)^G$ thus satisfies $H\circ
\rr_g = \Conj_{g^{-1}} \circ H$. The group structure is the naturally
induced via $(H_1H_2)(u)= H_1(u)H_2(u)$ for all $u\in P$.
Correspondingly, the \emph{gauge algebra} or \emph{(Lie) algebra of
  infinitesimal gauge transformations} is given by the algebra of
$G$-invariant $\lieg$-valued functions
\begin{equation}
    \label{eq:gauge-algebra}
    C^\infty(P,\lieg)^G
\end{equation}
with respect to the adjoint action. Thus one has $\Xi\circ \rr_g=
\Ad_{g^{-1}} \circ \Xi$ for such $\Xi\in C^\infty(P,\lieg)^G$. Again,
the Lie bracket is the one induced by the Lie bracket on $\lieg$,
$[\Xi_1,\Xi_2](u)= [\Xi_1(u),\Xi_2(u)]$.
\begin{remark}[Gauge group and gauge algebra]
    \label{remark:gauge-group-and-algebra}
    \begin{enumerate}
    \item There is a group isomorphism
        \begin{equation}
            \label{eq:iso-gauge-group}
            \Phi: C^\infty(P,G)^G \longrightarrow \Gau(P)= \{ \tau: 
            P\longrightarrow P\:|\: \tau \circ \rr_g = \rr_g\circ \tau 
            \textrm{ and } \pp \circ \tau = \pp \} 
        \end{equation}
        between the gauge group and the gauge transformations
        $\Gau(P)$ of the principal fibre bundle which is denoted and
        defined by
        \begin{equation}
            \label{eq:iso-gauge-group-explicit}
            \Phi_H(u)= (\Phi(H))(u) = \rr_{H(u)}u= u.H(u)
        \end{equation}
        for all $H\in C^\infty(P,G)^G$ and $u\in P$.
    \item There is a vector space isomorphism
        \begin{equation}
            \label{eq:iso-gauge-algebra}
            \phi: C^\infty(P,\lieg)^G \longrightarrow \gau(P)=
            \Gamma^\infty(P,VP)^G = \{ V\in \Gamma^\infty(P,TP)\:|\:
            \rr_g^*V=V \textrm{ and } T \pp \circ V=0\}
        \end{equation}
        between the gauge algebra and the infinitesimal gauge
        transformations of the principal fibre bundle which is given
        by
        \begin{equation}
            \label{eq:iso-gauge-algebra-explicit}
            \phi(\Xi)= \Xi_P,  \quad \quad 
            \Xi_P(u)= \ddto u.\exp(t\Xi(u)) 
        \end{equation}
        for all $\Xi\in C^\infty(P,\lieg)^G$. In addition, the
        isomorphism is an anti-homomorphism of Lie algebras,
        \begin{equation}
            \label{eq:iso-gauge-algebra-bracket}
            [\Xi_1,\Xi_2]_P= - [(\Xi_1)_P,(\Xi_2)_P].
        \end{equation}
    \item The notion \emph{infinitesimal} comes from the fact that a
        vector field $V\in \Gamma^\infty(P,TP)$ is an infinitesimal
        gauge transformation in $\gau(P)$ if and only if its flow
        $\Fl^V_t\in \Gau(P)$ is a gauge transformation for all
        possible $t\in I$ with some $I\subseteq \mathbb{R}$. Under the
        above isomorphisms this relation is encoded in the generalized
        exponential map $\exp: C^\infty(P,\lieg)^G\longrightarrow
        C^\infty(P,G)^G$ from
        Appendix~\ref{subsec:connection-principal} since
        \begin{equation}
            \label{eq:flow-infinitesimal-gauge-transformation}
            \Fl^{\Xi_P}_t= \Phi_{\exp(t\Xi)}  
        \end{equation}
        and $\ddto \exp(t\Xi(u))= \Xi(u)$ for the one-parameter group
        $\exp(t\Xi)$.
    \end{enumerate}    
\end{remark}

The \emph{(local) gauge transformations} of the fields are given by
the natural action of the gauge group $\Gau(P)$ via pullback. Then
all structures yield transformed ones of the same type. For $H\in
C^\infty(P,G)^G$ the tangent map of $\Phi_H\in \Gau(P)$ is given
by
\begin{equation}
    \label{eq:tangent-map-gauge-transformation}
    T_u\Phi_H Z_u= T_u \rr_{H(u)} Z_u+ \left( T_{H(u)}l_{H^{-1}(u)}
        T_uH Z_u \right)_P(\Phi_H(u))
\end{equation}
for all $Z_u\in T_uP$ and with the left multiplication $l_gh=gh$ in
$G$. Thus one finds the following gauge transformations
\begin{equation}
    \label{eq:gauge-transformation-fields}
    H:
    \begin{array}{rcl}
        f &\longmapsto & f'= \Phi_{H^{-1}}^* f= \rep(H)f\\
        \omega &\longmapsto& \omega'= \Phi_{H^{-1}}^* \omega = 
        \Ad_H \omega + \delta H^{-1}\\
        \Omega & \longmapsto& \Omega'= \Phi_{H^{-1}}^* \Omega = 
        \Ad_H \omega
    \end{array}
%     \label{eq:gauge-transformation-matter-field}
%     \label{eq:gauge-transformation-gauge-potential}
\end{equation}
for all matter fields $f\in C^\infty(P,V)^G$, gauge potentials
$\omega\in \mathcal{C}$ and force fields $\Omega\in
(\Gamma^\infty_\hor(P,\Dach{2}T^*P))^G$. The expression $\delta
H^{-1}: TP\longrightarrow \lieg$ is the \emph{left logarithmic
  derivative} of $H^{-1}$ defined by $\delta
H^{-1}(Z_u)=T_{H^{-1}(u)}l_{H(u)}
T_uH^{-1}Z_u$. Further, it is $(\rep(H)\alpha)_u=
\rep(H(u))\circ \alpha_u$ for each vector-valued form $\alpha\in
\Gamma^\infty(P,\Dach{\bullet}T^*P)\otimes V$ with a corresponding
representation. The pullback of a connection one-form $\omega\in
\mathcal{C}$ indeed yields $\omega'\in \mathcal{C}$ with curvature
$\Omega'=\D_{\omega'}\omega'$. The \emph{gauge covariance} of the
covariant derivative $\D_\omega$ is encoded in the fact that
\begin{equation}
    \label{eq:gauge-covariance}
    (\D_\omega\alpha)'= \D_{\omega'}\alpha'.
\end{equation}
The \emph{infinitesimal (local) gauge transformations} with respect to
a function $\Xi\in C^\infty(P,\lieg)^G$ are obtained by considering
the gauge transformation with respect to the one-parameter group
$H_t= \exp(t\Xi)$ and to take the derivative $\ddto$. The
infinitesimal version of \eqref{eq:gauge-transformation-fields} then
is
\begin{equation}
    \label{eq:infinitesimal-gauge-transformation-fields}
    \Xi:
    \begin{array}{rcl}
        f&\longmapsto & \delta_\Xi f= -\Lie_{\Xi_P} f = \rep'(\Xi) f \\
        \omega&\longmapsto& \delta_\Xi\omega = \ad_{\Xi}\omega-\D \omega
        = [\Xi,\omega]-\D \omega\\
        \Omega&\longmapsto& \delta_\Xi\Omega =  -\Lie_{\Xi_P}\Omega =
        \ad(\Xi) \Omega. 
    \end{array}
\end{equation}
There, $\Lie_{\Xi_P}= \ddto (\Fl_t^{\Xi_P})^*$ denotes the Lie
derivative in direction $\Xi_P$ and $\rep':
\lieg\longrightarrow \End(V)$ is the induced Lie algebra
representation. For all vector-valued forms one again sets
$(\rep'(\Xi) \alpha)_u= \rep'(\Xi(u)) \alpha_u$. For the horizontal
forms $\alpha\in \Gamma^\infty_\hor(P,\Dach{\bullet}T^*P)\otimes V$
one additionally has
\begin{equation}
    \label{eq:inf-gauge-trafo-commutator}
    [\delta_{\Xi_1}, \delta_{\Xi_2}]\alpha 
    = \delta_{[\Xi_1,\Xi_2]} \alpha. 
\end{equation}
The infinitesimal gauge transformations $\delta_\Xi$ only yield
elements $\delta_{\Xi}\alpha$ of the same type if the corresponding
fields give rise to a vector space. This is not true for the gauge
potentials. The principal connections are only an affine space over
the invariant and horizontal $\lieg$-valued one forms. Thus one has
\begin{equation}
    \label{eq:infinitesimal-action-gauge-potential}
    \mathcal{C}\ni \omega \longrightarrow \delta_{\Xi}\omega\in
    T_\omega \mathcal{C}= (\Gamma^\infty_\hor(P,T^*P)\otimes \lieg)^G.
\end{equation}
In order to reproduce the formulas physicists are familiar with, one
makes use of a \emph{local gauge} which is nothing but a local section
\begin{equation}
    \label{eq:local-gauging}
    \sigma\in \Gamma^\infty(U,P|_U)
\end{equation}
of the principal fibre bundle over an open set $U\subseteq M$. Then
one defines
\begin{eqnarray}
    \label{eq:pull-back-fields}
    \phi &=& \sigma^*f \in C^\infty(U,V)\\
    A &=& \rep' \circ \sigma^* \omega \in \Gamma^\infty(U,
    T^*U)\otimes \End(V)\\
    F &=& \rep'\circ \sigma^*\Omega= \rep' \circ (\D \sigma^*\omega+
    \frac{1}{2}[\sigma^*\omega,\sigma^*\omega]_\wedge) \in
    \Gamma^\infty(U,\Dach{2}T^*U)\otimes \End(V)  \\ 
    \mathcal{U} &=& \rep \circ \sigma^*H \in C^\infty(U,\Aut(V))\\
    \mathfrak{U}&=& \rep'\circ \sigma^*\Xi \in C^\infty(U,\End(V))
\end{eqnarray}    
for all $f\in C^\infty(P,V)^G$, $\omega\in \mathcal{C}$, $H\in
C^\infty(P,G)^G$, and $\Xi\in C^\infty(P,\lieg)^G$. Then the structure
equation is $F=\D A+ A\wedge A$ and the Bianchi identity reads $\D F=
F \wedge A- A\wedge F$
% \begin{eqnarray}
%     \label{eq:structure-equation-local}
%     F&=& \D A+ A\wedge A\\
%     \label{eq:Bianchi-local}
%     \D F&=& F \wedge A- A\wedge F  
% \end{eqnarray}
where, for instance, $(A\wedge A)(X,Y)= A(X)\circ A(Y)-A(Y)\circ A(X)$
for all $X,Y\in \Gamma^\infty(U,TU)$. The local and infinitesimal
local gauge transformations \eqref{eq:gauge-transformation-fields} and
\eqref{eq:infinitesimal-gauge-transformation-fields} now have the form
\begin{equation}    
    \label{eq:local-and-infinitesimal-local-gauge-transformations}
    \mathcal{U}:
    \begin{array}{rcl}
        \phi &\longmapsto & \mathcal{U} \phi\\
        A & \longmapsto & \mathcal{U} A \mathcal{U}^{-1}+ \mathcal{U} \D
        \mathcal{U}^{-1}\\
        F & \longmapsto & \mathcal{U} F \mathcal{U}^{-1}
    \end{array}
    \quad \textrm{and} \quad 
    \mathfrak{U}:
    \begin{array}{rcl}
        \phi &\longmapsto & \mathfrak{U} \phi\\
        A &\longmapsto & [\mathfrak{U},A]- \D \mathfrak{U}\\
        F &\longmapsto & [\mathfrak{U},F].
    \end{array}
\end{equation}
Equation \eqref{eq:covariant-exterior-derivative-matter-fields} yields
$\sigma^* \D_\omega f = (\D +A) \sigma^* f= (\D+A)\phi$ where
$(A\phi)_p(X_p)=A_p(X_p)(\phi(p))$. Usually one chooses a basis
$\{e_a\}_{a=1,\dots, \dim G}$ of the Lie algebra $\lieg$ such that
$[e_a,e_b]= C_{ab}^ce_c$ with the structure constants $C_{ab}^c$ and
sets $T_a= \rep'(e_a)$. Here and in the following we use Einstein's
summation convention. Then one obviously has $\omega= \omega^a e_a$
and $\Omega= \Omega^a e_a$ with simple differential forms $\omega^a\in
\Gamma^\infty(P,T^*P)$ and $\Omega^a\in \Gamma^\infty(P,\Dach{2}T^*P)$
and one defines $A^a= \sigma^*\omega^a\in \Gamma^\infty(U,T^*U)$. If
$U\subseteq M$ is the domain of a chart $x: U\longrightarrow
x(U)\subseteq \mathbb{R}^n$ the forms have a $C^\infty(U)$-module
basis generated by the $\D x^\mu$ for $\mu=1,\dots n= \dim M$. This
yields
\begin{equation}
    \label{eq:gauge-field-local}
    A= (\sigma^*\omega^a)\rep'(e_a)=A^aT_a= A^a_\mu\D x^\mu T_a 
\end{equation}
with $A^a_\mu\in C^\infty(U)$. Analogously,
\begin{equation}
    \label{eq:force-field-local}
    F=\frac{1}{2} F^a_{\mu\nu} \D x^\mu \wedge \D x^\nu T_a 
\end{equation}
and \eqref{eq:structure-equation} yields the well-known formula
\begin{equation}
    \label{eq:structure-equation-local-chart}
    F^a_{\mu\nu}= \partial_\mu A^a_\nu - \partial_\nu A^a_\mu +
    C^a_{bc} A^b_\mu A^c_\nu. 
\end{equation}
The fields $\mathcal{A}_\mu= A^a_\mu T_a$ and
$\mathcal{F}_{\mu\nu}=F^a_{\mu\nu}T_a\in C^\infty(U, \End(V))$
which have the same transformation behaviour as $A$ and $F$ then are
the gauge and force fields physicists work with. Moreover, with the
definition $\sigma^*\D_\omega f= D_\mu(\sigma^*f) \D x^\mu$ and
\eqref{eq:covariant-exterior-derivative-matter-fields} one finds the
usual covariant derivative 
\begin{equation}
    \label{eq:covariant-derivative-local-matter-field}
    D_\mu= \partial_\mu + A^a_\mu T_a= \partial_\mu + \mathcal{A}_\mu 
\end{equation}
for the matter fields $\phi$. Analogously, for $\lieg$-valued forms
$\alpha\in (\Gamma^\infty(P,\Dach{k}T^*P)\otimes \lieg)^G$ with
$\rep'\circ \sigma^* \alpha= \alpha^a_{\mu_1\dots \mu_k} \D
x^{\mu_1}\wedge \dots \wedge \D x^{\mu_k} T_a$ one sets $\rep' \circ
\sigma^*(\D_\omega \alpha)= D_{b\mu}^a \alpha^b_{\mu_1 \dots \mu_k} \D
x^\mu\wedge \D x^{\mu_1}\wedge \dots \wedge \D x^{\mu_k} T_a$ and
obtains
\begin{equation}
    \label{eq:covariant-derivative-local-force-field}
    D^a_{b\mu}= \partial_\mu \delta^a_b - C^a_{bc} A^c_\mu. 
\end{equation}
With \eqref{eq:covariant-exterior-derivative-inv-hor-forms} the
Bianchi identity \eqref{eq:Bianchi-global} implies the
\emph{homogeneous field equations}
\begin{equation}
    \label{eq:Bianchi-local}
    \sum_{\mathrm{zykl \{ \rho\mu\nu \}}} D^a_{b\rho} F^b_{\mu\nu}=0.
    % \sum_{\mathrm{zykl \{ \rho\mu\nu \}}} \partial_\rho
    % F^b_{\mu\nu}- C^a_{bc} A^c_\rho F^b_{\mu\nu} 
\end{equation}
for all $a=1,\dots,\dim G$.

\begin{remark}[Introduction of charges]
    In physical applications one has to take care of the
    dimensions. By convention one thus introduces the coupling
    constants which have the interpretation of charges. In
    electrodynamics one sets $D_\mu= \partial_\mu + ie\mathcal{A}_\mu$
    where $e$ is the electric charge.
\end{remark}

\begin{remark}[Lagrangians, field equations and interactions]
\label{remark:lagrangians-and-all-that}
    \begin{enumerate}
    \item In order to find the dynamics of such a gauge theory one
        needs some more geometrical data. All this can be found in
        \cite{bleecker:2005} and shall only be outlined briefly. A
        $G$-invariant \emph{Lagrangian} $L$ yields a function
        $\mathcal{L}_0: C^\infty(P,V)^G\longrightarrow C^\infty(M)$ of
        the form $\mathcal{L}_0(f)=L(f,\D f)$ with $L(\rep_g f,\rep_g
        \D f)=L(f,\D f)$ which in physics literature is referred to as
        \emph{invariance under global gauge transformations}. The
        demand for invariance under local gauge transformations is
        satisfied by 
        \begin{equation}
            \label{eq:Lagrangian-with-cov-der}
            \mathcal{L}(f,\omega)=L(f,\D_\omega f)
        \end{equation}
        involving the gauge potential $\omega$. The metric $h$ of the
        space-time $M$ and a $G$-invariant metric $k$ on the Lie
        algebra $\lieg$, $k(\Ad_g\xi,\Ad_g \eta)=k(\xi,\eta)$, induce
        a metric $hk$ on the horizontal $\lieg$-valued forms on
        $P$. This gives rise to the \emph{self-action density}
        \begin{equation}
            \label{eq:self-action-density}
            \mathcal{S}(\omega)=-\frac{1}{2} hk(\D_\omega
            \omega,\D_\omega\omega)= -\frac{1}{2} \norm{\Omega}^2.
        \end{equation}
        Taking a volume density $\mu$ with respect to $h$ one can
        define the \emph{action}
        $\int_{U}(\mathcal{L}+\mathcal{S})(f,\omega)\mu$ over any open
        subset $U\subseteq M$ with compact closure. The
        \emph{principle of least action} which claims that $f$ and
        $\omega$ are \emph{stationary} in the sense that
        \begin{equation}
            \label{eq:stationarity}
            \ddto \int_U(\mathcal{L}+\mathcal{S})(f+tf',\omega
            +t\alpha) \mu=0  
        \end{equation}
        for all $f'\in C^\infty(P,V)^G$ and $\alpha\in
        (\Gamma^\infty_\hor(P,T^*P)\otimes \lieg)^G$, then is
        equivalent to the \emph{Lagrange equation} for the matter
        field $f$ and the \emph{inhomogeneous field equation} for the
        gauge field $\omega$.
    \item In order to find physically relevant examples of matter
        fields and corresponding Lagrangians one needs the theory of
        \emph{spin structures, spinors and Dirac operators}, confer
        \cite{Cartan:1966,Chevalley:1997,friedrich:1997} for the
        mathematical background and \cite{naber:2000} for some
        physical considerations. On the Minkowski space, the
        Lagrangian typically is of the form
        \begin{equation}
            \label{eq:general-Lagrangian}
            \mathcal{L}_0 (\phi)=
            \langle \partial_\mu\phi, \partial^\mu\phi 
            \rangle -m^2 \langle \phi,\phi \rangle - \mathcal{V}
            (\langle \phi,\phi \rangle)
        \end{equation}
        where $\langle\cdot,\cdot\rangle$ is some metric on the
        considered vector space $V$. The first summand is seen as
        kinetic energy, the constant $m$ is interpreted as mass and
        $\mathcal{V}$ is some interaction potential. The structure
        group of the considered principal bundle then is an
        appropriate subgroup $G$ of the orthogonal group of
        $\langle\cdot, \cdot\rangle$. For matrix Lie groups with the
        defining representation the self-interaction then is given by
        the trace,
        \begin{equation}
            \label{eq:self-interaction-matrix-group}
            \mathcal{S}(\mathcal{A})= -\frac{1}{2}
            \tr(\mathcal{F}_{\mu\nu} \mathcal{F}^{\mu\nu}). 
        \end{equation}
        In the physical interpretation the demand of local gauge
        invariance and the introduction of the covariant derivative in
        \eqref{eq:Lagrangian-with-cov-der}, usually referred to as
        \emph{minimal coupling}, determine the interaction of the
        gauge potentials with the matter fields.
    \end{enumerate}
\end{remark}

\begin{remark}[Associated vector bundles]
    \label{remark:associated-vector-bundles}
    It is well-known that the representation $\rep: G\longrightarrow
    \Aut(V)$ of the structure group leads to an associated vector
    bundle $\pe: E=P\times_G V \longrightarrow M$, confer
    Appendix~\ref{sec:associated-vector-bundles}. Then all structures
    and operations from above have a corresponding description in
    terms of the associated vector bundle. Some of the well-known
    facts are summarized in
    Proposition~\ref{proposition:associated-vector-bundles}. The
    matter fields can be seen as sections of this associated bundle
    due to the isomorphism $s: C^\infty(P,V)^G \longrightarrow
    \Gamma^\infty(M,E)$ and the gauge potentials and force fields give
    rise to covariant derivatives $\nabla^E$ on $E$ with curvatures
    $\mathcal{R}^E$. Further, an element $H\in C^\infty(P,G)^G$ of the
    gauge group induces a vector bundle automorphism $ \hat{H}\in
    \Aut(E)$, this means a fibrewise linear diffeomorphism
    $\hat{H}:E\longrightarrow E$ with $\pe\circ \hat{H}=\pe$, by
    $\hat{H}[u,v]=[u,\rep (H(u))v]=[\Phi_H(u),v]$ for all equivalence
    classes $[u,v]\in P\times_G V$. Analogously, an element $\Xi\in
    C^\infty(P,\lieg)^G$ of the gauge algebra induces a section
    $\hat{\Xi}\in \Gamma^\infty(M,\End(E))$ in the endomorphism bundle
    of $E$ by $\hat{\Xi}[u,v]=[u,\rep'(\Xi(u))v]$. With the
    isomorphism $s$ one finds that $\hat{H}\circ s(f) = s(\rep(H)f)$.
\end{remark}

\section{Module structures in classical gauge theories}
\label{subsec:algebraic-properties-gauge-theory}

In order to investigate the impact of a noncommutative space-time
structure on the above formulation in the aspired framework of
deformation quantization, it is necessary to identify the crucial
algebraic structures of a gauge theory which are in contact with the
algebra of functions on the space-time. Considering a star product
algebra $(C^\infty(M)[[\lambda]],\star)$ as the new noncommutative
version of the space-time, the idea then is to adapt the related
structures by an appropriate deformation in order to maintain the
algebraic structures which are assumed to be the fundamental content
of the theory. From the purely algebraic point of view it is clear
that the simplest structures related to an algebra are corresponding
modules. Indeed, the above formulation points out that certain right
modules have an important physical interpretation.

The matter fields, given by invariant vector valued functions
$C^\infty(P,V)^G$ on the principal fibre bundle are obviously a module
with respect to the functions $C^\infty(M)$ on the base with the
pointwise multiplication making use of the pullback with the bundle
projection. For $f\in C^\infty(P,V)^G$ and $a\in C^\infty(M)$ the new
matter field $fa$ is defined by
\begin{equation}
    \label{eq:right-module-matter-fields}
    (fa)(u)= (\pp^*a \cdot f)(u)= a(\pp(u)) f(u)
\end{equation}
for all $u\in P$ where the multiplication on the right side is the
scalar one of the vector space $V$. The required $G$-invariance is
clear by the property $\pp\circ \rr_g=\pp$ and the defining
property of the representation $\rep$. Due to the commutativity of the
pointwise product of functions this is as well a left as a right
module structure but as indicated by the notation it will be
considered as a right module structure.

As explained in Remark~\ref{remark:associated-vector-bundles} and
Proposition~\ref{proposition:associated-vector-bundles} the matter
fields can also be seen as the sections $\Gamma^\infty(M,P\times_G V)$
of an associated vector bundle and the considered isomorphism is a
module isomorphism. The natural right module structure of the sections
is the crucial algebraic property and any such finitely generated and
projective module can be seen as the module of sections on a vector
bundle. This well-known statement, confer the more detailed discussion
in Section~\ref{sec:quantization-vector-bundles}, points out the
importance of the above right module structure which is obviously
induced by the corresponding $G$-invariant right module structure of
the functions $C^\infty(P)$ itself. For $f\in C^\infty(P)$ and $a\in
C^\infty(M)$ one defines $fa= f\cdot \pp^*a \in C^\infty(P)$ and finds
the property

\begin{equation}
    \label{eq:invariant-right-module-functions}
    \rr_g^* (fa)=(\rr_g^* f)a.
\end{equation}

Without regard to the $G$-invariance this basic module structure of
functions also occurs in a much more general framework and only makes
use of the projection $\pp: P\longrightarrow M$ which is a surjective
submersion.

The local gauge transformations $\Phi\in \Gau(P)$ and the
infinitesimal counterparts $\Xi_P\in \Gamma^\infty(VP)^G$ which are
applied to matter fields are obviously seen to respect the considered
module structure and thus are elements of the so-called commutant. For
the functions $f\in C^\infty(P)$ one easily finds that
\begin{equation}
    \label{eq:gauge-trafo-commutant}
    \Phi^*(fa)= (\Phi^*f)a
\end{equation}
and
\begin{equation}
    \label{eq:inf-gauge-trafo-commutant}
    \delta_{\Xi_P}(fa) = -\Lie_{\Xi_P}(fa)= (-\Lie_{\Xi_P}f)a= (\delta_{\Xi_P}f)a
\end{equation}
for all $a\in C^\infty(M)$.

The set $\mathcal{C}$ of gauge potentials $\omega$ has no such right
module structure. It is an affine vector space over the horizontal and
$G$-invariant $\lieg$-valued one-forms $\alpha\in
(\Gamma^\infty_\hor(P, T^*P)\otimes \lieg)^G$. This vector space and
all higher differential forms of the same type are again modules with
respect to $C^\infty(M)$.

\section{Deformation quantization and classical gauge theories}
\label{subsec:definition-def-quant-sursub-principal}

The aim to adapt the crucial $C^\infty(M)$-module structure of
$C^\infty(P)$ to a given star product $\star$ on the manifold $M$
gives rise to the following definitions which were already given in
the publication \cite{bordemann.neumaier.waldmann.weiss:2007a:pre}.

\begin{definition}[Deformation quantization of surjective submersions]
    \label{definition:deformation-quantization-sursub}
    Let $\pp: P\longrightarrow M$ be a surjective submersion and
    $\star$ be a differential star product on $M$.
    \begin{enumerate}
    \item A \emph{deformation quantization of the surjective
          submersion} is a right
        $(C^\infty(M)[[\lambda]],\star)$-module structure $\bullet$ of
        $C^\infty(P)[[\lambda]]$ meaning that
        \begin{equation}
            \label{eq:deformed-right-module}
            f \bullet (a\star b) = (f\bullet a) \bullet b
        \end{equation}
        for all $f\in C^\infty(P)[[\lambda]]$ and $a,b \in
        C^\infty(M)[[\lambda]]$ such that
        \begin{equation}
            \label{eq:deformation-surjective-submersion}
            f \bullet a 
            = f \cdot \mathsf{p}^*a 
            + \sum_{r=1}^\infty \lambda^r \rho_r(f,a)
        \end{equation}
        with $\mathbb{C}$-bilinear maps $\rho_r:C^\infty(P)\times
        C^\infty(M)\longrightarrow C^\infty(M)$ which are
        bidifferential operators.
    \item Two such deformations $\bullet$ and $\tilde{\bullet}$ are
        said to be \emph{equivalent} if and only if there exists a
        formal series $T = \id_{C^\infty(P)} + \sum_{r=1}^\infty
        \lambda^r T_r$ of differential operators $T_r \in
        \Diffop(C^\infty(P))$ such that for all $f \in
        C^\infty(P)[[\lambda]]$ and $a \in C^\infty(M)[[\lambda]]$
        \begin{equation}
            \label{eq:equivalence-surj-sub}
            T(f \bullet a)=
            T(f) \mathbin{\tilde{\bullet}} a.
        \end{equation}
    \end{enumerate}
\end{definition}

The notion of bidifferential operators used for the definition means
that the local expression of these operators shall always exists and
that with respect to local charts they are bidifferential operators in
the common sense of calculus. The global description of such maps is a
subtle point and will be clarified in the course of the further
investigations of these structures. Taking the $G$-invariance into
account the following definition is natural.

\begin{definition}[Deformation quantization of principal fibre bundles]
    \label{definition:deformation-quantization-principal}
    Let $\mathsf{p}:P\longrightarrow M$ be a principal fibre bundle
    with structure group $G$ and principal right action $\mathsf{r}: P
    \times G \longrightarrow P$, $ \mathsf{r}(u,g)=\mathsf{r}_gu$, and
    $\star$ a differential star product on $M$.
    \begin{enumerate}
    \item A \emph{deformation quantization of the principal fibre
          bundle} is a $G$-invariant deformation quantization of the
        surjective submersion $\mathsf{p}: P \longrightarrow M$ with
        respect to $\star$, this means a right module structure
        $\bullet$ as in \eqref{eq:deformed-right-module} and
        \eqref{eq:deformation-surjective-submersion} with the
        additional property
        \begin{equation}
            \label{eq:deformation-right-module-principal}
            \mathsf{r}_g^* (f \bullet a)
            = (\mathsf{r}_g^*f) \bullet a
        \end{equation}
        for all $f \in C^\infty(P)[[\lambda]]$, $a \in
        C^\infty(M)[[\lambda]]$, and $g \in G$.
    \item Two such deformations $\bullet$ and $\tilde{\bullet}$ are
        said to be equivalent if they are equivalent in the sense of
        Definition~\ref{definition:deformation-quantization-sursub}
        with $G$-invariant operators $T_r$ which means that in
        addition to \eqref{eq:deformation-right-module-principal} for
        all $g \in G$ one has
        \begin{equation}
            \label{equivalence-principal-G-invariance}
            \mathsf{r}_g^* \circ T_r
            = T_r \circ \mathsf{r}_g^*.
        \end{equation}
    \end{enumerate}
\end{definition}

In the diploma thesis \cite{weiss:2006a} the existence of deformation
quantizations of principal fibre bundles was already shown for the
special case of a symplectic base manifold $M$. There the star product
is assumed to be a Fedosov star product \cite{fedosov:1994a} depending
on a symplectic and torsion-free covariant derivative $\nabla^M$ on
$M$ and a formal series $\Omega^M = \sum_{r=1}^\infty \lambda^r
\Omega^M_r \in \lambda \Gamma^\infty(M,\Dach{2} T^*M)[[\lambda]]$ of
closed two-forms, $\D \Omega^M_r=0$. With an adapted version of the
Fedosov construction for bimodules it is shown that there always exist
corresponding deformation quantizations of any principal fibre bundle
over $M$. The construction and thus the right module structure itself
further depend on the choice of a principal connection one-form
$\omega$ and an always existing $G$-invariant, torsion-free covariant
derivative $\nabla^P$ on $P$ which respects the vertical bundle and is
related to $\nabla^M$. Although it was shown that deformation
quantizations for different choices of $\nabla^P$ are equivalent with
explicitly given equivalence transformations
\cite[Satz~3.1.16]{weiss:2006a} it was not possible to clarify the
classification in this framework. This and the fact that there is no
evidence for a symplectic structure on a space-time gives rise to
investigate the above notion of deformation quantization for arbitrary
star products on Poisson manifolds. In the present work we prove the
following two fundamental theorems.

\begin{theorem}
    \label{theorem:existence-uniqueness-sursub}
    Every surjective submersion $\pp:P\longrightarrow M$ with a star
    product $\star$ on $M$ admits a deformation quantization which is
    unique up to equivalence.
\end{theorem} 

\begin{theorem}
    \label{theorem:existence-uniqueness-principal}
    Every principal fibre bundle $\pp:P\longrightarrow M$ with a star
    product $\star$ on $M$ admits a deformation quantization which is
    unique up to equivalence.
\end{theorem}

Although the Fedosov construction in \cite{weiss:2006a} only works in
the symplectic setting, it has further properties that might be useful
for a physical interpretation. In particular, it yields a functorial
dependence of the module structure $\bullet$ on the connection
$\omega$ and $\nabla^P$. Under a gauge transformation $\Phi\in
\Gau(P)$ the crucial dependence on $\omega$, denoted by
$\bullet_\omega$, has the behaviour
\begin{equation}
    \label{eq:functorial-dependence-fedosov}
    \Phi^*(f \mathbin{\bullet_{\omega}} a) = \Phi^*f
    \mathbin{\bullet_{\Phi^* \omega}} a. 
\end{equation}
Although $\bullet_\omega$ and $\bullet_{\Phi^*\omega}$ are equivalent
the pullback $\Phi^*$ is no equivalence transformation. 

Equation \eqref{eq:functorial-dependence-fedosov} in particular shows,
that the local gauge transformations are no longer in the
commutant. If this algebraic property is seen to be fundamental and
shall be maintained, the actions of the gauge transformations have to
be replaced by the natural action of the endomorphisms in the
commutants. For the infinitesimal versions one obviously makes the
same observation. Since they are classically given by Lie derivatives,
this means by differential operators, it is an interesting task to
compute the commutant of the deformation quantization within all
differential operators of the algebra $C^\infty(P)$ which will be done
in this work.

The situation for gauge potentials and gauge fields has not been
studied so far and is a difficult problem since it is not clear which
structures should be adapted. Moreover, the gauge potentials already
appear as a necessary choice for the Fedosov construction which gives
evidence that the matter and gauge fields already are strongly related
in this setting. Nevertheless, it is a first step of understanding
also to investigate the deformation problem of the vector space of
horizontal and $G$-invariant forms together with their commutants. All
this will now be done in a very general way which yields a further
contribution to the final aim of adapting the whole geometry of
principal fibre bundles to a given star product.

\chapter{Deformation theory of algebras and modules}
\label{cha:deformation-algebras-modules}

% The questions about the existence of star products on arbitrary
% Poisson manifolds and the characterization of the corresponding
% equivalence classes are already answered, confer the famous
% publications of Kontsevich \cite{kontsevich:1997:pre,kontsevich:1997a}
% and other important works dealing with this topic. Although the
% considerations there are of a different kind,
The basic ideas of deformation quantization as introduced in
\cite{bayen.et.al:1977a,bayen.et.al:1978a} have their origin in the
early works \cite{gerstenhaber:1963a, gerstenhaber:1964a,
  gerstenhaber:1966a, gerstenhaber:1968a, gerstenhaber:1974a} of
Gerstenhaber on algebraic deformation theory. Concretely, they provide
the algebraic setting for the deformation theory of algebras. The main
point there is that the obstruction for an order by order construction
of an associative deformation of an algebra structure is encoded in a
certain Hochschild cohomology of the algebra. In the same way the
obstruction for a construction of an equivalence transformation
between two associative deformations is given by another Hochschild
cohomology. This fundamental ideas can be reformulated for the
deformation theory of other kinds of algebraic structures. Donald and
Flanigan \cite{donald.flanigan:1974a} have first considered the
deformation theory of modules. Later this topic has been discussed in
\cite{yau:2005} by Yau who has also studied
various other structures.

In most cases, however, it is not possible to compute the crucial
cohomologies explicitly or they are not vanishing and thus this
approach does not help for the decision whether it is possible to
execute the constructions of the relevant structures or not. However,
the deformation theory of Gerstenhaber is one of the fundamental
guidelines of this work.
% The reason for this is twofold. First of all the framework for
% algebras can be extended in a natural way to find a corresponding
% and adapted obstruction theory for modules over an algebra which
% should be deformed.
With the further assumption that the Hochschild cochains are
differential operators it will be possible to compute the cohomologies
for the right modules occurring in the framework of surjective
submersions and principal fibre bundles.

This chapter presents the general framework for all this. Besides
giving a motivation by the well-known concepts and definitions there
are presented the modifications and adaptions of the ideas and
approaches for the investigation of deformed modules. The first
section is dedicated to the well-known notion of Hochschild
complexes. It is introduced together with a short summary of the
naturally occurring operations like the insertions and the
cup product. After this very general statements we concentrate on
complexes which are induced by algebras and right modules. For the
first case the algebraic structure of the Hochschild cohomology is
well-known and given by a so-called Gerstenhaber algebra. Leaving the
purely algebraic description it is explained in
Section~\ref{sec:algebras-modules-type} what shall be understood by
algebra and module structures of particular types. This new definition
specifies in an axiomatic way which further properties of the
algebraic structures can be implemented in the purely algebraic
setting. In comparison to 
%the first publication
\cite{bordemann.neumaier.waldmann.weiss:2007a:pre} the crucial
conditions are formulated in a slightly more general way in the
present work. The Sections~\ref{sec:deformation-algebras} and
\ref{sec:deformation-modules} then show that the obstruction theory
for deformations of algebras and modules can be reformulated within
this refined context. Moreover,
Section~\ref{sec:commutant-module-structures} shows that the
considered Hochschild cohomologies even play a crucial role for the
investigation of the commutant of a deformed right module
structure. With respect to the applications
Section~\ref{sec:G-invariant-types} contains a discussion of algebraic
structures which are invariant under group actions. Finally this
chapter is closed with a first simple example given by the projective
modules. There it is possible to compute the purely algebraic
cohomology by an explicit homotopy. In the following presentation it
is assumed that the reader is familiar with the basic concepts of
homological algebra. The crucial definitions can be found in
Appendix~\ref{cha:homological-algebra} or in the respective
literature, for instance \cite{jacobson:1985a,jacobson:1989a,
  lang:1997a}.

From now on let $\mathbb{K}$ be a field of characteristic $0$ like
$\mathbb{K}=\mathbb{Q},\mathbb{R},\mathbb{C}$. This will be sufficient
for our purposes. However, all considerations of this chapter are
still valid for a commutative ring $\mathbb{K}$ with at least $1\neq
0$ and $0 \neq 2,\frac{1}{2} \in \mathbb{K}$.

\section{The Hochschild complex}
\label{sec:Hochschild}

In order to build up an obstruction theory for the deformations of
algebras and modules one needs the notion of Hochschild complexes and
Hochschild cohomologies which shall be introduced in the
following. All presented structures of this section are well-known and
already appear in the early works
\cite{hochschild:1945,hochschild:1946} of Hochschild and
\cite{gerstenhaber:1963a} of Gerstenhaber. In our presentation, the
very general considerations and definitions in the literature are
already adapted to the algebraic framework of our later examples. This
is convenient in order to introduce the more specific framework which
will be used later on, and to emphasize the crucial additional
aspects.

\subsection{The Hochschild complex with values in a bimodule}
\label{subsec:Hochschild-bimodule}

Let $\mathcal{A}$ and $\mathcal{M}$ be $\mathbb{K}$-vector spaces.
For all integers $k\in \mathbb{Z}$ consider the derived
$\mathbb{K}$-vector spaces
\begin{equation}
    \label{eq:Hochschild-complex-bimodule}
    \HC^k(\mathcal{A},\mathcal{M})=
    \left\{
        \begin{array}{ccc}
            \{0 \}
            && k<0\\
            \mathcal{M} 
            & \textrm{for} & k=0\\
            \Hom_{\mathbb{K}} (\underbrace{\mathcal{A}\times
              \dots \times \mathcal{A}}_{k \textrm{ times}}, \mathcal{M})
            &&  k \geq 1,
        \end{array}
    \right.
\end{equation}
this means the \emph{$\mathbb{K}$-multilinear maps} of $\mathcal{A}$ with
values in $\mathcal{M}$ for $k\ge 1$. The direct sum $\HC^\bullet
(\mathcal{A},\mathcal{M})=\bigoplus_{k\in \mathbb{Z}} \HC^k
(\mathcal{A},\mathcal{M})$ yields a $\mathbb{Z}$-graded space and the
corresponding degree of homogeneous elements is often referred to as
\emph{dimension} or \emph{tensor degree}, since
$\HC^k(\mathcal{A},\mathcal{M})$ for $k\ge 1$ can also be seen as the
linear maps $\Hom_{\mathbb{K}} (\mathcal{A}^{\otimes k}, \mathcal{M})$
of the $k$-fold tensor product $\mathcal{A}^{\otimes k}$ of
$\mathcal{A}$ over $\mathbb{K}$. Besides this, it will be convenient
for further definitions and computations to shift the grading by one
which means to consider
\begin{equation}
    \label{eq:grading-Hochschild-algebra}
    \HC[1]^\bullet(\mathcal{A},\mathcal{M}) = \bigoplus_{k\in
      \mathbb{Z}} \HC^{k+1}(\mathcal{A}, \mathcal{M})
\end{equation}
where the \emph{degree} of homogeneous elements is defined by
\begin{equation}
    \label{eq:degree-Hochschild-algebra}
    \deg \phi = k \quad \Leftrightarrow \quad \phi\in
    \HC^{k+1}(\mathcal{A},\mathcal{M}).
\end{equation}
So, an element of $\HC^{k}(\mathcal{A},\mathcal{M})$ is said to have
dimension or tensor degree $k$ but degree $k-1$. Within this framework
one can now define the following basic operations on
$\HC^\bullet(\mathcal{A},\mathcal{M})$.

\begin{definition}[The insertion maps and the cup product]
    Let $\mathcal{A},\mathcal{M}$ be $\mathbb{K}$-vector spaces.
    \begin{enumerate}
    \item For homogeneous elements
        $\phi\in\HC^{k+1}(\mathcal{A},\mathcal{M})$ and $\psi\in
        \HC^{m+1}(\mathcal{A},\mathcal{A})$ with $k,m\in \mathbb{N}_0$
        and for $i=0,\dots, k=\deg \phi$ the \emph{insertion $\phi
          \circi \psi\in \HC^{k+m+1} (\mathcal{A},\mathcal{M})$ of
          $\psi$ in $\phi$ after the $i$-th position} is defined by
        \begin{equation}
            \label{eq:insertion-after-ith-position}
            (\phi\circi \psi) (a_1,\dots, a_{k+m+1}) = \phi(a_1,\dots , a_i,
            \psi(a_{i+1}, \dots, a_{i+m+1}), a_{i+m+2}, \dots, a_{k+m+1})
        \end{equation}
        for all $a_1,\dots, a_{k+m+1}\in \mathcal{A}$.
    \item Linear combination of all possible insertions yields a new
        multiplication
        \begin{equation}
            \label{eq:total-insertion-algebra}
            \phi\circ \psi = \sum_{i=0}^{\deg \phi} (-1)^{i \deg \psi}
            \phi\circi \psi
        \end{equation}
    \item Let $\mathcal{M}$ be equipped with the additional structure
        of an associative $\mathbb{K}$-algebra. Then for all $k,m \in
        \mathbb{N}_0$ the \emph{cup product}
        \begin{equation}
            \label{eq:cup-product-map}
            \cup: \HC^{k}(\mathcal{A},\mathcal{M}) \times
            \HC^{m}(\mathcal{A},\mathcal{M}) \longrightarrow
            \HC^{k+m}(\mathcal{A},\mathcal{M}) 
        \end{equation}
        is defined by
        \begin{equation}
            \label{eq:cup-product}
            (\phi\cup \psi) (a_1,\dots, a_{k+m})= \phi(a_1, \dots ,
            a_k) \psi(a_{k+1},\dots, a_{k+m})
        \end{equation}     
        for all $\phi \in \HC^{k}(\mathcal{A},\mathcal{M})$, $\psi\in
        \HC^{m}(\mathcal{A},\mathcal{M})$ and $a_1,\dots, a_{k+m}\in
        \mathcal{A}$.
    \end{enumerate}
\end{definition}

All maps allow bilinear extensions defined on all elements in
$\HC^\bullet(\mathcal{A},\mathcal{M})$ and
$\HC^\bullet(\mathcal{A},\mathcal{A})$. For the insertion one therefore
simply defines $\phi\circi \psi=0$ if $i>\deg \phi$. If one of the
involved tensor degrees is less than zero, all maps are defined to
be zero, too.  By \eqref{eq:cup-product-map}, the cup product is
homogeneous with respect to the tensor degree. Instead, the maps
$\circi, \circ: \HC^{k+1}(\mathcal{A},\mathcal{M}) \times
\HC^{m+1}(\mathcal{A},\mathcal{M}) \longrightarrow
\HC^{k+m+1}(\mathcal{A},\mathcal{M}) $ are homogeneous with respect to
the shifted degree,
\begin{equation}
    \label{eq:insertion-algebra-homogeneous}
    \deg (\phi\circ \psi)= \deg \phi + \deg \psi.
\end{equation}

\begin{remark}[Associativity of the cup product]
    The cup product $\cup$ is obviously associative. Thus,
    $(\HC^\bullet (\mathcal{A},\mathcal{M}),\cup)$ is a graded,
    associative algebra with respect to the tensor degree.
\end{remark}

In the case $\mathcal{M}=\mathcal{A}$ the insertions can be
composed. But as stated in the following proposition this does not
lead to an associative multiplication.

\begin{proposition}
    \label{proposition:properties-insertions}
    Let $\mathcal{A}$, $\mathcal{M}$ be $\mathbb{K}$-vector spaces and
    let $\phi\in \HC^\bullet(\mathcal{A},\mathcal{M})$, $\psi,\chi\in
    \HC^\bullet(\mathcal{A},\mathcal{A})$ be homogeneous elements.
    \begin{enumerate}
    \item With respect to their composition the corresponding
        insertions satisfy
        \begin{equation}
            \label{eq:associativity-circi}
            (\phi \circi \psi)\circp{j} \chi= 
            \left\{ 
                \begin{array}{lcl}
                    (\phi\circp{j} \chi) \circp{i+ \deg \chi} \psi &&
                    j < i\\
                    \phi \circi (\psi \circp{j-i} \chi) & \textrm{for}
                    & i\leq 
                    j \leq i+\deg \psi\\
                    (\phi \circp{j-\deg \psi} \chi) \circi \psi &&
                    j > i+\deg \psi.
                \end{array}
            \right.
        \end{equation}
    \item In consequence, the map $\circ$ satisfies
        \begin{eqnarray}
            \label{eq:associativity-insertion}
            (\phi\circ \psi)\circ \chi - \phi\circ (\psi\circ \chi)
            &=& \sum_{i=1}^{\deg \phi}\sum_{j=0}^{i-1}
            (-1)^{i \deg \psi +j \deg \chi} (\phi \circi
            \psi)\circp{j} \chi \nonumber \\
            &&+
            \sum_{i=0}^{\deg \phi-1}\sum_{j=i+\deg \psi+1}^{\deg
              \phi+\deg \psi}
            (-1)^{i \deg \psi +j \deg \chi} (\phi \circi
            \psi)\circp{j} \chi
            \nonumber \\
            &=& (-1)^{\deg \psi \deg \chi} ((\phi\circ \chi) \circ
            \psi - \phi\circ (\chi\circ \psi)). 
        \end{eqnarray}
    \end{enumerate}
\end{proposition}

A complete and detailed version of the combinatorial proof can be
found in \cite{waldmann:2007a}. There, the case
$\mathcal{M}=\mathcal{A}$ is treated, but the considerations are
exactly the same in the slightly more general case of the present
situation.

\begin{proposition}
    \label{proposition:Gerstenhaber-bracket}
    In the case $\mathcal{M}=\mathcal{A}$ the supercommutator
    \begin{equation}
        \label{eq:super-commutator-insertion}
        [\phi,\psi]= \phi\circ \psi -(-1)^{\deg \phi \deg\psi}
        \psi \circ \phi 
    \end{equation}
    is superantisymmetric and satisfies the super Jacobi
    identity
    \begin{equation}
        \label{eq:super-jacobi-gerstenhaber}
        [\phi,[\psi,\chi]]= [[\phi,\psi],\chi]+ (-1)^{\deg \phi
          \deg \psi} [\psi, [\phi,\chi]]
    \end{equation}
    with respect to the degree $\deg$. Thus,
    $(\HC^\bullet(\mathcal{A},\mathcal{A}),\deg, [\cdot, \cdot])$
    becomes a super Lie algebra.
\end{proposition}

Again, a detailed proof can be found in \cite{waldmann:2007a}.

\begin{remark}
    As explained in \cite{gerstenhaber:1963a},
    Equation~\eqref{eq:associativity-circi} shows that $\{
    \HC[1]^\bullet(\mathcal{A},\mathcal{A}), \circi\}$ defines a
    so-called \emph{right pre-Lie system} and that $\HC[1]^\bullet
    (\mathcal{A},\mathcal{M})$ is a \emph{right module over} $\{
    \HC[1]^\bullet(\mathcal{A},\mathcal{A}), \circi\}$. Then, for
    $\mathcal{M}=\mathcal{A}$,
    Equation~\eqref{eq:associativity-insertion} is a general
    consequence of the fact that
    $(\HC[1]^\bullet(\mathcal{A},\mathcal{A}), \circ)$ becomes a
    \emph{graded right pre-Lie ring}. Further, this right pre-Lie ring
    gives rise to the \emph{graded Lie ring} or super-Lie algebra
    $(\HC[1]^\bullet(\mathcal{A},\mathcal{A}), [\cdot, \cdot])$
    satisfying \eqref{eq:super-jacobi-gerstenhaber}. For more details
    and the explicit definitions of the mentioned notions, confer
    \cite{gerstenhaber:1963a}.
\end{remark}

\begin{definition}[Gerstenhaber bracket]
    \label{definition:Gerstenhaber-bracket}
    Let $\mathcal{A}$ be a $\mathbb{K}$-vector space. The
    nonassociative product $\circ$ of $\HC^\bullet
    (\mathcal{A},\mathcal{A})$ is called \emph{Gerstenhaber product}
    and the super Lie bracket $[\cdot,\cdot]$ is called the
    \emph{Gerstenhaber bracket}.
\end{definition}

The Gerstenhaber bracket yields an important characterization of
associative multiplications.

\begin{lemma}[Associative multiplications]
    \label{lemma:condition-associative-multiplication}
    Let $\mathcal{A}$ be a $\mathbb{K}$-vector space. Then, a bilinear
    map $\mu: \mathcal{A}\times \mathcal{A}\longrightarrow
    \mathcal{A}$, seen as an element $\mu\in \HC^2(\mathcal{A},
    \mathcal{A})$, is an associative multiplication, if and only if
    \begin{equation}
        \label{eq:condition-associative-multiplication}
        [\mu,\mu]=0.
    \end{equation}
\end{lemma}

\begin{proof}
    By definition, $\deg \mu=1$ and thus $[\mu,\mu]= \mu\circ \mu+
    \mu\circ \mu= 2\mu\circ \mu$. Then, 
    \begin{equation*}
        (\mu\circ \mu) (a,b,c)= \sum_{i=0}^1 (-1)^i (\mu\circi \mu)
        (a,b,c) = \mu(\mu(a,b),c)- \mu(a,\mu(b,c))
    \end{equation*}
    and the Lemma is clear since $\frac{1}{2}\in \mathbb{K}$.
\end{proof}

After this general considerations we now consider the case where
$\mathcal{A}$ is an associative $\mathbb{K}$-algebra and $\mathcal{M}$
is an $(\mathcal{A},\mathcal{A})$-bimodule. Here and in the following
$\mathcal{A}$ will always be a unital algebra with unit $1$. For
convenience's sake the multiplication in the algebra and the module
multiplication are simply denoted by $xy$ with corresponding elements
$x,y$ if the meaning is clear by the context. Then, consider the
$\mathbb{K}$-linear maps
\begin{equation}
    \label{eq:Hochschild-differentials}
    \delta^k: \HC^k (\mathcal{A},\mathcal{M})
    \longrightarrow \HC^{k+1} (\mathcal{A},\mathcal{M})
\end{equation}
which for $k\in \mathbb{N}$ are defined by
\begin{eqnarray}
    \label{eq:Hochschild-bimodules}
    (\delta^k \phi)(a_1,\dots, a_{k+1})
    &=& a_1 \phi(a_2, \dots ,
    a_{k+1}) \nonumber \\
    && + \sum_{i=1}^k (-1)^i \phi(a_1, \dots, a_i a_{i+1},
    \dots, a_{k+1})\\
    && + (-1)^{k+1}\phi(a_1,\dots, a_k) a_{k+1}. \nonumber
\end{eqnarray}
For $k=0$ one has $(\delta^0\phi)(a)=a\phi-\phi a$ for $\phi\in
\mathcal{M}$ and $a\in \mathcal{A}$. The second term of the sum in
\eqref{eq:Hochschild-bimodules} is given by $\sum_{i=1}^k (-1)^i
\phi(a_1, \dots, a_i a_{i+1}, \dots, a_{k+1})= -(\phi\circ \mu)
(a_1,\dots, a_{k+1})$. It is a straightforward but lengthy computation
using the bimodule properties of $\mathcal{M}$ which shows that
\begin{equation}
        \label{eq:Hochschild-differential-square}
        \delta^{k+1} \circ \delta^k=0.
\end{equation}
By setting the maps $\delta^k$ to be trivial for $k<0$ all this yields
a $\mathbb{K}$-linear map $\delta:
\HC^\bullet(\mathcal{A},\mathcal{M}) \longrightarrow
\HC^{\bullet+1}(\mathcal{A},\mathcal{M})$ which increases
the %(tensor)
degree by one and has square zero, $\delta^2=0$. Thus, these
structures define a cochain complex, confer
Appendix~\ref{cha:homological-algebra}, which leads to the following
definition.

\begin{definition}[The Hochschild complex of an algebra with values in
    a bimodule]
    \label{definition:Hochschild-complex}
    Let $\mathcal{A}$ be an associative $\mathbb{K}$-algebra and
    $\mathcal{M}$ be a $\mathbb{K}$-vector space with an
    $(\mathcal{A},\mathcal{A})$-bimodule structure. Then, the cochain
    complex $(\HC^\bullet(\mathcal{A},\mathcal{M}),\delta)$ is called
    the \emph{Hochschild complex of $\mathcal{A}$ with values} or
    \emph{coefficients in $\mathcal{M}$} and the \emph{Hochschild
      differential $\delta$}. For all $k\in \mathbb{N}_0$ the
    corresponding $k$-th cohomology group
    \begin{equation}
        \label{eq:Hochschild-cohomology}
        \HH^k(\mathcal{A},\mathcal{M})= \frac{\ker \delta^k} {\image
          \delta^{k-1}}
    \end{equation}
    is called the $k$-th \emph{Hochschild cohomology group} and one sets
    \begin{equation}
        \label{eq:Hochschild-cohomology-graded-space}
        \HH^\bullet(\mathcal{A},\mathcal{M})= \bigoplus_{k\in
          \mathbb{N}_0} \HH^k (\mathcal{A},\mathcal{M}).
    \end{equation}
\end{definition}

\subsection{The Hochschild complex of an associative algebra }
\label{subsec:Hochschild-complex-algebra}

It is clear from the Definition \ref{definition:Hochschild-complex}
that an associative algebra, seen as a bimodule over itself, yields a
Hochschild complex where the associativity is crucial for
\eqref{eq:Hochschild-differential-square}.

\begin{definition}[The Hochschild complex of an associative algebra]
    \label{definition:Hochschild-algebra}
    Let $(\mathcal{A},\mu)$ be an associative algebra. Then
    $(\HC^\bullet(\mathcal{A},\mathcal{A}),\delta)$ is called
    the \emph{Hochschild complex of the algebra $\mathcal{A}$}.
\end{definition}

Clearly, the algebra multiplication is then seen as a cochain
\begin{equation}
    \label{eq:algebra-multiplication-cochain}
    \mu \in \HC^2(\mathcal{A},\mathcal{A}).
\end{equation}
%One should note that if $\mathcal{A}$ is commutative there exist two
%ways to define the bimodule structure which yield slightly different
%differentials $\delta$ and thus different complexes. %%%Nein, ist dasselbe

In this section the well-known additional features of Hochschild
complexes of associative algebras are presented. As already seen in
Proposition~\ref{proposition:Gerstenhaber-bracket} and
Lemma~\ref{lemma:condition-associative-multiplication} the case
$\mathcal{M}=\mathcal{A}$ yields further structures and different ways
to describe the algebra and its Hochschild complex. Moreover, the
Hochschild cohomology $\HH^\bullet(\mathcal{A},\mathcal{A})$ is one of
the most important structures encoding crucial properties of the
associative algebra $\mathcal{A}$. For completeness' sake all this
will be recalled in the following such that finally one can state the
well-known algebraic structure of
$\HH^\bullet(\mathcal{A},\mathcal{A})$.

\begin{lemma}[The Hochschild complex of an algebra]
    \label{lemma:Hochschild-algebra}
    Let $(\mathcal{A},\mu)$ be an associative algebra.
    \begin{enumerate}
    \item The associativity of the multiplication $\mu$ yields that
        $\mu\in \HC^2(\mathcal{A},\mathcal{A})$ is a cocycle, $\delta
        \mu=0$. It is even a coboundary $\mu =\delta \iota$ with
        respect to the identity cochain $\iota\in
        \HC^1(\mathcal{A},\mathcal{A})$, $\iota(a)=a$.
    \item The Hochschild differential can be expressed by
        \begin{equation}
            \label{eq:Hochschild-differential-algebra}
            \delta \phi= (-1)^{\deg \phi}[\mu,\phi]= -[\phi,\mu]
        \end{equation}
        with the Gerstenhaber bracket $[\cdot,\cdot]$.
    \item The cup product can be expressed by
        \begin{equation}
            \label{eq:cup-product-as-insertions}
            \phi\cup \psi= (\mu\circp{0}\phi)\circp{k}\psi
        \end{equation}
        for homogeneous $\phi\in \HC^k(\mathcal{A},\mathcal{A})$ and
        arbitrary $\psi\in \HC^\bullet(\mathcal{A},\mathcal{A})$.
    \end{enumerate}
\end{lemma}

\begin{proof}
    All assertions are clear by the definitions and are shown in
    \cite{gerstenhaber:1963a}.
\end{proof}

The expression \eqref{eq:Hochschild-differential-algebra} for the
Hochschild differential $\delta$, the super Jacobi identity
\eqref{eq:super-jacobi-gerstenhaber} and Lemma
\ref{lemma:condition-associative-multiplication} yield a different
proof for $\delta^2=0$ since
\begin{equation*}
    \delta^2\phi = [[\phi,\mu],\mu]= [\phi, [\mu,\mu]]-
    [[\phi,\mu],\mu]= -\delta^2\phi.
\end{equation*}
Further, Equation~\eqref{eq:Hochschild-differential-algebra} makes it
easy to investigate the behaviour of the Hochschild differential with
respect to the Gerstenhaber product $\circ$, the Gerstenhaber bracket
$[\cdot,\cdot]$ and the cup product $\cup$.

\begin{proposition}[Compatibilities of the Hochschild differential]
    \label{proposition:Hochschild-differential-operations}
    Let $(\mathcal{A},\mu)$ be an associative
    $\mathbb{K}$-algebra. Then, the Hochschild differential $\delta$
    of the corresponding Hochschild complex has the following
    properties.
    \begin{enumerate}
    \item For $\phi,\psi \in \HC^\bullet(\mathcal{A},\mathcal{A})$
        with homogeneous $\psi$ it is
        \begin{equation}
            \label{eq:Hochschild-differential-bracket}
            \delta[\phi,\psi]= [\phi, \delta\psi] + (-1)^{\deg
              \psi}[\delta\phi, \psi].
        \end{equation}
    \item For $\phi\in \HC^k(\mathcal{A},\mathcal{A})$ and $\psi\in
        \HC^\bullet(\mathcal{A},\mathcal{A})$ it is
        \begin{equation}
            \label{eq:Hochschild-differential-cup-product}
            \delta(\phi\cup \psi)= \delta \phi \cup \psi +
            (-1)^k\phi\cup \delta\psi.
        \end{equation}
    \item For $\phi\in \HC^k(\mathcal{A},\mathcal{A})$ and $\psi\in
        \HC^m(\mathcal{A},\mathcal{A})$ it is
        \begin{equation}
            \label{eq:Hochschild-differential-circ}
            \phi\circ \delta\psi -\delta(\phi\circ \psi) +
            (-1)^{m-1}\delta\phi\circ \psi = (-1)^{m-1}(\psi \cup
            \phi- (-1)^{km} \phi\cup \psi).
        \end{equation}
    \end{enumerate}
\end{proposition}

\begin{proof}
    The first part is a direct consequence of Lemma
    \ref{lemma:Hochschild-algebra} and
    \eqref{eq:super-jacobi-gerstenhaber}. Explicitly, one computes
    $\delta[\phi,\psi]= -[[\phi,\psi],\mu]= -[\phi,[\psi,\mu]]-
    (-1)^{\deg \psi} [[\phi,\mu],\psi]= [\phi,\delta\psi]+ (-1)^{\deg
      \psi}[\delta\phi, \psi]$. The proofs of the other statements are
    basic computations using the definitions, Proposition
    \ref{proposition:properties-insertions} and Lemma
    \ref{lemma:Hochschild-algebra}. All this can be found in
    \cite{waldmann:2007a}.
\end{proof}

\begin{remark}
    \begin{enumerate}
    \item The second statement
        \eqref{eq:Hochschild-differential-cup-product} still holds in
        the general framework of Hoch\-schild complexes of algebras
        $\mathcal{A}$ with coefficients in bimodules $\mathcal{M}$
        which are associative algebras themselves if the additional
        multiplication of $\mathcal{M}$ is compatible with the module
        structure, this means if $a(m_1m_2)=(am_1)m_2$,
        $(m_1m_2)a=m_1(m_2 a)$ and $(m_1 a)m_2= m_1(a m_2)$ for all
        $m_1,m_2\in \mathcal{M}$ and $a\in \mathcal{A}$.
    \item Equation~\eqref{eq:Hochschild-differential-bracket} means
        that $\delta$ is a super derivation of
        $(\HC[1]^\bullet(\mathcal{A},\mathcal{A}), [\cdot,\cdot])$
        from the right of degree one. In the same way, the Leibniz
        rule \eqref{eq:Hochschild-differential-cup-product} shows that
        $\delta$ is a super derivation of
        $(\HC^\bullet(\mathcal{A},\mathcal{M}), \cup)$ from the left.
    \end{enumerate}
\end{remark}

The Equations~\eqref{eq:Hochschild-differential-bracket} and
\eqref{eq:Hochschild-differential-cup-product} assure on one hand that
the products of cocycles are again cocycles and on the other hand that
the coboundaries build an ideal within the cocycles. Since
Equation~\eqref{eq:Hochschild-differential-circ} further shows that
for cocycles $\phi,\psi$ the right hand side is a coboundary the
following corollary is obvious.

\begin{corollary}
    With Proposition
    \ref{proposition:Hochschild-differential-operations} it is clear
    that the Gerstenhaber bracket $[\cdot,\cdot]$ and the cup product
    $\cup$ can be well-defined on the Hochschild cohomology groups
    $\HH^\bullet(\mathcal{A},\mathcal{A})$ where the bracket yields a
    super Lie algebra
    $(\HH[1]^\bullet(\mathcal{A},\mathcal{A}),[\cdot,\cdot])$ and the
    cup product an associative and super commutative algebra
    $(\HH^\bullet(\mathcal{A},\mathcal{A}), \cup)$.
\end{corollary}

The assertions of \cite[Thm.~5, Cor.~2]{gerstenhaber:1963a} which are
the content of the next proposition contributes the last part for the
characterization of the algebraic structure of the Hochschild
cohomology.

\begin{proposition}
    Let $\mathcal{A}$ be an associative algebra and $\xi\in
    \HH^k(\mathcal{A},\mathcal{A})$, $\eta \in
    \HH^m(\mathcal{A},\mathcal{A})$, $\zeta \in
    \HH^\bullet(\mathcal{A},\mathcal{A})$. Then, the Leibniz rule
    \begin{equation}
        \label{eq:Leibniz-rule-cohomology}
        [\eta\cup \zeta, \xi] = [\eta,\xi]\cup \zeta + (-1)^{(k-1)m}
        \eta\cup [\zeta,\xi]
    \end{equation}
    holds. This means that $[\cdot,\xi]$ is a graded derivation of
    $\cup$ of degree $\deg \xi$.
\end{proposition}

\begin{corollary}[The Gerstenhaber algebra
    $\HH^\bullet(\mathcal{A},\mathcal{A})$]
    Let $\mathcal{A}$ be an associative algebra. Then, the
    corresponding Hochschild cohomology has the structure of a
    so-called \emph{right Gerstenhaber algebra}
    $(\HH^\bullet(\mathcal{A},\mathcal{A}), \cup, [\cdot, \cdot])$,
    this means of a $\mathbb{Z}$-graded, associative and super
    commutative algebra with multiplication $\cup$ and a super Lie
    algebra structure $[\cdot,\cdot]$ with respect to the shifted
    grading $\HH[1]^\bullet(\mathcal{A},\mathcal{A})$ satisfying
    \eqref{eq:Leibniz-rule-cohomology}.  Moreover, if $\mathcal{A}$ is
    unital with unit $\mathbb{1}\in \mathcal{A}$ the class
    $[\mathbb{1}]\in \HH^0(\mathcal{A},\mathcal{A})$ is a unit with
    respect to the cup product.
\end{corollary}

As stated in the beginning of the chapter the Hochschild cohomology of
an algebra describes crucial properties. The zeroth Hochschild
cohomology $\HH^0(\mathcal{A},\mathcal{A})= \ker \delta^0$ for
instance is nothing but the center of the algebra since for $a,b\in
\mathcal{A}=\HC^0(\mathcal{A},\mathcal{A})$ one has $(\delta a)(b)=
ba-ab=-\ad(a)b$. This further shows that the coboundaries $\delta a\in
\HC^1(\mathcal{A},\mathcal{A})=\End_{\mathbb{K}}(\mathcal{A},\mathcal{A})$
are inner derivations. For an arbitrary $D\in
\HC^1(\mathcal{A},\mathcal{A})$ one finds $(\delta D) (a,b)= aD(b)-
D(ab)+D(a)b$ which shows that the cocycles in $\ker \delta^1$ are the
derivations of $\mathcal{A}$. Thus the first cohomology is the space
of outer derivations. As it will be explained in Section
\ref{sec:deformation-algebras} the second and third Hochschild
cohomology are crucial for the deformation theory of the algebra.

\subsection{The  Hochschild complex of a right module}
\label{subsec:Hochschild-complex-module}

Similar to associative algebras, every module over an associative
algebra gives rise to a particular Hochschild complex. Without loss of
generality the situation is investigated for right modules.

So let $(\mathcal{A},\mu)$ be an associative $\mathbb{K}$-algebra and
$\mathcal{E}$ be a vector space over $\mathbb{K}$ with a right
$\mathcal{A}$-module structure $\rho: \mathcal{E}\times \mathcal{A}
\longrightarrow \mathcal{E}$. This can be seen as a
$\mathbb{K}$-linear map
\begin{equation}
    \label{eq:right-module-as-ring-morphism}
    \rho: \mathcal{A}
    \longrightarrow \End_{\mathbb{K}}(\mathcal{E},\mathcal{E})
\end{equation}
with values in the $\mathbb{K}$-linear endomorphisms of $\mathcal{E}$
and thus as an element $\rho\in \Hom_{\mathbb{K}}
(\mathcal{A}, \End_{\mathbb{K}}(\mathcal{E}, \mathcal{E}))$. Then, of
course, one has the identification $\rho(e,a)=\rho(a)(e)$ for all
$e\in \mathcal{E}$ and $a\in \mathcal{A}$ and the condition for $\rho$
to be a right module structure with respect to $\mu$ reads
\begin{equation}
    \label{eq:right-module-new-description}
    \rho(b)\circ \rho(a) = \rho(\mu(a,b))
\end{equation}
for all $a,b \in \mathcal{A}$. 
%This means that a right module
%structure $\rho$ can be seen as an algebra homomorphisms from
%$(\mathcal{A},\mu)$ into the opposite $(\End_{\mathbb{K}}(\mathcal{E},
%\mathcal{E}),\circ)^{\mathrm{opp}}$ of the unital $\mathbb{K}$-algebra
%$\End_{\mathbb{K}}(\mathcal{E}, \mathcal{E})$ with the composition
%$\circ$ of maps as multiplication. 
The $\mathbb{K}$-vector space $\End_{\mathbb{K}}(\mathcal{E},
\mathcal{E})$ can be equipped with the canonical
$(\mathcal{A},\mathcal{A})$-bimodule structure
\begin{equation}
    \label{eq:bimodule-of-endomorphisms}
    a \cdot D \cdot b = \rho(b) \circ D \circ \rho(a)
\end{equation}
for $a,b \in \mathcal{A}$ and $D\in \End_{\mathbb{K}}(\mathcal{E},
\mathcal{E})$ where $\circ$ is the usual composition of maps. This
gives rise to a particular Hochschild complex.

\begin{definition}[The Hochschild complex of a right module]
    \label{definition:Hochschild-module}
    Let $\mathcal{E}$ be a right module over an associative algebra
    $\mathcal{A}$. Then
    $(\HC^\bullet(\mathcal{A},\End_{\mathbb{K}}(\mathcal{E},
    \mathcal{E})),\delta)$ is called the \emph{Hochschild complex of
      the right module $\mathcal{E}$}.
\end{definition}

\begin{remark}[The Hochschild complex of a right module]
    \label{remark:Hochschild-module}
    Since the endomorphisms $\End_{\mathbb{K}}(\mathcal{E},
    \mathcal{E})$ build an algebra the Hochschild complex
    $\HC^\bullet(\mathcal{A},\End_{\mathbb{K}}(\mathcal{E},
    \mathcal{E}))$ is equipped with a cup product given by
    \begin{equation}
        \label{eq:cup-product-module}
        \phi\cup \psi(a_1,\dots, a_{k+m})= \phi(a_1, \dots, a_k) \circ
        \psi(a_{k+1},\dots ,a_{k+m}).
    \end{equation}
    Since the right module structure $\rho$ can be seen as a cochain
    \begin{equation}
        \label{eq:right-module-structure-cochain}
        \rho \in
        \HC^1(\mathcal{A}, \End_{\mathbb{K}} (\mathcal{E}, \mathcal{E}))
    \end{equation}
    the Hochschild differential has the form
    \begin{eqnarray}
        \label{eq:Hochschild-differential-module}
        (\delta \phi)(a_1,\dots, a_{k+1}) &=&
        (\phi\cup \rho) (a_2,\dots,a_{k+1},a_1) \nonumber \\
        && +\sum_{i=1}^k (-1)^i (\phi \circp{i-1} \mu) (a_1,\dots, a_{k+1})\\
        && +(-1)^{k+1} (\rho \cup \phi)(a_{k+1},a_1, \dots, a_k). \nonumber
    \end{eqnarray}
\end{remark}

\begin{remark}[Right modules over commutative algebras]
    \label{remark:special-case-commutative-algebra}
    If $(\mathcal{A},\mu)$ is a commutative associative algebra
    equation \eqref{eq:right-module-new-description} can be written in
    the form
    \begin{equation}
        \label{eq:right-module-cup-insertion}
        \rho\circ \mu= \rho\cup\rho.       
    \end{equation}
    Moreover, there exists an additional choice of the bimodule
    structure of $\End_{\mathbb{K}}(\mathcal{E}, \mathcal{E})$ since
    \eqref{eq:bimodule-of-endomorphisms} then can be changed into
    \begin{equation}
        \label{eq:new-bimodule-of-endomorphism}
        a\cdot D\cdot b= \rho(a)\circ D\circ \rho(b).
    \end{equation}
    This yields a slightly different Hochschild differential with the
    simple form
    \begin{equation}
        \label{eq:new-Hochschild-differential}
        \delta \phi= \rho\cup \phi- \phi\circ \mu +(-1)^{k+1}\phi\cup \rho.
    \end{equation}
    %In fact the difference is only a sign.
    In this chapter the general situation will be discussed. In the
    later applications, however, it is convenient to use
    \eqref{eq:new-bimodule-of-endomorphism} which is possible without
    any restrictions.
\end{remark}

Again the Hochschild cohomology encodes crucial information of the
right module structure. The zeroth Hochschild cohomology
$\HH^0(\mathcal{A}, \End_{\mathbb{K}} (\mathcal{E}, \mathcal{E}))=
\ker \delta^0$ for instance is simply the set of module endomorphisms
\begin{equation}
    \label{eq:commutant-general}
    \End_{\mathcal{A}}(\mathcal{E},\mathcal{E})= \{D
    \in \End_{\mathbb{K}} (\mathcal{E}, \mathcal{E}) \; \big| \; D 
    \circ \rho(a) =  \rho(a) \circ D \; \textrm{for all} \; a \in
    \mathcal{A} \},
\end{equation}
which is called the \emph{commutant} of the given right module
structure. The interpretation of the first and second Hochschild
cohomology will arise in the framework of the deformation theory
discussed in Section \ref{sec:deformation-modules}.

\section{Algebras and modules of a particular  type}
\label{sec:algebras-modules-type}

With the purely algebraic framework built up so far one already could
formulate the well-known definitions and results concerning
deformations of associative algebras and their generalizations to
right modules. There, of course, the corresponding Hochschild
cohomologies play the crucial role. However, since the considered
algebra multiplication $\mu$ and the right module structure $\rho$
typically have additional properties which are also required for the
higher orders of corresponding deformations one has to introduce a
slightly more specific framework containing these aspects. In order to
provide the implementation of the additional properties into the
considerations it will be crucial that these properties can be
expressed in an adequate way and satisfy certain compatibility
conditions. With respect to the further applications and since all
considerations are based on the framework of Hochschild complexes it
is necessary to regard the given structures as corresponding cochains
$\mu\in \HC^2(\mathcal{A},\mathcal{A})$ and $\rho\in
\HC^1(\mathcal{A}, \End_{\mathbb{K}}
(\mathcal{E},\mathcal{E}))$. Then, the particular properties are
assumed to be expressed by particular Hochschild cochains. This reads
as
\begin{equation}
    \label{eq:HC-type-multiplication}
    \mu \in \HCtype^2(\mathcal{A},\mathcal{A})\subseteq
    \HC^2(\mathcal{A},\mathcal{A}) 
\end{equation}
and
\begin{equation}
    \label{eq:HC-type-module}
    \rho \in
    \HCtype^1(\mathcal{A},\mathcal{D}) \subseteq
    \HC^1(\mathcal{A}, \End_{\mathbb{K}}(\mathcal{E},\mathcal{E}))  
\end{equation}
where, in principle, `$\mathrm{type}$' describes the properties
concerning the arguments of $\mathcal{A}$ and $\mathcal{D}
\subseteq \End_{\mathbb{K}} (\mathcal{E},\mathcal{E})$ is the subset
of maps with the correct behaviour with respect to arguments of
$\mathcal{E}$. This separation is of course not strict. In general,
the characterization `$\mathrm{type}$' and $\mathcal{D}$ describe the
behaviour of the cochains in a whole.

\begin{remark}
    \label{remark:particular-types}
    \begin{enumerate}
    \item It should be noted that in general the properties for the
        multiplication $\mu$ and the right module structure $\rho$
        could be different. Since we have to ensure that these two are
        compatible in the sense explained later and with respect to
        the examples in the following chapters we do not distinguish
        the different kinds of `$\mathrm{type}$' in the notation.
    \item The characterization of the right module structure $\rho$,
        for example the demands to be continuous with respect to some
        topology or to be a differential operator, really depends on
        whether it is seen as a map $\rho: \mathcal{E}\times
        \mathcal{A}\longrightarrow \mathcal{E}$ or as a map $\rho:
        \mathcal{A} \longrightarrow \End_{\mathbb{K}} (\mathcal{E},
        \mathcal{E})$ as in \eqref{eq:HC-type-module}, although these
        two perspectives are equivalent in the purely algebraic
        setting.  In
        Remark~\ref{remark:notion-differential-right-module} this
        difference is pointed out for the notion of differential
        operators. For $\mu_r$ this complication does not appear with
        \eqref{eq:HC-type-multiplication}. Of course there could occur
        further subtleties if one would use the universal property of
        the tensor product and define the cochains as linear maps of
        the tensor product instead of multilinear maps. But
        nevertheless, if the types are fixed in the explained way one
        can afterwards always identify the multilinear maps with the
        corresponding linear maps of the tensor product if this view
        is more convenient.
    \end{enumerate}
\end{remark}

The above mentioned implementation of the additional properties will
be practicable if they are given in the form
\eqref{eq:HC-type-multiplication} and \eqref{eq:HC-type-module} and if
these specifications can be extended to the whole Hochschild complexes
yielding subcomplexes which are closed under the existing
operations. In the following, the structures of such particular types
will be the ones of interest. The necessary requirements are
summarized in the following two definitions.

\begin{definition}[Associative algebras of particular types]
    \label{definition:algebra-type}
    An associative $\mathbb{K}$-algebra $(\mathcal{A},\mu)$ is said to
    be of a particular \emph{type} if for all $k\in \mathbb{N}$ there
    are vector subspaces $\HCtype^k(\mathcal{A},\mathcal{A})\subseteq
    \HC^k(\mathcal{A},\mathcal{A})$ of cochains such that with
    $\HCtype^0(\mathcal{A},\mathcal{A})= \mathcal{A}$ the following
    assertions hold.
    \begin{enumerate}
    \item $\mu \in \HCtype^2(\mathcal{A},\mathcal{A})$.
    \item $\HCtype^{\bullet}(\mathcal{A},\mathcal{A})$ is closed under
        the insertions $\circi$ after the $i$-th position and the
        natural action of the symmetric group on the arguments in
        $\mathcal{A}$.
    \end{enumerate}
\end{definition}

Regarding the structures defined in Section \ref{sec:Hochschild} it is
clear that the above assumptions assure that
$\HCtype^{\bullet}(\mathcal{A},\mathcal{A})$ is a subcomplex of
$\HC^{\bullet} (\mathcal{A},\mathcal{A})$ with the Hochschild
differential $\delta$ which only depends on $\mu$. Further, the
related Hochschild cohomology $\HHtype^{\bullet} (\mathcal{A},
\mathcal{A})$ is a Gerstenhaber algebra with the same multiplication
and bracket.

\begin{definition}[Module structures of particular types]
    \label{definition:module-type}
    Let $(\mathcal{A},\mu)$ be an associative $\mathbb{K}$-algebra of
    a particular type. An $\mathcal{A}$-module structure $\rho$ on a
    $\mathbb{K}$-vector space $\mathcal{E}$ is said to be of a
    particular \emph{type} which is \emph{compatible} with
    $\HCtype^{\bullet}(\mathcal{A},\mathcal{A})$ if for all $k\in
    \mathbb{N}_0$ there are vector subspaces $\HCtype^k
    (\mathcal{A},\mathcal{D})\subseteq
    \HC^k(\mathcal{A},\End_{\mathbb{K}}(\mathcal{E}, \mathcal{E}))$ of
    cochains with values in a subalgebra $\mathcal{D}
    \subseteq \End_{\mathbb{K}} (\mathcal{E}, \mathcal{E})$ such that
    the following assertions hold.
    \begin{enumerate}
    \item $\rho\in \HCtype^1(\mathcal{A}, \mathcal{D})$.
    \item $\id_\mathcal{E} \in \HCtype^0 (\mathcal{A},
        \mathcal{D})\subseteq \mathcal{D}$.
    \item $\HCtype^{\bullet} (\mathcal{A}, \mathcal{D})$ is closed
        under the insertions $\circi$ after the $i$-th position with
        respect to the complex $\HCtype^{\bullet} (\mathcal{A},
        \mathcal{A})$, the cup product $\cup$ in
        \eqref{eq:cup-product-module}, and the natural action of the
        symmetric group on the arguments in $\mathcal{A}$.
    \end{enumerate}
\end{definition}

The third condition for the cup product implies that
$\HCtype^0(\mathcal{A},\mathcal{D})\subseteq \mathcal{D}$ is a
subalgebra. With \eqref{eq:bimodule-of-endomorphisms} it is obvious
that $\mathcal{D}$ is an $(\mathcal{A},\mathcal{A})$-subbimodule of
$\End_{\mathbb{K}}(\mathcal{E},\mathcal{E})$ which in general is not
true for $\HCtype^0(\mathcal{A},\mathcal{D})$. The above conditions
and \eqref{eq:Hochschild-differential-module}, or
\eqref{eq:new-Hochschild-differential} respectively, assure that
$\HCtype^{\bullet} (\mathcal{A},\mathcal{D})$ is a subcomplex of
$\HC^{\bullet} (\mathcal{A},\End_{\mathbb{K}}(\mathcal{E},
\mathcal{E}))$. The corresponding Hochschild cohomology then is
denoted by $\HHtype^{\bullet}(\mathcal{A},\mathcal{D})$.
% With respect to later applications we make the following obvious but
% useful observation.
%%% 
%%% folgendes Lemma ist nur richtig, wenn D' mit type verträglich
%%% 
% \begin{lemma}[Module structures of refined types]
%     Let $\rho\in \HCtype^1(\mathcal{A}, \mathcal{D})$ be a module
%     structure of a particular type with subalgebra $\mathcal{D}
%     \subseteq \End_{\mathbb{K}} (\mathcal{E}, \mathcal{E})$ as in
%     Definition~\ref{definition:module-type}. Further, let $\rho$ in
%     fact take values in a subalgebra $\mathcal{D}'\subseteq
%     \mathcal{D}$ with $\id\in \mathcal{D}'$.\\
%     Then, $\rho\in \HCtype^1(\mathcal{A}, \mathcal{D}')$ is again a
%     module structure of a particular type by considering $\HCtype^k
%     (\mathcal{A},\mathcal{D}')$. Then, this type is referred to as the
%     \emph{refined type} of the initial one.
% \end{lemma}

% \begin{remark}
%     \label{remark:Hochschild-type-practice}
%     It is clear that in practice one first has to find adequate
%     subsets of the algebraic Hochschild complexes which satisfy the
%     stated three conditions, in particular the closedness under the
%     insertions and the cup product. Then, algebra multiplications
%     and right module structures which are elements in such subsets
%     by definition are of the given type and finally yield the
%     corresponding subcomplexes.
% \end{remark}

%\begin{remark}[Deformations of algebras of a particular type]
%    \label{remark:deformations-algebra-type}

%\end{remark}

\section{Deformations of associative algebras}
\label{sec:deformation-algebras}

This section is a summary of the basic definitions and results
concerning Gerstenhaber's deformation theory of associative algebras
which later are easily extended in the setting of modules. The
slightly more specific framework used here does not affect the proofs
of the statements which can all be found in the original work
\cite{gerstenhaber:1964a} of Gerstenhaber. Given a particular type of
associative algebras as in Definition~\ref{definition:algebra-type} it
is obvious how the notion of deformations has to be reformulated in
this framework. Since all relevant structures, in particular the
deformations itself and the equivalence transformations between them,
can be expressed in terms of the Hochschild complex the restriction to
the subcomplex of the given type yields the desired refined
notion. The used notation in the following presentation closely
follows \cite{waldmann:2007a}. Using the notion of formal power series
in a formal parameter $\lambda$ as explained in
Section~\ref{sec:noncommutative-space-times} the crucial definition
now is the following.

\begin{definition}[Formal associative deformation]
    \label{definition:deformation-algebra}
    Let $(\mathcal{A},\mu_0)$ be an associative $\mathbb{K}$-algebra
    of a particular type as in Definition \ref{definition:algebra-type}.
    \begin{enumerate}
    \item A formal associative deformation of $(\mathcal{A},\mu_0)$ is
        a $\mathbb{K}[[\lambda]]$-bilinear associative multiplication
        $\mu$ for $\mathcal{A}[[\lambda]]$ of the form
        \begin{equation}
            \label{eq:formal-associative-deformation}
            \mu= \sum_{r=0}^\infty \lambda^r \mu_r:
            \mathcal{A}[[\lambda]] \times
            \mathcal{A}[[\lambda]] \longrightarrow
            \mathcal{A}[[\lambda]] 
        \end{equation}
        with $\mathbb{K}$-bilinear maps $\mu_r\in
        \HCtype^2(\mathcal{A},\mathcal{A})$.
    \item Two such deformations $\mu$ and $\tilde{\mu}$ are said to be
        equivalent if there exists a $\mathbb{K}[[\lambda]]$-linear
        algebra isomorphism $S= \id+ \sum_{r=1}^\infty \lambda^r S_r:
        (\mathcal{A}[[\lambda]],\mu) \longrightarrow
        (\mathcal{A}[[\lambda]], \tilde{\mu})$ with $S_r\in
        \HCtype^1(\mathcal{A},\mathcal{A})$.
    \item Associative deformations and their equivalence up to an
        order $r\in \mathbb{N}$ are defined in the corresponding way.
    \end{enumerate}
\end{definition}

With the same argumentation as in Lemma
\ref{lemma:condition-associative-multiplication} a formal power series
$\mu= \sum_{r=0}^\infty \lambda^r \mu_r$ of $\mathbb{K}$-linear maps
$\mu_r:\mathcal{A}\times \mathcal{A} \longrightarrow \mathcal{A}$
defines an associative multiplication if the correspondingly defined
Gerstenhaber bracket $[\mu,\mu]=0$ vanishes. This condition has to be
fulfilled in all orders of $\lambda$. In order $r\in \mathbb{N}_0$ it
reads 
\begin{equation}
    \label{eq:condition-associativity-orders}
    \sum_{s=0}^r [\mu_s, \mu_{r-s}]=0.
\end{equation}
In lowest order $\lambda^0$ this is simply the associativity condition
$[\mu_0,\mu_0]=0$ for the undeformed product. By Lemma
\ref{lemma:Hochschild-algebra} the Hochschild differential can be
expressed in terms of the Gerstenhaber bracket with $\mu_0$. Since
$[\mu_0,\mu_r]=[\mu_r,\mu_0]= -\delta\mu_r$, condition
\eqref{eq:condition-associativity-orders} can be read as
\begin{equation}
    \label{eq:condition-associativity-with-differential-order-one}
    \delta \mu_1=0
\end{equation}
for the first order and for all higher orders $r\ge 2$ one gets
\begin{equation}
    \label{eq:condition-associativity-with-differential}
    \delta \mu_r= \frac{1}{2} \sum_{s=1}^{r-1} [\mu_s, \mu_{r-s}].
\end{equation}
The last equation points out a possibility to construct such
deformations order by order and leads to the following well-known
result, confer \cite[Chap.~1, Prop.~3]{gerstenhaber:1964a}.

\begin{proposition}
    \label{proposition:obstruction-existence-algebra}
    Let $(\mathcal{A},\mu_0)$ be an associative $\mathbb{K}$-algebra
    of a particular type. Further, let $\mu^{(r)}= \mu_0+ \dots
    +\lambda^r \mu_r$ be an associative deformation up to order
    $r$. Then the condition for $\mu_{r+1}\in
    \HCtype^2(\mathcal{A},\mathcal{A})$ to define an associative
    deformation $\mu^{(r+1)}= \mu^{(r)}+\lambda^{r+1} \mu_{r+1}$ up to
    order $r+1$ is
    \begin{equation}
        \label{eq:condition-mu-order-by-order}
        \delta \mu_{r+1}= R_r 
    \end{equation}
    with $R_r \in \HCtype^3(\mathcal{A},\mathcal{A})$ explicitly given
    by $R_0=0$ and
    \begin{equation}
        \label{eq:obstruction-existence-algebra}
        R_r= \frac{1}{2} \sum_{s=1}^r
        [\mu_s,\mu_{r+1-s}] \quad \textrm{for } r\ge 1.
    \end{equation}
    Moreover, $\delta R_r=0$ whence the obstruction in order
    $r+1$ is the class $[R_r]\in
    \HHtype^3(\mathcal{A},\mathcal{A})$.
\end{proposition}

Concerning equivalence one finds the following well-known statement,
confer \cite[Chap.~1, Prop.~1]{gerstenhaber:1964a}.

\begin{proposition}
 \label{proposition:equivalence-algebra}
 Let $(\mathcal{A},\mu_0)$ be an associative $\mathbb{K}$-algebra of a
 particular type.
    \begin{enumerate}
    \item An associative deformation of $\mu$ is always equivalent to
        a deformation of the form $\tilde{\mu}= \mu_0+
        \sum_{s=r}^\infty \lambda^s \tilde{\mu}_s$ where the first
        non-vanishing cochain $\tilde{\mu}_r$ is a cocycle, $\delta
        \tilde{\mu}_r=0$ but no coboundary.
    \item If $\HHtype^2(\mathcal{A},\mathcal{A})=0$ is trivial, all
        associative deformations are equivalent.
    \end{enumerate}
\end{proposition}

The two propositions show that the second and the third Hochschild
cohomology groups of an algebra contain crucial information about its
deformation theory. $\HHtype^3(\mathcal{A},\mathcal{A})$ encodes the
obstruction for a continuation of a given deformation up to an
arbitrary order $r$ to a deformation up to order $r+1$. So, if this
cohomology group is trivial there exists a simple way to construct
associative deformations order by order. The elements in
$\HHtype^2(\mathcal{A}, \mathcal{A})$ are exactly the equivalence
classes of deformations up to order one. This is the case since $\mu=
\sum_{s=0}^\infty \lambda^s \mu_s$ is an associative deformation up to
order one if and only if $\delta \mu_1=0$ and two such deformations
$\mu$ and $\tilde{\mu}$ are equivalent up to order one if and only if
$\mu_1-\tilde{\mu}_1=\delta \phi$ is a coboundary. The general
classification of associative deformations, however, is more
difficult. But nevertheless, if the second Hochschild cohomology group
vanishes all deformations are equivalent.

\section{Deformations of right modules}
\label{sec:deformation-modules}

Motivated by Gerstenhaber's deformation theory of associative algebras
and algebras with other properties, confer
\cite[Chap.~1]{gerstenhaber:1964a}, it is now a straightforward
generalization which yields the notion of a deformation theory of
modules over algebras in the slightly modified framework. A first
purely algebraic definition for finite dimensional modules was given
in the work of Donald and Flanigan \cite{donald.flanigan:1974a}. There
and in many other works on the topic, confer \cite{yau:2005} and the
references therein, the module is seen as a ring homomorphism as in
\eqref{eq:right-module-as-ring-morphism} which then is the structure
to be deformed. In the case which is relevant in this work, though, it
is not sufficient only to perform this approach but one has to look
for a more general notion where the deformation is a right module with
respect to a given deformation of the underlying algebra. Without loss
of generality this new concept is only considered for right
modules. For left modules all results can be derived similarly.

\begin{definition}[Deformation of a right module structure]
    \label{definition:deformation-right-module}
    Let $(\mathcal{A},\mu_0)$ be an associative $\mathbb{K}$-algebra
    and $(\mathcal{E},\rho_0)$ be a right $\mathcal{A}$-module, both
    of particular types as in the Definitions
    \ref{definition:algebra-type} and \ref{definition:module-type}.
    Further let $\mu=\sum_{r=0}^\infty \lambda^r \mu_r:
    \mathcal{A}[[\lambda]] \times \mathcal{A}[[\lambda]]
    \longrightarrow \mathcal{A}[[\lambda]]$ be a formal associative
    deformation of $\mu_0$ with $\mu_r \in
    \HCtype^2(\mathcal{A},\mathcal{A})$ as in Definition
    \ref{definition:deformation-algebra}.
    \begin{enumerate}
    \item A \emph{deformation of the right module structure $\rho_0$
          with respect to $\mu$ of the given type} is a
        $\mathbb{K}[[\lambda]]$-bilinear right
        $(\mathcal{A}[[\lambda]],\mu)$-module structure $\rho$ of
        $\mathcal{E}[[\lambda]]$ of the form
        \begin{equation}
            \label{eq:module-deformation}
            \rho=\sum_{r=0}^\infty \lambda^r \rho_r: \mathcal{E}[[\lambda]]
            \times \mathcal{A}[[\lambda]] \longrightarrow
            \mathcal{E}[[\lambda]] 
        \end{equation}
        with $K[[\lambda]]$-bilinear maps $\rho_r\in
        \HCtype^1(\mathcal{A},\mathcal{D})$.
    \item Two such deformations $\rho$ and $\tilde{\rho}$ are said to
        be equivalent if there exists a formal series
        \begin{equation}
            \label{eq:equivalence-trafo-module-deformation}
            T= \id_{\mathcal{E}}+ \sum_{r=1}^\infty \lambda^r T_r \in
            \HCtype^0(\mathcal{A},\mathcal{D})[[\lambda]]
        \end{equation}
        such that $T(\rho(e,a))= \tilde{\rho}(T e,a)$, or equivalently
        \begin{equation}
            \label{eq:equivalence-module-deformation-HC}
            T\circ \rho(a) = \tilde{\rho}(a)\circ T,
        \end{equation}
        for all $a\in \mathcal{A}[[\lambda]]$ and $e\in
        \mathcal{E}[[\lambda]]$.
    \item The deformations up to an order $r\in \mathbb{N}$ as well as
        their equivalences are defined in a corresponding way.
    \end{enumerate}
\end{definition}

\begin{remark}
    Note that the condition
    \eqref{eq:equivalence-module-deformation-HC} for equivalence can
    be formulated in terms of the cup product
    \eqref{eq:cup-product-module} and then reads
    \begin{equation}
        \label{eq:equivalence-module-deformation-cup}
        T\cup \rho= \tilde{\rho}\cup T.
    \end{equation}
    Due to \eqref{eq:inverse-element} it is clear that a formal series
    $T$ of operators as in
    \eqref{eq:equivalence-trafo-module-deformation} is invertible with
    $T^{-1}\in \HCtype^0(\mathcal{A},\mathcal{D})[[\lambda]]$. Thus,
    such a map $T$ and a deformation $\rho$ of the considered type
    define a new one via
    \begin{equation}
        \label{eq:defining-new-module-structure}
        \tilde{\rho}=T\cup \rho\cup T^{-1}.
    \end{equation}
    Note that analogue assertions already hold for deformed algebra
    structures.
\end{remark}

A comparison of the Hochschild complex of an associative algebra
$\mathcal{A}$ with the one obtained in Section
\ref{subsec:Hochschild-complex-module} when regarding the algebra as a
right module over itself already yields an indication which Hochschild
cohomologies are crucial for the deformations of right modules. It is
evident that for all $k\geq 1$ there is a natural way to identify the
cochains in $\HC^k(\mathcal{A},\mathcal{A})$ with the ones in $\HC^{k
  -1}(\mathcal{A}, \Hom_{\mathbb{K}} (\mathcal{A},
\mathcal{A}))$. This means that the transition to the module point of
view induces a simple shift in the grading of the considered
Hochschild complexes. Since the obstructions discussed in Section
\ref{sec:deformation-algebras} lie in the second and third cohomology
one could thus already guess that the obstruction theory for the
deformation of general right modules is closely related to the first
and second Hochschild cohomology $\HHtype^{\bullet}
(\mathcal{A},\mathcal{D})$. That this is really the case is the
content of the following two propositions, confer \cite[Lemma~2.1 and
2.2]{bordemann.neumaier.waldmann.weiss:2007a:pre}.

\begin{proposition}
    \label{proposition:obstruction-deformation-module}
    In the setting of
    Definition~\ref{definition:deformation-right-module}, let
    $\rho^{(r)}=\rho_0+ \dots + \lambda^r \rho_r$ be a right module
    structure with respect to $\mu$ up to order $r$ with $\rho_s\in
    \HCtype^1(\mathcal{A},\mathcal{D})$ for all $s=0, \dots, r$. Then
    the condition for $\rho_{r+1}\in
    \HCtype^1(\mathcal{A},\mathcal{D})$ to define a right module
    structure $\rho^{(r+1)}= \rho^{(r)}+ \lambda^{r+1} \rho_{r+1}$ up
    to order $r+1$ is
    \begin{equation}
        \label{eq:condition-rho-order-by-order}
        \delta \rho_{r+1} = R_r
    \end{equation}
    with $R_r \in \HCtype^2 (\mathcal{A}, \mathcal{D})$ explicitly
    given by $R_0=\rho_0\circ \mu_1$ and
    \begin{equation}
        \label{eq:obstruction-existence}
        R_r (a, b) 
        =  
        \sum_{s=0}^r \rho_s (\mu_{r+1-s}(a, b))
        - \sum_{s=1}^r \rho_s(b) \circ \rho_{r+1-s}(a) \quad
        \textrm{for } r\ge 1.
    \end{equation}
    Moreover, $\delta R_r = 0$ whence the recursive obstruction in
    order $r+1$ is the class $[R_r] \in \HHtype^2(\mathcal{A},
    \mathcal{D})$.
\end{proposition}

\begin{proof}
    The guideline of the proof is of course given by the well-known
    considerations in the basic references \cite{gerstenhaber:1964a,
      donald.flanigan:1974a, yau:2005}. Nevertheless, it is carried
    out completely in order to show that all the slight modifications
    match together in more specific framework.

    The main goal in general is of course to find $\rho_r\in
    \HCtype^1(\mathcal{A},\mathcal{D})$ for all $r\in \mathbb{N}$ such
    that $\rho=\sum_{r=0}^{\infty} \lambda^r \rho_r$ is a right module
    structure with respect to $\mu$. The failure of an arbitrary
    $\rho$ to be such a right module structure is encoded in the
    $\mathbb{K}[[\lambda]]$-bilinear map $C: \mathcal{A}[[\lambda]]
    \times \mathcal{A}[[\lambda]]
    \longrightarrow \End_{\mathbb{K}[[\lambda]]}
    (\mathcal{E}[[\lambda]], \mathcal{E}[[\lambda]])$ with
    $C(a,b)=\rho(\mu(a,b))-\rho(b) \circ \rho(a)$. With
    $\End_{\mathbb{K}[[\lambda]]} (\mathcal{E}[[\lambda]],
    \mathcal{E}[[\lambda]])= \End_{\mathbb{K}}
    (\mathcal{E},\mathcal{E})[[\lambda]]$ and the conditions in
    Definition \ref{definition:module-type} it is clear that
    $C=\sum_{r=0}^{\infty} \lambda^r C_r$ with elements $C_r \in
    \HCtype^2(\mathcal{A},\mathcal{D})$ given by $C_r(a,b)=
    \sum_{s=0}^r [\rho_s (\mu_{r-s}(a,b)) -\rho_s(b) \circ
    \rho_{r-s}(a)]$. Due to the associativity of $\mu$ the map $C$
    satisfies the equation
    \begin{eqnarray}
        \label{eq:lemma-for-failure-module}
        \lefteqn{C(\mu(a,b),c) - C(a,\mu(b,c))} \nonumber\\
        &=& \rho(\mu(\mu(a,b),c))-
        \rho(c) \circ \rho(\mu(a,b)) - \rho(\mu(a,\mu(b,c))) +
        \rho(\mu(b,c))\circ \rho(a) \nonumber\\
        &=&  \rho(\mu(b,c))\circ \rho(a) -\rho(c) \circ \rho(b) \circ
        \rho(a) + \rho(c) \circ \rho(b) \circ
        \rho(a) -\rho(c) \circ \rho(\mu(a,b))\\
        &=& C(b,c) \circ \rho(a) - \rho(c) \circ C(a,b). \nonumber
    \end{eqnarray}
    After these preliminary considerations the actual proof is the
    following. By assumption, there are given $\rho_0, \dots, \rho_r$
    such that $C_0=\dots = C_r=0$. For an arbitrary $\rho_{r+1}\in
    \HCtype^1(\mathcal{A},\mathcal{D})$ one gets
    \begin{eqnarray*}
        (\delta \rho_{r+1})(a,b) &=& a\cdot \rho_{r+1}(b)
        -\rho_{r+1}(\mu_0(a,b)) +\rho_{r+1} (a)\cdot b\\
        &=& \rho_{r+1}(b) \circ \rho_0(a) - \rho_{r+1}(\mu_0(a,b)) +
        \rho_0(b) \circ \rho_{r+1}(a)
    \end{eqnarray*}
    and thus the failure of the so defined $\rho^{(r+1)}$ to be a
    right module structure up to order $r+1$ is given by
    \begin{eqnarray*}
        C_{r+1}(a,b) &=& \rho_{r+1}(\mu_0(a,b)) + \sum_{s=0}^r
        \rho_s(\mu_{r+1-s}(a,b)) \\
        &&- \rho_{r+1}(b) \circ \rho_0(a) -
        \rho_0(b) \circ \rho_{r+1}(a) - \sum_{s=1}^r \rho_s(b) \circ
        \rho_{r+1-s}(a)\\
        &=& -(\delta \rho_{r+1})(a,b) + R_r(a,b)
    \end{eqnarray*}
    with $R_r \in \HCtype^2(\mathcal{A},\mathcal{D})$ as in
    \eqref{eq:obstruction-existence}. Clearly, $R_0=\rho_0\circ
    \mu_1$. Then, condition \eqref{eq:condition-rho-order-by-order}
    follows immediately. Since $C_0=\dots =C_r=0$
    Equation~\eqref{eq:lemma-for-failure-module} yields in order $r+1$
    \begin{eqnarray*}
        0 &=& C_{r+1}(b,c) \circ \rho_0(a) - C_{r+1}(\mu_0(a,b),c) +
        C_{r+1}(a, \mu_0(b,c)) - \rho_0(c) \circ C_{r+1}(a,b)\\
        &=& (\delta C_{r+1})(a,b,c)\\
        &=& (\delta R_r)(a,b,c)
    \end{eqnarray*}
    where one makes use of $\delta^2=0$. So $\delta R_r=0$ and the
    remaining assertions follow.
\end{proof}

In particular, the proposition states that if the obstruction vanishes
an order by order construction can be done which yields a deformation
$\rho$ of the desired type. So in this case the existence of a
deformation is given.

\begin{corollary}[Existence of module deformations]
    \label{corollary:existence-module-deformations}
    If the second Hochschild cohomology of a right module structure
    $\rho_0$ of a particular type is trivial,
    $\HHtype^2(\mathcal{A},\mathcal{D})= \{0\}$, all deformations up
    to a finite order can be extended to deformations.
\end{corollary}

In a similar way one gets the following statement concerning the
equivalence of two given deformations.

\begin{proposition}
    \label{proposition:obstruction-equivalence-module}
    In the setting of Definition
    \ref{definition:deformation-right-module}, let $\rho$ and
    $\tilde{\rho}$ be two deformations of $\rho_0$ of a particular
    type and let $T^{(r)}= \id + \dots +\lambda^r T_r$ be an
    equivalence transformation between them up to order $r$ with $T_s
    \in \HCtype^0(\mathcal{A},\mathcal{D})$ for $s=0,\dots,r$. Then
    the condition for $T_{r+1} \in \HCtype^0(\mathcal{A},\mathcal{D})$
    to define an equivalence transformation $T^{(r+1)}= T^{(r)} +
    \lambda^{r+1} T_{r+1}$ up to order $r+1$ is
    \begin{equation}
        \label{eq:condition-equivalence-order-by-order}
        \delta T_{r+1} = E_r
    \end{equation}
    with $E_r\in \HCtype^1(\mathcal{A}, \mathcal{D})$ explicitly given by
    \begin{equation}
        \label{eq:obstruction-equivalence}
        E_r(a)= \sum_{s=0}^r (\tilde{\rho}_{r+1-s} (a) \circ T_s - T_s
        \circ \rho_{r+1-s}(a)).
    \end{equation}
    Moreover, $\delta E_r=0$ whence the recursive obstruction in order
    $r+1$ is the class $[E_r]\in
    \HHtype^1(\mathcal{A},\mathcal{D})$.
\end{proposition}

\begin{proof}
    Again, the proof starts with some well-known preliminary
    considerations. The failure of arbitrary $T_r\in
    \HCtype^0(\mathcal{A},\mathcal{D})$, $r\in \mathbb{N}$, to yield
    an equivalence transformation $T=\id + \sum_{r=1}^{\infty}
    \lambda^r T_r$ is now encoded in the
    $\mathbb{K}[[\lambda]]$-linear map $C: \mathcal{A}[[\lambda]]
    \longrightarrow \End_{\mathbb{K}[[\lambda]]}
    (\mathcal{E}[[\lambda]], \mathcal{E}[[\lambda]])$ with $C(a) = T
    \circ \rho(a) - \tilde{\rho}(a) \circ T$. Again,
    $C=\sum_{r=0}^{\infty} \lambda^r C_r$ with $C_r \in
    \HCtype^1(\mathcal{A},\mathcal{D})$ given by $C_r(a)= \sum_{s=0}^r
    [T_s \circ \rho_{r-s}(a) -\tilde{\rho}_s(a) \circ T_{r-s}]$.  Due
    to the right module property of $\rho$ and $\tilde{\rho}$ the map
    $C$ satisfies the equation
    \begin{eqnarray}
        \label{eq:lemma-for-failure-equivalence}
        C(\mu(a,b)) &=& T\circ \rho(\mu(a,b)) - \tilde{\rho}(\mu(a,b))
        \circ T \nonumber \\
        &=& T\circ \rho(b)\circ \rho(a) - \tilde{\rho}(b) \circ
        \tilde{\rho}(a) \circ T  \\
        &=& C(b) \circ \rho(a) + \tilde{\rho}(b)\circ C(a). \nonumber
    \end{eqnarray}
    Now, by assumption there are given $T_0=\id, \dots, T_r$ such that
    $C_0=\dots = C_r=0$. Since $\rho$ and $\tilde{\rho}$ coincide in
    zeroth order the failure of the extension $T^{(r+1)}$ for an
    arbitrary $T_{r+1}\in \HCtype^0(\mathcal{A},\mathcal{D})$ to be an
    equivalence transformation up to order $r+1$ reads
    \begin{eqnarray*}
        C_{r+1}(a) 
        &=& T_{r+1} \circ \rho_0(a) - \rho_0(a) \circ T_{r+1} +
        \sum_{s=0}^r ( T_s \circ \rho_{r+1-s}(a) -
        \tilde{\rho}_{r+1-s}(a) \circ T_s) \\ 
        &=& (\delta T_{r+1})(a) - E_r(a)
    \end{eqnarray*}
    with $E_r \in \HCtype^1(\mathcal{A},\mathcal{D})$ as in
    \eqref{eq:obstruction-equivalence}. Then, condition
    \eqref{eq:condition-equivalence-order-by-order} follows
    immediately. Finally,
    Equation~\eqref{eq:lemma-for-failure-equivalence} yields in order
    $r+1$
    \begin{eqnarray*}
        0 &=& C_{r+1}(b) \circ \rho_0(a) - C_{r+1}(\mu_0(a,b)) +
        \rho_0(b) \circ C_{r+1}(a)\\
        &=& (\delta C_{r+1})(a,b)\\
        &=& - (\delta E_r)(a,b),
    \end{eqnarray*}
    and thus $\delta E_r=0$.
\end{proof}

Again, the proposition assures that if the obstruction vanishes the
orderwise construction of equivalence transformations can be carried
out.

\begin{corollary}[Equivalence of module deformations]
    \label{corollary:equivalence-module-deformations}
    If the first Hochschild cohomology of a right module structure
    $\rho_0$ of a particular type is trivial,
    $\HHtype^1(\mathcal{A},\mathcal{D})= \{0\}$, the module structure
    is \emph{rigid}. This means that all deformations of $\rho$ with
    respect to the same deformation of the underlying algebra, are
    equivalent.
\end{corollary}

Note that even if there occur obstructions in higher orders, this
means $[E_r]\neq 0$ in $\HHtype^1(\mathcal{A},\mathcal{D})$, it could
still be possible to construct equivalence transformations order by
order in a more complicated way than in the presented one. For example
it could still be possible to \emph{change} the already found
$T_1,\dots, T_r$ in order to find an adequate $T_{r+1}$.

\begin{remark}
    If the algebra $(\mathcal{A},\mu_0)$ is commutative, it is obvious
    that all the results concerning module deformations do not depend
    on the choice between the two possible bimodule structures
    \eqref{eq:bimodule-of-endomorphisms} and
    \eqref{eq:new-bimodule-of-endomorphism} which yield different
    Hochschild complexes. In both cases the statements are exactly the
    same and which complex finally will be chosen is just a matter of
    convenience.
\end{remark}

Before we come to a more concrete example in the subsequent subsection
we make the following interesting general observation.

\begin{proposition}[Rigidity of the unit]
    \label{proposition:unit-acts-as-unit}
    Let $(\mathcal{A},\mu_0)$ be a unital algebra with unit $1$
    together with a deformation $\mu=\sum_{r=0}^\infty\lambda^r \mu_r$
    such that $\mu(a,1)=a=\mu(1,a)$ for all $a\in
    \mathcal{A}[[\lambda]]$. Further, let $(\mathcal{E},\rho_0)$ be a
    right $\mathcal{A}$-module with $\rho_0(1)=\id_\mathcal{E}$.
    Then, every deformation $\rho=\sum_{r=0}^\infty\lambda^r\rho_r$ of
    $\rho_0$ with respect to $\mu$ satisfies
    \begin{equation}
        \label{eq:unit-acts-as-unit}
        \rho(1)=\id_{\mathcal{E}[[\lambda]]}.
    \end{equation}
\end{proposition}

\begin{proof}
    Due to the assumptions, $\rho(1)= \id_{\mathcal{E}}+
    \sum_{r=1}^\infty\lambda^r\rho_r(1)$ is an invertible
    map. Then \eqref{eq:unit-acts-as-unit} follows from
    $\rho(1)\circ \rho(1)= \rho(1)$ which is a consequence
    of the given right module property. 
\end{proof}

Note that the assertion in
Proposition~\ref{proposition:unit-acts-as-unit} does not depend on a
particular type and is thus valid in general.

\section{The commutant of a deformed module structure}
\label{sec:commutant-module-structures}

The set of module endomorphisms was already introduced in the end of
Section \ref{subsec:Hochschild-complex-module} and can be considered
for every module structure. Now we come back to this topic and
investigate the relation between the commutants of a given right
module structure and its deformations. Whenever one is given a right
module structure $\rho_0:\mathcal{E}\times \mathcal{A}\longrightarrow
\mathcal{E}$ of a particular type as in Definition
\ref{definition:deformation-right-module} one is interested in the
module endomorphisms which are operators with the correct behaviour as
well. In general, the commutant of interest is the one within the
subalgebra $\mathcal{D}\subseteq \End_{\mathbb{K}}
(\mathcal{E},\mathcal{E})$, explicitly given by
\begin{equation}
    \label{eq:commutant-undeformed}
    \mathcal{B}_0=  
    \left\{ 
        A\in \mathcal{D}\; \big| \; A \circ \rho_0(a) = \rho_0(a)
        \circ A \; \textrm{for all} \; a \in \mathcal{A} 
    \right\}.
\end{equation}
It is clear that $\mathcal{B}_0$ is a subalgebra of $\mathcal{D}$ with
the usual composition $\circ$ of maps as multiplication. It will be
referred to as \emph{classical} or \emph{undeformed commutant}.

In analogy to this the commutant of a deformation $\rho$ of $\rho_0$
with respect to a deformation $\mu$ of $\mu_0$ as in Definition
\ref{definition:deformation-right-module} is defined by
\begin{equation}
    \label{eq:commutant-of-deformation}
    \mathcal{K}_0 = 
    \left\{
        A \in \mathcal{D}[[\lambda]] 
        \; \big| \;
        A \circ \rho(a) = \rho(a) \circ A
        \; \textrm{for all} \; a \in \mathcal{A}[[\lambda]]
    \right\}.
\end{equation}
If a formal series $A= \sum_{r=0}^\infty \lambda^r A_r \in
\mathcal{D}[[\lambda]]$ is an element of $\mathcal{K}_0$ the condition
in the lowest order of $\lambda$ shows that $A_0 \in \mathcal{B}_0$ is
an element of the classical commutant. This is of course a motivation
to investigate the relation between the commutant $\mathcal{K}_0$ and
the formal power series $\mathcal{B}_0[[\lambda]]$.

This investigation is possible with homological techniques for the
commutants within the subalgebra
$\HCtype^0(\mathcal{A},\mathcal{D})\subseteq \mathcal{D}$. For the
undeformed situation this means to consider the zeroth Hochschild
cohomology
\begin{equation}
    \label{eq:commutant-undeformed-special}
    \mathcal{B}= \HHtype^0(\mathcal{A},\mathcal{D})= \ker
    \delta^0\subseteq \HCtype^0(\mathcal{A},\mathcal{D})
\end{equation} 
which is again a subalgebra $\mathcal{B}\subseteq
\HCtype^0(\mathcal{A},\mathcal{D})$ and a subalgebra
$\mathcal{B}\subseteq \mathcal{B}_0$. The new commutant for the
deformed structure is now given by
\begin{equation}
    \label{eq:commutant-of-deformation-special}
    \mathcal{K} = 
    \left\{
        A \in \HCtype^0(\mathcal{A},\mathcal{D})[[\lambda]] 
        \; \big| \;
        A \circ \rho(a) = \rho(a) \circ A
        \; \textrm{for all} \; a \in \mathcal{A}[[\lambda]]
    \right\}.
\end{equation}

Note that one really has 
\begin{equation}
    \label{eq:notion-commutant-unique}
    \mathcal{B}=\mathcal{B}_0 \quad \textrm{and} \quad
    \mathcal{K}=\mathcal{K}_0 
\end{equation}
if the right module structure is of a particular type with the
additional property that
\begin{equation}
    \label{eq:convenient-type}
    \HCtype^0(\mathcal{A},\mathcal{D})=\mathcal{D}.
\end{equation}
A comparison of \eqref{eq:commutant-of-deformation-special} with the
notion of an equivalence transformation between deformed right module
structures as in Definition \ref{definition:deformation-right-module}
shows that the elements in $\mathcal{K}$ have the same property as the
self equivalence transformations of $\rho$ but do not necessarily
start with the identity in the lowest order. Considering the proof of
Proposition~\ref{proposition:obstruction-equivalence-module} this
yields a guideline how to construct an element of $\mathcal{K}$ from
one of $\mathcal{B}$. For the following proposition we will use the
fact that it is always possible to decompose
$\HCtype^0(\mathcal{A},\mathcal{D})$ into a direct sum
\begin{equation}
    \label{eq:Complement}
    \HCtype^0(\mathcal{A},\mathcal{D}) 
    = \mathcal{B} \oplus \overline{\mathcal{B}}
    = \ker \delta^0 \oplus \overline{\mathcal{B}}
\end{equation}
of $\mathcal{B}$ and a complementary linear subspace
$\overline{\mathcal{B}} \subseteq \HCtype^0(\mathcal{A},\mathcal{D})$
since we work over a field $\mathbb{K}$.

\begin{proposition}[The commutant of a deformed right module structure]
    \label{proposition:description-commutant}
    Let $(\mathcal{E},\rho_0)_{(\mathcal{A},\mu_0)}$ be a right module
    structure of a particular type and let $\rho$ be a corresponding
    deformation of $\rho_0$ with respect to a deformation $\mu$ as in
    Definition \ref{definition:deformation-right-module}. Further, let
    the first Hochschild cohomology be trivial,
    \begin{equation}
        \label{eq:first-cohomology-vanishing}
        \HHtype^1(\mathcal{A},\mathcal{D})= \{ 0 \}.
    \end{equation}
    Then, every choice of a complementary subspace
    $\overline{\mathcal{B}}$ induces a unique
    $\mathbb{K}[[\lambda]]$-linear bijection
    \begin{equation}
        \label{eq:bijection-commutant}
        \rho': \mathcal{B}[[\lambda]]
        \longrightarrow
        \mathcal{K} \subseteq
        \HCtype^0(\mathcal{A},\mathcal{D})[[\lambda]] 
    \end{equation}
    of the form $\rho' = \id + \sum_{r=1}^\infty \lambda^r \rho'_r$
    with $\mathbb{K}$-linear maps
    \begin{equation}
        \label{eq:corrections-in-complement}
        \rho'_r: \mathcal{B} \longrightarrow
        \overline{\mathcal{B}} \subseteq
        \HCtype^0(\mathcal{A},\mathcal{D}) 
    \end{equation}
    for $r\geq 1$.
\end{proposition}

\begin{proof}
    As stated above, an element $A\in \mathcal{B}$ can be seen as the
    zeroth order of an element $\rho'(A)\in \mathcal{K}$ which shall
    be constructed recursively order by order in a unique way so that
    it is the value of $A$ under a map $\rho'=\id+ \sum_{r=1}^\infty
    \lambda^r \rho'_r$. To this end the results of
    Proposition~\ref{proposition:obstruction-equivalence-module} can
    be extended to the present situation. Assume to have already found
    $\rho'^{(r)}(A)=A+ \dots+ \lambda^r \rho'_r(A)$ satisfying the
    requested condition to be in $\mathcal{K}$ up to order $r$, this
    means that $\rho'^{(r)}(A)\circ \rho(a)- \rho(a)\circ
    \rho'^{(r)}(A)= \lambda^{r+1}C_{r+1}+ \cdots$ for all $a\in
    \mathcal{A}$. The condition for extending this to an element
    $\rho'^{(r+1)}(A)=\rho'^{(r)}(A)+ \lambda^{r+1}\rho'_{r+1}(A)$
    satisfying the crucial condition up to order $r+1$ is
    \begin{equation}
        \label{eq:condition-correction-commutant}
        (\delta (\rho'_{r+1}(A)))(a)=E_r(a)= \sum_{s=0}^r (\rho_{r+1-s}
        (a) \circ \rho'_s(A) - \rho'_s(A) \circ \rho_{r+1-s}(a))
    \end{equation}
    with $\delta E_r=0$ in analogy to
    Proposition~\ref{proposition:obstruction-equivalence-module}. Since
    the first Hochschild cohomology is trivial such an element
    $\rho'_{r+1}(A)\in \HCtype^0(\mathcal{A},\mathcal{D})$ always
    exists and it is even unique with the further condition
    \begin{equation}
        \label{eq:condition-correction-in-complement}
        \rho'_{r+1}(A)\in \overline{\mathcal{B}}
    \end{equation}
    for a choice $\HCtype^0(\mathcal{A},\mathcal{D})= \ker \delta^0
    \oplus \overline{\mathcal{B}}$.  This procedure in fact yields an
    injective map
    \begin{equation}
        \label{eq:embedding-classical-commutant}
        \rho'=\id+ \sum_{r=0}^\infty \lambda^r \rho'_r: \mathcal{B}
        \longrightarrow (\mathcal{B}+ \lambda
        \overline{\mathcal{B}}[[\lambda]]) \cap \mathcal{K}.
    \end{equation}
    The maps $\rho'_r: \mathcal{B}\longrightarrow
    \overline{\mathcal{B}}$ are $\mathbb{K}$-linear which is seen by
    an easy induction over $r$ using the linearity of the defining
    conditions \eqref{eq:condition-correction-commutant} and
    \eqref{eq:condition-correction-in-complement}. The
    $\mathbb{K}[[\lambda]]$-linear extension of $\rho'$ has again
    values in $\mathcal{K}$ and is not only injective but also
    surjective which is seen as follows. Let
    $A=\sum_{r=0}^\infty\lambda^r A_r\in \mathcal{K}$. Then
    $A_0\in\mathcal{B}$ and $A-\rho'(A_0)=\sum_{r=1}^\infty \lambda^r
    B_r\in \mathcal{K}$ starts in order one with $B_1\in
    \mathcal{B}$. This way, induction finally yields a preimage
    $A_0+\lambda B_1 +\cdots$ of $A$.
\end{proof}

By definition, the commutant $\mathcal{B}$ yields a
$(\mathcal{B},\mathcal{A})$-bimodule structure of $\mathcal{E}$,
denoted by
\begin{equation}
    \label{eq:notation-bimodule}
    {}_{(\mathcal{B},\circ)} (\id_{\mathcal{B}},
    \mathcal{E},\rho_0)_{(\mathcal{A},\mu_0)}. 
\end{equation}
This notation is reasonable, since $\id_{\mathcal{B}}:
\mathcal{B}\longrightarrow \End_{\mathbb{K}}(\mathcal{E},\mathcal{E})$
can be seen as a left module structure in the same way as $\rho_0$ is
a right module structure. The above proposition implies that the
bimodule \eqref{eq:notation-bimodule} has a corresponding counterpart.

\begin{corollary}[Induced deformations of the classical commutant]
    \label{corollary:deformation-commutant}
    In the situation of Proposition
    \ref{proposition:description-commutant} the choice of a
    complementary subspace $\cc{\mathcal{B}}$ immediately induces the
    following additional structures.
    \begin{enumerate}
    \item The isomorphism $\rho'$ yields an associative deformation
        \begin{equation}
            \label{eq:deformation-commutant}
            \mu'= \circ + \sum_{r=1}^{\infty} \lambda^r \mu'_r:
            \mathcal{B}[[\lambda]] \times
            \mathcal{B}[[\lambda]] \longrightarrow
            \mathcal{B}[[\lambda]]  
        \end{equation}
        of the classical commutant by
        \begin{equation}
            \label{eq:deformation-commutant-explicit}
            \mu'(A,B)= \rho'^{-1}(\rho'(A) \circ \rho'(B)).
        \end{equation}
    \item $\rho'$ further induces a deformed left
        $(\mathcal{B}[[\lambda]], \mu')$-module structure on
        $\mathcal{E}[[\lambda]]$ of the usual action of
        $\mathcal{B}\subseteq \End_{\mathbb{K}}(\mathcal{E},
        \mathcal{E})$ on $\mathcal{E}$, such that
        $\mathcal{E}[[\lambda]]$ becomes a bimodule with respect to
        the two deformed algebras $(\mathcal{B}[[\lambda]], \mu')$ and
        $(\mathcal{A}[[\lambda]], \mu)$.
    \end{enumerate}
\end{corollary}

\begin{proof}
    The map $\mu'$ defined in
    \eqref{eq:deformation-commutant-explicit} starts with $\circ$
    since the isomorphism $\rho'$ and its inverse start with the
    identity in the lowest order in $\lambda$. The associativity is
    clear with the one of the composition $\circ$ of maps. The left
    module structure in the second part is obvious by
    \eqref{eq:deformation-commutant-explicit}. The bimodule structure
    is a direct consequence of the fact that $\rho'$ only takes values
    in the commutant.
\end{proof}

The next definition presents a natural notion of deformations of
bimodule structures.
 
\begin{definition}[Deformation of bimodules]
    \label{definition:deformation-bimodules}
    Let $(\mathcal{A},\mu_0)$ and $(\mathcal{B},\mu'_0)$ be
    associative $\mathbb{K}$-algebras of particular types. Further,
    let $\mathcal{E}$ be a $\mathbb{K}$-vector space endowed with a
    right module structure $\rho_0\in
    \HCtype^1(\mathcal{A},\mathcal{D})$ compatible with
    $\HCtype^{\bullet}(\mathcal{A},\mathcal{A})$ and a left module
    structure $\rho'_0\in \HCtypep^1(\mathcal{B},\mathcal{D}')$
    compatible with $\HCtypep^{\bullet}(\mathcal{B},\mathcal{B})$ such
    that $\mathcal{E}$ is a $(\mathcal{B},\mathcal{A})$-bimodule,
    denoted by
    \begin{equation}
        \label{eq:undeformed-bimodule}
        {}_{(\mathcal{B},\mu'_0)} (\rho'_0,
        \mathcal{E},\rho_0)_{(\mathcal{A},\mu_0)}.
    \end{equation}
    \begin{enumerate}
    \item A deformation of this bimodule
        \eqref{eq:undeformed-bimodule} is a bimodule
        \begin{equation}
            \label{eq:deformation-bimodule}
            {}_{(\mathcal{B}[[\lambda]],\mu')} (\rho',
            \mathcal{E}[[\lambda]],\rho)_{(\mathcal{A}[[\lambda]],\mu)} 
        \end{equation}
        where the left and right module structures $\rho'$ and $\rho$
        are deformations with respect to corresponding algebra
        deformations $\mu'$ and $\mu$ as in Definition
        \ref{definition:deformation-right-module}.
    \item Two such deformations $(\mu',\rho',\rho,\mu)$ and
        $(\tilde{\mu}',\tilde{\rho}',\tilde{\rho},\tilde{\mu})$ are
        said to be equivalent if there exist formal series
        \begin{equation}
            \label{eq:equivalence-transformations-bimodule-left}
            S'= \id_{\mathcal{B}}+ \sum_{r=1}^\infty \lambda^r S'_r, 
            \quad \quad T'= \id_{\mathcal{E}} + \sum_{r=1}^\infty
            \lambda^r T'_r 
        \end{equation}
        with $S'_r\in \HCtypep^1(\mathcal{B},\mathcal{B})$, $T'_r\in
        \HCtypep^0(\mathcal{B},\mathcal{D'})$ and
        \begin{equation}
            \label{eq:equivalence-transformation-bimodule-right}
            S= \id_{\mathcal{A}}+ \sum_{r=1}^\infty \lambda^r S_r, 
            \quad \quad T= \id_{\mathcal{E}} + \sum_{r=1}^\infty
            \lambda^r T_r 
        \end{equation}
        with $S_r\in \HCtype^1(\mathcal{A},\mathcal{A})$, $T_r\in
        \HCtype^0(\mathcal{A},\mathcal{D})$ such that
        \begin{equation}
            \label{eq:equivalence-algebras-explicit}
            S'\circ \mu' = \tilde{\mu}' \circ (S'\otimes S'), \quad \quad 
            S\circ \mu = \tilde{\mu} \circ (S\otimes S)
        \end{equation}
        and
        \begin{equation}
            \label{eq:equivalence-modules-explicit}
            T'\circ \rho'(A)= \tilde{\rho}'(S'A) \circ T', \quad \quad 
            T\circ \rho (a)= \tilde{\rho} (S a) \circ T
        \end{equation}
        for all $A\in \mathcal{B}[[\lambda]]$ and $a\in
        \mathcal{A}[[\lambda]]$.

        In the special situation where
        \begin{equation}
            \label{eq:special-situation-bimodule}
            \HCtype^0(\mathcal{A},\mathcal{D})=
            \HCtypep^0(\mathcal{B},\mathcal{D}')
        \end{equation}
        the two deformations are said to be equivalent if the
        equivalence transformations on $\mathcal{E}[[\lambda]]$
        coincide, this means if
        \begin{equation}
            \label{eq:special-situation-bimodule-equivalence}
            T=T'.
        \end{equation}
    \item The deformations up to order $r\in \mathbb{N}$ as well as
        their equivalence are defined in a corresponding way.
    \end{enumerate}
\end{definition}

\begin{remark}
    \begin{enumerate}
    \item The equations in \eqref{eq:equivalence-algebras-explicit}
        state nothing but the usual equivalence of deformed algebras
        where the algebra multiplications are identified with the
        corresponding linear maps of the tensor product as explained
        in Remark \ref{remark:particular-types}.
    \item The equivalence of deformed bimodules in general is just the
        separate equivalence of the left and the right module
        structures, now with respect to equivalent algebras and not to
        fixed ones as in
        Definition~\ref{definition:deformation-right-module}.
    \item Note that the necessary condition for
        \eqref{eq:special-situation-bimodule} is that either
        $\mathcal{D}\subseteq \mathcal{D}'$ or $\mathcal{D}'\subseteq
        \mathcal{D}$.
    \end{enumerate}
\end{remark}

It is obvious that the bimodule structure in
Corollary~\ref{corollary:deformation-commutant} is such a deformation
of \eqref{eq:notation-bimodule} where the types of the algebra
$(\mathcal{B},\circ)$ and the corresponding left module structure
$\id_{\mathcal{B}}$ of $\mathcal{E}$ are the purely algebraic
ones. This means one considers
$\HCtypep^{\bullet}(\mathcal{B},\mathcal{B}) =
\HC^{\bullet}(\mathcal{B}, \mathcal{B})$ and
\begin{equation}
    \label{eq:type-for-commutant-algebra-left-module}
    \HCtypep^{\bullet} (\mathcal{B},\mathcal{D}') = 
    \HC^{\bullet}(\mathcal{B},\mathcal{D}') \quad \textrm{with} \quad
    \mathcal{D}'= \HCtype^0(\mathcal{A},\mathcal{D})\subseteq
    \mathcal{D} 
\end{equation}
where $\mathcal{D}'$ has the natural $(\mathcal{B},
\mathcal{B})$-bimodule structure given by $B_1\cdot D'\cdot
B_2=B_1\circ D'\circ B_2 $ for $B_1,B_2\in \mathcal{B}$ and $D'\in
\mathcal{D}'$.

The necessary definition
\eqref{eq:type-for-commutant-algebra-left-module} in particular shows
that $\HCtypep^0(\mathcal{B},\mathcal{D}')
=\HCtype^0(\mathcal{A},\mathcal{D})$. Thus the condition
\eqref{eq:special-situation-bimodule} for the stronger notion of
equivalence is always satisfied. One can show that the deformations in
Corollary~\ref{corollary:deformation-commutant} are indeed unique up
to equivalence in this sense.

\begin{proposition}[The uniqueness of the induced deformations]
    \label{proposition:classes-of-deformations}
    Let $(\mathcal{E},\rho_0)_{(\mathcal{A},\mu_0)}$ be a right module
    structure of a particular type.
    \begin{enumerate}
    \item If the first Hochschild cohomology
        $\HHtype^1(\mathcal{A},\mathcal{D})= \{ 0 \}$ is trivial and
        $\mu$ is a given deformation of $\mu_0$ different choices of a
        corresponding deformation $\rho$ of $\rho_0$ and a
        complementary subspace $\overline{\mathcal{B}}$ in Proposition
        \ref{proposition:description-commutant} yield equivalent
        deformations $\mu'$ of the commutant $(\mathcal{B},\circ)$.
    \item Furthermore, if the first Hochschild cohomology
        $\HHtype^1(\mathcal{A},\mathcal{D})= \{ 0 \}$ is trivial,
        equivalent deformations of $\mu_0$ and corresponding choices
        as above yield equivalent deformations of the bimodule
        structure \eqref{eq:notation-bimodule} in the sense of
        Definition \ref{definition:deformation-bimodules}.
    \item If the first two Hochschild cohomologies are trivial, this
        means if $\HHtype^1(\mathcal{A},\mathcal{D})= \{ 0 \}$ and
        $\HHtype^2(\mathcal{A}, \mathcal{D}) = \{0\}$, there is a map
        \begin{equation}
            \label{eq:DefTheorien}
            \Def_{\mathrm{type}}(\mathcal{A}) \longrightarrow
            \Def(\HHtype^0(\mathcal{A}, \mathcal{D})),
        \end{equation}
        where $\Def_{\mathrm{type}}$ denotes the set of equivalence
        classes of associative deformations of a particular type.
    \end{enumerate}
\end{proposition}

\begin{proof}
    For the first part let $\rho$ and $\tilde{\rho}$ be two
    deformations of $\rho_0$. Due to the assumption
    $\HHtype^1(\mathcal{A},\mathcal{D})= \{ 0 \}$ and Corollary
    \ref{corollary:equivalence-module-deformations} these are
    equivalent in the sense of Definition
    \ref{definition:deformation-right-module}. So there exists an
    equivalence transformation $T= \id_{\mathcal{E}}+
    \sum_{r=1}^{\infty} \lambda^r T_r$ with $T\circ \rho(a) =
    \tilde{\rho}(a) \circ T$ for all $a\in
    \mathcal{A}[[\lambda]]$. Obviously, the invertible map $T$ gives
    rise to an isomorphism
    \begin{equation}
        \label{eq:isomorphism-commutants}
        \Conj_T: \mathcal{K}_\rho \longrightarrow
        \mathcal{K}_{\tilde{\rho}},\quad \quad \Conj_T A= T\circ A\circ
        T^{-1} \quad \textrm{for } A\in \mathcal{K}_{\rho}, 
    \end{equation}
    where $\mathcal{K}_\rho$ and $\mathcal{K}_{\tilde{\rho}}$ denote
    the corresponding commutants. For each of
    these deformations we now consider a choice of a complementary
    subspace $\overline{\mathcal{B}}$ and get the corresponding maps
    $\rho'$ and $\tilde{\rho}'$ as in Proposition
    \ref{proposition:description-commutant}. They are invertible and
    so it is possible to define the map
    \begin{equation}
        \label{eq:equivalence-transformation-commutant}
        S'= \tilde{\rho}'^{-1} \circ \Conj_T \circ \rho'=
        \id_{\mathcal{B}} + \sum_{r=1}^\infty \lambda^r S'_r:
        \mathcal{B}[[\lambda]] \longrightarrow \mathcal{B}[[\lambda]] 
    \end{equation}
    which obviously is of the stated form with $S'_r\in
    \HC^1(\mathcal{B},\mathcal{B})$. Then this $S'$ is an equivalence
    transformation between the two deformations $\mu'$ and
    $\tilde{\mu}'$ induced by $\rho'$ and $\tilde{\rho}'$ since by
    definition
    \begin{eqnarray}
        \label{eq:equivalence-commutant-check}
        (S'\circ \mu')(A,B)&=& (\tilde{\rho}'^{-1} \circ \Conj_T \circ 
        \rho'\circ \rho'^{-1}) (\rho'(A)\circ \rho'(B))\nonumber \\
        &=& \tilde{\rho}'^{-1}(T\circ \rho'(A)\circ  T^{-1} \circ T
        \circ \rho'(B) \circ T^{-1}) \nonumber \\
        &=& \tilde{\rho}'^{-1}(\tilde{\rho}'(S'A) \circ
        \tilde{\rho}'(S'B))\\
        &=& (\tilde{\mu}'\circ (S'\otimes S'))(A,B). \nonumber
    \end{eqnarray}

    The statement of the second part can be easily traced back to the
    first part. So let $S$ be an equivalence transformation between
    the two deformations $\mu$ and $\tilde{\mu}$ as in Definition
    \ref{definition:deformation-bimodules}, this means $S\circ \mu=
    \tilde{\mu}\circ (S\otimes S)$. For two corresponding
    deformations $\rho$ and $\tilde{\rho}$ of $\rho_0$ one now defines
    \begin{equation}
        \label{eq:helping-right-module}
        \hat{\rho}=\rho \circ S^{-1}
    \end{equation}
    which is a deformation of $\rho_0$ with respect to $\tilde{\mu}$
    because of $\hat{\rho}(\tilde{\mu}(a,b))= (\rho\circ
    S^{-1})(\tilde{\mu}(a,b))= \rho(\mu(S^{-1}a,S^{-1}b))=
    \rho(S^{-1}b) \circ \rho(S^{-1}a)= \hat{\rho}(b)\circ
    \hat{\rho}(a)$ for all $a,b\in \mathcal{A}[[\lambda]]$. Thus
    $\hat{\rho}$ is equivalent to $\tilde{\rho}$ via some $T$ as
    before, this means $T\circ \hat{\rho}(a)= \tilde{\rho}(a)\circ T$.
    Then, the maps $S$, $T$ and the therewith constructed map $S'$ as
    in \eqref{eq:equivalence-transformation-commutant} are the
    equivalence transformations between $(\mu',\rho',\rho,\mu)$ and
    $(\tilde{\mu}', \tilde{\rho}',\tilde{\rho},\tilde{\mu})$. The
    computation \eqref{eq:equivalence-commutant-check} is exactly the
    same and the new $T$ immediately leads to $T\circ \rho(a)=T\circ
    \hat{\rho}(Sa)= \tilde{\rho}(Sa) \circ T$ and $T\circ \rho'(A)=
    T\circ \rho'(A)\circ T^{-1}\circ T = (\Conj_T\circ \rho')(A) \circ
    T= \tilde{\rho}'(S'(A)) \circ T$.

    The third part is a direct consequence of the second one. Since in
    addition the second Hochschild cohomology is trivial the
    deformations $\rho$ always exist and all the above structures can
    be constructed which finally yields the stated map due to the
    given uniqueness up to equivalence.
\end{proof}

\section{Invariant algebra and module structures}
\label{sec:G-invariant-types}

In this section we investigate the particular situation where the
algebra and module structures of a particular type are additionally
invariant under the action of a group $G$. Although it is not used in
this work it is remarkable that the discussion of $G$-invariant
deformations of algebra structures, especially $G$-invariant star
products, is still an area of vivid research. The basic notion was
given in \cite{bertelson.bieliavsky.gutt:1998a}. Moreover, if $G$ is a
Lie group it is also interesting to discuss the infinitesimal version
of $G$-invariance, this means the invariance under the induced
representation of the Lie algebra $\lie{g}$ of $G$. For more details
and the physical applications of invariant star products, in
particular in the framework of phase space reduction and quantum
momentum mappings, the reader is referred to the references
\cite{bordemann.herbig.waldmann:2000a, mueller-bahns.neumaier:2004a,
  mueller-bahns.neumaier:2004b}.

For our purpose we concentrate on $G$-invariant module structures and
if necessary the $G$-invariant deformation of the underlying algebra
is assumed to be given. This general case is discussed but later we
will be in the easier situation where the action on the algebra is
trivial. In the following, it will be pointed out in which cases the
$G$-invariant algebra and module structures define a new particular
type. As it will be seen, the crucial condition for that is the
invariance of the involved vector spaces. Besides this we can make
interesting additional observations concerning the commutant of
invariant deformations. This will play an important role in the
applications.

First of all we introduce the used notation and recall some well-known
definitions. Given an arbitrary \emph{action} of a group $G$ on a set
$V$, this means a group homomorphism $G\longrightarrow \Aut(V)$ into
the automorphisms of $V$, we simply write
\begin{equation}
    \label{eq:notation-action-of-element}
    g\acts v 
\end{equation}
for the action of a group element $g\in G$ on an element $v\in V$. If
$V$ is a vector space over some field $\mathbb{K}$ the action is
called a representation if all group elements act by
$\mathbb{K}$-linear maps. The space of all $G$-invariant elements is
denoted by
\begin{equation}
    \label{eq:set-G-inv-elements}
    V^G=\{ v\in V\:|\: g\acts v=v \: \textrm{for all} \: g\in G\}. 
\end{equation}
The \emph{orbit} of an element $v\in V$ is defined as the set
\begin{equation}
    \label{eq:orbit}
    G\acts v= \{g\acts v\:|\: g\in G\}
\end{equation}
Analogously, one defines the sets $g\acts W=\{g\acts w\:|\: w\in W\}$
for linear subspaces $W\subseteq V$. Such a $W\subseteq V$ is called
invariant if $G\acts W= \{g\acts w\:|\: g\in G, w\in W\} \subseteq
W$. Having representations on $\mathbb{K}$-vector spaces $V_1,\dots,
V_k$, $k\in \mathbb{N}$, and $W$ one always has a corresponding
representation on the $\mathbb{K}$-vector space
$\Hom_{\mathbb{K}}(V_1\times \dots \times V_k, W)$ of multilinear
maps. Explicitly, this representation is defined by
\begin{equation}
    \label{eq:representation-multi-linear-maps}
    (g\acts \phi) (v_1,\dots, v_k)= g\acts (\phi(g^{-1}\acts v_1,\dots,
    g^{-1}\acts v_k))
\end{equation}
for all $g\in G$ with inverse $g^{-1}$, $v_i\in V_i$ and $\phi\in
\Hom_{\mathbb{K}}(V_1\times \dots \times V_k, W)$. In the literature
the $G$-invariant elements $g\acts \phi=\phi$ of this representation
are often called \emph{$G$-equivariant} maps since they mediate the
actions,
\begin{equation}
    \label{eq:equivariant-multi-linear-maps}
    g\acts \phi(v_1,\dots,v_k)= \phi(g\acts v_1,\dots, g\acts v_k).
\end{equation}
By the above definition the composition of maps is always
$G$-invariant, this means
\begin{equation}
    \label{eq:action-composition-of-maps}
    g\acts (\phi\circ \psi) = (g\acts \phi) \circ (g\acts \psi)
\end{equation}
for all accordingly given maps $\phi$ and $\psi$ between sets with a
$G$-action. The following definition is obvious.

\begin{definition}[$G$-invariant algebra and module structures]
    \label{definition:G-inv-algebra-module-structures}
    An associative $\mathbb{K}$-algebra $(\mathcal{A},\mu)$ is said to
    be \emph{$G$-invariant} with respect to a representation of a
    group $G$ on the vector space $\mathcal{A}$ if the group acts by
    algebra automorphisms, this means if
    \begin{equation}
        \label{eq:algebra-G-inv}
        g\acts \mu(a, b)= \mu(g\acts a, g\acts b) \quad \textrm{for
          all } \quad a,b\in \mathcal{A}, g\in G. 
    \end{equation}
    If this is the case a right $(\mathcal{A},\mu)$-module structure
    $\rho$ of a $\mathbb{K}$-vector space $\mathcal{E}$ is said to be
    \emph{$G$-invariant} with respect to a further representation of
    $G$ on $\mathcal{E}$ if 
    \begin{equation}
        \label{eq:module-G-inv}
        g\acts \rho(e,a)= \rho(g\acts e, g\acts a) \quad \textrm{for
          all } e\in \mathcal{E}, a\in \mathcal{A}, g\in G.
    \end{equation}
\end{definition}

The Hochschild complexes of algebras $\mathcal{A}$ and modules
$\mathcal{E}$ as defined in the
Sections~\ref{subsec:Hochschild-complex-algebra} and
\ref{subsec:Hochschild-complex-module} basically consist of maps
between those vector spaces. Thus it is clear that representations of
a group $G$ on $\mathcal{A}$ and $\mathcal{E}$ induce representations
on the vector spaces of these complexes. Then, the above definition
contains nothing but the $G$-invariance $g\acts \mu=\mu$ and $g\acts
\rho=\rho$ of the corresponding cochains.

\begin{definition}[Algebra and module structures of a particular
    $G$-invariant type]
    \label{definition:G-invariant-type}
    If $\mu$ and $\rho$ as in
    Definition~\ref{definition:G-inv-algebra-module-structures} are
    $G$-invariant structures of a particular type as in
    Section~\ref{sec:algebras-modules-type} they are said to be of a
    corresponding \emph{$G$-invariant type} if all vector spaces
    occurring in the complexes
    $\HCtype^\bullet(\mathcal{A},\mathcal{A})$ and
    $\HCtype^\bullet(\mathcal{A},\mathcal{D})$ are $G$-invariant, this
    means if
    \begin{eqnarray}
        \label{eq:G-compatible-types}
        G\acts \HCtype^k(\mathcal{A},\mathcal{A}) &\subseteq&
        \HCtype^k(\mathcal{A},\mathcal{A}), \nonumber \\
        G\acts \HCtype^k(\mathcal{A},\mathcal{D}) &\subseteq&
        \HCtype^k(\mathcal{A},\mathcal{D}), \quad \textrm{and} \\
        G\acts \mathcal{D}&\subseteq& \mathcal{D}. \nonumber
    \end{eqnarray}
\end{definition}

Since all operations on the considered complexes are defined by the
composition of maps and the given $G$-invariant algebra and module
structures, they are $G$-invariant as well. This means that the
insertions $\circi$ after the $i$-th position, the cup product $\cup$
and the Hochschild differential $\delta$ have the properties
\begin{eqnarray}
    \label{eq:insertions-G-invariant}
    g\acts (\phi \circi \psi)&=&(g\acts \phi) \circi (g\acts \psi),\\ 
    \label{eq:cup-product-G-invariant}
    g\acts (\phi \cup \psi)&=& (g\acts \phi) \cup (g\acts \psi),\\
    \label{eq:Hochschild-delta-G-invariant}
    g\acts (\delta \phi) &=& \delta (g\acts \phi),
\end{eqnarray}
for all possible cochains $\phi,\psi\in
\HCtype^\bullet(\mathcal{A},\mathcal{A})$ or
$\HCtype^\bullet(\mathcal{A},\mathcal{D})$ and $g\in G$. Further it
is clear that $\mathcal{D}$ is a $G$-invariant
$(\mathcal{A},\mathcal{A})$-bimodule, this means
\begin{equation}
    \label{eq:Diffop-G-invariant-bimodule}
    g\acts (a\cdot D\cdot b)= (g\acts a)\cdot (g\acts D)\cdot
    (g\acts b)
\end{equation}
for all $g\in G$, $a,b \in \mathcal{A}$, $D\in \mathcal{D}$.

These consequences justify
Definition~\ref{definition:G-invariant-type} since structures of
$G$-invariant types in fact define new particular types in the sense
of the Definitions~\ref{definition:algebra-type} and
\ref{definition:module-type}. One simply has to consider the complexes
$\HCtype^\bullet(\mathcal{A},\mathcal{A})^G$ and
$\HCtype^\bullet(\mathcal{A},\mathcal{D})^G$ of $G$-invariant
cochains. For degree zero one of course has
$\HC^0(\mathcal{A},\mathcal{A})^G=\mathcal{A}$. Note, that by the
definition \eqref{eq:representation-multi-linear-maps} all
$G$-invariant $k$-cochains $\phi\in
\HCtype^k(\mathcal{A},\mathcal{D})$ indeed obey the equation
\begin{equation}
    \label{eq:G-invariant-cochains}
    g\acts (\phi(a_1,\dots, a_k)(e))= \phi(g\acts a_1,\dots, g\acts
    a_k) ( g\acts e).
\end{equation}

\begin{remark}[Representations and cohomology]
    \label{remark:representations-and-cohomology}
    For $G$-invariant types
    Equation~\eqref{eq:Hochschild-delta-G-invariant} implies that the
    cocycles and coboundaries of the initial Hochschild complex
    $\HCtype^\bullet(\mathcal{A},\mathcal{D})$ are invariant vector
    subspaces. Thus one has a representation on the cohomologies as
    well. But note that the cohomology of the $G$-invariant type in
    general is not the same as the $G$-invariant cohomology of the
    initial type. In formulas this means that for all $k\ge 1$
    \begin{equation}
        \label{eq:cohomology-and-G-invariance}
        H^k(\HCtype^\bullet(\mathcal{A},\mathcal{D})^G) \neq
        \HHtype^k(\mathcal{A},\mathcal{D})^G. 
    \end{equation}
    The only exception is given by the zeroth cohomology. There one
    has
    \begin{equation}
        \label{eq:zeroth-cohomology-and-G-invariance}
        H^0(\HCtype^\bullet(\mathcal{A},\mathcal{D})^G)
        =\HHtype^0(\mathcal{A},\mathcal{D})^G
    \end{equation}
    since both sides are given by the elements $\phi\in
    \HCtype^0(\mathcal{A},\mathcal{D})$ with $g\acts \phi=\phi$ and
    $\delta\phi=0$.
\end{remark}

\begin{remark}[Modules with a trivial representation on the algebra]
    \label{remark:action-only-on-module}
    If the representation of $G$ on the algebra $\mathcal{A}$ is
    trivial, meaning that all group elements act by the identity map,
    the Hochschild complex of the $G$-invariant right module structure
    $\rho$ is given by
    \begin{equation}
        \label{eq:G-inv-complex-action-trivial}
        \HCtype^\bullet(\mathcal{A},\mathcal{D})^G=
        \HCtype^\bullet(\mathcal{A},\mathcal{D}^G). 
    \end{equation}
\end{remark}

With the above statements all results from the
Sections~\ref{sec:deformation-algebras}, \ref{sec:deformation-modules}
and \ref{sec:commutant-module-structures} can be applied in order to
investigate deformations of right modules that inherit a given
$G$-invariance. The crucial cohomology then is of course
$H^\bullet(\HCtype(\mathcal{A},\mathcal{D})^G)$. If the first one is
trivial, $H^1(\HCtype(\mathcal{A},\mathcal{D})^G)=\{0\}$,
Proposition~\ref{proposition:description-commutant} states that with
the notation $\mathcal{B}= \HHtype^0(\mathcal{A},\mathcal{D})$ any
decomposition
\begin{equation}
    \label{eq:decomposition-G-invariance-trivial-case}
    \HCtype^1(\mathcal{A},\mathcal{D})^G= \mathcal{B}^G\oplus
    \cc{\mathcal{B}^G} 
\end{equation}
leads to a bijection $\rho':\mathcal{B}^G[[\lambda]]\longrightarrow
\mathcal{K}^G$ between the commutants within
$\HCtype^0(\mathcal{A},\mathcal{D})^G[[\lambda]]$. Thus one has the
property
\begin{equation}
    \label{eq:action-on-commutant-trivial}
    g\acts \rho'(A)(e)=\rho'(A)(g\acts e).
\end{equation}
In most applications, however, one is rather interested in the
commutants within $\HCtype^0(\mathcal{A},\mathcal{D})[[\lambda]]$, or
even $\mathcal{D}[[\lambda]]$, and wants to find a bijection
$\rho':\mathcal{B}[[\lambda]]\longrightarrow \mathcal{K}$ such that
the induced deformed bimodule is $G$-invariant like the undeformed one
in \eqref{eq:notation-bimodule}. This is the case if the complement of
$\mathcal{B}$ is an invariant vector space.

\begin{proposition}[The commutant of a $G$-invariant module structure] 
    \label{proposition:commutant-G-invariance}
    Let $\rho$ be a $G$-invariant deformation of a right module
    structure $\rho_0$ of a particular $G$-invariant type and let
    $\HHtype^1(\mathcal{A},\mathcal{D})=\{0\}$ as in
    Proposition~\ref{proposition:description-commutant}. If the
    complementary subspace $\cc{\mathcal{B}}$ is $G$-invariant, this
    means
    \begin{equation}
        \label{eq:complement-invariant}
        G\acts \cc{\mathcal{B}} \subseteq \cc{\mathcal{B}},
    \end{equation}
    the induced bijection $\rho': \mathcal{B}[[\lambda]] \cong
    \mathcal{K}$ and the deformed algebra structure
    $(\mathcal{B}[[\lambda]],\mu')$ as in
    Corollary~\ref{corollary:deformation-commutant} are
    $G$-invariant. Explicitly, this means that
    \begin{eqnarray}
        \label{eq:action-commutant-module}
        g\acts (\rho'(A)(e))&=&\rho'(g \acts A)(g \acts e),\\
        \label{eq:action-commutant-algebra}
        g\acts (\mu'(A,B))&=& \mu'(g \acts A, g \acts B)
    \end{eqnarray}
    for all $A, B\in \mathcal{B}[[\lambda]]$, $e\in
    \mathcal{E}[[\lambda]]$ and $g\in G$.

    Moreover, if additionally the first $G$-invariant cohomology is
    trivial,
    \begin{equation}
        \label{eq:trivial-first-G-inv-cohomology}
        H^1(\HCtype(\mathcal{A},\mathcal{D})^G)=\{0\},
    \end{equation}
    the so derived $G$-invariant deformations $\rho'$ and $\mu'$ and
    the induced $G$-invariant bimodule are unique up to $G$-invariant
    equivalence in analogy to
    Proposition~\ref{proposition:classes-of-deformations}.
\end{proposition}

\begin{proof}
    First one notes that $\mathcal{B}= \HHtype^0(\mathcal{A},
    \mathcal{D})$ always is an invariant subspace. Then the assumption
    \eqref{eq:complement-invariant} implies the $G$-invariance of the
    decomposition $\HCtype^0(\mathcal{A},\mathcal{D})
    =\mathcal{B}\oplus \cc{\mathcal{B}}$. The recursive construction
    of $\rho'= \id+ \sum_{r=1}^\infty\lambda^r \rho'_r$ from the proof
    of Proposition~\ref{proposition:description-commutant} shows that
    each $\rho'_r$ and thus $\rho'$ is $G$-invariant. With the given
    assumptions the results of this section imply that for all
    $\mathcal{A}\in \mathcal{B}$ and $r\ge 0$ the elements $g \acts
    \rho'_{r+1}(A)$ satisfy the same defining
    Equation~\eqref{eq:condition-correction-commutant} as
    $\rho'_{r+1}(g\acts A)$. The $G$-invariance of the complementary
    subspace $\cc{\mathcal{B}}$ and the condition
    \eqref{eq:condition-correction-in-complement} finally imply that
    these elements are equal by uniqueness. Thus one has the stated
    $G$-invariance of $\rho'$. By definition, this implies
    \eqref{eq:action-commutant-module} and
    \eqref{eq:action-commutant-algebra}.

    The assertion concerning $G$-invariant equivalence is clear with
    Proposition~\ref{proposition:classes-of-deformations} since
    \eqref{eq:trivial-first-G-inv-cohomology} assures the
    $G$-invariance of all maps occurring in the proof, in particular of
    the relevant equivalence transformations.
\end{proof}

% It is remarkable that the $G$-invariant Hochschild complexes in
% general do not have the property \eqref{eq:convenient-type} of
% Remark~\ref{remark:special-case-zero-Hochschild}. Even if one
% considers a particular type with
% $\HCtype^0(\mathcal{A},\mathcal{D})=\mathcal{D}$ one gets
% \begin{equation}
%     \label{eq:convenient-G-invariant}
%     \HCtype^0(\mathcal{A},\mathcal{D})^G=\mathcal{D}^G.
% \end{equation}

\section{Projective modules}
\label{sec:projective-modules}

As a first example we consider projective right modules. The
well-known definition of a projective module can be found in
Appendix~\ref{sec:resolutions-derived-functors}. It turns out that
they are rigid and that deformations thereof always exist since the
crucial Hochschild cohomologies vanish. In fact it is possible to
compute all cohomologies using an explicit homotopy.

So let $\mathcal{A}$ be an associative $\mathbb{K}$-algebra and
$\mathcal{E}$ be a projective right $\mathcal{A}$-module where the
module multiplication is written as $e\cdot a$ for $e\in \mathcal{E}$
and $a\in \mathcal{A}$. By one of the equivalent characterizations of
projective modules stated in Remark~\ref{remark:projective-module} we
can choose a dual basis, this means families
\begin{equation}
    \label{eq:dual-basis}
    \{e_i\}_{i\in I}\subseteq \mathcal{E} \quad \textrm{and} \quad
    \{\epsilon^i\}_{i\in I} \subseteq 
    \Hom_{\mathcal{A}}(\mathcal{E},\mathcal{A}) 
\end{equation}
for some index set $I$ where $\Hom_{\mathcal{A}}
(\mathcal{E},\mathcal{A})$ denotes the set of all right
$\mathcal{A}$-linear homomorphisms. Then, by definition, all $e\in
\mathcal{E}$ can be written as
\begin{equation}
    \label{eq:dual-basis-property}
    e = \sum_{i\in I} e_i \cdot \epsilon^i(e)
\end{equation}
where, respectively, only finitely many $\epsilon^i(e)$ are different
from zero even if there should be infinitely many $e_i$ and
$\epsilon^i$. Now consider the maps
\begin{equation}
    \label{eq:homotopy-maps-projective}
    h^k: \HC^k(\mathcal{A}, \End_{\mathbb{K}}
    (\mathcal{E},\mathcal{E})) 
    \longrightarrow
    \HC^{k-1}(\mathcal{A},  \End_{\mathbb{K}}
    (\mathcal{E},\mathcal{E}))  
\end{equation}
defined by $h^0=0$ and
\begin{equation}
    \label{eq:homotopy-maps-projective-explicit}
    (h^k\phi)(a_1, \dots, a_{k-1})e = \sum_{i\in I} \phi(
    \epsilon^i(e), a_1, \dots, a_{k-1})e_i
\end{equation}
for $k\ge 1$. Using the summation convention for $\sum_{i\in I}$ and
with \eqref{eq:bimodule-of-endomorphisms} we compute for $k\ge 1$
\begin{eqnarray*}
    \label{eq:delta-h-projective}
    \lefteqn{(\delta^{k-1} h^k \phi) (a_1, \dots , a_k)e}\\
    &=& (h^k\phi) (a_2, \dots, a_k) (e\cdot a_1)\\
    && + \sum_{s=1}^{k-1} (-1)^s (h^k\phi) (a_1, \dots, a_s a_{s+1},
    \dots, a_k) e + (-1)^k \left((h^k\phi) (a_1,\dots , a_{k-1})
        e \right) \cdot a_k\\
    &=& \phi(\epsilon^i(e\cdot a_1), a_2, \dots, a_k) e_i\\
    && + \sum_{s=1}^{k-1} (-1)^s \phi(\epsilon^i(e), a_1, \dots, a_s
    a_{s+1}, \dots, a_k) e_i + (-1)^k \left(\phi(\epsilon^i(e), a_1, \dots,
        a_{k-1}) e_i\right) \cdot a_k
\end{eqnarray*}
and analogously
\begin{eqnarray*}
    \label{eq:h-delta-projective}
    \lefteqn{(h^{k+1}\delta^k \phi) (a_1, \dots, a_k) e }\\
    &=& \left((\delta^k\phi) (\epsilon^i(e), a_1, \dots,
        a_k)\right)e_i\\ 
    &=& \phi(a_1, \dots, a_k)(e_i \cdot \epsilon^i(e)) -
    \phi(\epsilon^i(e) a_1, a_2, \dots, a_k) e_i\\
    && -\sum_{s=1}^{k-1} (-1)^s \phi( \epsilon^i(e), a_1, \dots, a_s
    a_{s+1}, \dots, a_k)e_i - (-1)^k \left(\phi(\epsilon^i(e), a_1, \dots,
        a_{k-1})e_i \right) \cdot a_k.
\end{eqnarray*}
With the properties of the dual basis this yields $\left(\delta^{k-1}
    h^k \phi + h^{k+1} \delta^k\phi \right) (a_1,\dots, a_k) e=
\phi(a_1,\dots, a_k)e$. Thus one has
\begin{equation}
    \label{eq:homotopy-projective}
    \delta^{k-1} \circ h^k + h^{k+1} \circ \delta^k= \id_{\HC^k},
\end{equation}
which shows that $\id_{\HC^k}$ is homotopic to zero for all $k\ge 1$,
confer Definition~\ref{definition:homotopy}. Consequently, the
Hochschild cohomology is trivial. For $h^1$ one further computes
\begin{eqnarray*}
    (h^1 \delta^0 D)e &=& (\delta^0 D) (\epsilon^i(e))e_i\\
    &=&  D(e_i\cdot \epsilon^i(e))- D(e_i)\cdot \epsilon^i(e)\\
    &=& D(e) - D(e_i) \cdot \epsilon^i(e).
\end{eqnarray*}
This shows that $\id- h^1 \circ \delta^0$ is a projection onto the
commutant $\ker \delta^0$ of the right module since $\epsilon^i$ is
right $\mathcal{A}$-linear. Altogether one gets the following
statement.

\begin{proposition}[The cohomology of a projective module]
    \label{proposition:cohomology-projective}
    Let $\mathcal{E}$ be a projective right module over an associative
    algebra $\mathcal{A}$. Then
    \begin{equation}
        \label{eq:cohomology-projective}
        \HH^k(\mathcal{A},\End_{\mathbb{K}}(\mathcal{E},\mathcal{E}))
        = \left\{
            \begin{array}{c}
                \End_{\mathcal{A}}(\mathcal{E},\mathcal{E})\\
                \{0\}
            \end{array}
            \textrm{for}
            \begin{array}{c}
                k=0\\
                k\ge 1.
            \end{array}\right.
    \end{equation}
    In addition, the choice of a dual basis yields an explicit
    homotopy and a projection onto the commutant.
\end{proposition}
The assertion of course includes the cases where $\mathcal{E}$ is a
free module or where $\mathcal{E}=\mathcal{A}$ is the associative
algebra itself with unit $1$. In the last case the maps $h^k$ can be
simply defined by $h^k(\phi)(a_1,\dots, a_{k-1})b= \phi(b, a_1,\dots,
a_{k-1}) 1$ as stated in
\cite[Remark~2.3]{bordemann.neumaier.waldmann.weiss:2007a:pre}.

\begin{remark}[Projective modules of particular types and deformed
    vector bundles]
    \begin{enumerate}
    \item In general it is not true that the map $h^k$ respects a
        given particular type of cochains. Nevertheless, the simple
        algebraic definition gives strong evidence that this is indeed
        the case in many typical situations. The type of differential
        module structures which will be introduced in the next chapter
        is such an example, confer
        Proposition~\ref{proposition:projective-differential-modules}.
    \item Proposition \ref{proposition:cohomology-projective} shows
        that projective modules are rigid and that deformations can be
        constructed order by order. As it will be explained in
        Section~\ref{sec:quantization-vector-bundles} an important
        application of this fact lies in the deformation theory of
        vector bundles over Poisson manifolds as investigated in
        \cite{bursztyn.waldmann:2000b, waldmann:2001b,
          waldmann:2002b}.
    \end{enumerate}
\end{remark}

\chapter{Differential algebra and module structures}  
\label{cha:differential-G-invariant-structures}

In this chapter we specify what shall be understood by algebra and
module structures acting by differential operators. To this end we use
the notion of multidifferential operators between modules with respect
to some associative and commutative algebra. The basic definition and
some well-known facts of these operators are presented in the first
section. After these preliminaries the second section contains a
discussion of the general framework which is necessary for the
definition of differential algebra and module structures. As it will
be shown these definitions determine a particular type of algebraic
structures as introduced in
Section~\ref{sec:algebras-modules-type}. It will further be shown that
one is always able to speak of $G$-invariant differential algebra and
module structures with respect to given representations of a group
$G$. Thus it is possible in the aspired applications to investigate
the deformations of such structures using the general results
established in Chapter~\ref{cha:deformation-algebras-modules}. Besides
presenting the basic framework we already make some general
observations concerning the Hochschild complexes induced by these
examples of particular types.
% This will be helpful for the aspired
% computation of the Hochschild cohomologies in
% Chapter~\ref{cha:deformation-on-bundles}.
Finally, it is shown that the induced cohomology of projective modules
of the differential type always vanishes since the homotopy from
Section~\ref{sec:projective-modules} still holds in the presented
differential framework.

\section{Algebraically defined multidifferential operators}
\label{sec:algebraic-multi-diffop}

The motivation for a general notion of differential operators
naturally arises from the simple and fundamental definition used
within the framework of ordinary analysis. There, a differential
operator $D$ of degree $l\in \mathbb{N}_0$ on the functions $f\in
C^\infty(\mathbb{R})$ is a map $D:f\longmapsto \sum_{i=0}^l c^i
\frac{\partial^i f} {\partial x^i}$ with functions $c^i\in
C^\infty(\mathbb{R})$. Considering $C^\infty(\mathbb{R})$ as a module
over itself the algebraic properties of such a map $D$ are the
guideline for the general algebraic definition. The most
characteristic property is derived from the Leibniz rule implying that
$f\longmapsto cD(f)-D(cf)$ for any $c\in C^\infty(\mathbb{R})$ is a
differential operator of degree $l-1$.

In the following summary of the well-known definitions we adapt the
notation used in \cite[App.~A]{waldmann:2007a} where one can find
further details. However, the basic definition of differential
operators was first given in the works of Grothendieck, confer
\cite[Def.~16.8.1]{grothendieck:1967}. The necessary structures for
the definition of multidifferential operators used in this work are an
associative and commutative $\mathbb{K}$-algebra $\mathcal{C}$ and
$\mathbb{K}$-vector spaces $\mathcal{E}_1,\dots, \mathcal{E}_N$, $N\in
\mathbb{N}$, and $\mathcal{F}$ with a compatible $\mathcal{C}$-module
structure.

For these structures one considers the $\mathbb{K}$-multilinear maps
$D: \mathcal{E}_1\times \dots \times \mathcal{E}_N\longrightarrow
\mathcal{F}$ which, as always, are identified with the
$\mathbb{K}$-linear maps $\Hom_{\mathbb{K}}(\mathcal{E}_1\otimes \dots
\otimes \mathcal{E}_N,\mathcal{F})$ due to the universal property of
the tensor product $\otimes=\otimes_{\mathbb{K}}$. Any of the given
$\mathcal{C}$-module structures induces corresponding ones on the
multilinear maps. One defines them as a left $\mathcal{C}$-module
structure
\begin{equation}
    \label{eq:left-module-multi-linear-maps}
    (c\cdot D)(e_1,\dots,e_N)= c(D(e_1,\dots, e_N)) 
\end{equation}
and right $\mathcal{C}$-module structures
\begin{equation}
    \label{eq:right-module-multi-linear-maps}
    (D \mathbin{\cdot^{(i)}} c) (e_1,\dots,e_N) = D(e_1,\dots, ce_i,
    \dots, e_N) 
\end{equation}
for $i=1,\dots, N$ which commute with each other and
\eqref{eq:left-module-multi-linear-maps}. Note that in general $D
\mathbin{\cdot^{(i)}} c \neq D \mathbin{\cdot^{(j)}} c$ for $i\neq j$.
Denoting by $\LL_c$ and $\LL_c^{(i)}$ the left multiplication with an
element $c\in \mathcal{C}$ in $\mathcal{F}$ and $\mathcal{E}_i$,
respectively, these commuting module structures are given by
\begin{equation}
    \label{eq:module-multi-linear-maps-short}
    c\cdot D= \LL_c \circ D \quad \textrm{and} \quad 
    D \mathbin{\cdot^{(i)}} c= D\circ \LL_c^{(i)} 
\end{equation}
where in the last case $\LL_c^{(i)}$ is the obvious extension to the
tensor product $\mathcal{E}_1\otimes \dots \otimes \mathcal{E}_N$.
Then it is possible to define the operations
\begin{equation}
    \label{eq:quasi-commutator}
    D\longmapsto \ad^{(i)}_{c} (D) = c\cdot D - D\mathbin{\cdot^{(i)}} c 
\end{equation}
for all $c\in \mathcal{C}$ and $i=1,\dots,N$. Note, that all
$\ad^{(i)}_{c}$ commute. According to the above stated motivation a
multilinear operator $D$ is called a differential operator if it
vanishes after applying a finite number of operations as in
\eqref{eq:quasi-commutator}.

For the description of the \emph{multiorders of differentiation} one
makes use of multiindices $L=(l_1,\dots,l_N),$ $K=(k_1,\dots, k_N)\in
\mathbb{Z}^N$. Their addition and subtraction is defined
componentwise. Further one sets $L\le K$ if $l_i\le k_i$ for all $i=
1,\dots ,N$. For $L\ge 0$, the absolute value of $L$ is given by
$|L|=l_1+ \dots + l_N$. Particular multiindices of absolute value one
are given by the canonical basis vectors $\mathsf{e}_i\in
\mathbb{Z}^N$ with $1$ at the $i$-th position and $0$ at the others.
Now a multilinear map $D$ is called a differential operator of order
$L=(l_1,\dots,l_N)$ if
\begin{equation}
    \label{eq:diffop-quasi-commutator}
    \ad^{(1)}_{c^{(1)}_1}\circ \dots \circ \ad^{(1)}_{c^{(1)}_{l_1}}\circ
    \dots \circ
    \ad^{(N)}_{c^{(N)}_1} \circ \dots \circ \ad^{(N)}_{c^{(N)}_{l_N}} (D)=0
\end{equation}
for all ${c^{(1)}_1}, \dots, {c^{(N)}_{l_N}}\in \mathcal{C}$. For
later applications it is convenient to reformulate this definition as
a recursive one. This allows to prove all later assertions by
induction over the absolute value $|L|$.

\begin{definition}[Multidifferential operator
    \protect{\cite[Def.~A.4.1]{waldmann:2007a}}] 
    \label{definition:multi-diffops}
    Let $\mathcal{C}$ be an associative and commutative
    $\mathbb{K}$-algebra and $\mathcal{E}_1,\dots, \mathcal{E}_N$ and
    $\mathcal{F}$ be $\mathbb{K}$-vector spaces with a compatible
    $\mathcal{C}$-module structure. The \emph{multidifferential
      operators $\Diffop^L_{\mathcal{C}} (\mathcal{E}_1,
      \dots,\mathcal{E}_N; \mathcal{F})$ with arguments in
      $\mathcal{E}_1,\dots, \mathcal{E}_N$ and values in $\mathcal{F}$
      of multiorder $L=(l_1,\dots, l_N)\in \mathbb{Z}^N$} are then
    recusively defined by
    \begin{equation}
        \label{eq:Diffop-zero}
        \Diffop^L_{\mathcal{C}} (\mathcal{E}_1,
        \dots,\mathcal{E}_N; \mathcal{F})= \{0\} \quad \textrm{if
          there is some } l_i< 0 
    \end{equation}
    and 
    \begin{equation}
        \label{eq:Diffop-definition}
        \begin{split}
            \Diffop^L_{\mathcal{C}} (\mathcal{E}_1, \dots,\mathcal{E}_N;
            \mathcal{F}) = & \big\{ D\in
                \Hom_{\mathbb{K}}(\mathcal{E}_1\otimes \dots \otimes
                \mathcal{E}_N,\mathcal{F}) \: \big| \\
            & \phantom{\big\{} 
                \LL_c\circ D - D\circ \LL_c^{(i)} \in
                \Diffop^{L-\mathsf{e}_i}_{\mathcal{C}} (\mathcal{E}_1,
            \dots,\mathcal{E}_N; \mathcal{F})\\
            &  \phantom{\big\{} \textrm{for all } c\in
            \mathcal{C}, i=1, \dots, N \big\}
        \end{split}
        \end{equation}
    for $L\ge 0$.
\end{definition}

A direct consequence of Definition~\ref{definition:multi-diffops} is
the following proposition.

\begin{proposition}[Filtration and module structures
    \protect{\cite[Prop.~A.4.3]{waldmann:2007a}}]
    \label{proposition:diffops-filtered-and-module}
    Let $\mathcal{E}_1,\dots, \mathcal{E}_N$ and $\mathcal{F}$ be
    $\mathcal{C}$-modules as in
    Definition~\ref{definition:multi-diffops}. Then,
    \begin{equation}
        \label{eq:diffop-filtration-explicitly}
        \Diffop^L_{\mathcal{C}} (\mathcal{E}_1,\dots,
        \mathcal{E}_N;\mathcal{F}) \subseteq
        \Diffop^K_{\mathcal{C}} (\mathcal{E}_1,\dots,
        \mathcal{E}_N;\mathcal{F}),
    \end{equation}
    if $L \le K$. This implies that 
    \begin{equation}
        \label{eq:diffop-filtration}
        \Diffop^\bullet_{\mathcal{C}} (\mathcal{E}_1,\dots,
        \mathcal{E}_N; \mathcal{F}) = \Union_{L\ge 0}
        \Diffop^L_{\mathcal{C}} (\mathcal{E}_1,\dots,
        \mathcal{E}_N;\mathcal{F}) \subseteq
        \Hom_{\mathbb{K}}(\mathcal{E}_1 \otimes \dots \otimes
        \mathcal{E}_N,\mathcal{F})  
    \end{equation}
    is a filtered subspace. Moreover, $\Diffop^L_{\mathcal{C}}
    (\mathcal{E}_1,\dots, \mathcal{E}_N;\mathcal{F})$ inherits the
    $\mathcal{C}$-module structures
    \eqref{eq:module-multi-linear-maps-short}.
\end{proposition}

As it is typical for algebraically defined multidifferential
operators, the proofs of the assertions are easy inductions over the
absolute value $|L|$. In \cite{waldmann:2007a} one finds all the
arguments in great detail. Note that it is crucial that $\mathcal{C}$
is commutative.  It is an important feature of multidifferential
operators that they can be composed.

\begin{proposition}[Composition of differential operators
    \protect{\cite[Prop.~A.4.4]{waldmann:2007a}}]
    \label{proposition:composition-diffop}
    Let $\mathcal{E}_1^{(1)}, \dots, \mathcal{E}_{N_1}^{(1)}, \dots,
    \mathcal{E}_1^{(M)}, \dots \mathcal{E}_{N_M}^{(M)}$ as well as
    $\mathcal{F}_1, \dots, \mathcal{F}_M$ and $\mathcal{G}$ be
    $\mathcal{C}$-modules. Then, for
    \begin{equation}
        \label{eq:diffops-for-composition}
        D_i \in \Diffop^{K_i}_{\mathcal{C}} (\mathcal{E}_1^{(i)},
        \dots, \mathcal{E}_{N_i}^{(i)}; \mathcal{F}_i),\: i=1,\dots, M
        \quad \textrm{and} \quad D\in
        \Diffop_{\mathcal{C}}^L(\mathcal{F}_1,\dots, \mathcal{F}_M;
        \mathcal{G}) 
     \end{equation}
     one has
     \begin{equation}
         \label{eq:composition-diffops}
         D\circ (D_1\otimes \dots \otimes D_M)\in
         \Diffop_{\mathcal{C}}^{(K_1+ \cc{l_1}, \dots, K_M+ \cc{l_M})}
         (\mathcal{E}_1^{(1)}, \dots, \mathcal{E}_{N_M}^{(M)};
         \mathcal{G}),
     \end{equation}
     where $\cc{l_i}=(l_i,\dots, l_i)\in \mathbb{Z}^{N_i}$ for all
     $i=1, \dots, M$.
\end{proposition}

The proof can again be found in \cite{waldmann:2007a} and is an
induction over $r= |K_1|+ \dots + |K_M|+ |L|$.

The algebraic definition of multidifferential operators yields
invariant subspaces if the relevant structures are compatible with the
representations of a group $G$. In detail, one finds the following
well-known assertion.

\begin{lemma}[Representations on multidifferential operators]
    \label{lemma:Diffops-and-representations}
    Let $\mathcal{C}$ a $G$-invariant, associative, and commutative
    algebra and let $\mathcal{E}_1,\ldots,\mathcal{E}_k, \mathcal{F}$
    be $G$-invariant $\mathcal{C}$-modules with respect corresponding
    representations of a group $G$. Then, the differential operators
    are invariant subspaces of $\Hom_{\mathbb{K}}(\mathcal{E}_1\otimes
    \dots \otimes \mathcal{E}_N,\mathcal{F})$ with respect to the
    induced representation
    \eqref{eq:representation-multi-linear-maps}. For all $L\in
    \mathbb{N}_0^k$ one has
    \begin{equation}
        \label{eq:Diffops-invariant-subspace}
        G\acts \Diffop^L_{\mathcal{C}}(\mathcal{E}_1,\dots,
        \mathcal{E}_N;\mathcal{F}) \subseteq
        \Diffop^L_{\mathcal{C}}(\mathcal{E}_1,\dots,
        \mathcal{E}_N;\mathcal{F}).
    \end{equation}
\end{lemma}

\begin{proof}
    With the given definition of multidifferential operators the
    proof is an easy induction over $|L|$ using the compatibility of
    the structures. In the usual notation the crucial equation is
    $\LL_{g\acts c} \circ (g\acts D)- (g\acts D) \circ
    \LL^{(i)}_{g\acts c}= g\acts (\LL_c\circ D - D\circ \LL_c)$.
\end{proof}

In the applications one sometimes uses the fact that the notion of 
multidifferential operators is functorial.

\begin{lemma}[Functoriality of $\Diffop$]
    \label{lemma:diffop-functorial}
    Let $N\in \mathbb{N}_0$, $\gamma: \mathcal{C}\longrightarrow
    \mathcal{C}'$ be an isomorphism of associative and commutative
    $\mathbb{K}$-algebras and $\epsilon_i:\mathcal{E}_i
    \longrightarrow \mathcal{E}_i'$, $i=1,\dots, N$, $\zeta:
    \mathcal{F}\longrightarrow \mathcal{F}'$ be isomorphisms along
    $\gamma$ of $\mathcal{C}$- and $\mathcal{C}'$-modules.  Then, for
    all $L\in \mathbb{N}_0^N$ there is an isomorphism
    \begin{equation}
        \label{eq:isomorphism-diffop}
        \Omega:
        \Diffop^L_{\mathcal{C}} (\mathcal{E}_1,\dots, \mathcal{E}_N;
        \mathcal{F}) 
        \longrightarrow 
        \Diffop^L_{\mathcal{C}'} (\mathcal{E}_1',\dots, \mathcal{E}_N';
        \mathcal{F}') 
    \end{equation}
    of $\mathcal{C}$- and $\mathcal{C}'$-modules along $\gamma$
    defined by
    \begin{equation}
        \label{eq:isomorphism-diffop-explicit}
        (\Omega D)(e_1',\dots, e_N')=
        \zeta(D(\epsilon_1^{-1}(e_1'), \dots, \epsilon_N^{-1}(e_N')))
    \end{equation}
    for all $D\in \Diffop^L_{\mathcal{C}}(\mathcal{E}_1,\dots,
    \mathcal{E}_N; \mathcal{F})$ and $e_i'\in \mathcal{E}_i'$,
    $i=1,\dots, N$. 

    Thus, $\Diffop^L$, this means the assignment of the space
    $\Diffop^L_{\mathcal{C}}(\mathcal{E}_1,\dots, \mathcal{E}_N;
    \mathcal{F})$ of differential operators to a tuple
    $(\mathcal{C},\mathcal{E}_1,\dots , \mathcal{E}_N, \mathcal{F})$
    consisting of an algebra and modules, can be seen as a functor
    between correspondingly defined categories.

    The functoriality is of course also given in the setting of given
    representations of a group $G$. If the initial structures and maps
    are $G$-invariant the same is true for the isomorphism $\Omega$.
\end{lemma}

\begin{proof}
    That $\Omega$ takes values in the stated multidifferential
    operators is a simple induction over $r=|L|$ making use of the
    equation $\LL_{c'}\circ (\Omega D) - (\Omega D)\circ \LL_{c'}^{(i)}=
    \Omega (\LL_{\gamma^{-1}(c')} \circ D- D\circ
    \LL_{\gamma^{-1}(c')}^{(i)})$. Then all further assertions are
    obvious.
\end{proof}

If one considers differential operators
$\Diffop^l_{\mathcal{C}}(\mathcal{E},\mathcal{E})$ of a free
$\mathcal{C}$-module $\mathcal{E}$ it is easy to see that they can be
identified with matrices of differential operators in
$\Diffop^l_{\mathcal{C}}(\mathcal{C},\mathcal{C})$.

\begin{lemma}[Differential operators of free modules]
    \label{lemma:Diffop-free-modules}
    Let $\mathcal{C}$ be a commutative, associative
    $\mathbb{K}$-algebra and $\mathcal{E}$ be a finitely generated
    free left $\mathcal{C}$-module with a module basis
    $\mathfrak{B}=\{e_1,\ldots, e_N\} \subseteq \mathcal{E}$. Then,
    for every differential operator $D\in \Diffop^l_\mathcal{C}
    (\mathcal{E}, \mathcal{E})$ with $l\in \mathbb{N}_0$ there exist
    $N^2$ uniquely defined differential operators $D^i_j \in
    \Diffop^l_\mathcal{C} (\mathcal{C},\mathcal{C})$,
    $i,j=1,\ldots,N$, with
    \begin{equation}
        \label{eq:Diffops-free-module}
        D(c^j e_j)= D^i_j(c^j) e_i
    \end{equation}
    for all $c^j \in \mathcal{C}$. The map
    \begin{equation}
        \label{eq:Diffop-matrix}
        \Phi_\mathfrak{B}: D \longmapsto (D^i_j)_{i,j=1,\ldots,N}
    \end{equation}
    is an isomorphism
    \begin{equation}
        \label{eq:Iso-Diffop-matrix}
        \Phi_\mathfrak{B}: \Diffop^l_\mathcal{C} (\mathcal{E},\mathcal{E}) 
        \cong 
        \Mat_{N\times N} (\Diffop^l_\mathcal{C} (\mathcal{C},\mathcal{C}))
    \end{equation}
    of $\mathcal{C}$-bimodules which depends on the choice of the
    module basis $\mathfrak{B}$ for $\mathcal{E}$. The module
    structures $c \cdot D= \LL_c \circ D$ and $D \cdot c= D\circ
    \LL_c$ for the $N\times N$-matrices on the right hand side of
    \eqref{eq:Iso-Diffop-matrix} are defined componentwise and
    obviously induced by the multiplication in
    $\mathcal{C}$.

    If $\mathcal{C}$ and $\mathcal{E}$ are $G$-invariant structures
    and if the module basis $\mathfrak{B}=\{e_1,\ldots, e_N\}$
    consists of $G$-invariant elements
    \begin{equation}
        \label{eq:G-inv-module-basis}
        g\acts e_i=e_i \quad \textrm{for all} \quad i\in \{1,\dots, N\},
    \end{equation}
    the map $\Phi_{\mathfrak{B}}$ is $G$-invariant with the
    componentwise representation of $G$ on the matrices.
\end{lemma}

\begin{proof}
    For arbitrary $j \in \{1,\ldots, r\}$ and $c\in \mathcal{C}$ one
    surely has $D(c e_j)=D^i_j(c) e_i$ with uniquely defined elements
    $D^i_j(c) \in \mathcal{C}$. The also uniquely defined maps $D^i_j:
    \mathcal{C}\longrightarrow \mathcal{C}$ are $\mathbb{K}$-linear.
    Thus they satisfy \eqref{eq:Diffops-free-module}. The assertion
    $D^i_j \in \Diffop^l_\mathcal{C}(\mathcal{C},\mathcal{C})$ can be
    shown by induction over $l$. For $l=0$ one computes $D(c e_j)= c
    D(e_j) = c (d^i_j e_i)$ with $d^i_j \in \mathcal{C}$. So $D^i_j$
    is a multiplication in $\mathcal{C}$ and thus
    $D^i_j=\LL^{\mathcal{C}}_{d^i_j} \in
    \Diffop^0_\mathcal{C}(\mathcal{C},\mathcal{C})$. For $l>0$
    evaluation of $(\LL_{c'}\circ D- D\circ \LL_{c'}) (c e_j)$ with
    $c,c'\in \mathcal{C}$ leads to $\LL^{\mathcal{C}}_{c'} \circ D^i_j
    - D^i_j\circ \LL^{\mathcal{C}}_{c'} = (\LL_{c'} \circ D - D \circ
    \LL_{c'})^i_j$ for all $i,j\in \{1,\ldots,N\}$ and the induction
    easily follows. The statement concerning $G$-invariance is clear.
\end{proof}

\section{Differential algebra and module structures}
\label{sec:differential-algebra-module-structures}

In this section we investigate algebra and module structures which can
be seen as differential operators in the sense of the previous
section. The main goal is to present the general conditions which are
necessary to define a notion of \emph{differential} such that the
properties of the Definitions~\ref{definition:algebra-type} and
\ref{definition:module-type} are achieved. For associative algebras
this definition is rather trivial.

\begin{definition}[Associative algebras of differential type]
    \label{definition:algebra-type-diff}
    An associative algebra $(\mathcal{A},\mu)$ is said to be of
    \emph{differential type} or simply \emph{differential} with
    respect to a commutative, associative algebra $\mathcal{C}$ if
    $\mathcal{A}$ is a (left) $\mathcal{C}$-module and if there exist
    $m,n\in \mathbb{N}_0$ such that the algebra multiplication is a
    bidifferential operator
    \begin{equation}
        \label{eq:algebra-type-diff}
        \mu\in \Diffop^{(m,n)}_{\mathcal{C}} (\mathcal{A},
        \mathcal{A}; \mathcal{A}).
    \end{equation}
\end{definition}

Due to the general properties of algebraically defined
multidifferential operators this definition yields a particular type
as in Definition~\ref{definition:algebra-type}.

\begin{corollary}[Associative algebras of differential type]
    \label{corollary:algebra-type-diff}
    A differential associative algebra $(\mathcal{A},\mu)$ as in
    Definition~\ref{definition:algebra-type-diff} satisfies the
    conditions of Definition~\ref{definition:algebra-type} with the
    differential Hochschild cochains
    $\HCdiff^\bullet(\mathcal{A},\mathcal{A})$ given by the
    multidifferential operators
    \begin{equation}
        \label{eq:Hochschild-complex-diff-algebra}
        \HCdiff^k(\mathcal{A},\mathcal{A}) =
        \Union_{L\in \mathbb{N}_0^k} \Diffop^L_{\mathcal{C}}
        (\mathcal{A},\dots, \mathcal{A} 
        ; \mathcal{A}) \quad \textrm{for } k\in \mathbb{N} 
    \end{equation}
    and $\HCdiff^0(\mathcal{A},\mathcal{A})=\mathcal{A}$.
\end{corollary}

\begin{proof}
    The only nontrivial point is the closedness under the insertions
    which is a well-known fact for differential operators, confer
    Proposition~\ref{proposition:composition-diffop}. The closedness
    under the cup product and all other operations discussed in
    Section \ref{subsec:Hochschild-complex-algebra} is clear since
    they can be expressed in terms of insertions.
\end{proof}

\begin{remark}
    \label{remark:differential-algebra-structure}
    If $\mathcal{A}$ is commutative and the $\mathcal{C}$-module
    structure is compatible with the algebra multiplication, the
    latter is always differential with $m=n=0$. This is in particular
    the case for $\mathcal{A}=\mathcal{C}$.
\end{remark}

For right module structures $\rho: \mathcal{E}\times \mathcal{A}
\longrightarrow \mathcal{E}$ and the derived Hochschild complexes the
situation is more difficult. There, the differential type has to be
defined in a more subtle way in order to guarantee the properties of
Definition \ref{definition:module-type}, in particular the closedness
under the cup product.

\begin{definition}[Right module structures of differential type]
    \label{definition:module-type-diff}
    Let $(\mathcal{A}, \mu)$ be a differential associative algebra
    with respect to a commutative associative algebra $\mathcal{C}$ as
    in Definition \ref{definition:algebra-type-diff}. A right
    $\mathcal{A}$-module structure $\rho$ of a $\mathbb{K}$-vector
    space $\mathcal{E}$ is said to be of \emph{differential type} or
    simply \emph{differential} with respect to
    $\HCdiff^\bullet(\mathcal{A},\mathcal{A})$ if the following
    conditions are satisfied.
    \begin{enumerate}
    \item $\mathcal{E}$ has a left module structure
        $l:\mathcal{B}\times \mathcal{E}\longrightarrow \mathcal{E}$
        with respect to a commutative associative algebra
        $\mathcal{B}$, also seen as
        $l:\mathcal{B}\longrightarrow \End_{\mathbb{K}}
        (\mathcal{E},\mathcal{E})$.
    \item There exists an algebra homomorphism $\gamma: \mathcal{C}
        \longrightarrow \mathcal{B}$.
    \item The right module structure is a differential operator
        \begin{equation}
            \label{eq:module-type-diff}
            \rho\in \Diffop^{L_\rho}_{\mathcal{C}} (\mathcal{A};
            \Diffop^{l_\rho}_{\mathcal{B}} (\mathcal{E},\mathcal{E}))
        \end{equation}
        with $L_\rho,l_\rho\in \mathbb{N}_0$ and where the spaces
        $\mathcal{D}^l=\Diffop^l_{\mathcal{B}}
        (\mathcal{E},\mathcal{E})$ for all $l\in \mathbb{N}_0$ are
        equipped with the natural left $\mathcal{C}$-module structure
        \begin{equation}
            \label{eq:induced-module-diffop}
            cD=\LL_cD= l(\gamma (c)) \circ D \quad \textrm{for all } c\in
            \mathcal{C}, D\in \mathcal{D}^l.
        \end{equation}
    \end{enumerate}
\end{definition}

The left multiplication $\LL_c$ with $c\in \mathcal{C}$ in
\eqref{eq:induced-module-diffop} in fact defines a left
$\mathcal{C}$-module structure since $l(\gamma (c))\in
\Diffop^0_\mathcal{B}(\mathcal{E}, \mathcal{E})$ and $\LL_c
(\LL_{c'}D)= l(\gamma (c)) \circ l(\gamma (c')) \circ D=
l(\gamma(cc')) \circ D= \LL_{cc'} D$. Using this module structure to
define corresponding multidifferential operators of $\mathcal{A}$ with
values in $\Diffop^l_{\mathcal{B}} (\mathcal{E},\mathcal{E})$ leads to
the following proposition.

\begin{proposition}
    \label{proposition:module-type-diff}
    A differential right module structure $\rho$ as in Definition
    \ref{definition:module-type-diff} satisfies the conditions of
    Definition \ref{definition:module-type} with the differential
    Hochschild cochains $\HCdiff^\bullet(\mathcal{A},\mathcal{D})$
    defined by
    \begin{equation}
        \label{eq:subalgebra-Hochschild-module-diff}
        \HCdiff^0(\mathcal{A},\mathcal{D})= \mathcal{D}=
        \Union_{l\in \mathbb{N}_0}
        \underbrace{\Diffop^l_{\mathcal{B}} (\mathcal{E};\mathcal{E})} 
        _{= \mathcal{D}^l} 
    \end{equation}
    and
    \begin{equation}
        \label{eq:Hochschild-complex-diff-module}
        \HCdiff^k(\mathcal{A},\mathcal{D}) =
        \Union_{L\in \mathbb{N}_0^k} \Union_{l\in \mathbb{N}_0}
        \underbrace{\Diffop^L_{\mathcal{C}} (\mathcal{A}, \dots, \mathcal{A};
          \Diffop^l_{\mathcal{B}}(\mathcal{E},\mathcal{E}))}
        _{=:\HC^{k,L,l}} \quad
        \textrm{for } k\in \mathbb{N}
    \end{equation}
    which are filtered spaces with $\mathcal{D}^l\subseteq
    \mathcal{D}^{l'}$ if $l\le l'$ and $\HC^{k,L,l}\subseteq
    \HC^{k,L',l'}$ if $L\le L'$ and $l\le l'$. By setting
    $\HC^{0,L,l}= \mathcal{D}^l$ the cup product $\cup$ as in
    \eqref{eq:cup-product-module} has the property
    \begin{equation}
        \label{eq:cup-product-differential}
        \HC^{k_1,L_1,l_1}\cup \HC^{k_2,L_2,l_2}\subseteq \HC^{k_1+k_2,
          (L_1, L_2+\overline{l_1}), l_1+l_2} \quad \textrm{with }
        \overline{l_1}=(l_1,\dots,l_1)\in \mathbb{N}_0^{k_2}
    \end{equation}
    for all $k_i\in\mathbb{N}_0$,
    $L_i\in\mathbb{N}_0^{k_i}$, $l_i\in \mathbb{N}_0$ and $i=1,2$.
\end{proposition}

\begin{proof}
    It is sufficient to prove \eqref{eq:cup-product-differential}
    which guarantees the closedness under the cup product. The
    remaining conditions in Definition \ref{definition:module-type}
    then are obvious or follow from
    Proposition~\ref{proposition:composition-diffop}. The stated
    filtration is clear with the filtrations of $\Diffop^\bullet
    (\mathcal{E}, \mathcal{E})$ and
    $\Diffop^\bullet(\mathcal{A},\dots, \mathcal{A};
    \Diffop^l(\mathcal{E}, \mathcal{E}))$ for all $l\in \mathbb{N}_0$.

    The left multiplication $l(b) \in \mathcal{D}^0$ in $\mathcal{E}$
    with an element $b\in \mathcal{B}$ obviously induces a left
    multiplication $\LL_b$ and a right multiplication $\mathsf{R}_b$
    in $\mathcal{D}^l$ for all $l\in \mathbb{N}_0$ which are given by
    $\LL_bD=l(b)\circ D$ and $\mathsf{R}_bD=D\circ l(b)$. Then one
    first notes that for $\phi\in \HC^{k,L,l}$ and $b\in \mathcal{B}$
    one has $\LL_b\circ \phi \in \HC^{k,L,l}$ which is seen by an easy
    induction over $|L|$ because $\LL_c \circ (\LL_b\circ \phi)-
    (\LL_b\circ \phi)\circ \LL_c^{(i)}=\LL_b \circ (\LL_c \circ \phi-
    \phi\circ \LL_c^{(i)})$ with the notation introduced in the
    previous section. Analogously one sees that $\mathsf{R}_b \circ
    \phi \in \HC^{k,L,l}$ and by the very definition of differential
    operators one further gets $\LL_b \circ \phi-
    \mathsf{R}_b\circ\phi \in \HC^{k,L,l-1}$.

    For fixed $k_1,k_2$ the proof of the statement
    \eqref{eq:cup-product-differential} is an induction over $r=
    |L_1|+|L_2|+l_1+l_2$. Let $\phi_1\in \HC^{k_1,L_1,l_1}$ and
    $\phi_2\in \HC^{k_2,L_2,l_2}$. With the definition
    \begin{equation*}
        (\phi_1\cup\phi_2)(a_1,\dots, a_{k_1+k_2})= \phi_1(a_1,\dots,
        a_{k_1})\circ \phi_2(a_{k_1+1},\dots, a_{k_1+k_2})
    \end{equation*}
    it is always clear that $\phi_1\circ \phi_2$ has values in
    $\mathcal{D}^{l_1+l_2}$ since the usual composition of maps
    satisfies $\mathcal{D}^{l_1}\circ \mathcal{D}^{l_2}\subseteq
    \mathcal{D}^{l_1+l_2}$. For $r=0$ the $\mathcal{C}$-linearity of
    $\phi_1\cup\phi_2$ is clear by the fact that
    $|L_1|=|L_2|=l_1=0$. Under the assumption that the statement holds
    for $r$ consider elements $\phi_s\in \HC^{k_s,L_s,l_s}$ for
    $s=1,2$ with $|L_1|+|L_2|+l_1+l_2=r+1$.  With $c\in \mathcal{C}$
    one then computes for $i=1,\dots,k_1$
    \begin{equation*}
        \LL_c \circ (\phi_1 \cup \phi_2) - (\phi_1 \cup
        \phi_2)\circ \LL_c^{(i)}
        = (\LL_c \circ \phi_1 - \phi_1\circ \LL_c^{(i)}) \cup
        \phi_2. 
    \end{equation*}
    Since $\LL_c \circ \phi_1 - \phi_1\circ \LL_c^{(i)} \in
    \HC^{k_1,L_1-\mathsf{e}_i,l_1}$ it follows by assumption that
    \begin{equation}
        \label{eq:first-part-induction}
        \LL_c \circ (\phi_1\cup \phi_2) - (\phi_1\cup \phi_2)
        \circ \LL_c^{(i)} \in 
        \HC^{k_1+k_2,(L_1, L_2+\overline{l_1})-\mathsf{e}_i,
          l_1+l_2}. 
    \end{equation}
    Since $\chi= \LL_c \circ \phi_2 - \phi_2 \circ \LL_c^{(j)} \in
    \HC^{k_2,L_2-\mathsf{e}_j,l_2}$ one obtains for $i=k_1+j$ with
    $j=1,\dots, k_2$
    \begin{eqnarray}
        \label{eq:crucial-point-cup}
        \LL_c \circ (\phi_1 \cup \phi_2) - (\phi_1 \cup
        \phi_2)\circ \LL_c^{(i)}
        &=& (\LL_c \circ \phi_1)\cup \phi_2 - \phi_1 \cup (\phi_2
        \circ \LL_c^{(j)}) \nonumber \\
        &=& (\LL_{\gamma (c)} \circ \phi_1- \mathsf{R}_{\gamma
          (c)}\circ \phi_1) \circ \phi_2 + \phi_1 \cup \chi. 
    \end{eqnarray}
    With $\LL_{\gamma (c)} \circ \phi_1- \mathsf{R}_{\gamma (c)} \circ
    \phi_1 \in \HC^{k_1,L_1, l_1-1}$ and the filtration,
    \eqref{eq:first-part-induction} is true for all
    $i=1,\dots,k_1+k_2$. By the definition of multidifferential
    operators this yields $\phi_1 \cup \phi_2 \in \HC^{k_1+k_2,
      (L_1,L_2+\overline{l_1}),l_1+l_2}$.
\end{proof}

\begin{remark}
\label{remark:notes-on-differential-Hochschild-complexes}
    \begin{enumerate}
    \item In the case $\mathcal{A}=\mathcal{C}$, what in particular
        means that $\mathcal{A}$ is commutative, $\mathcal{E}$
        inherits two different $\mathcal{A}$-module structures
        \begin{eqnarray}
            \label{eq:necessary-module}
            a \mathbin{\cdot_1} e &=& l(\gamma (a)) e \quad
            \textrm{and}\\
            a \mathbin{\cdot_2} e &=& \rho(a) e.
        \end{eqnarray}
        For the definition of $\HC^{k,L,l}$ it is important to use the
        left $\mathcal{A}$-module structure of
        $\Diffop^{l}(\mathcal{E},\mathcal{E})$ which is induced by
        \eqref{eq:necessary-module}. Otherwise, even in the case $\rho
        \in \HC^{1,L_\rho,0}$, the closedness under the cup product is
        not valid, confer the crucial point
        \eqref{eq:crucial-point-cup} in the proof.

        If, in addition, $\mathcal{E}=\mathcal{B}$ is a commutative
        and associative algebra and if all structures have to be
        traced back to the right module structure $\rho$, the
        assertions of Proposition \ref{proposition:module-type-diff}
        are only true, if $\rho(a)=l(\gamma (a))$ is a left
        multiplication in $\mathcal{E}$. This is a necessary condition
        for $\rho$ implying $L_\rho=l_\rho=0$ for the degrees of
        differentiation. For associative and commutative algebras
        $\mathcal{E}$ and $\mathcal{A}$ the differential setting is
        thus determined by an algebra homomorphism $\gamma:\mathcal{A}
        \longrightarrow \mathcal{E}$.
    \item Note that \eqref{eq:cup-product-differential} becomes false
        if one neglects the increase $\overline{l_1}$ in the
        multiorder of differentiation. Due to this fact which states
        that the order $(L_1,L_2+\overline{l_1})$ of differentiation
        of a corresponding cup product $\phi_1\cup \phi_2$ depends on
        the order $l_1$ of the values of $\phi_1$, one sees that in
        general the naive definition
        $\HCdiff^k(\mathcal{A},\mathcal{D}) = \Union_{L\in
          \mathbb{N}_0^k} \Diffop^L(\mathcal{A}, \dots,
        \mathcal{A};\mathcal{D})$ would not lead to a well-defined cup
        product since particular products $\phi_1\cup \phi_2$ possibly
        would not lead to differential operators of a certain degree
        of differentiation. Therefore it is essential to demand that
        the maps $\phi\in \HCdiff^\bullet (\mathcal{A},\mathcal{D})$
        have values in the differential operators $\mathcal{D}^l$ with
        some fixed order $l$ of differentiation which does not depend
        on the arguments of $\phi$.
    \end{enumerate}
\end{remark}

\begin{remark}[$G$-invariant differential algebra and module
    structures]
    \label{remark:G-invariant-differential-type}
    Due to Lemma~\ref{lemma:Diffops-and-representations} it is clear
    that the conditions \eqref{eq:G-compatible-types} in
    Definition~\ref{definition:G-invariant-type} are satisfied for
    $G$-invariant and differential algebras $\mathcal{A}$ and
    $\mathcal{A}$-modules $\mathcal{E}$ if the additionally introduced
    structures in the Definitions~\ref{definition:algebra-type-diff}
    and \ref{definition:module-type-diff} are $G$-invariant. This
    means that the left $\mathcal{C}$-module structure of
    $\mathcal{A}$, the left $\mathcal{B}$-module structure of
    $\mathcal{E}$ and the algebra morphism $\gamma:\mathcal{C}
    \longrightarrow \mathcal{B}$ have to be $G$-invariant with respect
    to corresponding representations. In this situations we then can
    consider the particular type of $G$-invariant differential
    structures.
\end{remark}

\begin{corollary}
    \label{corollary:Differential-Hochschild-module}
    With $\mu$ as in \eqref{eq:algebra-type-diff}, $\rho$ as in
    \eqref{eq:module-type-diff} and the resulting
    $(\mathcal{A},\mathcal{A})$-bimodule structure of $\mathcal{D}$ as
    in \eqref{eq:bimodule-of-endomorphisms} the differential $\delta$
    of the corresponding differential Hochschild complex
    $(\HCdiff^\bullet(\mathcal{A},\mathcal{D}),\delta)$ has the
    restrictions
    \begin{equation}
        \label{eq:Hochschild-differential-general}
        \delta: \HC^{k,L, l} \longrightarrow \HC^{k+1,\tilde{L},
          l+l_\rho} 
    \end{equation} 
    with $\tilde{L}= (\tilde{l}_1,\dots, \tilde{l}_{k+1})$ where for
    $k\ge 2$
    \begin{eqnarray}
        \label{eq:grading-differential-general}
        \tilde{l}_1 &=& \max \{l_1+m, L_\rho+l, l_1+l_\rho \}
        \nonumber \\
        \tilde{l}_i &=& \max \{l_i+m,l_{i-1}+n, l_i+l_\rho \}
        \quad \textrm{for}\quad i=2,\dots,k\\
        \tilde{l}_{k+1} &=& \max \{l_k+n, l_k,L_\rho \}. \nonumber
    \end{eqnarray}
    If for commutative $\mathcal{A}$ the structure
    \eqref{eq:new-bimodule-of-endomorphism} is used,
    Equation~\eqref{eq:grading-differential-general} has to be
    replaced by
    \begin{eqnarray}
        \label{eq:grading-differential-commutative}
        \tilde{l}_1 &=& \max \{L_\rho, l_1+m \}
        \nonumber \\
        \tilde{l}_i &=& \max \{l_{i-1}+l_\rho, l_i+m, l_{i-1}+n \}
        \quad \textrm{for}\quad i=2,\dots,k\\
        \tilde{l}_{k+1} &=& \max \{l_k+l_\rho, l_k+n, L_\rho+l \}. \nonumber
    \end{eqnarray}
    In both cases every $\HC^{k,L,l}$ is a left $\mathcal{B}$-module
    via $b\cdot \phi= \LL_b \circ \phi$ with $\LL_b$ as in the proof
    of Proposition~\ref{proposition:module-type-diff}.  For $l_\rho =
    0$, this means if $\mathcal{E}$ is a
    $(\mathcal{B},\mathcal{A})$-bimodule, one has
    \begin{equation}
        \label{eq:Hochschild-delta-left-multiplication}
        \delta (b\cdot \phi) = b\cdot (\delta \phi).
    \end{equation}
\end{corollary}

\begin{proof}
    For $\phi\in \HC^{k,L, l}$ and $a_1,\dots, a_{k+1}\in \mathcal{A}$
    consider the expression \eqref{eq:Hochschild-differential-module}.
    Due to \eqref{eq:cup-product-differential}, the occurring cochains
    are $\phi\cup \rho \in \HC^{k+1, (L,L_\rho +l), l+l_\rho}$ and
    $\rho \cup \phi \in \HC^{k+1, (L_\rho, L+ \overline{l_\rho}),
      l+l_\rho}$. Further, $\phi\circp{i-1} \mu \in \HC^{k+1, (l_1,
      \dots,l_{i-1}, l_i+m, l_i+n, l_{i+1}, \dots, l_k),l}$ by the
    general property of multidifferential operators and the
    precondition concerning $\mu$. This shows
    \eqref{eq:Hochschild-differential-general} and a simple counting
    of orders of differentiation leads to
    \eqref{eq:grading-differential-general}. With the bimodule
    structure \eqref{eq:new-bimodule-of-endomorphism} the
    considerations are exactly the same using
    \eqref{eq:new-Hochschild-differential}. Equation
    \eqref{eq:Hochschild-delta-left-multiplication} is obvious with
    the definition of $\delta$.
\end{proof}

\begin{remark}
    The grading in \eqref{eq:Hochschild-differential-general} for the
    cases $k=0$ and $k=1$ can be treated in the same way.
\end{remark}

\begin{remark}
    Note that
    $\mathcal{D}^l$ with $l\in \mathbb{N}_0$ is an
    $(\mathcal{A},\mathcal{A})$-bimodule via
    \eqref{eq:bimodule-of-endomorphisms} or
    \eqref{eq:new-bimodule-of-endomorphism} only if $\rho\in
    \HC^{1,L_\rho,0}$.
\end{remark}

As shown, the considered differential structures yield corresponding
\emph{differential Hochschild complexes}. In general this can be
defined as follows.

\begin{definition}[Differential Hochschild complex]
    \label{definition:differential-Hochschild-complex}
    Let $\mathcal{C}$ be an associative and commutative
    $\mathbb{K}$-algebra. Further let $\mathcal{A}$ be an associative
    $\mathbb{K}$-algebra and $\mathcal{M}$ be an
    $(\mathcal{A},\mathcal{A})$-bimodule as in Definition
    \ref{definition:Hochschild-complex}, both endowed with a left
    module structure with respect to $\mathcal{C}$. If the
    algebraically defined multidifferential operators
    \begin{equation}
        \label{eq:differential-Hochschild-complex}
        \HCdiff^k(\mathcal{A},\mathcal{M})= \Union_{L\in
          \mathbb{N}_0^k} \Diffop^L_{\mathcal{C}}(\mathcal{A},\dots
        \mathcal{A};\mathcal{M}) \subseteq 
        \HC^k( \mathcal{A},\mathcal{M}), \quad k\in \mathbb{N},
    \end{equation}
    and $\HCdiff^0(\mathcal{A},\mathcal{M})=\mathcal{M}$ build a
    subcomplex $(\HCdiff^k(\mathcal{A},\mathcal{M}), \delta)$, this is
    called the \emph{(algebraic) differential Hochschild complex of
      $\mathcal{A}$ with values in $\mathcal{M}$ over $\mathcal{C}$}.
\end{definition}

In general, \eqref{eq:differential-Hochschild-complex} does not define
a subcomplex since the considered subspaces could not be closed under
the Hochschild differential $\delta$. In order to guarantee this one
has to demand adequate properties of the involved algebraic
structures.

\begin{lemma}
    \label{lemma:conditions-differential-Hochschild-complex}
    Let $\mathcal{A}$ and $\mathcal{M}$ as above satisfy the following
    conditions:
    \begin{enumerate}
    \item The algebra multiplication in $\mathcal{A}$ is a
        differential operator $\mu \in
        \Diffop^{(m,n)}_{\mathcal{C}}(\mathcal{A},\mathcal{A};\mathcal{A})$
        with $m,n\in \mathbb{N}_0$.
    \item  The left and right $\mathcal{A}$-module
    structures of $\mathcal{M}$ have the property that for all $k\in
    \mathbb{N}_0$, $L=(l_1,\dots, l_k)\in \mathbb{N_0}^k$ and $\phi\in
    \Diffop^L_{\mathcal{C}}(\mathcal{A},\mathcal{M})$ there exist
    $s,t,u_\phi,v_\phi\in \mathbb{N}_0$  such that the maps
    \begin{align}
        \label{eq:differential-left-module}
        &\left((a_1,\dots, a_{k+1})\longmapsto a_1 \phi(a_2, \dots,
            a_{k+1})\right) \in \Diffop^{L_1}_{\mathcal{C}}
        (\mathcal{A},\dots, \mathcal{A}; \mathcal{M}),\\
        &\left((a_1,\dots, a_{k+1})\longmapsto \phi(a_1, \dots, a_k)
            a_{k+1} \right) \in
        \Diffop^{L_2}_{\mathcal{C}}(\mathcal{A},\dots, \mathcal{A};
        \mathcal{M})
    \end{align}
    are differential operators of the multiorders $L_1=(u_\phi,
    l_1+s, \dots, l_k+s)\in \mathbb{N_0}^{k+1}$ and\\ $L_2=(l_1+t,
    \dots, l_k+t, v_\phi)\in \mathbb{N_0}^{k+1}$.
    \end{enumerate}
    Then \eqref{eq:differential-Hochschild-complex} defines a
    differential Hochschild complex since for $L=(l_1,\dots, l_k)\in
    \mathbb{N_0}^k$
    \begin{equation}
        \label{eq:differential-differential-Hochschild}
        \delta: \Diffop^L_{\mathcal{C}}(\mathcal{A}, \dots,
        \mathcal{A};\mathcal{M}) \longrightarrow
        \Diffop^{\tilde{L}}_{\mathcal{C}}(\mathcal{A}, \dots, 
        \mathcal{A};\mathcal{M}),
    \end{equation}
    where $\tilde{L}= (\tilde{l}_1, \dots, \tilde{l}_{k+1})\in
    \mathbb{N}_0^{k+1}$ with
    \begin{eqnarray}
        \label{eq:grading-Hochschild-differential}
        \tilde{l}_1 &=& \max \{ u_\phi, l_1+t, l_1+m \}
        \nonumber \\
        \tilde{l}_i &=& \max \{l_{i-1}+s , l_i+t, l_i+m, l_{i-1}+n \}
        \quad \textrm{for}\quad i=2,\dots,k\\
        \tilde{l}_{k+1} &=& \max \{l_k+s , v_\phi, l_k+n \}. \nonumber
    \end{eqnarray}
\end{lemma}

The proof is obvious and a direct consequence of the definition of
$\delta$ and the properties of differential operators. In the
subcomplexes derived from differential algebra and module structures
we have basically shown that Lemma
\ref{lemma:conditions-differential-Hochschild-complex} is
satisfied. In the case of algebras we have $s=0=t$ and
$u_\phi=0=v_\phi$ for all $\phi\in \HCdiff^\bullet (\mathcal{A},
\mathcal{A})$. Depending on the possibly given choice of the
$(\mathcal{A},\mathcal{A})$-bimodule structure of $\mathcal{D}$ we
have for right modules either $s=0$, $t=l_\rho$ and $v_\phi=L_\rho$,
$u_\phi=l$ for $\phi\in \HC^{k,L,l}$ in the case
\eqref{eq:bimodule-of-endomorphisms} or $s=l_\rho$, $t=0$ and
$u_\phi=L_\rho$, $v_\phi=l$ for $\phi\in \HC^{k,L,l}$ in the case
\eqref{eq:new-bimodule-of-endomorphism}.

As a simple but important example one can again consider the case of
projective modules. The results of
Section~\ref{sec:projective-modules} can be reformulated in the
setting of differential module structures.

\begin{proposition}[Projective modules of the differential type]
    \label{proposition:projective-differential-modules}
    Let $\mathcal{A}$ be a commutative, associative algebra and let
    $\mathcal{E}$ be a differential right module as in
    Definition~\ref{definition:module-type-diff} with
    $\mathcal{A}=\mathcal{C}=\mathcal{B}$ and
    $\gamma=\id_\mathcal{A}$. Further, let $\mathcal{E}$ be a
    projective $\mathcal{A}$-module.

    Then, the homotopy map $h^k$ from
    Section~\ref{sec:projective-modules} satisfies
    \begin{equation}
        \label{eq:homotopy-projective-differential-explicit}
        h^k: \Diffop^{(l_1,\dots,l_k)}(\mathcal{A},\dots, \mathcal{A};
        \Diffop^l(\mathcal{E},\mathcal{E})) \longrightarrow
        \Diffop^{(l_2,\dots,l_k)}(\mathcal{A},\dots,\mathcal{A};
        \Diffop^{l_1}(\mathcal{E},\mathcal{E})) 
    \end{equation}
    for all $l_1,\dots,l_k,l\in \mathbb{N}_0$ and $k\in \mathbb{N}$.
    This shows that the corresponding cohomology is trivial,
    \begin{equation}
        \label{eq:projective-differential-cohomology}
        \HHdiff^k(\mathcal{A},\mathcal{D})=
        \left\{
            \begin{array}{c}
                \Diffop^0_{\mathcal{A}}(\mathcal{E},\mathcal{E})\\
                \{0\}
            \end{array}
        \right.
        \textrm{for}
        \begin{array}{c}
            k=0\\
            k\ge 1.
        \end{array}
    \end{equation}
\end{proposition}

\begin{proof}
    The proof is a twofold induction. Using the defining
    Equation~\eqref{eq:homotopy-maps-projective-explicit}, a simple
    computation shows that for all $\phi\in
    \HCdiff^k(\mathcal{A},\mathcal{D})$ and $a,a_1,\dots, a_{k-1}\in
    \mathcal{A}$ one has
    \begin{equation}
        \label{eq:first-induction}
        \LL_a\circ h^k\phi(a_1,\dots, a_{k-1})-h^k\phi(a_1,\dots,
        a_{k-1}) \circ \LL_a= h^k(\LL_a\circ \phi- \phi\circ
        \LL_a^{(1)})(a_1,\dots,a_{k-1}).  
    \end{equation}
    Due to this equation the first induction over $l_1$ shows that
    $h^k\phi$ takes values in $\Diffop^{l_1}(\mathcal{E},\mathcal{E})$
    if $\phi$ is a differential operator of multiorder
    $(l_1,\dots,l_k)$ for arbitrary $l_2,\dots,l_k\in \mathbb{N}_0$. A
    second induction over $r=l_2+\dots+l_k$ making use of the fact
    that
    \begin{equation}
        \label{eq:second-induction}
        \LL_a \circ h^k\phi -h^k\phi \circ \LL_a^{(i)}= h^k(\LL_a
        \circ \phi -\phi \circ \LL_a^{(i+1)}) 
    \end{equation}
    for all $i=1,\dots,k_1$ then shows the remaining assertion
    contained in \eqref{eq:homotopy-projective-differential-explicit}.
    The statement concerning the differential cohomology then follows
    in the same way as the computation of the purely algebraic one in
    Section~\ref{sec:projective-modules}.
\end{proof}

The above proposition yields an easy proof for the well-known facts
concerning deformed vector bundles, confer
Section~\ref{sec:quantization-vector-bundles}.

We close this section with two easy observations concerning
differential Hochschild complexes which will be very useful for the
later computation of the corresponding cohomologies.

\begin{remark}[Functoriality of differential Hochschild complexes]
    \label{remark:functoriality-Hochschild-diff}
    It is well-known that Hochschild complexes have a functorial
    behaviour with respect to their initial data. With
    Lemma~\ref{lemma:diffop-functorial} it is a simple exercise to see
    that this is still true in the differential or $G$-invariant
    differential setting. Thus, different choices of structures as
    $\mathcal{A}$, $\mathcal{E}$ and $\gamma:\mathcal{C}
    \longrightarrow \mathcal{B}$ in
    Definition~\ref{definition:algebra-type-diff} and
    \ref{definition:module-type-diff} which are related by structure
    preserving isomorphisms yield isomorphic differential Hochschild
    complexes with exactly the same properties.
\end{remark}

The results of Lemma~\ref{lemma:Diffop-free-modules} induce the
following simple observation.

\begin{lemma}[Differential Hochschild complexes and free modules]
    \label{lemma:special-situation-free-module}
    Let there be given the structures in
    Definition~\ref{definition:module-type-diff} for the following
    specific situation: Let $\mathcal{A}= \mathcal{C}$ be commutative and
    let the right $\mathcal{A}$-module structure of $\mathcal{E}$ be given
    by
    \begin{equation}
        \label{eq:right-module-structure-left-multiplication}
        \rho(a)=l_{\gamma (a)},
    \end{equation}
    so that $\rho \in \HC^{1,0,0}$. Then, for all $l\in \mathbb{N}_0$
    let the space $\mathcal{D}^l=\Diffop^l_{\mathcal{B}}
    (\mathcal{E},\mathcal{E})$ be equipped with the
    $(\mathcal{A},\mathcal{A})$-bimodule structure
    \eqref{eq:new-bimodule-of-endomorphism}, this means $a\cdot D\cdot
    a'= \rho(a)\circ D\circ \rho(a')$ for all $a,a'\in \mathcal{A}$
    and $D\in \mathcal{D}^l$, so that the left module structure of
    $\mathcal{D}^l$ with respect to $\mathcal{C}=\mathcal{A}$ is
    unique. In addition, let $\mathcal{E}$ be a free
    $\mathcal{B}$-module with module basis
    $\mathfrak{B}=\{e_1,\ldots, e_N\}$.

    Then, the isomorphisms $\Phi_\mathfrak{B}:
    \Diffop^l_{\mathcal{B}}(\mathcal{E};\mathcal{E}) \cong
    \Mat_{N\times
      N}(\Diffop^l_{\mathcal{B}}(\mathcal{B};\mathcal{B}))$ as in
    Lemma~\ref{lemma:Diffop-free-modules} are isomorphisms of
    $(\mathcal{A},\mathcal{A})$-bimodules when using the obvious
    $\mathcal{A}$-module structure of $\mathcal{B}$ induced by
    $\gamma$. Further, this gives rise to an isomorphism
    \begin{eqnarray}
        \label{eq:Iso-complex-matrix}
        \left(\HCdiff^\bullet(\mathcal{A},
        \Diffop_{\mathcal{B}}(\mathcal{E},\mathcal{E})),\delta \right)  
        &\cong&
        \left(\HCdiff^\bullet(\mathcal{A},
        \Mat_{N\times N}(\Diffop_{\mathcal{B}} (\mathcal{B},
        \mathcal{B}))),\delta\right)\\
        &\cong&
        \label{eq:Iso-matrix-matrix}
        \left( \Mat_{N\times N} (\HCdiff^\bullet(\mathcal{A},
            \Diffop_{\mathcal{B}} (\mathcal{B},
            \mathcal{B}))),\delta \right)
    \end{eqnarray}
    of complexes where the differential in
    \eqref{eq:Iso-matrix-matrix} is defined componentwise,
    $\delta(\phi^i_j):=(\delta \phi^i_j)$. In particular, one has
    \begin{equation}
        \label{eq:Iso-complex-matrix-grading}
        \Diffop^L_\mathcal{A} (\mathcal{A},\ldots,\mathcal{A};
        \Diffop^l_{\mathcal{B}}(\mathcal{E},\mathcal{E})) 
        \cong
        \Diffop^L_\mathcal{A} (\mathcal{A},\ldots,\mathcal{A}; 
        \Mat_{N\times N}
        (\Diffop^l_{\mathcal{B}}(\mathcal{B},\mathcal{B}))).
    \end{equation}
    In the $G$-invariant setting the $G$-invariant module basis
    guarantees that the isomorphisms \eqref{eq:Iso-complex-matrix},
    \eqref{eq:Iso-matrix-matrix}, and
    \eqref{eq:Iso-complex-matrix-grading} are $G$-invariant.
\end{lemma}

\begin{proof}
    That $\Phi_{\mathfrak{B}}$ is an isomorphism of
    $(\mathcal{A},\mathcal{A})$-bimodules is obvious with the given
    module structures and the fact that $\Phi_{\mathfrak{B}}$ is
    already an isomorphism of
    $(\mathcal{B},\mathcal{B})$-bimodules. The isomorphisms
    \eqref{eq:Iso-complex-matrix-grading} and
    \eqref{eq:Iso-complex-matrix} then follow by
    Lemma~\ref{lemma:diffop-functorial} and
    Remark~\ref{remark:functoriality-Hochschild-diff}. The isomorphism
    \eqref{eq:Iso-matrix-matrix} is an isomorphism of complexes since
    the Hochschild differential $\delta$ is defined componentwise.
\end{proof}

Note that \eqref{eq:right-module-structure-left-multiplication} is
necessary in order to guarantee that $\mathcal{B}$ has a right
$\mathcal{A}$-module structure which is related to the one of
$\mathcal{E}$.

\begin{corollary}[The cohomology of free modules]
    Due to \eqref{eq:Iso-matrix-matrix} the knowledge of the
    cohomology $\HHdiff^\bullet(\mathcal{A},
    \Diffop_{\mathcal{B}}(\mathcal{E},\mathcal{E}))$ is equivalent to
    the know\-ledge of $\HHdiff^\bullet(\mathcal{A},
    \Diffop_{\mathcal{B}}(\mathcal{B},\mathcal{B}))$.
\end{corollary}

\chapter{Sheaf theory and Hochschild cohomologies}
\label{cha:sheaves}

It is a basic concept in differential geometry to investigate global
structures locally. The definition of a smooth manifold itself makes
use of local charts and many other fundamental notions are defined by
using local expressions of the involved data as well. In order to
compute Hochschild cohomologies that arise in a global geometric
context it is thus natural to ask if the problem can be solved by
investigating the local situation. This approach is of course only
possible if the global objects are related to corresponding local
expressions and if the global and local data determine each other in a
sufficient way. These aspired relations are the content of the notions
of presheaves and sheaves over topological spaces. As it will become
clear in this chapter this is the adequate framework to treat the
above problem. The basic definitions and general results of sheaf
theory that will be used in the subsequent considerations are
introduced in the first section. After that we concentrate on sheaves
over smooth manifolds and present a different point of view of the
sheaves of sheaf homomorphisms and particular subsheaves thereof. This
leads to the important result that all differential operators of a
finite order of differentiation between particular sheaves carry a
sheaf structure themselves. With this conclusion and the observations
concerning sheaves over different manifolds made in
Section~\ref{sec:sheaves-different-manifolds} it is possible to apply
these concepts to the discussion of differential Hochschild complexes
of modules as introduced in
Section~\ref{sec:differential-algebra-module-structures}. It will be
pointed out that in typical situations these complexes give rise to
corresponding presheaves. Although one has to abandon the desirable
properties of a sheaf it is shown that under certain assumptions one
is still able to compute the global cohomologies by computing the
local ones. Together with the investigation of the $G$-invariant case
which enforces a slightly different setting the chapter culminates in
the two important
Propositions~\ref{proposition:sheaves-and-cohomology} and
\ref{proposition:sheaves-and-cohomology-G-inv}.

\section{Sheaves and sheaf homomorphisms}
\label{sec:sheave-sheaf-homomorphisms}

The origin of sheaf theory can in principle be seen in the attempt to
find an axiomatic description of the convenient properties of
functions on a topological space with respect to the process of
restricting them to open subsets. It is obvious that functions have
the property to be determined by their restrictions and thus they
yield a generic example of a structure where the knowledge of the
global information is equivalent to the knowledge of the local one.
The notions of presheaves and sheaves developed from this idea have
important applications in many realms of mathematics, for instance in
algebraic geometry, complex analysis and many others. The formulation
of the essential properties can be performed in different equivalent
ways. An important reference for the basics of sheaf theory is the
early work of Godement \cite{godement:1964a}. Slightly different
approaches and applications can moreover be found in \cite[Chap.~II,
Sect.~1]{hartshorne:1977a}, \cite[Chap.~2]{vaisman:1973}, and
\cite[Chap.~II]{wells:1980a}. In the following presentation we
basically follow the latter two references.

\begin{definition}[Presheaves in terms of categories and functors]
    \label{definition:presheaves-functors}
    %\cite[Chap.~2, Sect.~1]{vaisman:1973}
    Let $M$ be a topological space and $\catopen{M}$ be the
    category of its open, nonempty subsets and inclusions. Further,
    let $\mathcal{C}$ be an arbitrary category.
    \begin{enumerate}
    \item A contravariant functor
        \begin{equation}
            \label{eq:presheaf-functor}
            \mathcal{F}: \catopen{M} \longrightarrow \mathcal{C}
        \end{equation}
        is called a \emph{presheaf} $\mathcal{F}$ \emph{over} or
        \emph{on} $M$ with values in $\mathcal{C}$.
    \item A \emph{homomorphism} between two presheaves
        $\mathcal{F},\mathcal{G}: \catopen{M} \longrightarrow
        \mathcal{C}$ over $M$ is a natural transformation between the
        two functors. Correspondingly, an \emph{isomorphism} is a
        natural equivalence.
    \end{enumerate}
\end{definition}

This precise but abstract definition obviously allows a more concrete
formulation which is also used to introduce a notion of
sheaves. Although we mainly work with the following description it is
often very useful to have the more basic definition in mind. With
respect to our further purpose we adjust the notation and restrict
ourselves to certain types of categories.

\begin{definition}[Presheaves and sheaves of sets]
    \label{definition:presheaves-sheaves}
      % \cite[Chap.~1, Sect.~1.9; Chap.~2, Sect.~1.1]{godement:1964a}, 
      % \cite[Chap.~II, Def.~1.1, 1.2]{wells:1980a}
    Let $M$ be a topological space and let $\underline{M}$ denote the
    set of its open, nonempty subsets. Further, let $\mathcal{C}$ be a
    category with certain sets as objects and certain maps as
    morphisms.
    \begin{enumerate}
    \item A \emph{presheaf $\mathcal{F}$ of objects of $\mathcal{C}$
          over $M$} is given by assigning to each $U\in \underline{M}$
        a set $\mathcal{F}(U)$ in the class of objects of
        $\mathcal{C}$,
        \begin{equation}
            \label{eq:assignment-presheaf}
            U \stackrel{\mathcal{F}}{\longmapsto} \mathcal{F}(U),
        \end{equation}
        and to each inclusion $V\subseteq U$ of open sets a
        \emph{restriction map} or \emph{restriction homomorphism}
        \begin{equation}
            \label{eq:restriction-morphisms}
            r^U_V: \mathcal{F}(U) \longrightarrow \mathcal{F}(V),
        \end{equation}
        such that
        \begin{enumerate}
        \item 
            \label{item:sheaf-1} 
            for all $U\in \underline{M}$ the map
            $r^U_U=\id_{\mathcal{F}(U)}$ is the identity morphism
            of $\mathcal{F}(U)$.
        \item 
            \label{item:sheaf-2}
            for all $W\subseteq V \subseteq U \subseteq M$ one
            has $r^V_W \circ r^U_V = r^U_W$.  
        \end{enumerate}
    \item A presheaf $\mathcal{F}$ is called a \emph{sheaf} if for
        every open covering $U= \Union_{i\in I} U_i$ with open
        $U\subseteq M$ the following further conditions are satisfied.
        \begin{enumerate}
            \setcounter{enumii}{2}
        \item 
            \label{item:sheaf-3}
            If $s,t \in \mathcal{F}(U)$ and $r^U_{U_i}(s)=
            r^U_{U_i}(t)$ for all $i\in I$, then $s=t$.
        \item 
            \label{item:sheaf-4}
            If $s_i\in \mathcal{F}(U_i)$ are given for all $i\in I$
            with the property that $r^{U_i}_{U_i\cap U_j}(s_i)=
            r^{U_j}_{U_i\cap U_j}(s_j)$ for all $U_i\cap U_j\neq
            \emptyset$, then there exists an $s\in \mathcal{F}(U)$
            such that $r^U_{U_i}(s)=s_i$ for all $i\in I$.            
        \end{enumerate}
    \end{enumerate}
\end{definition}

\begin{definition}[(Pre)sheaf homomorphisms]
    % \cite[Chap.~II, Def.~1.1, 1.2]{wells:1980a}
    A \emph{(homo)morphism} between two (pre)sheaves $\mathcal{F},
    \mathcal{G}$ (in the sense of Definition
    \ref{definition:presheaves-sheaves})
    \begin{equation}
        h: \mathcal{F} \longrightarrow \mathcal{G}
    \end{equation}
    is a collection $h=\{h_U\}_{U\in \underline{M}}$ of maps
    \begin{equation}
        h_U: \mathcal{F}(U) \longrightarrow \mathcal{G}(U)
    \end{equation}
    which is compatible with the restriction maps. This means that the
    diagram
    \begin{displaymath}
        \xymatrix{\mathcal{F}(U) \ar[r]^{h_{U}} \ar[d]_{r^U_V} &
          \mathcal{G}(U) \ar[d]^{r^U_V}\\ 
          \mathcal{F}(V) \ar[r]_{h_V} & \mathcal{G}(V)}
    \end{displaymath}
    commutes for all $V\subseteq U\subseteq M$. The set of all
    (pre)sheaf morphisms of $\mathcal{F}$
    into $\mathcal{G}$ is denoted by $\Hom(\mathcal{F},\mathcal{G})$.

    If all the maps $h_U$ are inclusions $\mathcal{F}$ is called a
    \emph{sub(pre)sheaf} of $\mathcal{G}$. In the case that all $h_U$
    are isomorphisms, $h$ is called a \emph{(pre)sheaf isomorphism}.
\end{definition}

\begin{remark}[Notation, (pre)sheaves with additional structures,
    restricted sheaves]
    \begin{enumerate}
    \item Usually, one denotes the value of a restriction map as
        $r^U_V(s)=s|_V$.
    \item In the most cases the considered categories $\mathcal{C}$
        consist of sets with an algebraic structure and structure
        preserving morphisms. Due to Definition
        \ref{definition:presheaves-functors} it is automatically clear
        that the restriction homomorphisms preserve these structures
        which therefore is always required in Definition
        \ref{definition:presheaves-sheaves}. Then, the same is
        demanded for corresponding homomorphisms. This way one is able
        to define sheaves of groups, vector spaces and all other kinds
        of algebraic structures, as well as the morphisms between
        them.  Further, this notion of a (pre)sheaf morphism can be
        extended to a notion of a (pre)sheaf morphism where the
        (pre)sheaves $\mathcal{F}$ and $\mathcal{G}$ have different
        algebraic structures. Definition
        \ref{definition:sheaf-modules} gives an example, confer also
        Remark \ref{remark:sheaf-modules}.
    \item It is clear that every (pre)sheaf $\mathcal{F}$ over $M$
        yields a (pre)sheaf $\mathcal{F}|_U$ over $U$ for all open
        subsets $U\in \underline{M}$ by restricting Definition
        \ref{definition:presheaves-sheaves} to open subsets of $U$.
    \end{enumerate}
\end{remark}

\begin{definition}[(Pre)sheaves of modules]
    \label{definition:sheaf-modules}
    %\cite[Chap.~II, Def.~1.8]{wells:1980a}
    Let $\mathcal{A}$ be a presheaf of rings and $\mathcal{E}$ be a
    (pre)sheaf of abelian groups over $M$. Then, $\mathcal{E}$ is
    called a \emph{(pre)sheaf of $\mathcal{A}$-modules} if for every
    $U\in\underline{M}$ the group $\mathcal{E}(U)$ is an
    $\mathcal{A}(U)$-module and for $V\subseteq U$ one has
    \begin{equation}
        \label{eq:sheaf-module-compatible}
        (a\cdot e)|_V = a|_V \cdot e|_V
    \end{equation}
    for all $a\in \mathcal{A}(U)$ and $e\in \mathcal{E}(U)$.
\end{definition}

\begin{remark}
    \label{remark:sheaf-modules}
    It is clear that the above notion of a sheaf of modules yields a
    sheaf homomorphism $\mathcal{A}\times \mathcal{E}\longrightarrow
    \mathcal{E}$ of (pre)sheaves of sets. But the converse is not true
    since the module structure is a nontrivial condition.
\end{remark}

\begin{proposition}[The sheaf of local morphisms]
    \label{proposition:sheaf-Hom}
    % \cite[Chap.~II, Sect.~1, Ex.~1.15]{hartshorne:1977a}
    Let $\mathcal{F}$ be a presheaf and $\mathcal{G}$ be a sheaf of
    sets over $M$ as in Definition
    \ref{definition:presheaves-sheaves}. Then, the presheaf of sets
    \begin{equation}
        \label{eq:sheaf-local-morphisms}
        U\longmapsto \Homs(\mathcal{F},\mathcal{G})(U)
        =\Hom(\mathcal{F}|_U, \mathcal{G}|_U) \quad 
        \textrm{for all } U\in \underline{M} 
    \end{equation}
    is a sheaf of sets which is called
    the \emph{sheaf of local morphisms} of $\mathcal{F}$ into
    $\mathcal{G}$ and denoted by $\Homs(\mathcal{F},\mathcal{G})$.
\end{proposition}

\begin{proof}
    By definition, the elements in $\Hom(\mathcal{F}|_U,
    \mathcal{G}|_U)$ are collections
    $\{h_W\}_{W\subseteq U}$ of maps. The restriction
    maps $r^U_V: \Hom(\mathcal{F}|_U, \mathcal{G}|_U) \longrightarrow
    \Hom(\mathcal{F}|_V, \mathcal{G}|_V)$ for $V\subseteq U$ are
    defined by
    \begin{equation}
        \label{eq:restriction-for-morphisms}
        r^U_V( \{ h_U \}_{W\subseteq U})= \{ h_W \}_{W\subseteq V}
    \end{equation}
    which obviously satisfy the conditions (a) %\ref{item:sheaf-1} 
    and (b). %\ref{item:sheaf-2}
    Now let $U=\Union_{i\in I} U_i$ as in
    Definition \ref{definition:presheaves-sheaves}. For (c)
    % \ref{item:sheaf-3} 
    we have to consider collections $\{h_W\}_{W\subseteq U}$ and
    $\{k_W\}_{W\subseteq U}$ with $\{h_W\}_{W\subseteq U_i}=
    \{k_W\}_{W\subseteq U_i}$ for all $i\in I$. Then, for $W\subseteq
    U$ and $f\in \mathcal{F}(W)$ one has for all $U_i$ with $U_i\cap
    W\neq \emptyset$
    \begin{equation*}
        h_W(f)|_{W\cap U_i}= h_{W\cap
          U_i}(f|_{W\cap U_i})= k_{W\cap
          U_i}(f|_{W\cap U_i}) = k_W(f)|_{W \cap
          U_i}. 
    \end{equation*}
    Since $\mathcal{G}$ is a sheaf and $f$ was arbitrary this yields
    $h_W= k_W$. Thus, $\{h_W\}_{W\subseteq U}= \{k_W\}_{W\subseteq U}$
    and (c) %\ref{item:sheaf-3}
    is shown. For (d) %\ref{item:sheaf-4}
    we have to consider collections $\{h^i_W\}_{W\subseteq U_i} \in
    \Hom(\mathcal{F}|_{U_i}, \mathcal{G}|_{U_i})$ for all $i\in I$
    with $\{h^i_W\}_{W\subseteq U_i\cap U_j}= \{h^j_W\}_{W\subseteq
      U_i\cap U_j}$. Then, for $W\subseteq U$ one can define a map
    $h_W: \mathcal{F}(W) \longrightarrow \mathcal{G}(W)$ by
    \begin{equation}
        \label{eq:sheaf-property-morphisms}
        h_W(f)|_{W\cap U_i}= h^i_{W\cap
          U_i}(f|_{W\cap U_i}) 
    \end{equation}
    for all $f\in \mathcal{F}(W)$ and $W\cap U_i\neq \emptyset$ which
    is possible since $\mathcal{G}$ is a sheaf and $h^i_{W\cap U_i\cap
      U_j}(f|_{W\cap U_i\cap U_j})=h^j_{W\cap U_i\cap U_j}(f|_{W\cap
      U_i\cap U_j})$.  Altogether, this leads to a collection of maps
    $\{h_W\}_{W\subseteq U}$ such that by definition $h_W=h^i_W$ for
    all $W\subseteq U_i$ and $i\in I$. Further, for $Z\subseteq W
    \subseteq U$ it follows
    \begin{equation*} 
        h_W(f)|_Z|_{Z\cap U_i}= h^i_{Z\cap
          U_i}(f|_{Z\cap U_i})=
        h_Z(f|_Z)|_{Z\cap U_i}
    \end{equation*}
    for $f\in \mathcal{F}(W)$ and $i\in I$. Again, the sheaf property
    of $G$ leads to $h_W(f)|_Z= h_Z(f|_Z)$, so $\{h_W\}_{W\subseteq
      U}$ is a sheaf morphism which restricts to the given local ones.
\end{proof}

\begin{remark}[Subsheaves of $\Homs(\mathcal{F},\mathcal{G})$]
    \label{remark:subsheaves-hom}
    If $\mathcal{F}$ and $\mathcal{G}$ have further algebraic
    properties that should be preserved by the maps of corresponding
    sheaf homomorphisms, this leads to a subsheaf of
    $\Homs(\mathcal{F},\mathcal{G})$, if the morphism defined in
    \eqref{eq:sheaf-property-morphisms} has the same properties as the
    locally given ones.
\end{remark}

\begin{example}[Subsheaves of $\Homs(\mathcal{F},\mathcal{G})$]
    \label{example:subsheaves-hom}
    Let $\mathcal{F}, \mathcal{G}$ be as in Proposition
    \ref{proposition:sheaf-Hom}.
    \begin{enumerate}
    \item If $\mathcal{F}$ is a presheaf and $\mathcal{G}$ is a sheaf
        of abelian groups the sheaf homomorphisms consisting of group
        homomorphisms build a subsheaf of
        $\Homs(\mathcal{F},\mathcal{G})$, since $h_W$ in
        \eqref{eq:sheaf-property-morphisms} is a group homomorphism.
    \item Let $\mathcal{E}_1,\ldots, \mathcal{E}_k$, $k\in
        \mathbb{N}$, be (pre)sheaves of abelian groups and let
        $\mathcal{F}= \mathcal{E}_1 \times \ldots \times
        \mathcal{E}_k$ be the (pre)sheaf of sets of the Cartesian
        product.  Completely analogously, the sheaf homomorphisms
        consisting of maps which are group morphisms in every entry
        build a subsheaf of $\Homs(\mathcal{F},\mathcal{G})$ which is
        denoted by $\Homab(\mathcal{E}_1\times \ldots \times
        \mathcal{E}_k,\mathcal{G})$.
    \end{enumerate}
\end{example}

\section{The sheaves of local and differential operators}
\label{sec:sheaves-local-differential-operators}

As seen above the notion of sheaves and the homomorphisms between them
provides the appropriate framework to understand the behaviour and the
relations between local and global information. In particular, one is
able to study an element $f\in \mathcal{F}(M)$ by considering the
restrictions $f|_U$ to open subsets of a covering of $M$. In many
typical situations, though, one is dealing with sheaves $\mathcal{F}$
and $\mathcal{G}$ over $M$ and is given a `global' map
$h_M:\mathcal{F}(M)\longrightarrow \mathcal{G}(M)$ of a particular
type which one would like to restrict as well and to realize as a
sheaf homomorphism. Thus, it is necessary to find criteria for the
possibility of reconstructing all the maps of a particular sheaf
morphism $\{h_U\}_{U\subseteq M}$ of a particular type from $h_M$.  As
we will see, this new point of view provides a very useful
characterization of the sheaf of local morphisms allowing a convenient
investigation of specific subsheaves. Basically, the observations and
proofs of this section are generalized versions of the considerations
in \cite[App.~A]{waldmann:2007a}.

In the following we always consider a smooth manifold $M$ as a
topological space and denote the sheaf of smooth functions with values
in $\mathbb{K}$ on it by $C^\infty_M$ which is a sheaf of associative
and commutative algebras.
 
\begin{lemma}[Open coverings]
     \label{lemma:open-coverings}
     Let $M$ be a smooth manifold and $V\subseteq M$ be an open subset.
     Then there exists an open covering $V=\Union_{i\in I} O_i$ with
     the property that for all $i\in I$
     \begin{enumerate}
     \item $O_i\subset O_i^\cl \subset V$,
     \end{enumerate}
     where $O_i^\cl$ denotes the closure of $O_i$. Additionally, one
     can achieve that
     \begin{enumerate}
         \setcounter{enumi}{1}
     \item there exist open sets $R_i\subset V$ with $O^\cl_i \subset R_i
         \subset R^\cl_i\subset V$.
     \end{enumerate}
 \end{lemma}
 
\begin{proof}
    Using centered charts $x_p:V\supset V_p\longrightarrow
    x_p(V_p) \subset \mathbb{R}^n$ of $V$ for every point $p\in V$,
    this means $x_p(p)=0$, the proof is obvious since there exist balls
    around $0\in \mathbb{R}^n$ such that the images under the
    homeomorphism $x_p^{-1}$ have the required properties.
\end{proof}

\begin{lemma}[Sheaves of left $C^\infty_M$-modules]
    \label{lemma:sheaves-and-functions}
    Let $\mathcal{E}$ be a sheaf of left $C^\infty_M$-modules over $M$
    and $V \subseteq U \subseteq M$ be open subsets.
    \begin{enumerate}
    \item Every smooth function $\chi \in C^\infty(U)$ with
        $\supp(\chi)\subseteq V$ defines a map
        \begin{equation}
            \label{eq:expansion-map}
            \chi: \mathcal{E}(V) \longrightarrow \mathcal{E}(U),
        \end{equation}
        which for all $e\in\mathcal{E}(V)$ is given by
        \begin{equation}
            \label{eq:definition-expansion}
            \begin{array}{rcl}
            (\chi e)|_V&=& \chi|_V \cdot e,\\
            (\chi e)|_{U\setminus \supp(\chi)}&=&0. 
            \end{array}
        \end{equation}
        In particular, for $e\in \mathcal{E}(U)$ we have
        \begin{equation}
            \label{eq:expansion-multiplication}
            \chi (e|_V) = \chi \cdot e.
        \end{equation}
    \item Let $U=\Union_{i\in I} U_i$ be an open covering,
        $\{\chi_i\in C^\infty(M)\}_{i\in I}$ be a subordinate
        partition of unity and $\{e_i\in \mathcal{E}(U_i)\}_{i\in I}$
        be a family of elements. Then there exists a unique element
        $e\in \mathcal{E}(U)$ with
        \begin{equation}
            \label{eq:glueing-together}
            e|_W=\sum_{
              \begin{subarray}{c}
                  i \in I\\
                  \chi_i|_W \neq 0
              \end{subarray}
            } 
            (\chi_i e_i)|_W
        \end{equation}
        for all open subsets $W\subseteq U$ of the open covering
        \begin{equation}
            \label{eq:covering-partition-of-unity}
            \{W\subseteq U \textrm{ open} \: | \: 
            \chi_i|_W = 0 \textrm{ for all 
              except at most finitely many } i \in I \}.   
        \end{equation} 
        This element is denoted by $e= \sum_{i \in I} \chi_i e_i$
        since this shows the relevant dependencies.\\
        In particular, for $e\in \mathcal{E}(U)$ this yields
        \begin{equation}
            \label{eq:sum-of-restriction}
            e= \sum_{i\in I}\chi_i (e|_{U_i}).
        \end{equation}
        % \item For $\chi\in C^\infty(U)$ with $\supp(\chi)\subseteq
        %     V$ and $\eta \in C^\infty(V)$ with $\supp(\eta)\in W$ we
        %     have
        %     \begin{equation}
        %         \label{eq:composition-expansion-map}
        %         \chi \circ \eta= \chi \eta : \mathcal{E}(W) \longrightarrow
        %         \mathcal{E}(U)
        %     \end{equation}
        %     where $\chi \eta\in C^\infty(U)$ with $\supp(\chi
        %     \eta)\subseteq \supp(\chi) \cap \supp(\eta) \subseteq W$ is
        %     the image of $\eta$ under $\chi$ in the case
        %     $\mathcal{E}=C^\infty_M$, that means
        %     \begin{eqnarray}
        %         \label{eq:comp-exp-map-1}
        %         \chi \eta|_V&=& \chi|_V \cdot \eta,\\
        %         \label{eq:comp-exp-map-2}
        %         \chi \eta|_{U\setminus \supp(\chi)}&=&0.
        %     \end{eqnarray}
    \end{enumerate}
    Now, let $O\subset R\subset V$ be open subsets as the
    ones in the covering of Lemma \ref{lemma:open-coverings}.
    \begin{enumerate}
        \setcounter{enumi}{2}
    \item There exists a smooth function $\chi\in C^\infty(U)$ with
        \begin{equation}
            \label{eq:testfunction-one}
            \supp(\chi)\subseteq V \quad \textrm{and} \quad
            \chi|_{O^\cl}=1.
        \end{equation}
        Then, for $e\in \mathcal{E}(V)$ one has
        \begin{equation}
            \label{eq:testfunction-one-sheaf}
            (\chi e)|_O=e|_O.
        \end{equation} 
    \item There exists a smooth function
        $\chi\in C^\infty(U)$ with
        \begin{equation}
            \label{eq:testfunction-zero}
            \chi|_{O^\cl}=0 \quad \textrm{and} \quad
            \chi|_{U\setminus R}=1.
        \end{equation}
        Then, for $e\in \mathcal{E}(U)$ with $e|_V=0$ one has
        \begin{equation}
            \label{eq:testfunction-zero-sheaf}
            \chi\cdot e=e.
        \end{equation}
    \end{enumerate}
\end{lemma}

\begin{proof}
    The first assertion holds due to the sheaf properties of
    $\mathcal{E}$. The last ones are a consequence of the Urysohn
    lemma, confer \cite[Prop.~2.26]{lee:2003a} or
    \cite[Cor.~A.1.5]{waldmann:2007a} and
    \cite[Chap.~7]{querenburg:2001a}, stating that two disjoint closed
    subsets $A_1,A_2\subseteq M$ can be separated by a smooth function
    $\chi\in C^\infty(M)$ with values in $[0,1]$ and the property that
    $\chi|_{A_1}=1$ and $\chi|_{A_2}=0$. The right hand side of
    \eqref{eq:glueing-together} is a finite sum and thus an element
    $e_W\in \mathcal{E}(W)$. Now consider $W\cap W'\neq
    \emptyset$. For $\chi_i$ with $\chi_i|_W\neq 0$ but
    $\chi_i|_{W'}=0$ one easily sees that $\chi_i e_i|_{W'}=0$ by
    checking it on $W'\cap U_i$ and $W'\cap (U\setminus \supp
    \chi_i)$. With this it follows
    \begin{equation*}
         e_W|_{W\cap W'}= \sum_{
          \begin{subarray}{c}
              i \in I \\
              \chi_i|_W \neq 0\\
              \chi_i|_{W'} \neq 0
          \end{subarray}
        }(\chi_i e_i)|_{W\cap W'} = e_{W'}|_{W\cap W'}.
    \end{equation*}
    The sheaf property of $\mathcal{E}$ then assures the existence of
    a unique element $e\in \mathcal{E}(U)$ with $e|_W=e_W$. Equation
    \eqref{eq:sum-of-restriction} is obvious. In the third part one
    considers the disjoint and closed sets $O^\cl$ and $U\setminus V$
    in the topology of the subspace $U$ and verifies that the
    functions $\chi\in C^\infty(U)$ with $\chi|_{O^\cl}=1$ and
    $\chi|_{U\setminus V}=0$ have the stated properties. Equation
    \eqref{eq:testfunction-one-sheaf} is obvious.  The existence of a
    function with \eqref{eq:testfunction-zero} follows because $O^\cl$
    and $U\setminus R$ are closed and disjoint sets. For
    \eqref{eq:testfunction-zero-sheaf} one considers the open cover of
    $U$ consisting of $V$ and $U\setminus R^\cl$ to check that the
    corresponding restrictions of $\chi\cdot e$ and $e$ coincide.
\end{proof}

\begin{definition}[Local map]
    \label{definition:local-maps}
    Let $\mathcal{E}_1, \ldots, \mathcal{E}_k$, $k\in\mathbb{N}$, and
    $\mathcal{G}$ be presheaves of abelian groups over $M$ and $U\in
    \underline{M}$. Then a map $D: \mathcal{E}_1(U) \times \ldots
    \times \mathcal{E}_k(U)\longrightarrow \mathcal{G}(U)$ which is a
    group morphism in every argument is said to be \emph{local} if for
    all $e_1\in \mathcal{E}_1(U), \ldots,e_k\in \mathcal{E}_k(U)$ and
    $V\subseteq U$ one has that
    \begin{equation}
        \label{eq:locality}
        D(e_1,\ldots, e_k)|_V= 0 \quad \textrm{if} \quad e_i|_V=0
        \quad  \textrm{for one } i \in \{1,\ldots, k\} .
    \end{equation}
    The set of all local maps on $U\in \underline{M}$ is denoted by
    \begin{equation}
        \Loc(\mathcal{E}_1, \ldots, \mathcal{E}_k; \mathcal{G})(U).
    \end{equation}
\end{definition}

\begin{remark}[Local map on functions]
    \label{remark:local-maps-on-functions}
    The presented notion of local maps is obviously consistent with
    the well-known one for operators $D:C^\infty(M)\longrightarrow
    C^\infty(M)$ on the functions of a manifold. There, $D$ is local if
    \begin{equation}
        \label{eq:local-operator-functions}
        \supp (D \chi) \subseteq \supp (\chi)
    \end{equation}
    for all $\chi\in C^\infty(M)$.
\end{remark}

\begin{corollary}
    \label{corollary:special-property-local-maps}
    Local maps $D$ as above have the useful property
    that for elements $e_i,e'_i \in \mathcal{E}_i(U)$ with
    $e_i|_V=e'_i|_V$ for all $i\in \{1,\ldots, k\}$ and $V\subseteq U$
    \begin{equation}
        \label{eq:local-maps-equal}
        D(e_1,\ldots, e_k)|_V= D(e'_1,\ldots, e'_k)|_V.
    \end{equation}
\end{corollary}

\begin{lemma}[Presheaf homomorphisms consist of local maps]
    \label{lemma:sheaf-homs-local}
    Let $\mathcal{E}_1, \ldots, \mathcal{E}_k$, $k\in\mathbb{N}$, and
    $\mathcal{G}$ be presheaves of abelian groups over $M$.  Then, the
    maps $h_U\in \{h_U\}_{U\subseteq M} \in \Homab(\mathcal{E}_1\times
    \ldots \times \mathcal{E}_k,\mathcal{G})(M)$ are local.
\end{lemma}

\begin{proof}
    The assertion is clear with $h_U(e_1,\ldots, e_k)|_V= h_V(e_1|_V,
    \ldots, e_k|_V)=0$.
\end{proof}

The following proposition shows that local maps and sheaf
homomorphisms are really the same for sheaves of $C^\infty_M$-modules.

\begin{proposition}[Sheaf homomorphisms and local maps]
    \label{proposition:sheaf-hom-loc}
    Let $\mathcal{E}_1, \ldots, \mathcal{E}_k$, $k\in\mathbb{N}$, be
    sheaves of $C^\infty_M$-modules and $\mathcal{G}$ be a sheaf of
    abelian groups over $M$. Then, every sheaf homomorphism
    $\{h_U\}_{U\subseteq M}\in \Homab(\mathcal{E}_1\times \ldots \times
    \mathcal{E}_k,\mathcal{G})(M)$ is already uniquely determined by the
    map $h_M: \mathcal{F}(M)\longrightarrow \mathcal{G}(M)$ and
    \begin{equation}
        \label{eq:iso-hom-loc-explizit}
        \{h_U\}_{U\subseteq M} \longmapsto h_M
    \end{equation}
    yields an isomorphism
    \begin{equation}
        \label{eq:iso-hom-loc}
        \Homab(\mathcal{E}_1\times
        \ldots \times \mathcal{E}_k,\mathcal{G})(M) 
        \cong \Loc(\mathcal{E}_1, \ldots, \mathcal{E}_k;
        \mathcal{G})(M)
    \end{equation}
    of abelian groups.
\end{proposition}

\begin{proof}
    Let $\{h_U\}_{U\subseteq M}$ and $\{h'_U\}_{U\subseteq M}$ be two
    sheaf homomorphisms in $\Homab(\mathcal{E}_1\times \ldots \times
    \mathcal{E}_k,\mathcal{G})(M)$ with $h_M=h'_M$.  For an open
    subset $V \subseteq M$ consider open subsets $O\subseteq V$ as the
    ones in the covering $\Union_{i \in I} O_i = V$ in
    Lemma~\ref{lemma:open-coverings}, this means with $O^\cl\subseteq
    V$, and corresponding functions $\chi \in C^\infty(M)$ as in the
    third part of Lemma~\ref{lemma:sheaves-and-functions}, this means
    with $\supp(\chi)\subseteq V$ and $\chi|_{O^\cl}=1$. Then,
    \begin{eqnarray*}
        h_V(e_1,\ldots,e_k)|_O &=&
        h_O((\chi e_1)|_O, \ldots, (\chi e_k)|_O)\\ 
        &=&h_M(\chi e_1,\ldots, \chi e_k)|_O\\
        &=& h'_V(e_1,\ldots,e_k)|_O 
    \end{eqnarray*}
    for all $e_1\in \mathcal{E}_1(V),\ldots, e_k\in \mathcal{E}_k(V)$.
    Since $\mathcal{G}$ is a sheaf, the two considered sheaf
    homomorphisms are equal. This shows that if a local map $D\in
    \Loc(\mathcal{E}_1, \ldots, \mathcal{E}_k; \mathcal{G})(M)$
    determines a sheaf homomorphism $\{D_U\}_{U\subseteq M}\in
    \Homab(\mathcal{E}_1\times \ldots \times
    \mathcal{E}_k,\mathcal{G})(M)$ with $D_M=D$, this is unique.
    Indeed, this is the case. For every $V\subseteq M$ a map $D_V$ can
    be defined by
    \begin{equation}
        \label{eq:recovering-of-sheaf-hom}
        D_V(e_1,\ldots, e_k)|_O = D(\chi e_1,\ldots,\chi e_k)|_O  
    \end{equation}
    for all elements $e_i\in \mathcal{E}_i(V)$ and all open subsets
    $O$ with corresponding functions $\chi$ as above.  The map $D_V:
    \mathcal{E}_1(V)\times \ldots \times \mathcal{E}_k(V)
    \longrightarrow \mathcal{G}(V)$ is well-defined by
    \eqref{eq:recovering-of-sheaf-hom} since $\mathcal{G}$ is a sheaf
    and since the right hand side does not depend on the choice
    $(O,\chi)$. If $(O',\chi')$ is another one it follows with $(\chi
    e_i)|_{O \cap O'} = (\chi' e_i)|_{O \cap O'}$ that $D(\chi
    e_1,\ldots, \chi e_k)|_{O \cap O'}=D(\chi' e_1,\ldots, \chi'
    e_k)|_{O \cap O'}$. In order to show that the so defined maps
    $D_V$ which are obviously group homomorphisms in every argument,
    define a sheaf homomorphism $\{D_V\}_{V\subseteq M} \in
    \Hom(\mathcal{F},\mathcal{G})$, let $V\subseteq U\subseteq M$ be
    fixed open subsets and $O\subseteq V$ as above. Further, let $\chi
    \in C^\infty(M)$ with $\supp(\chi)\subseteq V$ and $\chi|_O=1$.
    Then, the statement follows from
    \begin{eqnarray*}
        D_U(e_1,\ldots, e_k)|_V|_O 
        &=& D(\chi e_1, \ldots, \chi e_k)|_O\\
        &=& D(\chi(e_1|_V), \ldots, \chi(e_k|_V))|_O\\
        &=& D_V(e_1|_V, \ldots, e_k|_V)|_O.
    \end{eqnarray*}
    Due to the fact that $D_M=D$ and the proved uniqueness it is clear
    that this way every sheaf homomorphism $\{h_U\}_{U\subseteq M}$
    can be recovered from $h_M$. With the proved results and Lemma
    \ref{lemma:sheaf-homs-local} the map
    \eqref{eq:iso-hom-loc-explizit} apparently is an isomorphism.
\end{proof}

\begin{corollary}[The sheaf of local maps]
    \label{corollary:sheaf-local-maps}
    With the data of Proposition
    \ref{proposition:sheaf-hom-loc} the assignment
    \begin{equation}
        \label{eq:sheaf-local-maps}
        U \longmapsto
        \Loc(\mathcal{E}_1,\ldots,\mathcal{E}_k;\mathcal{G})(U)  
        \quad \textrm{for all} \quad U\in \underline{M}
    \end{equation}
    can be extended to a sheaf of abelian groups over $M$
    with unique restriction maps
    \begin{equation}
        \label{eq:rest-map-general}
        r^U_V:
        \Loc(\mathcal{E}_1,\ldots,\mathcal{E}_k;\mathcal{G})(U)  
        \longrightarrow 
        \Loc(\mathcal{E}_1,\ldots,\mathcal{E}_k;\mathcal{G})(V) 
    \end{equation}
    satisfying
    \begin{equation}
        \label{eq:property-rest-map}
        D(e_1,\ldots,e_k)|_V= (r^U_V D)
        (e_1|_V,\ldots,e_k|_V) 
    \end{equation}
    for all $D\in
    \Loc(\mathcal{E}_1,\ldots,\mathcal{E}_k;\mathcal{G})(U)$, $e_1\in
    \mathcal{E}_1(U),\ldots,e_k\in\mathcal{E}_k(U)$ and $V\subseteq
    U\subseteq M$.
\end{corollary}

\begin{proof}
    The statement is clear using the isomorphism
    \eqref{eq:iso-hom-loc} and Proposition \ref{proposition:sheaf-Hom}
    in order to define $\Loc(\mathcal{E}_1, \cdots, \mathcal{E}_k,
    \mathcal{G})$ as an isomorphic sheaf of $\Homs(\mathcal{E}_1\times
    \ldots \times \mathcal{E}_k, \mathcal{G})$. The induced restriction
    maps are directly seen to be the analogues of $D\longmapsto D_V$ in
    \eqref{eq:recovering-of-sheaf-hom}. Explicitly, for $V\subseteq U
    \subseteq M$ choose a covering of $V$ of open sets $O$
    satisfying $O^{\mathrm{cl}}\subseteq V$ and corresponding
    functions $\chi\in C^\infty(U)$ with $\supp(\chi)\subseteq V$ and
    $\chi|_O=1$. Then $r^U_V$ is defined by
     \begin{equation}
         \label{eq:definition-restriction-map}
         (r^U_V D)(e_1,\ldots,
         e_k)|_O= D(\chi e_1,\ldots,\chi 
         e_k)|_O 
     \end{equation} 
     for all $D\in \Loc(\mathcal{E}_1,
     \ldots,\mathcal{E}_k;\mathcal{G}) (U)$ and arbitrary elements
     $e_i\in \mathcal{E}_i(V)$.
\end{proof}

\begin{remark}
    \label{remark:local-maps-determined-by-restrictions}
    The first part in the proof of
    Proposition~\ref{proposition:sheaf-hom-loc} shows that the
    considered local operators are uniquely defined by their values on
    restricted elements.
\end{remark}

% \begin{proof}
%     Given two maps $D, D'\in \Loc(\mathcal{E}_1,\ldots,
%     \mathcal{E}_k;\mathcal{F}) (V)$ with $D(e_1|_V, \ldots, e_k|_V)=
%     D'(e_1|_V, \ldots, e_k|_V)$ for all $e_1\in
%     \mathcal{E}_1(U),\ldots, e_k\in \mathcal{E}_k(U)$ one has to
%     show that $D=D'$. Choose an open covering of $V$ with open sets
%     $O\subset V$ and corresponding functions $\chi\in C^\infty(U)$
%     as in the third part of Lemma \ref{lemma:sheaves-and-functions},
%     that means $\supp \chi\subseteq V$ and $\chi|_O=1$. Clearly,
%     $\chi|_V \in C^\infty(V)$ has the same properties. By Corollary
%     \ref{corollary:special-property-local-maps}, for $e_1\in
%     \mathcal{E}_1 (V), \ldots, e_k\in \mathcal{E}_k(V)$ one gets
%     $D(e_1,\ldots, e_k)|_O = D(\chi|_V \cdot e_1,\ldots, \chi|_V
%     \cdot e_k)|_O = D((\chi e_1)|_V, \ldots,(\chi e_k)|_V)|_O =
%     D'(e_1, \ldots, e_k)|_O$ and the assertion follows.
% \end{proof}

Corollary~\ref{corollary:sheaf-local-maps} is of course trivial and at
first sight it seems redundant to introduce the new notion of local
maps. But as we will see, this slightly different point of view is
very helpful in order to investigate sheaf homomorphisms consisting of
maps with further particular properties. This becomes clear with the
following important discussion of multidifferential maps and all later
applications thereof. First of all we need a generalized version of
the statement in \cite[Lemma~A.3.1]{waldmann:2007a}. 

\begin{lemma}
    \label{lemma:hom-local-diffops}
    Let $\mathcal{C}$ be a presheaf of associative and commutative
    algebras over a field $\mathbb{K}$ and $\mathcal{F}$ be a sheaf of
    $\mathcal{C}$-modules over $M$. Further let $V\subseteq M$ be an
    open subset, $\mathcal{E}_1,\ldots,\mathcal{E}_k$, $k\in
    \mathbb{N}$, be left $\mathcal{C}(V)$-modules and let $D\in
    \Hom_{\mathbb{K}} (\mathcal{E}_1,\ldots,\mathcal{E}_k;
    \mathcal{F}(V))$ be a $\mathbb{K}$-multilinear map. If there
    exist an open covering $V=\Union_{i\in I} O_i$ and differential
    operators $D_i\in \Diffop^L_{\mathcal{C}(V)} (\mathcal{E}_1,
    \ldots, \mathcal{E}_k; \mathcal{F}(V))$ of multiorder $L\in
    \mathbb{N}_0^k$ such that
    \begin{equation}
        \label{eq:hom-local-diffops}
        D(e_1,\ldots,e_k)|_{O_i}= D_i(e_1,\ldots, e_k)|_{O_i} \quad
        \textrm{for all} \quad i\in I, e_1\in \mathcal{E}_1, \ldots,
        e_k\in \mathcal{E}_k,
    \end{equation}
    then $D\in  \Diffop^L_{\mathcal{C}(V)} (\mathcal{E}_1, \ldots,
    \mathcal{E}_k; \mathcal{F}(V))$.
\end{lemma}

\begin{proof}
    The proof is an easy induction over $|L|$ using the formula
    $(D\circ \LL_c^{(l)}- \LL_c\circ D)(e_1,\ldots, e_k)|_{O_i}=
    (D_i\circ \LL_c^{(l)}- \LL_c\circ D_i)(e_1,\ldots, e_k)|_{O_i}$
    for all $c \in \mathcal{C}(V)$.
\end{proof}

% \begin{lemma}
%     \label{lemma:hom-local-diffops}
%     Let $\tilde{V}\subseteq P$ be an open subset,
%     $\mathcal{E}_1,\ldots,\mathcal{E}_k$ be $C^\infty(V)$ left modules
%     and $\mathcal{F}$ be a sheaf of $C^\infty_P$ left modules. Further
%     let $D\in \Hom_{\mathbb{K}} (\mathcal{E}_1,\ldots,\mathcal{E}_k;
%     \mathcal{F}(\tilde{V}))$ be a multilinear map. If there exist an
%     open covering $\tilde{V}=\Union_{i\in I} \tilde{O}_i$ and
%     differential operators $D_i\in
%     \Diffop^L_{C^\infty(V)}(\mathcal{E}_1, \ldots, \mathcal{E}_k;
%     \mathcal{F}(\tilde{V}))$ of multiorder $L\in \mathbb{N}_0^k$ such that
%     \begin{equation}
%         \label{eq:hom-local-diffops}
%         D(e_1,\ldots,e_k)|_{\tilde{O}_i}= D_i(e_1,\ldots, e_k)|_{\tilde{O}_i} \quad
%         \textrm{for all} \quad i\in I, e_1\in \mathcal{E}_1, \ldots,
%         e_k\in \mathcal{E}_k
%     \end{equation}
%     then $D\in  \Diffop^L_{C^\infty(V)} (\mathcal{E}_1, \ldots,
%     \mathcal{E}_k; \mathcal{F}(\tilde{V}))$
% \end{lemma}

% \begin{proof}
%     with the notation of (section or definition of algebraic
%     multidofferential operators) the proof is an easy induction over
%     $|L|$ using the formula $(D\circ L_a^{(l)}- L_a\circ
%     D)(e_1,\ldots, e_k)|_{\tilde{O}_i}= (D_i\circ L_a^{(l)}- L_a\circ
%     D_i)(e_1,\ldots, e_k)|_{\tilde{O}_i}$ for $a \in C^\infty(V)$.
% \end{proof}

\begin{proposition}[The sheaves of differential operators]
    \label{proposition:sheaf-differential-operators}
    Let $\mathcal{E}_1,\ldots, \mathcal{E}_k$, $k\in \mathbb{N}$, and
    $\mathcal{F}$ be sheaves of $C^\infty_M$-modules over $M$. For all
    open subsets $U\subseteq M$ and $L\in \mathbb{N}_0^k$ let
    \begin{equation}
        \label{eq:sheaf-diffop}
        \Diffop^L(\mathcal{E}_1,\ldots,\mathcal{E}_k;
        \mathcal{F})(U)=
        \Diffop^L_{C^\infty(U)}(\mathcal{E}_1(U),
        \ldots,\mathcal{E}_k(U);\mathcal{F}(U))   
    \end{equation}
    be the sets of algebraically defined multidifferential operators.
    \begin{enumerate}
    \item For all $U\subseteq M$ and $L\in \mathbb{N}_0^k$
        \begin{equation}
            \label{eq:Diffop-in-Loc}
            \Diffop^L(\mathcal{E}_1, \ldots, \mathcal{E}_k;
            \mathcal{F})(U) \subseteq
            \Loc(\mathcal{E}_1,\ldots,\mathcal{E}_k;
            \mathcal{F})(U).  
        \end{equation}
    \item For all $V\subseteq U \subseteq M$ and $L\in
        \mathbb{N}^k_0$
        \begin{equation}
            \label{eq:restriction-map-preserves-diffop}
            r^U_V: \Diffop^L(\mathcal{E}_1,
            \ldots, \mathcal{E}_k; \mathcal{F})(U)
            \longrightarrow \Diffop^L(\mathcal{E}_1, 
            \ldots, \mathcal{E}_k; \mathcal{F})(V).
        \end{equation}
    \item For every $L\in \mathbb{N}_0^k$ the assignment
        $\underline{M} \ni U\longmapsto \Diffop^L(\mathcal{E}_1,
        \ldots, \mathcal{E}_k; \mathcal{F})(U)$ together with the
        restriction maps $r^U_V$ in
        \eqref{eq:restriction-map-preserves-diffop}, yield a sheaf of
        $(C^\infty_M,C^\infty_M)$-bimodules, the so-called sheaf
        $\Diffop^L(\mathcal{E}_1, \ldots, \mathcal{E}_k;\mathcal{F})$
        of multidifferential operators from $\mathcal{E}_1\ldots,
        \mathcal{E}_k$ to $\mathcal{F}$ of multiorder $L$.
    \end{enumerate}
\end{proposition}

\begin{proof}
    As usual for algebraically defined multidifferential operators,
    the assertions are in the most cases proved by an induction over
    the absolute value $|L|\in \mathbb{N}_0$ of the multiindex $L$.
    \begin{enumerate}
    \item Let $D\in \Diffop^L(\mathcal{E}_1, \ldots, \mathcal{E}_k;
        \mathcal{F})(U)$, $V\subseteq U$ and $e_i\in
        \mathcal{E}_i(U)$, $i=1,\ldots,k$, with $e_l|_V=0$ for at
        least one certain $l\in\{1,\ldots,k\}$. Then we choose a
        covering of $V$ with open sets $O$ as in
        Lemma~\ref{lemma:open-coverings} and functions $\chi$ as in
        the fourth part of
        Lemma~\ref{lemma:sheaves-and-functions}. With the fact that
        $\chi\cdot e_l=e_l$ we have
        \begin{eqnarray*}
            D(e_1,\ldots, e_l, \ldots, e_k)|_O &=&  
            D(e_1,\ldots, \chi \cdot e_l, \ldots,
            e_k)|_O \\  
            &=&
            \label{eq:proof-diffop-in-local}
            \chi|_O \cdot D(e_1,\ldots, e_k)|_O + \left((D\circ 
                \LL_\chi^{(l)}-\LL_\chi\circ D) (e_1,\ldots,
                e_k)\right)|_O. 
        \end{eqnarray*}
        Due to the sheaf properties of $\mathcal{F}$ this leads to
        $D(e_1,\ldots,e_k)|_V=0$ by an induction over $|L|$ as
        mentioned above, since $(D \circ \LL_\chi^{(l)}-\LL_\chi \circ
        D)\in \Diffop^{L-\mathsf{e}_l} (\mathcal{E}_1, \ldots,
        \mathcal{E}_k; \mathcal{F}) (U)$ and $\chi|_O=0$.
    \item Let $V\subseteq U\subseteq M$, $D\in
        \Diffop^L(\mathcal{E}_1, \ldots, \mathcal{E}_k;
        \mathcal{F})(U)$, $e_i\in \mathcal{E}_i(V)$ for $i=1,\ldots,
        k$ and $a\in C^\infty(V)$. For the functions $\chi$ with
        respect to $O\subset V$ as in the definition
        \eqref{eq:definition-restriction-map} of the restriction map
        $r^U_V$ one has for $l\in \{1,\ldots,k\}$
        \begin{equation}
            \label{eq:quadrat-testfunction}
            \chi^2(a\cdot e_l)=(\chi a)\cdot (\chi e_l),
        \end{equation}
        since $(\chi^2(a\cdot e_l))|_V=\chi^2|_V\cdot a\cdot e_l=
        (\chi a)|_V \cdot (\chi e_l)|_V= ((\chi a)\cdot(\chi e_l))|_V$
        and the restrictions to $U\setminus \supp \chi$ vanish. With
        the fact that $D$ is local and $(\chi (a\cdot
        e_l))|_O=(\chi^2 (a \cdot e_l))|_O$ this leads to
        \begin{eqnarray}
            \label{eq:induktion-restriction-respects-diffop}
            \lefteqn{\left((r^U_V D) \circ
                  \LL_a^{(l)}-\LL_a\circ (r^U_V 
                  D) \right) 
              (e_1,\ldots, e_k)|_O} \nonumber \\
            &=& \left( D(\chi e_1, \ldots, \chi^2(a\cdot e_l),
                \ldots 
                , \chi e_k) - (\chi a) \cdot D(\chi e_1, \ldots, \chi
                e_k)\right)|_O \nonumber \\
            &=& \left( D(\chi e_1, \ldots, (\chi a)\cdot (\chi e_l),
                \ldots , \chi e_k) - (\chi a) \cdot D(\chi e_1,
                \ldots, \chi 
                e_k)\right)|_O\\
            &=& \left( D\circ \LL_{\chi a}^{(l)} - \LL_{\chi a}\circ
                D\right) (\chi e_1, \ldots, \chi
            e_k)|_O \nonumber \\ 
            &=& r^U_V \left( D\circ \LL_{\chi a}^{(l)} - \LL_{\chi a}  
                \circ D\right)(e_1, \ldots, e_k) |_O. \nonumber 
        \end{eqnarray}
        After checking the case $|L|=0$ the assertion $r^U_V D\in
        \Diffop^L(\mathcal{E}_1, \ldots, \mathcal{E}_k;
        \mathcal{F})(V)$ follows again by induction. If $r^U_V \left(
            D\circ \LL_{\chi a}^{(l)} - \LL_{\chi a} \circ D \right)
        \in \Diffop^{L-\mathsf{e}_l} (\mathcal{E}_1, \ldots,
        \mathcal{E}_k; \mathcal{F})(V)$,
        Equation~\eqref{eq:induktion-restriction-respects-diffop} and
        Lemma \ref{lemma:hom-local-diffops} show that $(r^U_V D) \circ
        \LL_a^{(l)} - \LL_a\circ (r^U_V D)\in \Diffop^{L-\mathsf{e}_l}
        (\mathcal{E}_1, \ldots, \mathcal{E}_k; \mathcal{F})(V)$.
    \item The sheaf axioms (a)-(c) follow from the fact that
        $\Loc(\mathcal{E}_1, \ldots, \mathcal{E}_k;\mathcal{F})$ is a
        sheaf and by the already proved statements. The only
        nontrivial consequence is the last sheaf axiom (d). Consider
        an open covering $U=\Union_{i\in I} U_i$ as in Definition
        \ref{definition:presheaves-sheaves} and $D_i \in
        \Diffop^L(\mathcal{E}_1, \ldots, \mathcal{E}_k;
        \mathcal{F})(U_i)$ with $D_i|_{U_i\cap U_j}=D_j|_{U_i\cap
          U_j}$. One has to show that the local map $D\in
        \Loc(\mathcal{E}_1, \ldots, \mathcal{E}_k;\mathcal{F})(U)$,
        defined by
        \begin{equation}
            \label{eq:last-axiom-local}
            D(e_1,\ldots,e_k)|_{U_i}=D_i(e_1|_{U_i},\ldots,
            e_k|_{U_i}) 
        \end{equation} 
        according to Equation~\eqref{eq:sheaf-property-morphisms},
        Proposition~\ref{proposition:sheaf-hom-loc} and
        Corollary~\ref{corollary:sheaf-local-maps}, is a differential
        operator $D\in \Diffop^L(\mathcal{E}_1, \ldots, \mathcal{E}_k;
        \mathcal{F})(U)$. Using the formula
        \begin{equation}
            \label{eq:last-axiom-local-diffops}
            (D\circ \LL_a^{(l)} - \LL_a \circ D) (e_1, \ldots,
            e_k)|_{U_i} = (\underbrace{D_i\circ \LL_{a|_{U_i}}^{(l)} -
              \LL_{a|_{U_i}}\circ D_i}_{H_i}) (e_1|_{U_i},\ldots,
            e_k|_{U_i}) 
        \end{equation}
        for $a\in C^\infty(U)$ this is again an easy induction over
        $|L|$ since one can show that the maps $H_i\in
        \Diffop^{L-\mathsf{e}_l} (\mathcal{E}_1,\ldots,
        \mathcal{E}_k;\mathcal{F}) (\tilde{V}_i)$ satisfy
        $H_i|_{\tilde{V}_i\cap \tilde{V}_j}=H_j|_{\tilde{V}_i\cap
          \tilde{V}_j}$.  Analogously to \eqref{eq:last-axiom-local},
        Equation~\eqref{eq:last-axiom-local-diffops} defines a unique
        element $H=D\circ \LL_a^{(l)} - \LL_a \circ D \in
        \Diffop^{L-\mathsf{e}_l} (\mathcal{E}_1,\ldots,
        \mathcal{E}_k;\mathcal{F}) (U)$ which is a multidifferential
        operator by the assumption of the induction. Finally, it is
        easy to verify that the left and right module structures $a
        \cdot D=\LL_a\circ D$ and $D \mathbin{\cdot^{(i)}} a= D \circ
        \LL_a^{(i)}$ for $i=1, \dots, N$ as in
        \eqref{eq:module-multi-linear-maps-short} with $a\in
        C^\infty_M(U)$ and $D\in \Diffop^L(\mathcal{E}_1, \ldots,
        \mathcal{E}_k; \mathcal{F})(U)$ for $U\subseteq M$ satisfy
        \begin{equation}
            \label{eq:module-structures-diffop-restrictions}
            (a\cdot D )|_{V} = a|_V \cdot D|_V \quad \textrm{and}
            \quad (D \mathbin{\cdot^{(i)}} a)|_V= D|_V
            \mathbin{\cdot^{(i)}} a|_V
        \end{equation}
        for all open subsets $V\subseteq U\subseteq M$.
    \end{enumerate}
\end{proof}

\begin{remark}
    \begin{enumerate}
    \item Proposition~\ref{proposition:sheaf-differential-operators}
        shows that the differential operators
        $\Diffop^L(\mathcal{E}_1, \ldots, \mathcal{E}_k;\mathcal{F})$
        of a finite order $L\in \mathbb{N}_0^k$ of differentiation are
        a subsheaf of $\Loc(\mathcal{E}_1, \ldots,
        \mathcal{E}_k;\mathcal{F})$. This way it is also clear that we
        have a new example of a subsheaf of $\Homs(\mathcal{E}_1\times
        \ldots \times \mathcal{E}_k,\mathcal{F})$ where all the maps
        of a sheaf homomorphism are corresponding differential
        operators.  
        % Note that this is not true in the general case
        % since one needs the explicit restriction map
        % \eqref{eq:restriction-map-preserves-diffop} to perform the
        % induction in the third part.
    \item Using the module structures
        \eqref{eq:module-multi-linear-maps-short}, the property
        \eqref{eq:module-structures-diffop-restrictions} implies that
        the sheaves of differential operators are sheaves of
        $C^\infty_M$-bimodules.
    \item Note that in general
        $\Diffop^\bullet(\mathcal{E}_1,\dots,\mathcal{E}_k,
        \mathcal{F})$ is only a presheaf since axiom (d) does not have
        to be satisfied. As a counterexample consider
        $\Diffop(C^\infty(\mathbb{R}), C^\infty(\mathbb{R}))$. Let
        $\chi_0\in C^\infty(\mathbb{R})$ be a function with
        $\supp(\chi_0) \subseteq (-1,1)$ and define $\chi_i\in
        C^\infty(\mathbb{R})$ by $\chi_i(x)=
        \chi_0(x+2i)+\chi_0(x-2i)$ for all $i\in \mathbb{N}$. Then let
        $U_i=(-2i-1,2i+1)$ and consider the differential operators
        $D_i\in \Diffop^i(C^\infty(U_i), C^\infty(U_i))$ given by
        $D_i=\sum_{l=0}^i \chi_l|_{U_i} \frac{\partial^l}{\partial
          x^l}$. Then there exists no $D\in
        \Diffop(C^\infty(\mathbb{R}), C^\infty(\mathbb{R}))$ with
        $D|_{U_i}=D_i$ since the order of $D$ would be infinite.
    \end{enumerate}
\end{remark}

In later applications we will make use of the following lemma
concerning sheaf endomorphisms consisting of differential maps.

\begin{lemma}
    For a sheaf $\mathcal{E}$ of left $C^\infty_M$-modules and $l\in
    \mathbb{N}_0$ let $\mathcal{D}^l=
    \Diffop^l(\mathcal{E};\mathcal{E})$ be the corresponding sheaf of
    differential operators. Then, for all $l,k \in \mathbb{N}_0$ the
    usual composition of maps $ \circ: \mathcal{D}^l(U) \times
    \mathcal{D}^k(U) \longrightarrow \mathcal{D}^{l+k}(U)$ for all
    $U\subseteq M$ defines a sheaf homomorphism
    \begin{equation}
        \label{eq:composition-of-maps-diffop-sheaf}
        \circ:
        \mathcal{D}^l \times \mathcal{D}^k \longrightarrow
        \mathcal{D}^{l+k}.
    \end{equation}
    This means that for $V\subseteq U\subseteq M$ and $D\in
    \mathcal{D}^l(U), D' \in \mathcal{D}^k(U)$ one has
    \begin{equation}
        \label{eq:composition-of-maps-diffop-sheaf-explicit}
        (D \circ D')|_V = D|_V \circ D'|_V.
    \end{equation}
\end{lemma}

\begin{proof}
    Equation~\eqref{eq:composition-of-maps-diffop-sheaf-explicit} is a
    direct consequence of \eqref{eq:property-rest-map} and
    Remark~\ref{remark:local-maps-determined-by-restrictions}.
\end{proof}

\section{Sheaves over different manifolds}
\label{sec:sheaves-different-manifolds}

Now we investigate the situation where the (pre)sheaves do not have
the same base manifold as underlying topological space. For this
purpose we introduce two simple definitions whose general versions can
be found in \cite[Chap.~II]{godement:1964a} and
\cite[Chap.~II]{hartshorne:1977a}. If $\pp: P \longrightarrow M$ is an
open map between two smooth manifolds we always use the notation
$U=\pp(\tilde{U})$ to denote open subsets $\tilde{U}\subseteq P$ and
the corresponding open images $U\subseteq M$.

\begin{lemma}[Pullback of sheaves]
    \label{lemma:pull-back-sheaf}
    Let $\pp: P\longrightarrow M$ be an open map between two smooth
    manifolds, this means $\pp(\tilde{U})\in \underline{M}$ for all
    $\tilde{U}\in \underline{P}$, and let $\mathcal{G}$ be a
    (pre)sheaf over $M$. Then the assignment
    \begin{equation}
        \label{eq:pull-back-sheaf}
        \underline{P} \ni \tilde{U} \longmapsto
        (\pp^*\mathcal{G}) (\tilde{U}):= \mathcal{G}(\pp(\tilde{U}))
    \end{equation}
    together with the restriction maps
    \begin{equation}
        \label{eq:pull-back-restriction}
        r^{\tilde{U}}_{\tilde{V}}:=
        r^{\pp(\tilde{U})}_{\pp(\tilde{V})}:
        (\pp^*\mathcal{G}) (\tilde{U}) \longrightarrow
        (\pp^*\mathcal{G}) (\tilde{V})
    \end{equation}
    yields a presheaf which we call \emph{pullback} $\pp^*\mathcal{G}$
    of $\mathcal{G}$ with $\pp$.
\end{lemma}

The proof of this lemma is obvious. 

\begin{remark}[Pullback and inverse image]
    \begin{enumerate}
    \item Note that, in general, the pullback of a sheaf as defined
        above is only a presheaf. An open covering
        $\tilde{U}=\Union_{i\in I} \tilde{U}_i$ induces an open
        covering $\pp(\tilde{U})=\Union_{i\in I}
        \pp(\tilde{U}_i)$. Then it is easy to verify that axiom (c) is
        still fulfilled. This is not always the case for axiom (d)
        since $\pp(\tilde{U}\cap \tilde{V})\subseteq
        \pp(\tilde{U})\cap \pp(\tilde{V})$ can be a proper subset, in
        which cases it is not possible to trace back property (d) of
        $\pp^*\mathcal{G}$ to the property (d) of $\mathcal{G}$.
    \item If $\pp$ is continuous instead of open there is a way to
        define a sheaf $\pp^{-1}\mathcal{G}$ over $P$ out of
        $\mathcal{G}$ which is called the \emph{inverse image} of
        $\mathcal{G}$, confer \cite[Chap.~II,
        Sect.~1.12]{godement:1964a} and \cite[Chap.~II,
        Sect.~1]{hartshorne:1977a}. This definition, though, requires
        some more technical preliminaries and since we do not need the
        sheaf property for our purposes we have only introduced the
        above notion of the pullback.
    \end{enumerate}
\end{remark}

\begin{lemma}[Direct image of sheaves]
     \label{lemma:direct-image-sheaf}
     Let $\pp: P\longrightarrow M$ be a continuous map between two
     smooth manifolds and let $\mathcal{G}$ be a
     (pre)sheaf over $P$. Then the assignment
    \begin{equation}
        \label{eq:direct-image-sheaf}
        \underline{M} \ni U \longmapsto
        (\pp_*\mathcal{G})(U):= \mathcal{G}(\pp^{-1}(U)) 
    \end{equation}
    together with the restriction maps
    \begin{equation}
        \label{eq:direct-image-restriction}
        r^U_V:=
        r^{\pp^{-1}(U)}_{\pp^{-1}(V)}:
        (\pp_*\mathcal{G})(U) \longrightarrow (\pp_*\mathcal{G})(V) 
    \end{equation}
    yields a (pre)sheaf, the so-called \emph{direct image}
    $\pp_*\mathcal{G}$ of $\mathcal{G}$ under $\pp$.
\end{lemma}

\begin{proof}
    Clearly, an open covering $U=\Union_{i\in I} U_i$ induces an open
    covering $\pp^{-1}(U)=\Union_{i\in I} \pp^{-1}(U_i)$. Then, the
    proof of (a), (b) and (c) is again obvious. For (d) one needs
    $\pp^{-1}(U\cap V)= \pp^{-1}(U)\cap \pp^{-1}(V)$.
\end{proof}

\begin{remark}
    If $\pp$ is open and continuous it follows that
    $\pp_*(\pp^*\mathcal{G})=\mathcal{G}$ since
    $\pp(\pp^{-1}(U))=U$. But since $\pp^{-1}(\pp(\tilde{U}))
    \supseteq \tilde{U}$ one has $\pp^*(\pp_* \mathcal{G})\neq
    \mathcal{G}$ in general.
\end{remark}

\begin{lemma}
\label{lemma:presheaf-homs-pull-back-direct-image}
    Let $\pp: P\longrightarrow M$ be an open and continuous map as
    above. Further let $\mathcal{F}$ be a (pre)sheaf of sets over $M$ and
    $\mathcal{G}$ be a sheaf of sets over $P$.
    Then there is a canonical isomorphism
    \begin{equation}
        \label{eq:isomorphism-sheaf-hom-pull-back-direct-image}
        \Hom(\mathcal{F}, \pp_*\mathcal{G})\cong \Hom(\pp^*
        \mathcal{F}, \mathcal{G})
    \end{equation}
    between the corresponding presheaf homomorphisms of presheaves
    over $M$ and $P$.
\end{lemma}

\begin{proof}
    For $\{h_U\}_{U\in \underline{M}}\in \Hom(\mathcal{F},
    \pp_*\mathcal{G})$ we define $\{h_{\tilde{U}}\}_{\tilde{U}\in
      \underline{P}}$ by $h_{\tilde{U}}:=
    r^{\pp^{-1}(\pp(\tilde{U}))}_{\tilde{U}}\circ
    h_{\pp(\tilde{U})}$. By assumption and the presheaf property of
    $\mathcal{G}$ this defines an element in $\Hom(\pp^*
    \mathcal{F}, \mathcal{G})$. Conversely, such a presheaf
    homomorphism defines an element in $\Hom(\mathcal{F},
    \pp_*\mathcal{G})$ with the definition
    $h_U:=h_{\pp^{-1}(U)}$. It is obvious that both constructions are
    inverse to each other.
\end{proof}

\begin{remark}
    Since $\mathcal{G}$ and $\pp_*\mathcal{G}$ are sheaves one is able
    to consider the sheaves of local morphisms
    $\Homs(\pp^*\mathcal{F},\mathcal{G})$ over $P$ or
    $\Homs(\mathcal{F},\pp_*\mathcal{G})$ over $M$.
    % Das ist kein Widerspruch zum gezeigten Lemma, da der Isomorphismus
    % je fuer eine feste Abbildung $\pp: \tilde{U}\longrightarrow U$ gilt. 
\end{remark}

\begin{example}
    \label{example:presheaf-morphism-functions}
    With $\pp: P \longrightarrow M$ as above there is a presheaf
    homomorphism $h=\{h_{\tilde{U}}\}_{\tilde{U}\in \underline{P}}:
    \pp^*C^\infty_M \longrightarrow C^\infty_P$ with
    \begin{equation}
        \label{eq:presheaf-morphism-functions}
        h_{\tilde{U}}:= (\pp|_{\tilde{U}})^*:
        (\pp^*C^\infty_M) (\tilde{U})= C^\infty(\pp(\tilde{U}))
        \longrightarrow C^\infty(\tilde{U}). 
    \end{equation}
    Since all $h_{\tilde{U}}$ are inclusions and with
    Lemma~\ref{lemma:presheaf-homs-pull-back-direct-image},
    $C^\infty_M$ is a subsheaf of $\pp_*C^\infty_P$.
\end{example}

\begin{remark}[Sheaves over different manifolds]
    \label{remark:generalization-sheaves-on-different-spaces}
    The assertions of Proposition~\ref{proposition:sheaf-hom-loc} and
    all derived results, now stated for sheaves over $P$, can be
    reformulated and extended to the situation where some
    $\mathcal{F}\in \{\mathcal{E}_1,\dots,\mathcal{E}_k\}$ are
    presheaves $\mathcal{F}= \pp^*\mathcal{E}$ of
    $\pp^*C^\infty_M$-modules over $P$ coming from sheaves
    $\mathcal{E}$ of $C^\infty_M$-modules over $M$.

    In order to find adequate coverings of $\tilde{V}\subseteq
    \tilde{U}\subseteq P$ with subsets $\tilde{O}$ and functions
    $\chi\in C^\infty(U)$ as in the crucial
    Equations~\eqref{eq:recovering-of-sheaf-hom},
    \eqref{eq:definition-restriction-map} and
    \eqref{eq:proof-diffop-in-local} one considers the corresponding
    open coverings with subsets $O$ of $V=\pp(\tilde{V})\subseteq M$
    and functions $\chi$ and uses the induced covering of
    $\pp(\tilde{V})$ with the subsets
    \begin{equation}
        \label{eq:induced-covering-total-space}
        \tilde{O}= \pp^{-1}(O)\cap \tilde{V}.
    \end{equation}
    Then, $\pp(\tilde{O})=O$ and the proofs and definitions can be
    made in a completely analogous way.

    Proposition~\ref{proposition:sheaf-differential-operators} has an
    analogous reformulation. In this case of multidifferential
    operators of multiorder $L\in\mathbb{N}_0^k$, $k\in \mathbb{N}$,
    let $\mathcal{E}_1,\ldots, \mathcal{E}_k$ be sheaves of
    $C^\infty_M$-modules over $M$ and $\mathcal{F}$ be a sheaf of
    $C^\infty_P$-modules over $P$. Using the presheaf homomorphism of
    Example~\ref{example:presheaf-morphism-functions}, $\mathcal{F}$
    inherits a $\pp^*C^\infty_M$-module structure. Then the sheaf
    \begin{equation}
        \label{eq:sheaf-diffop-over-P}
        \Diffop^L_{\pp^*C^\infty_M} (\pp^*\mathcal{E}_1, \ldots,
        \pp^*\mathcal{E}_k; \mathcal{F})
    \end{equation}
    of $(C^\infty_P, \pp^* C^\infty_M)$-bimodules over $P$ is
    well-defined in the obvious way. By use of Lemma
    \ref{lemma:presheaf-homs-pull-back-direct-image} the given
    structures also yield a sheaf 
    \begin{equation}
        \label{eq:sheaf-diffop-over-M}
        \Diffop^L_{C^\infty_M}
        (\mathcal{E}_1, \ldots, \mathcal{E}_k; \pp_*\mathcal{F})
    \end{equation}
    of $(\pp_*C^\infty_P, C^\infty_M)$-bimodules over $M$ which is the
    direct image of the afore mentioned sheaf.
\end{remark}

% \begin{lemma}
%     \label{lemma:local-form-of-Diffops}
%     Let $\tilde{U}\subseteq P$ be an open subset such that for
%     $U=\pp(\tilde{U})$ there is a local chart $(U,x)$ of $M$. Further,
%     let there be given $\phi\in \Diffop^L(C^\infty(M),
%     \Diffop^l (C^\infty(P)))$ with multiorders of differentiation $L
%     = (l_1, \ldots, l_k)\in \mathbb{N}_0^k$ and $l\in
%     \mathbb{N}_0$. Then there exist uniquely defined differential
%     operators $\phi_{\tilde{U}}^{\alpha_1 \cdots \alpha_k}\in \Diffop^l
%     (C^\infty(\tilde{U}))$ with multiindices $\alpha_i\in
%     \mathbb{N}_0^n$, $i=1,\ldots, k$, such that
%     \begin{equation}
%         \label{eq:local-form-of-Diffop}
%         \phi|_{\tilde{U}} (a_1,\ldots,a_k) =
%         \sum_{|\alpha_1|\le l_1,\ldots, 
%           |\alpha_k|\le l_k}
%         \left(\frac{\partial^{|\alpha_1|}a_1} {\partial x^{\alpha_1}}
%             \cdots 
%             \frac{\partial^{|\alpha_k|}a_k}{\partial
%               x^{\alpha_k}}\right) \cdot \phi_{\tilde{U}}^{\alpha_1 \cdots \alpha_k}
%     \end{equation}
%     for all $a_1,\ldots, a_k\in C^\infty(U)$, where
%     $|_{\tilde{U}}$ denotes the restriction of differential
%     operators as in Appendix~\ref{cha:sheaves}.
% \end{lemma}
% 
% The proof is a straightforward generalization of the well-known
% considerations for ordinary differential operators and can be found in
% Appendix ???.
% 
% \begin{proof}
%     The proof is an induction over 
% \end{proof}

\section{Sheaves and invariant differential Hochschild complexes}
\label{sec:sheaves-Hochschild}

After the previous general considerations we now come to the
applications with respect to the deformation theory of algebras and
modules. In the present section we investigate under which
circumstances the concepts of sheaf theory can be used for the task of
computing Hochschild cohomologies in the ($G$-invariant) differential
setting introduced in the Sections~\ref{sec:G-invariant-types} and
\ref{sec:differential-algebra-module-structures}. Besides discussing
and presenting the necessary framework the main goal is to formulate
the important Propositions~\ref{proposition:sheaves-and-cohomology}
and \ref{proposition:sheaves-and-cohomology-G-inv}.

The basic setting is the following. Let $\mathcal{\mathcal{A}}$ be a
sheaf of $\mathbb{K}$-algebras and $\mathcal{E}$ be a sheaf of
$\mathbb{K}$-vector spaces and right $\mathcal{A}$-modules over a
topological space $M$. This assignment of right modules to all open
subsets $U\subseteq M$ then leads to the induced assignment
\begin{equation}
    \label{eq:assignment-Hochschild-complex}
    U\longmapsto
    \HC^\bullet(\mathcal{A}(U), \End_{\mathbb{K}}
    (\mathcal{E}(U),\mathcal{E}(U)))     
\end{equation}
of the corresponding Hochschild complexes. According to the basic
intention of sheaf theory the natural question then is if these
Hochschild complexes are related with each other and if the know\-ledge
of the global cohomology is equivalent to that of the local ones. If
this is the case the assertion clearly reveals a natural approach of
computing the global cohomology by considering the possibly simpler
local problem.

Having this aim in mind a first step is to figure out when the
assignment \eqref{eq:assignment-Hochschild-complex} can be extended to
a presheaf of vector spaces such that the family of Hochschild
differentials is a presheaf homomorphism. For the ($G$-invariant)
differential Hochschild complexes this will turn out to be the
appropriate guideline to find the necessary general frameworks where
all aspired issues are given.

\subsection{Sheaves and differential Hochschild complexes}
\label{subsec:sheaves-diff}

First of all we investigate Hochschild complexes
$\HCdiff^\bullet(\mathcal{A},\mathcal{D})$ of the differential type as
in Corollary~\ref{corollary:Differential-Hochschild-module}. In order
to define the vector spaces of this complex we need a set
$\{\mathcal{A},\mathcal{E}, \mathcal{C},\mathcal{B}, \gamma\}$ of
\emph{purely algebraic structures} as occurring in
Definition~\ref{definition:module-type-diff}. Now assume that
$\mathcal{A},\mathcal{E},\mathcal{C}$, and $\mathcal{B}$ are sheaves
of corresponding structures over a topological space $P$ and that
$\gamma: \mathcal{C}\longrightarrow \mathcal{B}$ is a sheaf
homomorphisms consisting of algebra homomorphisms. Due to the results
of Section~\ref{sec:sheaves-local-differential-operators}, especially
Proposition~\ref{proposition:sheaf-differential-operators}, the
induced assignments of differential operators to open subsets are
presheaves if the corresponding presheaves of commutative algebras, in
our case $\mathcal{C}$ and $\mathcal{B}$, are presheaves of functions
over some smooth manifold.

\begin{remark}[A convenient framework for differential Hochschild
    complexes]
    \label{remark:structures-sheaves-differential}
    With respect to the later applications where we have to work with
    sheaves over different topological spaces we assume from now on to
    be given the following structures.
    \begin{enumerate}
    \item A surjective submersion $\pp:P\longrightarrow M$ between
        smooth manifolds $P$ and $M$ inducing
        \begin{enumerate}
        \item the sheaf $\mathcal{C}=C^\infty_M$ of smooth functions
            on $M$ which is a sheaf of associative and commutative
            algebras over $M$.
        \item the sheaf $\mathcal{B}=C^\infty_P$ of smooth functions
            on $P$ which is a sheaf of associative and commutative
            algebras over $P$.
        \item the presheaf morphism $\gamma:
            \pp^*\mathcal{C}\longrightarrow \mathcal{B}$ which is
            given by the pullback of functions, this means
            $\gamma_{\tilde{U}}=(\pp|_{\tilde{U}})^*$ for all open
            $\tilde{U}\subseteq P$, confer
            Example~\ref{example:presheaf-morphism-functions}.
        \end{enumerate}
    \item A sheaf $\mathcal{A}$ of $\mathbb{K}$-algebras with
        multiplication $\mu$ and left $\mathcal{C}$-modules over $M$.
    \item A sheaf $\mathcal{E}$ of $\mathbb{K}$-vector spaces with a
        left $\mathcal{B}$-module structure $l$ and a right
        $\pp^*\mathcal{A}$-module structure $\rho$ over $P$.
    \end{enumerate}
\end{remark}

The structures $\mathcal{E}$ and $\mathcal{B}$ immediately induce the
sheaves
%\begin{equation}
%    \label{eq:sheaf-differential-Hochschild-algebra}
$\mathcal{D}^l=
\Diffop^l_{\mathcal{B}}(\mathcal{\mathcal{E}},\mathcal{E}) $
%\end{equation}
of $\mathbb{K}$-vector spaces and $\mathcal{B}$-modules over $P$ for
all $l\in \mathbb{N}_0$.  Since the map $\pp:P\longrightarrow M$ is
open it is possible to consider the pullback $\pp^*\mathcal{C}$. Using
the sheaf morphism $\gamma$ in the sense of
\eqref{eq:induced-module-diffop},
$\mathcal{D}^l=\Diffop^l_{\mathcal{B}}(\mathcal{E},\mathcal{E})$ is a
sheaf of $\pp^*\mathcal{C}$-modules over $P$. The pullbacks
$\pp^*\mathcal{C}$ and $\pp^*\mathcal{A}$ only are presheaves over $P$
but as stated in
Remark~\ref{remark:generalization-sheaves-on-different-spaces} this
still guarantees that one has the sheaves
\begin{equation}
    \label{eq:sheaf-differential-Hochschild}
    \HC^{k,L,l} =
    \Diffop^L_{\pp^*\mathcal{C}} 
    (\pp^*\mathcal{A},\dots, \pp^*\mathcal{A};
    \Diffop^l_{\mathcal{B}}(\mathcal{E},\mathcal{E}))
\end{equation}
of $\mathcal{B}$-modules over $P$ for all $L\in\mathbb{N}_0^k$, $k \in
\mathbb{N}$, and $l\in \mathbb{N}_0$.

Like $\tilde{U}\longmapsto \mathcal{D}(\tilde{U})=\Union_{l\in
  \mathbb{N}_0}\mathcal{D}^l(\tilde{U})$, the assignment
\begin{equation}
    \label{eq:assignment-Hochschild-diff}
    \tilde{U} \longmapsto 
    \HCdiff^k(\pp^*\mathcal{A},\mathcal{D}) (\tilde{U})= \Union_{L\in
      \mathbb{N}_0^k} \Union_{l\in \mathbb{N}_0}
    \HC^{k,L,l}(\tilde{U}) 
\end{equation}
of differential operators to open subsets $\tilde{U}\subseteq P$ only
induce presheaves of $\mathbb{K}$-vector spaces and
$\mathcal{B}$-modules. The restriction maps with respect to open
subsets $\tilde{V}\subseteq \tilde{U}\subseteq P$ satisfy
\begin{equation}
    \label{eq:restrictions-Hochschild-values}
    [\phi(a_1,\dots, a_k)(e)]|_{\tilde{V}}
    = \phi(a_1,\dots, a_k)|_{\tilde{V}} (e|_{\tilde{V}})
    = \phi|_{\tilde{V}}(a_1|_{\pp(\tilde{V})},\dots,
    a_k|_{\pp(\tilde{V})}) (e|_{\tilde{V}})
\end{equation}
for all $\phi\in
\HC^\bullet(\pp^*\mathcal{A},\mathcal{D})(\tilde{U})$,
$a_1,\dots,a_k\in \mathcal{A}(\pp(\tilde{U}))$, and $e\in
\mathcal{E}(\tilde{U})$.

In the same way one finds the presheaf
%\begin{equation}
%    \label{eq:Hochschild-presheaf-algebra}
$\HCdiff^\bullet(\mathcal{A},\mathcal{A})$
%\end{equation}
of $\mathcal{C}$-modules over $M$ which is derived from the sheaf
$\mathcal{A}$ of algebras over $M$. For this presheaf it is clear that
the local cup products, insertions after the $i$-th position and
Hochschild differentials are well-defined on the differential
subcomplexes only if the family $\mu= \{\mu_U\}_{U\subseteq M}$ of
algebra structures consists of differential structures.

Due to Proposition~\ref{proposition:module-type-diff} the cup products
as in \eqref{eq:cup-product-module} are always well-defined for
\eqref{eq:assignment-Hochschild-diff}. But with the defining
Equation~\eqref{eq:Hochschild-differential-module} one further sees
that in this case the local Hochschild differentials are well-defined
endomorphisms of the differential operators only if both the family
$\mu$ and the family $\rho=\{ \rho_{\tilde{U}} \}_{\tilde{U}\subseteq
  P}$ of right module structures consist of differential
elements. This can of course be checked globally since the families
$\mu$ and $\rho$ are sheaf homomorphisms $\mu: \mathcal{A}\times
\mathcal{A} \longrightarrow \mathcal{A}$ and $\rho: \mathcal{E}\times
\pp^*\mathcal{A}\longrightarrow \mathcal{E}$.

\begin{lemma}[Sheaves of differential algebras and modules]
    \label{lemma:differential-structures-sheaf}
    Let there be given structures as in
    Remark~\ref{remark:structures-sheaves-differential} such that the
    algebra multiplication $\mu_M$ of $\mathcal{A}(M)$ and the right
    $(\pp^*\mathcal{A})(P)=\mathcal{A}(M)$-module structure $\rho_P$
    of $\mathcal{E}(P)$ are of the differential type as in the
    Definitions~\ref{definition:algebra-type-diff} and
    \ref{definition:module-type-diff}.

    Then, all algebra structures of the family
    $\mu=\{\mu_U\}_{U\subseteq M}$ and all right module structures of
    the family $\rho=\{\rho_{\tilde{U}}\}_{\tilde{U}\subseteq P}$ are
    differential with the same degrees of differentiation as the
    global ones.
\end{lemma}

\begin{proof}
    As a consequence of Proposition~\ref{proposition:sheaf-hom-loc},
    Corollary~\ref{corollary:sheaf-local-maps}, and
    Remark~\ref{remark:generalization-sheaves-on-different-spaces} one
    finds $\mu_U=\mu_M|_U$ and $\rho_{\tilde{U}}=\rho_P|_{\tilde{U}}$.
    Then, the presheaf property of the considered differential
    cochains yields the statement.
%%%
%%% ähnlicher Beweis: 
%%%
%     The algebra and module structures are local maps in the sense of
%     Definition~\ref{definition:local-maps} and further uniquely defined
%     by the values on restrictions according to
%     Lemma~\ref{lemma:local-maps-determined-by-restrictions}. Due to
%     the assumption and the presheaf property of
%     \eqref{eq:sheaf-differential-Hochschild-algebra} and
%     \eqref{eq:sheaf-differential-Hochschild} the global maps can be
%     restricted. This yields $\mu_U=\mu_M|_U$ and
%     $\rho_{\tilde{U}}=\rho_P|_{\tilde{U}}$ and the assertion follows.
%     Note that, in principle, the assertion could have been deduced in
%     terms of presheaf homomorphisms using the analogues of
%     Proposition~\ref{proposition:sheaf-hom-loc} and
%     Corollary~\ref{corollary:sheaf-local-maps}.
\end{proof}

\begin{remark}
    The proof in particular includes the following general
    observation. If the Hochschild complexes of a particular type
    induce presheaves of vector spaces and if the global structures
    are of this type, the same is true for the local ones. The
    converse assertion of
    Lemma~\ref{lemma:differential-structures-sheaf} is also true since
    cochains of a fixed order of differentiation are sheaves.
\end{remark}

The general property \eqref{eq:property-rest-map} implying
\eqref{eq:restrictions-Hochschild-values} now shows that
\eqref{eq:assignment-Hochschild-diff} really defines presheaves of
Hochschild complexes if $\mu_M$ and $\rho_P$ are differential. The
concrete meaning thereof becomes clear in the following lemma.

\begin{lemma}[Presheaf of complexes]
    \label{lemma:presheaf-of-complexes}
    Let there be given structures as in
    Remark~\ref{remark:structures-sheaves-differential} such that
    $\mu_M$ and $\rho_P$ are of the differential type. Then, the
    presheaf $\HCdiff^\bullet(\pp^*\mathcal{A},\mathcal{D})$ is a
    presheaf of Hochschild complexes over $P$ which means that the
    family $\{\delta_{\tilde{U}}\}_{\tilde{U}\subseteq P}$ of
    Hochschild differentials is a presheaf endomorphism. Thus one has
    \begin{equation}
        \label{eq:restriction-Hochschild-differential}
        \delta_{\tilde{V}}\circ r^{\tilde{U}}_{\tilde{V}} =
        r^{\tilde{U}}_{\tilde{V}}\circ  \delta_{\tilde{U}}
    \end{equation}
    for all open subsets $\tilde{V}\subseteq \tilde{U}\subseteq P$,
    where $r^{\tilde{U}}_{\tilde{V}}$ denotes the restriction map.  In
    particular this induces the sheaf homomorphisms $\delta:
    \HC^{k,L,l}\longrightarrow \HC^{k,L,l}$.

    Moreover, one has corresponding (pre)sheaf homomorphisms given by
    the cup product $\cup$ and the insertions $\circi$ after the
    $i$-th position with respect to the presheaf
    $\pp^*\HCdiff^\bullet(\mathcal{A},\mathcal{A})$ having analogue
    morphisms.
\end{lemma}

\begin{proof}
    Using the properties \eqref{eq:property-rest-map} and
    \eqref{eq:composition-of-maps-diffop-sheaf-explicit} of the
    restriction maps, all assertions can be proved by restricting the
    defining equations of all structures since the considered
    differential cochains are uniquely defined by the values on
    restrictions.
\end{proof}

Note the different terminology used in some references, confer
\cite[Chap.~2]{wells:1980a}. What here is called a sheaf of complexes
is there referred to as differential sheaf. We avoid the latter since
it would lead to some confusion with the notion of differential
operators.

\begin{corollary}[The sheaves of cocycles]
    \label{corollary:sheaf-cocycles}
    In the situation of Lemma~\ref{lemma:presheaf-of-complexes}
    equation \eqref{eq:restriction-Hochschild-differential}
    immediately implies that the Hochschild cocycles build subsheaves
    $Z^{k,L,l}\subseteq \HC^{k,L,l}$.
\end{corollary}

The fact that the presheaves $\HC^{k,L,l}$ really are sheaves of
$\mathcal{B}=C^\infty_P$-modules has an important consequence for the
later applications since it is possible to obtain global cochains out
of local ones as explained in the two first parts of
Lemma~\ref{lemma:sheaves-and-functions}. So, for all open subsets
$\tilde{V}\subseteq \tilde{U}\subseteq P$ every smooth function
$\tilde{\chi}\in C^\infty(\tilde{U})$ with $\supp
\tilde{\chi}\subseteq \tilde{V}$ induces a map $\tilde{\chi}:
\HC^{k,L,l}(\tilde{V})\longrightarrow \HC^{k,L,l}(\tilde{U})$. By its
defining Equation~\eqref{eq:definition-expansion} and the properties
\eqref{eq:restriction-Hochschild-differential} and
\eqref{eq:Hochschild-delta-left-multiplication} of the Hochschild
differentials it is easy to see that
\begin{equation}
    \label{eq:Hochschild-differential-expansion}
    \delta_{\tilde{U}} \circ \tilde{\chi} = \tilde{\chi}\circ
    \delta_{\tilde{V}} \quad 
    \textrm{if} \quad l_\rho=0,
\end{equation}
this means if $\mathcal{E}$ is a sheaf of
$(\mathcal{B},\pp^*\mathcal{A})$-bimodules.
Equation~\eqref{eq:Hochschild-differential-expansion} is clearly
verified by checking it on the subsets $\tilde{V}$ and
$\tilde{U}\setminus \supp \tilde{\chi}$. In the case of $l_\rho=0$ the
cocycles $Z^{k,L,l}$ are even subsheaves of
$C^\infty_P$-modules. Further, this condition has an important
consequence for the coboundaries.

\begin{lemma}[The sheaves of coboundaries]
    \label{lemma:coboundaries-sheaf}
    Let $\mu_M$ and $\rho_P$ be differential with $l_\rho=0$. Then,
    the presheaves of Hochschild coboundaries $B^{k,L,l}\subseteq
    \HC^{k,L,l}$ are subsheaves of $C^\infty_P$-modules if there exist
    $\tilde{L}\in \mathbb{N}_0^{k-1}$ and $\tilde{l}\in \mathbb{N}_0$
    such that
    \begin{equation}
        \label{eq:coboundary-image-of-finite-order}
        B^{k,L,l}=\delta \HC^{k-1,\tilde{L},\tilde{l}}.
    \end{equation}
\end{lemma}

\begin{proof}
    The only nontrivial point to check is axiom (d) in Definition
    \ref{definition:presheaves-sheaves}. So let $\Union_{i\in I}
    \tilde{U}_i= \tilde{U}\subseteq P$ be an open covering and let
    $\phi_i= \delta_{\tilde{U}_i} \Theta_i \in B^{k,L,l}(\tilde{U}_i)$
    be coboundaries with $\phi_i|_{\tilde{U}_i\cap \tilde{U}_j}=
    \phi_j|_{\tilde{U}_i\cap \tilde{U}_j}$ and where $\Theta_i \in
    \HC^{k-1,\tilde{L},\tilde{l}} (\tilde{U}_i)$. Since $\HC^{k,L,l}$
    is a sheaf there exists a global cochain $\phi\in
    \HC^{k,L,l}(\tilde{U})$ with $\phi|_{\tilde{U}_i}= \phi_i$ and it
    remains to show that $\phi$ is a coboundary. Now, choose a
    subordinate partition of unity $\{\tilde{\chi}\}_{i\in I}$ and
    consider the global object $\Theta=\sum_{i\in I} \tilde{\chi}_i
    \Theta_i\in \HC^{k-1,\tilde{L},\tilde{l}} (\tilde{U})$ as
    explained in Lemma~\ref{lemma:sheaves-and-functions}. Then the
    assertion follows from
    \begin{equation*}
        \delta_{\tilde{U}}\Theta= \delta_{\tilde{U}} \sum_{i\in I}
        \tilde{\chi}_i \Theta_i
        = \sum_{i\in I} \tilde{\chi}_i(\delta_{\tilde{U}_i} \Theta_i)
        = \sum_{i\in I} \tilde{\chi}_i (\phi|_{\tilde{U}_i})
        =\phi,
    \end{equation*}
    where the last equation is nothing but
    \eqref{eq:sum-of-restriction} and
    \begin{equation}
        \label{eq:Hochschild-differential-global-sum}
        \delta_{\tilde{U}} \sum_{i\in I} \tilde{\chi}_i \Theta_i
        = \sum_{i\in I} \tilde{\chi}_i(\delta_{\tilde{U}_i} \Theta_i)
    \end{equation}
    follows from the defining Equation~\eqref{eq:glueing-together} and
    \eqref{eq:Hochschild-differential-expansion}.
\end{proof}

Lemma~\ref{lemma:coboundaries-sheaf} in particular states that in any
case one has the sheaves $\delta \HC^{k,\tilde{L},\tilde{l}}$. The
given proof provides an important tool for the computation of
Hochschild cohomologies. If the coboundaries $B^{k,L,l}$ are a sheaf
two global cocycles are in the same cohomology class if and only if
the same is true for the restrictions to any open subset. This in
particular means that the global cohomology can be determined if the
knowledge of the local ones is given. With regard to the applications
we formulate these result in the following proposition.

\begin{proposition}[Global and local differential Hochschild
    cohomology]
    \label{proposition:sheaves-and-cohomology}
    Let there be given structures as in
    Remark~\ref{remark:structures-sheaves-differential} such that
    $\mu_M$ and $\rho_P$ are differential with $l_\rho=0$ which means
    that $\mathcal{E}$ is a $(\mathcal{B}, \pp^*
    \mathcal{A})$-bimodule. Further, let $\phi,\psi\in \HC^{k,L,l}(P)$
    be global cocycles, $\delta_P \phi=0=\delta_P \psi$, and let
    $\Union_{i\in I}\tilde{U}_i=P$ be an open covering.  Then the
    following two assertions are equivalent:
    \begin{enumerate}
    \item $\phi$ and $\psi$ are in the same cohomology class,
        $\phi-\psi\in B^{k,L,l}(P)$.
    \item All local cocycles $\phi|_{\tilde{U}_i}$ and
        $\psi|_{\tilde{U}_i}$ are equivalent,
        $(\phi-\psi)|_{\tilde{U}_i}\in B^{k,L,l}(\tilde{U}_i)$, and
        there exist $\tilde{L}\in \mathbb{N}_0^{k-1}$, $\tilde{l}\in
        \mathbb{N}_0$ such that $(\phi-\psi)|_{\tilde{U}_i} =
        \delta_{\tilde{U}_i} \Theta_i$ with cochains $\Theta_i\in
        \HC^{k-1,\tilde{L},\tilde{l}}(\tilde{U}_i)$ for all $i\in I$.
    \end{enumerate}
    In particular, if the second statement is true every subordinate
    partition of unity $\{\tilde{\chi}_i\}_{i\in I}$ with
    $\tilde{\chi}\in C^\infty(P)$ induces an element
    \begin{equation}
        \label{eq:global-potential}
        \Theta=\sum_{i\in I} \tilde{\chi}_i \Theta_i \in
        \HC^{k-1,\tilde{L},\tilde{l}}(P) \quad \textrm{with} \quad  
        \phi-\psi= \delta_P \Theta. 
    \end{equation}
\end{proposition}

\begin{proof}
    The proof is obvious and a special case of the proof of
    Lemma~\ref{lemma:coboundaries-sheaf} for the family
    $\{(\phi-\psi)|_{\tilde{U}_i}\}_{i\in I}$.
\end{proof}

\begin{remark}[Quotients of sheaves]
    \label{remark:sheaves-quotients}
    If \eqref{eq:coboundary-image-of-finite-order} holds the above
    equivalence of global and local statements shows that the
    equivalence relation defining the cohomology is an equivalence
    relation on the sheaves $Z^{k,L,l}$ in the sense of \cite[Chap.~II,
    Sect.~1.9]{godement:1964a} and \cite[Chap.~II,
    Sect.~1]{hartshorne:1977a}. The equivalence classes
    $Z^{k,L,l}/\delta(\HC^{k-1,\tilde{L},\tilde{l}})$ with appropriate
    orders of differentiation and the full Hochschild cohomology
    $\HHdiff^\bullet(\mathcal{A},\mathcal{D})$ can always be equipped
    with a canonical presheaf structure. But in general, these are no
    sheaf structures. Axiom (c) in Definition
    \ref{definition:presheaves-sheaves} still holds but Axiom (d) does
    not have to be satisfied.
\end{remark}

\begin{remark}
    Note that one has analogous statements for complexes and
    cohomologies induced by algebra structures.
\end{remark}

\subsection{Sheaves and $G$-invariant differential Hochschild complexes}
\label{subsec:sheaves-diff-G-inv}

Now we consider (sub)sheaves $\mathcal{F}$ of algebraic structures
\emph{with a compatible action of a Lie group $G$} over some manifold
$M$. By this we clearly mean that for any open subset $U\subseteq M$
the group $G$ acts on $\mathcal{F}(U)$ by structure preserving maps
and that this action is compatible with the restriction maps. Thus one
could define the action of $G$ on a sheaf $\mathcal{F}$ as a group
homomorphism $G\longrightarrow \Aut(\mathcal{F})$ into the structure
preserving sheaf automorphisms. Note that without loss of generality
one is able to consider left actions.

The situation for functions on a manifold $P$ shows how typical
examples look like and how the previously discussed situation has to
be adapted to the $G$-invariant case. Usually, the left action on the
functions comes from a right action $\rr:P\times G \longrightarrow P$
on the underlying space $P$ and is given by pullback. For all subsets
$\tilde{U}\subseteq P$ consisting of orbits one has a well-defined
action on the functions $f \in C^\infty(\tilde{U})$ by
\begin{equation}
    \label{eq:action-functions-pullback}
    g \acts f= f \circ \rr_g
\end{equation}
for all $g \in G$ where $\rr_g(u)=\rr(u,g)$ for all $u\in P$. In the
aspired case of a principal fibre bundle $\pp:P\longrightarrow M$ with
structure group $G$ these invariant subsets are nothing but the
preimages $\pp^{-1}(U)$ of the projection with $U\subseteq M$. This
shows that the direct image
\begin{equation}
    \label{eq:direct-image-functions}
    \pp_* C^\infty_P
\end{equation}
is a sheaf of associative and commutative algebras over the base
manifold $M$ with a compatible action of $G$. The $G$-invariant
functions are nothing but the pullbacks of functions on the base
space $M$. This leads to the following important lemma.

\begin{lemma}[Principal fibre bundles and sheaves with $G$-actions]
    \label{lemma:sheaves-principal-G-action}
    Let $\pp : P \longrightarrow M$ be a principal fibre bundle with
    structure group $G$. Further let $\mathcal{E}_1, \ldots,
    \mathcal{E}_k, \mathcal{F}$ be sheaves of $C^\infty_P$-modules
    over $P$ such that the sheaves $\pp_* \mathcal{E}_1, \ldots,
    \pp_*\mathcal{E}_k, \pp_*\mathcal{F}$ of $\pp_*
    C^\infty_P$-modules over $M$ are equipped with a compatible
    $G$-action.

    Then, for all $L\in \mathbb{N}_0^k$ the sheaf
    \begin{equation*}
        \pp_* \mathcal{D}^L= \pp_* \Diffop^L(\mathcal{E}_1,
        \ldots, \mathcal{E}_k; \mathcal{F}) 
    \end{equation*}
    has a compatible $G$-action given by
    \eqref{eq:representation-multi-linear-maps} and the $G$-invariant
    multidifferential operators build a subsheaf
    $(\pp_*\Diffop^L(\mathcal{E}_1, \ldots, \mathcal{E}_k;
    \mathcal{F}))^G$.
\end{lemma}

\begin{proof}
    One has to show the compatibility of the restriction maps $r^U_V$
    of $\pp_*\mathcal{D}^L$ this means of the restrictions
    $r^{\pp^{-1}(U)}_{\pp^{-1}(V)}$ of $\mathcal{D}^L$ with the
    induced $G$-action, where $V\subseteq U\subseteq M$ are open
    subsets. To this end, choose an open cover of $V$ with sets
    $O\subseteq V$ as in Lemma \ref{lemma:open-coverings} and
    corresponding functions $\chi$ as in the third part of Lemma
    \ref{lemma:sheaves-and-functions}. Then it is clear that
    $\pp^{-1}(V)\subseteq \pp^{-1}(U)$ is covered by the sets
    $\pp^{-1}(O)\subseteq \pp^{-1}(V)$ and that the functions $\pp^*
    \chi \in C^\infty(\pp^{-1}(U))= \pp_*C^\infty_P(U)$ are
    $G$-invariant and satisfy $\supp \pp^*\chi \subseteq
    \pp^{-1}(\supp \chi) \subseteq \pp^{-1}(V)$ and
    $(\pp^*\chi)|_{\pp^{-1}(O)}=1$. For these functions, the induced
    maps $\pp^*\chi: \mathcal{E}_i(\pp^{-1}(V))\longrightarrow
    \mathcal{E}_i(\pp^{-1}(U))$, $i=1,\dots,k$, are $G$-invariant
    since by assumption $(g \acts ((\pp^*\chi) e_i))|_{\pp^{-1}(V)} =
    g \acts ((\pp^*\chi)|_{\pp^{-1}(V)} \cdot e_i) = (g \acts
    \pp^*\chi)|_{\pp^{-1}(V)} \cdot (g\acts e)= ((\pp^*\chi)(g \acts
    e))|_{\pp^{-1}(V)}$, and since the restriction to
    $\pp^{-1}(U)\setminus \supp \pp^*\chi$ vanishes. With this
    property the defining
    Equation~\eqref{eq:definition-restriction-map} for the restriction
    maps of $\mathcal{D}^L$ leads to
    \begin{eqnarray*}
        (g \acts D|_{\pp^{-1}(V)})(e_1,\ldots, e_k)|_{\pp^{-1}(O)} 
        &=& \left( g \acts (D|_{\pp^{-1}(V)}(g^{-1} \acts e_1, \ldots, g^{-1} \acts
            e_k)) \right)|_{\pp^{-1}(O)}\\
        &=& \left( g \acts (D ((\pp^*\chi)(g^{-1} \acts e_1), \ldots, (\pp^*\chi)
            (g^{-1} \acts e_k))) \right)|_{\pp^{-1}(O)}\\
        &=& (g \acts D)|_{\pp^{-1}(V)} (e_1, \ldots, e_k)|_{\pp^{-1}(O)}
    \end{eqnarray*}
    for all $D\in \Diffop^L(\mathcal{E}_1,\dots, \mathcal{E}_k;
    \mathcal{F})$ from which $g \acts D|_{\pp^{-1}(V)}= (g \acts
    D)|_{\pp^{-1}(V)}$ follows. This compatibility also shows that the
    invariant differential operators build a subsheaf. The last sheaf
    axiom is satisfied since the global element $D$ with
    $D|_{\pp^{-1}(V)}= D_{\pp^{-1}(V)}$ for the locally given and
    $G$-invariant $D_{\pp^{-1}(V)}$ is again $G$-invariant due to $(g
    \acts D)|_{\pp^{-1}(V)}= g \acts D_{\pp^{-1}(V)}=
    D|_{\pp^{-1}(V)}$.
\end{proof}

With the above results it is obvious how the assumptions of the
previous subsection have to be specified in the $G$-invariant case.

\begin{remark}[A framework for $G$-invariant differential Hochschild
    complexes]
    \label{remark:structures-sheaves-differential-G-inv}
    We now assume to be given the following structures.
    \begin{enumerate}
    \item A principal fibre bundle $\pp:P\longrightarrow M$ with structure
        group $G$ inducing
        \begin{enumerate}
        \item the sheaf $\mathcal{C}=C^\infty_M$ of smooth functions over
            $M$ with trivial $G$-action.
        \item the compatible action of $G$ on the sheaf
            $\pp_*\mathcal{B}=\pp_*C^\infty_P$ of associative, commutative
            algebras over $M$ given by the pullback of functions with
            respect to the principal action.
        \item the sheaf morphism $\gamma: \mathcal{C}\longrightarrow
            (\pp_*\mathcal{B})^G$ onto the sheaf of $G$-invariant
            functions given by the pullback of functions, this means
            $\gamma_{U}=(\pp|_{\pp^{-1}(U)})^*$ for all open subsets
            $U\subseteq M$.
        \end{enumerate}
    \item A sheaf $(\mathcal{A},\mu)$ of $\mathbb{K}$-algebras over
        $M$ which is a sheaf of $\mathbb{K}$-vector spaces and
        $\mathcal{C}$-modules with a compatible $G$-action.
    \item A sheaf $\mathcal{E}$ of $\mathbb{K}$-vector spaces with a
        left $\mathcal{B}$-module structure $l$ and a right
        $\pp^*\mathcal{A}$-module structure $\rho$ over $P$ such that
        $\pp^*\mathcal{E}$ is a sheaf of vector spaces and
        $\mathcal{B}$-modules with a compatible $G$-action.
    \end{enumerate}
\end{remark}

Since the $G$-action shall be taken into account we only work with
(pre)sheaves over $M$ this time, if necessary by considering the
direct images with respect to the projection $\pp$. With
$\mathcal{A}=\pp_*(\pp^*\mathcal{A})$ we can always use
Lemma~\ref{lemma:sheaves-principal-G-action}. Thus one gets that
$\HCdiff^\bullet(\mathcal{A},\mathcal{A})$ and the direct image
$\pp_*\HCdiff^\bullet(\pp^*\mathcal{A},\mathcal{D})=
\HCdiff^\bullet(\mathcal{A},\pp_*\mathcal{D})$ are presheaves of
Hochschild complexes and $C^\infty_M$-modules over $M$ with compatible
$G$-actions. In particular, one has the sheaves of $G$-invariant
multidifferential operators
\begin{equation}
    \label{eq:sheaf-Hochschild-algebra-G-inv}
    \Diffop^L_{\mathcal{C}}(\mathcal{A},\dots,\mathcal{A};\mathcal{A})^G
\end{equation}
and
\begin{equation}
    \label{eq:sheaf-Hochschild-module-G-inv}
    \left(\pp_*\HC^{k,L,l}\right)^G = \Diffop^L_{\mathcal{C}}
    (\mathcal{A},\dots,\mathcal{A};
    \pp_*\Diffop^l_{\mathcal{B}}(\mathcal{E},\mathcal{E}))^G.
\end{equation}
If the action on $\mathcal{A}$ is trivial one gets presheaves
$\HCdiff^\bullet(\mathcal{A},(\pp_*\mathcal{D})^G)$ consisting of sheaves
\begin{equation}
    \label{eq:sheaf-Hochschild-module-G-inv-special}
    \left(\pp_*\HC^{k,L,l}\right)^G = \Diffop^L_{\mathcal{C}}
    (\mathcal{A},\dots,\mathcal{A};
    (\pp_*\Diffop^l_{\mathcal{B}}(\mathcal{E},\mathcal{E}))^G).
\end{equation}

Note that the $G$-invariance of the families $\mu$ and $\rho$ is given
by assumption. If $\mu_M$ and $\rho_P$ are differential one thus has
analogue results to those of Lemma~\ref{lemma:presheaf-of-complexes}
and Corollary~\ref{corollary:sheaf-cocycles} for the $G$-invariant
cochains and cocycles. The $G$-invariant coboundaries of the above
complexes are given by the presheaves
\begin{equation}
    \label{eq:coboundaries-G-inv}
    (\pp_*B^k)^G= \delta\left(
        \HCdiff^{k-1}(\mathcal{A},\pp_*\mathcal{D})^G \right).
    % \quad \textrm{or} \quad
    % (\pp_*B^k)^G= \delta \HC^{k-1}(\mathcal{A},\pp_*\mathcal{D}^G).
\end{equation}
The adapted version of Lemma~\ref{lemma:coboundaries-sheaf} then
states that all $\delta((\pp_*\HC^{k,\tilde{L},\tilde{l}})^G)
\subseteq \HCdiff^k(\mathcal{A},\pp_*\mathcal{D})$ are subsheaves. In
the corresponding proof one considers open coverings $\Union_{i\in
  I}U_i= M$ and the structures assigned to $\pp^{-1}(U_i)$. Choosing a
subordinate partition of unity $\{\chi_i\}_{i\in I}$ with functions
$\chi_i\in C^\infty(M)$ the global object $\Theta=\sum_{i\in
  I}(\pp^{-1}\chi_i) \Theta_i$ then is $G$-invariant if the local
$\Theta_i$ are $G$-invariant since $\pp^{-1}\chi_i\in
C^\infty(P)^G$. These considerations finally yield the following
refined version of
Proposition~\ref{proposition:sheaves-and-cohomology}.

\begin{proposition}[Global and local $G$-invariant cohomology]
    \label{proposition:sheaves-and-cohomology-G-inv}
    Let there be given structures as in
    Remark~\ref{remark:structures-sheaves-differential-G-inv} such
    that $\mu_M$ and $\rho_P$ are differential with
    $l_\rho=0$. Further, let $\phi,\psi\in (\HC^{k,L,l}(P))^G$ be
    global $G$-invariant cocycles, $\delta_P \phi=0=\delta_P \psi$,
    and let $\Union_{i\in I} U_i=M$ be an open covering of the base
    manifold $M$. Then, the following two assertions are equivalent:
    \begin{enumerate}
    \item $\phi$ and $\psi$ are in the same $G$-invariant cohomology
        class, $\phi-\psi\in (B^{k,L,l}(P))^G$.
    \item All local cocycles $\phi|_{\pp^{-1}(U_i)}$ and
        $\psi|_{\pp^{-1}(U_i)}$ are equivalent,
        $(\phi-\psi)|_{\pp^{-1}(U_i)}\in (B^{k,L,l}
        (\pp^{-1}(U_i)))^G$, and there exist $\tilde{L}\in
        \mathbb{N}_0^{k-1}$, $\tilde{l}\in \mathbb{N}_0$ such that
        $(\phi-\psi)|_{\pp^{-1}(U_i)} = \delta_{\pp^{-1}(U_i)}
        \Theta_i$ with cochains $\Theta_i\in
        (\HC^{k-1,\tilde{L},\tilde{l}} (\pp^{-1}(U_i)))^G$ for all
        $i\in I$.
    \end{enumerate}
    In particular, if the second statement is true every subordinate
    partition of unity $\{ \chi_i \}_{i\in I}$ with $\chi_i\in
    C^\infty(M)$ induces an element
    \begin{equation}
        \label{eq:global-potential-G-inv}
        \Theta=\sum_{i\in I} (\pp^{-1}\chi_i) \Theta_i \in
        (\HC^{k-1,\tilde{L},\tilde{l}}(P))^G \quad \textrm{with} \quad 
        \phi-\psi= \delta_P \Theta. 
    \end{equation}
\end{proposition}

\chapter{Cohomology and projective resolutions}
\label{cha:bar-koszul}

After the rather simple preparations of the last chapter the aim of
the present one is to develop the crucial techniques of cohomology
computation which are motivated by general results of homological
algebra, in particular the theory of projective resolutions of modules
$\mathcal{A}$ over rings. The fundamental definitions and basic
results from homological algebra can be found in
\cite[Chap.~6]{jacobson:1989a} and are summarized in
Appendix~\ref{cha:homological-algebra}. The basic idea which is the
motivation for all subsequent considerations of this chapter is the
following observation concerning derived functors.

Consider a ring $\mathsf{R}$, an $\mathsf{R}$-module $\mathcal{A}$ and
an additive functor $F$ from the category $\mathsf{R}$-$\mathbf{mod}$
of $\mathsf{R}$-modules into the category $\mathbf{Ab}$ of abelian
groups. Without loss of generality we consider a contravariant
functor. Then, the application of the $k$-th (right) derived functor
$R^kF$ to $\mathcal{A}$ by definition is a cohomology group $H^k(FC)$
with an arbitrary projective resolution $(C,\epsilon)$ of
$\mathcal{A}$, confer Definition~\ref{definition:resolution}. The
functor is well-defined since this definition does not depend on the
choice of the resolution. This has an immediate consequence for the
task of computing the cohomology $H^\bullet$ of a given complex. If it
is possible to find a ring $\mathsf{R}$, a functor $F$ and a
resolution $(C,\epsilon)$ of a module $\mathcal{A}$ as above such that
the cohomology groups of interest are isomorphic to the images of
$\mathcal{A}$ under the derived functors, $H^k\cong R^kF(\mathcal{A})=
H^k(F(C))$, the freedom in the choice of the resolution can be used to
find a different complex $F(C')$, clearly with the same cohomology,
but where the problem is easier to handle.

This simple observation yields the guideline for the whole chapter.
The first section describes the well-known fact that the algebraic
Hochschild complexes as in
Definition~\ref{definition:Hochschild-complex} are of the above form
which is seen with the algebraic bar resolution. This algebraic
situation and the general ideas explained above will then serve as a
guideline for the more specific situation where
$\mathcal{A}=C^\infty(V)$ is an algebra of functions and where the
Hochschild complexes are of a particular type, namely the continuous
and differential versions. For them we have to introduce the
topological bar resolution in order to see that the new cohomologies
are given in a similar way as the purely algebraic ones. After these
preliminary steps we present the corresponding topological version of
the Koszul resolution and prove that this new resolution really
provides a different way to compute the considered differential
Hochschild cohomologies. To this end, the necessary chain
homomorphisms and homotopy maps will be given explicitly in order to
ensure that the additional property of the cochains to be differential
is respected. The obtained results are already stated in the
publication \cite{bordemann.neumaier.waldmann.weiss:2007a:pre} and
have their origin in \cite{bordemann.et.al:2005a:pre}. In the present
work all this is presented in a self-contained way containing the
explicit and refined proofs.

\section{The Hochschild cohomology as an Ext group}
\label{sec:hochschild-Ext}

The purely algebraic Hochschild cohomology defined in
Definition~\ref{definition:Hochschild-complex} is a well-known example
of a cohomology, which is coming from a derived functor as explained
above. This will now be explained in some detail in order to point out
the subtleties of the later modifications. Basically, we follow
the considerations in \cite[Sect.~6.11]{jacobson:1989a}.

As in Definition~\ref{definition:Hochschild-complex}, let $\mathbb{K}$
be a field of characteristic zero and let $(\mathcal{A},\mu)$ be a
unital $\mathbb{K}$-algebra with unit $1$. Considering the
\emph{opposite algebra} $(\mathcal{A}^\opp,\mu^\opp)$ where
$\mathcal{A}=\mathcal{A}^\opp$ as a vector space and where the
multiplication is given by $\mu^\opp (a,b)= \mu(b,a)$ for all
$a,b\in\mathcal{A}$ one can consider the extended algebra
\begin{equation}
    \label{eq:extended-algebra-general}
    \mathcal{A}^e= \mathcal{A}\otimes \mathcal{A}^\opp.
\end{equation}
Here and in the following $\otimes=\otimes_{\mathbb{K}}$ is the tensor
product of $\mathbb{K}$-vector spaces. Then, every $\mathbb{K}$-vector
space $\mathcal{M}$ which is an $(\mathcal{A},\mathcal{A})$-bimodule
can be seen as an $\mathcal{A}^e$-module via $(a\otimes b)m= amb$ and
vice versa via $am=(a\otimes 1)m$ and $mb=(1\otimes b)m$. With the
usual left and right multiplication the same is true for the algebra
$\mathcal{A}$ itself. With the ring structure of $\mathcal{A}^e$ one
now comes to the following well-known result, confer
\cite[Thm.~6.17]{jacobson:1989a}.

\begin{proposition}
    \label{proposition:Hochschild-Ext}
    Let $\mathcal{A}$ be an algebra and $\mathcal{M}$ be an
    $(\mathcal{A},\mathcal{A})$-bimodule as in
    Definition~\ref{definition:Hochschild-complex}. Then, the
    Hochschild cohomology is given by a derived functor. In detail,
    for all $k\in\mathbb{N}_0$ one has
    \begin{equation}
        \label{eq:Ext-Hochschild}
        \HH^k(\mathcal{A},\mathcal{M})=
        \Ext^k_{\mathcal{A}^e}(\mathcal{A},\mathcal{M})= 
        R^k \mathrm{hom}(\cdot, \mathcal{M})\mathcal{A}= 
        H^k( \Hom_{\mathcal{A}^e}(C, \mathcal{M})),
    \end{equation}
    where in the last expression $(C,\epsilon)$ is an arbitrary
    projective resolution of $\mathcal{A}$ as $\mathcal{A}^e$-module.
\end{proposition}

\begin{proof}
    The proof makes use of a particular projective resolution, namely
    the so-called \emph{bar resolution}, showing that $H^k(
    \Hom_{\mathcal{A}^e}(C, \mathcal{M}))$ is isomorphic to the
    Hochschild cohomology defined in
    Definition~\ref{definition:Hochschild-complex}. For $k\ge 0$
    consider the vector spaces
    \begin{equation}
        \label{eq:bar-modules-general}
        X_k= \mathcal{A}\otimes \underbrace{\mathcal{A}\otimes \dots \otimes
          \mathcal{A}}_{k \textrm{ times}}\otimes \mathcal{A},
    \end{equation}
    which are $(\mathcal{A},\mathcal{A})$-bimodules via
    \begin{equation}
        \label{eq:bar-module-structure-general}
        a(x_0\otimes \dots \otimes x_{k+1})b= ax_0 \otimes x_1\otimes
        \dots x_k \otimes x_{k+1}b.
    \end{equation}
    Moreover, there are isomorphisms
    \begin{equation}
         \label{eq:bar-modules-general-extended-algebra}
%         \begin{array}{c}
%             \label{eq:bar-modules-general-extended-algebra}
%             X_0 =  \mathcal{A} \otimes \mathcal{A} \cong \mathcal{A}^e,
%             \quad a\otimes b\longmapsto a\otimes b, \quad \textrm{with }
%             a,b\in\mathcal{A} \nonumber \\
%             X_1= \mathcal{A}\otimes \mathcal{A}\otimes \mathcal{A}\cong
%             \mathcal{A}^e\otimes \mathcal{A},\quad  a\otimes x\otimes b\longmapsto
%             (a\otimes b)\otimes x \quad \textrm{with } a,b,x\in \mathcal{A}
%             \textrm{ and} \\
%             X_k\cong \mathcal{A}^e \otimes X_{k-2},\quad  a\otimes x\otimes b
%             \longmapsto (a\otimes b)\otimes x \quad \textrm{with } a,b\in
%             \mathcal{A}, x\in X_{k-2} \textrm{ for } k\ge 1 \nonumber
%         \end{array}
        \begin{array}{rclrcll}
            X_0 =  \mathcal{A} \otimes \mathcal{A} &\cong& \mathcal{A}^e,
            &a\otimes b&\longmapsto& a\otimes b, &\textrm{with }
            a,b\in\mathcal{A}  \\
            X_1= \mathcal{A}\otimes \mathcal{A}\otimes \mathcal{A}&\cong&
            \mathcal{A}^e\otimes \mathcal{A}, & a\otimes x\otimes b&\longmapsto&
            (a\otimes b)\otimes x &\textrm{with } a,b,x\in \mathcal{A}
            \textrm{ and} \\
            X_k&\cong& \mathcal{A}^e \otimes X_{k-2},& a\otimes x\otimes b
            &\longmapsto& (a\otimes b)\otimes x & \textrm{with } a,b\in
            \mathcal{A}, x\in X_{k-2}, k\ge 1 
        \end{array}
    \end{equation}
    between $\mathcal{A}^e$-modules. Using
    \eqref{eq:bar-modules-general-extended-algebra} it is clear that
    the $X_k$ are $\mathcal{A}^e$-free since they are tensor products
    of $\mathcal{A}^e$ and $\mathbb{K}$-vector spaces.

    The $\mathbb{K}$-linear maps $d_k: X_k\longrightarrow X_{k-1}$,
    $k\ge 1$, which are defined by
    \begin{equation}
        \label{eq:differential-bar-general}
        d_k(x_0\otimes\cdots \otimes x_{k+1})=\sum_{i=0}^k (-1)^i
        x_0\otimes \dots \otimes x_i x_{i+1} \otimes \dots \otimes
        x_{k+1}, 
    \end{equation}
    are also $\mathcal{A}^e$-linear and satisfy $d_{k-1}\circ
    d_k=0$. Together with the $\mathcal{A}^e$-linear augmentation
    $\epsilon: X_0\longrightarrow \mathcal{A}$, given by
    \begin{equation}
        \label{eq:augmentation-bar-general}
        \epsilon(a\otimes b) = ab,
    \end{equation}
    and satisfying $\epsilon \circ d_1=0$, one gets an
    $\mathcal{A}^e$-free complex $(X,\epsilon)$ over
    $\mathcal{A}$. The $\mathbb{K}$-homomorphisms $h_1:
    \mathcal{A}\longrightarrow X_0$ and $h_k: X_k\longrightarrow
    X_{k+1}$, given by
    \begin{equation}
        \label{eq:bar-homotopy-map-general}
        h_k(x_0\otimes \dots \otimes x_{k+1})= 1\otimes x_0\otimes \dots
        \otimes x_{k+1},
    \end{equation}
    satisfy
    \begin{equation}
        \label{eq:bar-homotopy-general}
        \epsilon \circ h_{-1} = \id_{\mathcal{A}}, \quad d_1 \circ h_0
        + h_{-1} \circ \epsilon=
        \id_{X_0}, \quad \textrm{and} \quad d_{k+1}\circ h_k + h_{k-1}
        \circ d_k
        =\id_{X_k}, \textrm{ for } k\ge 1. 
    \end{equation}
    Thus, $(C,\epsilon)$ is indeed a projective, even free resolution
    of $\mathcal{A}$ as an $\mathcal{A}^e$-module. Now, for all $k\in
    \mathbb{N}_0$ there is an isomorphism
    \begin{eqnarray}
        \label{eq:isomorphism-bar-Hochschild-general}
        \Xi^k: \Hom_{\mathcal{A}^e}(X_k,\mathcal{M})&\longrightarrow&
        \HC^k(\mathcal{A},\mathcal{M}),\nonumber \\
        (\Xi^k \psi) (x_1,\dots, x_k)&=& \psi(1\otimes x_1\otimes \dots
        \otimes x_k\otimes 1), 
    \end{eqnarray}
    onto the modules of the Hochschild complex
    $(\HC^\bullet(\mathcal{A},\mathcal{M}),\delta)$ as in
    Definition~\ref{definition:Hochschild-complex}. These are chain
    homomorphisms
    \begin{equation}
        \label{eq:bar-Hochschild-chain-map}
        \Xi^{k+1} \circ d_{k+1}^*= \delta^k \circ \Xi^k
    \end{equation}
    and so the cohomology groups defined by the derived functor are
    the Hochschild cohomologies from
    Definition~\ref{definition:Hochschild-complex}.
\end{proof}

\begin{remark}
    In principle, the Hochschild cohomology can be defined by
    \eqref{eq:Ext-Hochschild} in an even more general algebraic
    framework, confer \cite[Sect.~6.11]{jacobson:1989a}.
\end{remark}

\section{Continuous Hochschild cohomologies for  $C^\infty(V)$}
\label{sec:Continuous-Hochschild}

% Proposition \ref{proposition:Hochschild-Ext} shows that one can use
% any suitable projective resolution of $\mathcal{A}$ as an
% $\mathcal{A}^e$-module in order to compute the purely algebraic
% Hochschild cohomology in a given situation.  However, if one wants
% to consider the cohomology in the framework of particular types of
% Hochschild complexes as explained in Section \ref{subsec:HCtype},
% especially differential ones as in Section
% \ref{subsec:differential-Hochschild-complexes}, the general
% considerations are not valid anymore. They only give a guideline and
% one has to check in each case both if the Hochschild complex is of
% the form \eqref{eq:Ext-Hochschild} and if different choices of
% resolutions
% still yield the same cohomology.\\

In this section we give a concrete example for Hochschild complexes
with \emph{continuous} cochains and reformulate
Proposition~\ref{proposition:Hochschild-Ext} in the topological
setting.

From now on let $\mathbb{K}$ be $\mathbb{R}$ or $\mathbb{C}$. Then we
consider the $\mathbb{K}$-algebra $\mathcal{A}=C^\infty(V)$ of smooth
functions on an arbitrary open subset $V\subseteq
\mathbb{R}^n$. Equipped with the usual Fr\'{e}chet topology of smooth
functions, $\mathcal{A}$ is a topological
$\mathbb{K}$-algebra. Further, let $\mathcal{M}$ be a locally convex
$\mathbb{K}$-vector space and a topological
$(\mathcal{A},\mathcal{A})$-bimodule which means that the trilinear
map $(a,m,b)\longmapsto a\cdot m\cdot b$ is continuous with respect to the
induced product topology. Then, one can consider the Hochschild
complex $(\HC(\mathcal{A},\mathcal{M}),\delta)$ of $\mathcal{A}$ with
values in $\mathcal{M}$ as in
Definition~\ref{definition:Hochschild-complex} and therein the
cochains which are continuous maps with respect to the induced
topology on the $k$-fold products $\mathcal{A}\times \dots \times
\mathcal{A}$ for $k\in \mathbb{N}$. By the well-known properties of
locally convex spaces, confer \cite{jarchow:1981a,treves:1967a}, and
the defining Equation~\eqref{eq:Hochschild-bimodules}, it is clear
that this yields a subcomplex
\begin{equation}
    \label{eq:continuous-Hochschild-complex}
    (\HCcont^\bullet (\mathcal{A},\mathcal{M}), \delta),
\end{equation}
the so-called \emph{continuous Hochschild complex of $\mathcal{A}$
  with values in $\mathcal{M}$}. The only thing to show is that the
Hochschild differential maps continuous cochains into continuous
ones. To do this, one basically uses the definition and the fact that
a linear map $f:E\longrightarrow F$ between two locally convex spaces
is continuous, if to every continuous seminorm $p_F$ on $F$ there
exists a continuous seminorm $p_E$ on $E$ such that for all
$e\in E$ one has $p_F(f(e))\le p_E(e)$,
\cite[Prop.~7.7]{treves:1967a}.

For the algebra $\mathcal{A}$ we now investigate the topological
versions of the bar complex, confer
\cite[Sect.~III.2.$\alpha$]{connes:1994a} and
\cite{bordemann.et.al:2005a:pre, pflaum:1998a}. In principle, one just
considers the completions of the spaces and continuous continuations
of the continuous maps occurring in the proof of
Proposition~\ref{proposition:Hochschild-Ext} with respect to the given
locally convex topologies. The extended algebra $\mathcal{A}\otimes
\mathcal{A}^\opp$ in \eqref{eq:extended-algebra-general} is replaced
by its topological counterpart
\begin{equation}
    \label{eq:extended-algebra}
    \mathcal{A}^e= C^\infty(V\times V),
\end{equation}
where the tensor product in \eqref{eq:extended-algebra-general} has
been completed with respect to the canonical Fr\'{e}chet topology of
smooth functions. Since these are \emph{nuclear} there are no
ambiguities in the definition of the tensor product and its
completion, confer \cite{jarchow:1981a,treves:1967a}. In the very same
way all the other structures of the \emph{topological bar complex} are
now derived from the purely algebraic ones. So one redefines the
spaces
\begin{equation}
    \label{eq:bar-complex}
    X_0 = \mathcal{A}^e=C^\infty(V\times V)
    \quad
    \textrm{and}
    \quad
    X_k = C^\infty(V\times V^{k} \times V)
\end{equation}
for $k \in \mathbb{N}$ with the $\mathcal{A}^e$-module structure
\begin{equation}
    \label{eq:module-structure-bar-complex}
    (\hat{a}\chi)(v,q_1,\ldots,q_k,w)=\hat{a}(v,w)\chi(v,q_1,\ldots,q_k,w)
\end{equation}
for $\hat{a}\in \mathcal{A}^e$, $\chi\in X_k$ and $v, w, q_1, \ldots,
q_k \in V$. This is the adapted version of
\eqref{eq:bar-module-structure-general}. In the same way the boundary
operators $d^X_k: X_k \longrightarrow X_{k-1}$ and the augmentation
$\epsilon: X_0 \longrightarrow \mathcal{A}$ are defined by
\begin{align}
    \label{eq:boundary-operator-bar-complex}
    &(d^X_k \chi)(v,q_1,\ldots,q_{k-1},w)
    =
    \chi(v ,v,q_1,\ldots,q_{k-1},w)
    \nonumber\\
    &\quad 
    +\sum_{i=1}^{k-1}
    (-1)^i\chi(v,q_1,\ldots,q_i,q_i,\ldots,q_{k-1},w)
    +(-1)^k \chi(v,q_1,\ldots,q_{k-1},w,w)
\end{align}
and
\begin{equation}
    \label{eq:augmentation}
    (\epsilon\hat{a})(v) =  \hat{a}(v, v).
\end{equation}
It is obvious that $d^X_k$ and $\epsilon$ are homomorphisms of
$\mathcal{A}^e$-modules and an easy computation in analogy to the
algebraic case yields $d^X_{k-1}\circ d^X_k=0$ for all $k\ge 2$ and
$\epsilon\circ d^X_1=0$. The homotopies $h^X_{-1}:
\mathcal{A}\longrightarrow X_0$ and $h^X_k: X_k\longrightarrow
X_{k+1}$ now are given by
\begin{equation}
    \label{eq:homotopy-maps-bar-complex}
    (h^X_{-1}a)(v,w) = a(w) 
    \quad
    \textrm{and}
    \quad
    (h^X_k \chi)(v,q_1, \ldots, q_{k+1},w)
    = \chi(q_1,\ldots,q_{k+1},w) \quad \textrm{for } 
    k \ge 0,  
\end{equation}
and again satisfy
\begin{eqnarray}
    \label{eq:homotopy-bar-complex}
    \epsilon \circ h^X_{-1}&=&\id_\mathcal{A}, \nonumber\\
    h^X_{-1}\circ \epsilon + d^X_1 \circ h^X_0 &=& \id_{X_0}
    \quad \textrm{and}\\
    h^X_{k-1} \circ d^X_k +d^X_{k+1} \circ h^X_k &=&
    \id_{X_k} \quad \textrm{for all } k\ge 1. \nonumber
\end{eqnarray}
Hence the sequence
\begin{equation}
    \label{eq:diagram-bar-resolution}
      \xymatrix{{0} 
        & \mathcal{A} \ar[l] 
        & X_0 \ar[l]_\epsilon 
        & X_1 \ar[l]_{d^X_1} 
        & \ldots \ar[l]_{d^X_2} 
        & X_k \ar[l]_{d^X_k} 
        & \ldots \ar[l]_{d^X_{k+1}} 
      },
\end{equation}
is exact and the bar complex $(X,d^X)$ defines a resolution
$((X,d^X),\epsilon)$ of $\mathcal{A}$. Note that the modules $X_k$ are
\emph{topologically free} as $\mathcal{A}^e$-modules, confer
\cite{connes:1994a} for a more general version of this.

With these structures we now reformulate Proposition
\ref{proposition:Hochschild-Ext} for the considered specific
situation.

\begin{proposition}
    \label{proposition:Continuous-Hochschild}
    For $\mathbb{K}=\mathbb{R}$ or $\mathbb{C}$ and $V\subseteq
    \mathbb{R}^n$ let $\mathcal{A}=C^\infty(V)$ be the Fr\'{e}chet
    algebra of smooth functions on $V$. Further, let $\mathcal{M}$ be
    a locally convex $\mathbb{K}$-vector space with a topological
    $(\mathcal{A},\mathcal{A})$-bimodule structure. Then, $\mathcal{M}$
    has a unique $\mathcal{A}^e= C^\infty(V\times V)$-module structure
    with
    \begin{equation}
        \label{eq:regularity-condition}
        (a\otimes b) \cdot m 
        = a\cdot m \cdot b \quad \textrm{for all } a,b\in \mathcal{A},
        m\in \mathcal{M}. 
    \end{equation}
    For all $k\in \mathbb{N}$ there is an isomorphism
    \begin{equation}
        \label{eq:iso-bar-cont-Hochschild-cont}
        \Xi^k: \Hom_{\mathcal{A}^e}^{\mathrm{cont}}(X_k,\mathcal{M})
        \longrightarrow \HCcont^k(\mathcal{A},\mathcal{M})
    \end{equation}
    between the continuous and $\mathcal{A}^e$-linear maps from the
    topological bar modules $X_k$ to $\mathcal{M}$ and the
    continuous Hochschild modules, given by
    \begin{equation}
        \label{eq:iso-bar-cont-Hochschild-cont-explicitly}
        \left(\Xi^k\psi\right)(a_1,\dots, a_k) = \psi(1\otimes a_1
        \otimes \dots \otimes a_k \otimes 1).
    \end{equation}
    Further, $\Xi$ is a chain map, this means
    \begin{equation}
        \label{eq:iso-cont-differentials}
        \Xi^{k+1} \circ (d^X_{k+1})^*= \delta^k\circ \Xi^k \quad
        \textrm{for all } 
        k\in \mathbb{N}_0, 
    \end{equation}
    so $\Xi$ is even an isomorphism of complexes and thus the
    cohomologies are the same,
    \begin{equation}
        \label{eq:cohomology-cont}
        \HHcont^\bullet(\mathcal{A},\mathcal{M}) = H^\bullet
        (\Hom_{\mathcal{A}^e}^{\mathrm{cont}}(X,\mathcal{M})).
    \end{equation}
\end{proposition}

\begin{proof}
    By the continuity conditions, the $\mathcal{A}^e$-module structure
    of $\mathcal{M}$ is the unique extension of
    \eqref{eq:regularity-condition}. Further, it is clear that the
    $\Xi^k$ defined by
    \eqref{eq:iso-bar-cont-Hochschild-cont-explicitly} take values in
    the continuous Hochschild cochains. This is seen by using the
    seminorms which induce the topologies, and the given continuity
    properties. The injectivity of $\Xi^k$ follows from the fact that
    a continuous and $\mathcal{A}^e$-linear map in
    $\Hom_{\mathcal{A}^e}^{\mathrm{cont}} (X_k,\mathcal{M})$ is
    already determined by the values on elements of the form $1\otimes
    a_1 \otimes \dots \otimes a_k \otimes 1$ since the products of
    those elements with factorising elements $a\otimes b\in
    \mathcal{A}^e$ build a dense subset in $X_k$. The surjectivity is
    given by of the same reason because every map defined by the right
    hand side of \eqref{eq:iso-bar-cont-Hochschild-cont-explicitly}
    has a unique $\mathcal{A}^e$-linear and continuous continuation on
    $X_k$. The boundary operators $d^X_k$ are continuous maps. This is
    easily seen for \eqref{eq:differential-bar-general} and thus it is
    true for its continuation $d^X_k$. Equation
    \eqref{eq:iso-cont-differentials} is a direct computation using
    \eqref{eq:differential-bar-general} and the $\mathcal{A}^e$-module
    structure for elements $a\otimes b$ induced by
    \eqref{eq:bar-module-structure-general}. Then, the rest of the
    proposition is clear.
\end{proof}

\section{Differential Hochschild cohomologies for $C^\infty(V)$}
\label{sec:Differential-Hochschild}

The topological setting described above turns out to be the
appropriate starting point in order to specify the aspired
characterization of Hochschild complexes and cohomologies in the
framework of resolutions for the differential case. As already seen in
Section~\ref{sec:differential-algebra-module-structures} it is rather
difficult to define a general setting of Hochschild complexes with
algebraic multidifferential operators as cochains. However, for smooth
functions $\mathcal{A}=C^\infty(V)$ as above there exists a slightly
different but very natural way to define differential operators. As we
will see now this can be used to define corresponding
\emph{differential} subcomplexes of the ones in
\eqref{eq:iso-bar-cont-Hochschild-cont} which are isomorphic
again. This and the subsequent considerations ending up in
Theorem~\ref{theorem:isomorphism-Hochschild-bar} then will
provide a strong tool to compute differential Hochschild cohomologies
since the different notions of differential operators
coincide in the crucial examples.

In the following it will be necessary to demand that the module
$\mathcal{M}$ in Proposition~\ref{proposition:Continuous-Hochschild}
has more specific properties. So from now on, let $\mathcal{M}$ be a
locally convex, complete, topological Hausdorff space which is a
topological left $\mathcal{A}$-module with respect to the Fr\'{e}chet
topology of $\mathcal{A}=C^\infty(V)$. This means that the bilinear
multiplication
\begin{equation}
    \label{eq:left-module-continuous}
    (\mathcal{A} \times \mathcal{M}) \ni (a,m) \longmapsto a\cdot m\in
    \mathcal{M} 
\end{equation}
is continuous. Furthermore, we demand that $\mathcal{M} $ has a right
$\mathcal{A}$-module structure with the property that there exists an
$l\in \mathbb{N}$ such that the right module multiplication can be
expressed in terms of the left module multiplication by
\begin{equation}
    \label{eq:right-module-structure-differential}
    m \cdot b
    = \sum_{|\beta|\le l}
    \frac{\partial^{|\beta|}b} {\partial v^\beta} \cdot m^\beta 
    \quad \textrm{for all } b \in \mathcal{A}, m \in \mathcal{M}
\end{equation}
with elements $m^\beta\in \mathcal{M}$ depending continuously on $m$.
Such a module $\mathcal{M}$ which shall be referred to as
\emph{$(\mathcal{A},\mathcal{A})$-bimodule of order $l$} clearly
satisfies the conditions of
Proposition~\ref{proposition:Continuous-Hochschild}.

\begin{lemma}
    \label{lemma:differential-left-right-module}
    Let $\mathcal{M}$ be an $(\mathcal{A},\mathcal{A})$-bimodule of
    order $l$. Then $\mathcal{M}$ is a topological
    $(\mathcal{A},\mathcal{A})$-bimodule and the induced left
    $\mathcal{A}^e$-module structure is explicitly given by
    \begin{equation}
        \label{eq:extended-structure-and-left-module-structure}
        \hat{a} \cdot m
        = \sum_{|\beta|\le l} 
        \left( 
            \Delta_0^*
            \frac{\partial^{|\beta|}\hat{a}}
            {\partial w^\beta}
        \right) 
        \cdot m^\beta
        \quad \textrm{for all } \hat{a}\in \mathcal{A}^e,
    \end{equation}
    where $\Delta_k^*$ denotes the pullback with the total diagonal
    map $\Delta_k: V \longrightarrow V^{k+2}$ for $k\in \mathbb{N}$
    which is defined by $\Delta_k(v)=(v,\dots,v)$ and where the
    differentiation is the one with respect to the second argument of
    $\hat{a}$.
\end{lemma}

\begin{proof}
    Since $\mathcal{A}$ is commutative the right and left module
    structures combine to a bimodule structure. Using the common
    seminorms $p_{K,r}$ for $C^\infty(V)$ defined by
    \begin{equation}
        \label{eq:seminorms-functions}
        p_{K,r}(a)=\max_{
          \begin{subarray}{c}
              w\in K\\
              |\beta|\le r
          \end{subarray}
        } \frac{\partial^{|\beta|} a}{\partial v^\beta}(w)
    \end{equation}
    for all compact subsets $K\subseteq V$ and $r\in \mathbb{N}_0$,
    the fact that $l$ in
    \eqref{eq:right-module-structure-differential} is uniform for
    $\mathcal{M}$ implies the continuity of the trilinear map $(a, m,
    b) \longmapsto a\cdot m\cdot b$. Thus $\mathcal{M}$ is a topological
    bimodule as in
    Proposition~\ref{proposition:Continuous-Hochschild}. That the
    induced right $\mathcal{A}^e$-module has the form
    \eqref{eq:extended-structure-and-left-module-structure} follows by
    continuity from
    \begin{equation*}
        (a\otimes b)\cdot m = a \cdot m \cdot b = a\cdot
        \sum_{|\beta|\le l} \frac{\partial^{|\beta|}b} {\partial
          v^\beta} \cdot m^\beta
        = \sum_{|\beta|\le l} 
        \left( 
            \Delta_0^*
            \frac{\partial^{|\beta|} (a\otimes b)}
            {\partial w^\beta}
        \right) 
        \cdot m^\beta.
    \end{equation*}
\end{proof}
For the rest of this and the following section we use the definition
of differential operators which is natural for the algebra of smooth
functions and, at first sight, slightly different to the purely
algebraic definition of Section~\ref{sec:algebraic-multi-diffop}.

\begin{definition}[Multidifferential maps of
    $\mathcal{A}=C^\infty(V)$]
    \label{definition:differential-maps}
    Let $k \in \mathbb{N}$, $\mathcal{A}=C^\infty(V)$ be as above and
    let $\mathcal{M}$ be a left $\mathcal{A}$-module. An
    $\mathbb{R}$-multilinear map $\phi: \mathcal{A}\times \ldots
    \times \mathcal{A} \longrightarrow \mathcal{M}$ with $k$ arguments
    is said to be differential of multiorder $L=(l_1, \ldots, l_k)
    \in \mathbb{N}_0^k$, if it has the form
    \begin{equation}
        \label{eq:differential-map}
        \phi(a_1, \ldots, a_k)
        = 
        \sum_{|\alpha_1| \le l_1,\ldots,|\alpha_k| \le l_k}
        \left( 
            \frac{\partial^{|\alpha_1|} a_1}{\partial v^{\alpha_1}}
            \cdots
            \frac{\partial^{|\alpha_k|} a_k}{\partial v^{\alpha_k}}
        \right)
        \cdot \phi^{\alpha_1\cdots \alpha_k}
    \end{equation}
    with multiindices $\alpha_1, \ldots, \alpha_k \in \mathbb{N}_0^n$
    and $\phi^{\alpha_1\cdots \alpha_k} \in \mathcal{M}$.  These
    $L$-differential maps are denoted by $\Diffop^L(\mathcal{A},\dots,
    \mathcal{A}; \mathcal{M})$.
\end{definition}

\begin{remark}
    \begin{enumerate}
    \item If $\mathcal{M}$ is a topological $\mathcal{A}$-module any
        differential map in the sense of
        Definition~\ref{definition:differential-maps} is continuous.
    \item As already stated in the beginning,
        Definition~\ref{definition:differential-maps} is consistent
        with the purely algebraic definition of multidifferential
        operators in Section~\ref{sec:algebraic-multi-diffop} for the
        crucial examples treated in this work as we will see later in
        Chapter~\ref{cha:deformation-on-bundles}.
    \item In this sense,
        \eqref{eq:right-module-structure-differential} means that the
        right module structure is differential with respect to the
        left one, so $\phi_m: b \longmapsto m \cdot b$ is a differential
        operator depending on $m$ of order $l$.
    \end{enumerate}
\end{remark}

Using Definition~\ref{definition:differential-maps} for a module
$\mathcal{M}$ as in Lemma~\ref{lemma:differential-left-right-module}
we can consider the Hochschild cochains in
$\HCcont^\bullet(\mathcal{A},\mathcal{M})$ consisting of
multidifferential operators. Due to the additional structures this
yields a subcomplex.

\begin{proposition}
    Let $\mathcal{A}$ and $\mathcal{M}$ be as in
    Proposition~\ref{proposition:Continuous-Hochschild} where
    $\mathcal{M}$ has the additional property
    \eqref{eq:right-module-structure-differential} with the fixed
    $l\in \mathbb{N}$. Then the definition
    \begin{equation}
        \label{eq:differential-complex}
        \HCdiff^k(\mathcal{A},\mathcal{M})=
        \Union_{L\in \mathbb{N}_0^k} \Diffop^L(\mathcal{A}, \dots,
        \mathcal{A}; \mathcal{M}) \quad \textrm{for all } k\in \mathbb{N}
    \end{equation}
    and $\HCdiff^0(\mathcal{A},\mathcal{M})=\mathcal{M}$ yields a
    subcomplex $(\HCdiff^\bullet(\mathcal{A},\mathcal{M}),\delta)$ of
    $(\HCcont^\bullet(\mathcal{A},\mathcal{M}),\delta)$, the so-called
    \emph{differential Hochschild complex of $\mathcal{A}$ with values
      in $\mathcal{M}$}. In particular, the Hochschild differential
    $\delta$ restricts to maps
    \begin{equation}
        \label{eq:differential-Hochschild-differential}
        \delta:
        \Diffop^L(\mathcal{A},\dots,\mathcal{A};\mathcal{M})
        \longrightarrow
        \Diffop^{\tilde{L}}
        (\mathcal{A},\dots,\mathcal{A};\mathcal{M})   
    \end{equation}
    for all $k\in\mathbb{N}_0$ and $L=(l_1,\dots, l_k)\in
    \mathbb{N}_0^k$ where $\tilde{L}=(\tilde{l}_1, \dots,
    \tilde{l}_{k+1})\in \mathbb{N}_0^{k+1}$ is given by
    \begin{eqnarray}
        \label{eq:grading-differential}
        \tilde{l}_1&=& l_1,\nonumber \\
        \tilde{l}_i &=& \max \{ l_{i-1}, l_i \} \quad \textrm{for }
        i=2,\dots, k,\\
        \tilde{l}_{k+1}&=& \max \{ l_k, l\}. \nonumber        
    \end{eqnarray}
\end{proposition}

\begin{proof}
    The continuity properties concerning $\mathcal{M}$ and
    Definition~\ref{definition:differential-maps} assure that
    $\HCdiff^k(\mathcal{A},\mathcal{M})\subseteq
    \HCcont^k(\mathcal{A},\mathcal{M})$ is a subset for all $k\in
    \mathbb{N}_0$. Then one only has to show
    \eqref{eq:differential-Hochschild-differential} and
    \eqref{eq:grading-differential}. For $L$ as in
    \eqref{eq:differential-Hochschild-differential} and $\phi\in
    \Diffop^L(\mathcal{A},\dots, \mathcal{A};\mathcal{M})$ as in
    \eqref{eq:differential-map} one computes with
    \eqref{eq:right-module-structure-differential}
    \begin{eqnarray*}
        \label{eq:differential-Hochschild-differential-proof}
        \lefteqn{ (\delta \phi) (a_1,\dots, a_{k+1})}\\
        &=& \sum_{|\alpha_1|\le l_1,\ldots,|\alpha_k|\le l_k}
        a_1 \cdot \left( 
            \frac{\partial^{|\alpha_1|} a_2}{\partial v^{\alpha_1}}
            \cdots
            \frac{\partial^{|\alpha_k|} a_{k+1}}{\partial
              v^{\alpha_k}} 
        \right) \cdot \phi^{\alpha_1\dots \alpha_k}\\
        && + \sum_{i=1}^k  (-1)^i \sum_{|\alpha_1|\le l_1,\ldots,
          |\alpha_k|\le l_k}
        \left( 
            \frac{\partial^{|\alpha_1|} a_1}{\partial v^{\alpha_1}}
            \dots \frac{\partial^{|\alpha_i|} (a_i a_{i+1})}{\partial 
              v^{\alpha_i}} \dots 
            \frac{\partial^{|\alpha_k|} a_{k+1}}{\partial
              v^{\alpha_k}} 
        \right)
        \cdot \phi^{\alpha_1\cdots \alpha_k}\\
        && +(-1)^{k+1} \sum_{|\alpha_1|\le l_1,\ldots, |\alpha_k|\le 
          l_k} 
        \left( 
            \frac{\partial^{|\alpha_1|} a_1}{\partial v^{\alpha_1}}
            \cdots
            \frac{\partial^{|\alpha_k|} a_k} {\partial v^{\alpha_k}}  
        \right)\cdot \phi^{\alpha_1\dots \alpha_k} \cdot a_{k+1}\\
        &=& \sum_{|\alpha_1|\le l_1,\ldots,|\alpha_k|\le l_k}
        \left(a_1 
            \frac{\partial^{|\alpha_1|} a_2}{\partial v^{\alpha_1}}
            \cdots
            \frac{\partial^{|\alpha_k|} a_{k+1}}{\partial
              v^{\alpha_k}} 
        \right) \cdot \phi^{\alpha_1\cdots \alpha_k}\\
        && + \sum_{i=1}^k (-1)^i
        \sum_{
          \begin{subarray}{c}
              |\alpha_1|\le l_1,\ldots,|\alpha_k|\le l_k\\
              0\le \beta\le \alpha_i
          \end{subarray}
        } \binom{\alpha_i}{\beta}
        \left( 
            \frac{\partial^{|\alpha_1|} a_1}{\partial v^{\alpha_1}}
            \dots \frac{\partial^{|\beta|} a_i}{\partial 
              v^{\beta}}
            \frac{\partial^{|\alpha_i-\beta|} a_{i+1}} {\partial 
              v^{\alpha_i-\beta}}
            \dots 
            \frac{\partial^{|\alpha_k|} a_{k+1}}{\partial
              v^{\alpha_k}} 
        \right)
        \cdot \phi^{\alpha_1\cdots \alpha_k}\\
        && +(-1)^{k+1} \sum_{
          \begin{subarray}{c}
              |\alpha_1|\le l_1,\ldots,|\alpha_k|\le l_k\\
              |\beta|\le l
          \end{subarray}
        }
        \left( 
            \frac{\partial^{|\alpha_1|} a_1}{\partial v^{\alpha_1}}
            \cdots
            \frac{\partial^{|\alpha_k|} a_k} {\partial v^{\alpha_k}}   
            \frac{\partial^{|\beta|} a_{k+1}} {\partial v^{\beta}}
        \right) \cdot \phi^{\alpha_1\cdots \alpha_k \beta},
      \end{eqnarray*}
      where the well-known Leibniz rule was used, confer
      \cite[Chap.~VII, Sect.~5, Ex.~21]{amann.escher:2006b}. Given two
      multiindices $\alpha,\beta$ of length $n$ their binomial
      coefficient is defined by
      $\binom{\alpha}{\beta}=\prod_{s=1}^n\frac{\alpha_s!}{\beta_s
        !(\alpha_s-\beta_s)!}$. This shows that $\delta \phi$ is again
      of the form \eqref{eq:differential-map} and a simple counting of
      orders of differentiation yields
      \eqref{eq:grading-differential}.
\end{proof}

\begin{remark}
    In comparison with the definitions made in
    Section~\ref{sec:differential-algebra-module-structures} the new
    definition \eqref{eq:differential-complex} contains a slight abuse
    of notation. In the later examples, however, these notions
    coincide.
\end{remark}

By use of the isomorphism $\Xi$ the differential elements in
$\Hom_{\mathcal{A}^e} (X_k,\mathcal{M})$ of multiorder
$L=(l_1,\ldots,l_k) \in \mathbb{N}_0^k$ for all $k\in \mathbb{N}$ can
be defined as
\begin{equation}
    \label{eq:Hom-bar-differential}
    \Hom_{\mathcal{A}^e}^{\mathrm{diff},L}
    (X_k,\mathcal{M})= (\Xi^k)^{-1}(\Diffop^L(\mathcal{A},\dots, 
    \mathcal{A};\mathcal{M})).  
\end{equation}
Equation~\eqref{eq:iso-cont-differentials} shows that
$\Hom_{\mathcal{A}^e}^{\mathrm{diff}} (X_\bullet, \mathcal{M})$ is a
subcomplex of $\Hom_{\mathcal{A}^e}^{\mathrm{cont}} (X_\bullet,
\mathcal{M})$. By the very construction,
\eqref{eq:iso-bar-cont-Hochschild-cont-explicitly} restricts to an
isomorphism of complexes
\begin{equation}
    \label{eq:bijection-Hom-bar-Hochschild-diff}
    \Xi: \left(
        \Hom_{\mathcal{A}^e}^{\mathrm{diff}} (X_\bullet, \mathcal{M}),
        (d^X)^*
    \right)
    \longrightarrow
    \left(
        \mathrm{HC}^\bullet_{\mathrm{diff}} (\mathcal{A},\mathcal{M}),
        \delta
    \right).
\end{equation}

The concrete form of the elements in \eqref{eq:Hom-bar-differential}
is clarified in the following lemma.

\begin{lemma}
    An element $\psi \in \Hom_{\mathcal{A}^e}^{\mathrm{diff},L}
    (X_k,\mathcal{M})$ has the form
    \begin{equation}
        \label{eq:Hom-modules-elements-differential}
        \psi(\chi) =
        \sum_{
          \begin{subarray}{c}
              |\alpha_1|\le l_1,\ldots,|\alpha_k|\le l_k\\
              |\beta|\le l
          \end{subarray}
              }
              \left(
              \Delta_k^*
          \frac{\partial^{|\alpha_1|+ \cdots + |\alpha_k| +|\beta|} \chi}
          {\partial q_1^{\alpha_1} \cdots \partial
          q_k^{\alpha_k} \partial w^\beta}
               \right)
          \cdot \psi^{\alpha_1 \cdots \alpha_k\beta}
    \end{equation}
    with multiindices $\alpha_1, \ldots, \alpha_k, \beta \in
    \mathbb{N}_0^n$ and $\psi^{\alpha_1 \cdots \alpha_k \beta}\in
    \mathcal{M}$.
\end{lemma}

\begin{proof}
    The stated form \eqref{eq:Hom-modules-elements-differential}
    follows again by continuity arguments. Let $\psi \in
    \Hom_{\mathcal{A}^e}^{\mathrm{diff},L} (X_k,\mathcal{M})$ be such
    that $\Xi \psi \in \Diffop^L(\mathcal{A},\dots,
    \mathcal{A};\mathcal{M})$ as in \eqref{eq:differential-map}. Then
    one computes
    \begin{eqnarray*}
        \lefteqn{\psi(a\otimes a_1\otimes \dots \otimes a_k \otimes b)}\\
        &=& \psi((a\otimes b)(1\otimes a_1\otimes \dots \otimes
        a_k\otimes 1))\\
        &=& a (\Xi \psi) (a_1,\dots, a_k) b\\
        &=& \sum_{
          \begin{subarray}{c}
              |\alpha_1|\le l_1,\ldots,|\alpha_k|\le l_k\\
              |\beta|\le l
          \end{subarray}
        } 
        a \frac{\partial^{|\alpha_1|} a_1}{\partial v^{\alpha_1}}
        \dots \frac{\partial^{|\alpha_k|} a_k}{\partial 
          v^{\alpha_k}}
        \frac{\partial^{|\beta|} b} {\partial
          v^{\beta}} 
        \cdot \psi^{\alpha_1\cdots \alpha_k\beta}\\
        &=&  \sum_{
          \begin{subarray}{c}
              |\alpha_1|\le l_1,\ldots,|\alpha_k|\le l_k\\
              |\beta|\le l
          \end{subarray}
        }
        \left(
            \Delta_k^*
            \frac{\partial^{|\alpha_1|+ \cdots + |\alpha_k| +|\beta|}
              (a\otimes a_1\otimes \dots \otimes a_k\otimes b)}
            {\partial q_1^{\alpha_1} \cdots \partial
              q_k^{\alpha_k} \partial w^\beta}
        \right)
        \cdot \psi^{\alpha_1 \cdots \alpha_k\beta}
    \end{eqnarray*}
    and \eqref{eq:Hom-modules-elements-differential} follows.
\end{proof}

\section{The topological Koszul resolution of $C^\infty(V)$}
\label{sec:Koszul}

Following the ideas of homological algebra explained in the
introduction of this chapter, we now want to find another resolution
of $\mathcal{A}=C^\infty(V)$ such that the cohomology derived by the
functor $\Hom_{\mathcal{A}^e}(\cdot,\mathcal{M})$ is the same as the
continuous and, in particular, the differential Hochschild cohomology
defined above. In the following it will be shown that the topological
version of the Koszul resolution $(K,\epsilon)$ of $\mathcal{A}$
achieves this aim. It will turn out that the values
$\Hom_{\mathcal{A}^e}(K_k,\mathcal{M})$ of the functor
$\Hom_{\mathcal{A}^e}(\cdot,\mathcal{M})$ will not have to be
restricted to homomorphisms of a particular form. This and the fact
that the Koszul resolution is finite once more point out the
importance of the presented approach since the framework of the
aspired computation of the cohomologies becomes easier in two ways. On
one hand there only remain finitely many cohomology groups to compute
and on the other hand one does not have to care about the rather
technical properties of continuous or differential cochains. In
addition, this last statement in particular shows that the continuous
and the differential cohomologies considered above coincide.

After this very motivating consideration we now recall the
\emph{topological Koszul complex over $C^\infty(V)$}. In order to
guarantee that all occurring maps are well-defined we need the further
assumption that
\begin{equation}
    \label{eq:choose-convex-subset}
    V\subseteq \mathbb{R}^n \quad \textrm{is a \emph{convex} open
      subset.} 
\end{equation}

Then one considers the finitely many (topologically free)
$\mathcal{A}^e$-modules
\begin{equation}
    \label{eq:Koszul-complex}
    K_0  
    = \mathcal{A}^e
    \quad
    \textrm{and}
    \quad
    K_k  = 
    \mathcal{A}^e\otimes_{\mathbb{R}}
    \Dach{k}(\mathbb{R}^n)^* 
    \cong C^\infty(V\times V, \Dach{k}(\mathbb{R}^n)^*).
\end{equation}
There $\Dach{k}(\mathbb{R}^n)^*$ denotes the antisymmetric tensor
product of the dual space $(\mathbb{R}^n)^*$ of $\mathbb{R}^n$. The
boundary operators $d^K_k: K_k\longrightarrow K_{k-1}$ for
$k=1,\dots,n$ are defined by
\begin{eqnarray}
    \label{eq:boundary-operator-Koszul-complex}
    \left((d^K_k\omega)(v, w)\right)(x_1, \dots, x_{k-1})
    &=&
    \left(\omega(v, w)\right)(v-w, x_1, \dots, x_{k-1}) \\
    &=& (\insa(v-w)\omega(v,w))(x_1,\dots, x_{k-1}) \nonumber
\end{eqnarray}
for $v, w \in V$ and $x_1, \ldots, x_{k-1} \in \mathbb{R}^n$ where
$\insa$ denotes the usual insertion maps
$\insa(x):\Dach{k}(\mathbb{R}^n)^*\longrightarrow
\Dach{k-1}(\mathbb{R}^n)^*$ of vectors $x\in \mathbb{R}^n$ into forms.
With a basis $\{e_i\}_{i = 1, \ldots, n}$ of $\mathbb{R}^n$ and the
corresponding coordinate functions $\{x^j\}_{j=1,\dots,n}$ defined by
$x^j(v^i e_i)=v^j$ the boundary operators are simply given by
\begin{equation}
    \label{eq:boundary-operator-Koszul-complex-short}
    d^K_k= \sum_{j=1}^n \xi^j \insa(e_j)
\end{equation}
with the functions
\begin{equation}
    \label{eq:functions-xi}
    \xi^j=x^j\otimes 1- 1\otimes x^j\in
    \mathcal{A}^e,
\end{equation}
this means $\xi^j(v,w)=(v^j-w^j)$ for all $v ,w\in V$ with $v=v^i e_i$
and $w=w^i e_i$.

Again, the maps $d^K_k$ are $\mathcal{A}^e$-module homomorphisms with
$d^K_{k-1} \circ d^K_k = 0 $ for all $k\ge 2$ since
$\insa(v-w)\insa(v-w)=0$ and $\epsilon \circ d^K_1 = 0$. Then consider
the maps $h^K_{-1} = h^X_{-1}$ and $h^K_k: K_k \longrightarrow
K_{k+1}$ defined by
\begin{equation}
    \label{eq:homotopy-map-Koszul-complex-k}
    (h^K_k\omega)(v, w)
    = \sum_{j=1}^n  e^j \wedge \int\limits_0^1 \D t \: t^k 
    \frac{\partial \omega}{\partial v^j} (tv+(1-t)w, w)
\end{equation}
for $k \ge 0$ where $\frac{\partial \omega}{\partial v^j}$ denotes the
element in $K_k$ obtained from $\omega$ by taking the partial
derivative of its function part in $\mathcal{A}^e$ with respect to the
$j$-th component of the first argument.  With the basis $\{e_i\}_{i =
  1, \ldots, n}$ of $\mathbb{R}^n$ and the corresponding dual basis
$\{e^i\}_{i = 1, \ldots, n}$ of $(\mathbb{R}^n)^*$ the elements
$\omega \in K_k$ can be written as
\begin{equation}
    \label{eq:element-Koszul-complex}
    \omega
    = \frac{1}{k!} \sum_{i_1,\ldots,i_k=1}^n 
    \omega_{i_1\ldots i_k} e^{i_1} \wedge \ldots \wedge e^{i_k}
\end{equation}
with $\omega_{i_1 \ldots i_k}\in \mathcal{A}^e$. Then,
\eqref{eq:homotopy-map-Koszul-complex-k} reads
\begin{equation}
    \label{eq:homotopy-map-koszul-explicitly}
    (h_K^k\omega)(v, w)
    = \sum_{i_1, \dots , i_k, j=1}^n \int\limits_0^1 \D t \: t^k 
    \frac{\partial \omega_{i_1\dots i_k}}{\partial v^j} (tv+(1-t)w, w)
    e^j\wedge e^{i_1}\wedge \dots \wedge e^{i_k}.
\end{equation}

\begin{lemma}
    With the above maps and $\epsilon: K_0\longrightarrow \mathcal{A}$
    as in \eqref{eq:augmentation} the following homotopy identities
    hold.
    \begin{eqnarray}
        \label{eq:homotopy-Koszul-complex}
        \epsilon \circ h^K_{-1}&=&\id_\mathcal{A}, \nonumber\\
        h^K_{-1}\circ \epsilon +d^K_1 \circ h^K_0 &=& \id_{K_0}
        \quad \textrm{and}\\
        h^K_{k-1} \circ d^K_k +d^K_{k+1} \circ h^K_k &=&
        \id_{K_k} \quad \textrm{for all } k\ge 1 \nonumber.
    \end{eqnarray}
\end{lemma}

\begin{proof}
    For the following computation it is convenient to use the
    abbreviation $\omega= \omega_{I_k}e^{I_k}\in K_k$ instead of
    \eqref{eq:element-Koszul-complex} with the obvious meaning. With
    the summation convention and $\insa(v-w)= (v^i-w^i)\insa(e_i)$ we
    compute for $k\ge 1$
    \begin{equation}
        \label{eq:differential-after-homotopy}
        (d^K_{k+1} (h^K_k \omega))(v,w)
        = 
        \int\limits_0^1 \D t\: t^k \frac{\partial \omega_{I_k}}{\partial
          v^j} (tv +(1-t)w,w) (v^i-w^i) \insa(e_i)(e^j\wedge e^{I_k})
    \end{equation}
    and 
    \begin{eqnarray}
        \label{eq:homotopy-after-differential}
        (h^k_{k-1} (d^K_k \omega))(v,w)
        &=& 
        \int\limits_0^1 \D t\: t^{k-1} \frac{\partial}{\partial
          v^j} \left( \omega_{I_k}(v,w)(v^i-w^i)\right)|_{(tv+(1-t)w,w)}
        e^j\wedge \insa(e_i) e^{I_k}\nonumber \\
        &=&
        \int\limits_0^1 \D t \: t^k \frac{\partial \omega_{I_k}}{\partial
          v^j}(tv+(1-t)w,w) (v^i-w^i)  e^j\wedge \insa(e_i)
        e^{I_k} \\
        && + \int\limits_0^1 \D t \: t^{k-1} \omega_{I_k}(tv+(1-t)w,w)
        e^i\wedge \insa(e_i) e^{I_k}. \nonumber
    \end{eqnarray}
    Since $e^i\wedge \insa(e_i) e^{I_k}=k e^{I_k}$, the last summand
    can be written as
    \begin{eqnarray*}
        \int\limits_0^1 \D t \: t^{k-1} \omega_{I_k}(tv+(1-t)w,w)
        e^i\wedge \insa(e_i) e^{I_k}
        &=& \int\limits_0^1 \D t \: \frac{\D} {\D t} \left( t^k
            \omega_{I_k} (tv+(1-t)w,w) \right) e^{I_k} \\ 
        && - \int\limits_0^1 \D t \: t^k \frac{\partial \omega_{I_k}}
        {\partial v^j}(tv+(1-t)w,w) (v^j-w^j) e^{I_k}. 
    \end{eqnarray*}
    This and the derivation property of $\insa(e_i)$ show that the sum
    of \eqref{eq:differential-after-homotopy} and
    \eqref{eq:homotopy-after-differential} is
    \begin{eqnarray*}
        \label{eq:homotopy-Koszul-proof}
        ((d^K_{k+1} \circ h^K_k+ h^K_{k-1} \circ d^K_k) \omega)(v,w)=
        \omega_{I_k}(v,w) e^{I_k}= \omega(v,w).
    \end{eqnarray*}
    This shows the last equation in
    \eqref{eq:homotopy-Koszul-complex}. The second follows completely
    analogously and the first is the same as for the bar resolution.
\end{proof}

\begin{remark}
    \begin{enumerate}
    \item Note that in \eqref{eq:homotopy-map-Koszul-complex-k} the
        convexity of $V$ is crucial to get a well-defined map. There
        are of course different ways to define the homotopies in
        \eqref{eq:homotopy-maps-bar-complex} and
        \eqref{eq:homotopy-map-Koszul-complex-k}. In
        \cite{bordemann.neumaier.waldmann.weiss:2007a:pre}, for
        example, the slightly more complicated formulas
        $(h^X_{-1}a)(v,w) = a(v)$, $(h^X_k \chi)(v,q_1, \ldots,
        q_{k+1},w) = (-1)^{k+1}\chi(v,q_1,\ldots,q_{k+1})$ and
        $(h_K^k\omega)(v, w) = -\sum_{j=1}^n e^j \wedge
        \int\limits_0^1 \D t \: t^k \frac{\partial \omega}{\partial
          w^j} (v, tw+(1-t)v)$ have been used.
    \item The explicit formulas
        \eqref{eq:homotopy-map-Koszul-complex-k} and
        \eqref{eq:homotopy-map-koszul-explicitly} of the maps $h^k$
        show that it is necessary to use the topological version
        $\mathcal{A}^e= C^\infty(V\times V)$ of the extended
        algebra. If $\omega\in K_k$ as in
        \eqref{eq:element-Koszul-complex} has factorising coefficient
        functions $\omega_{i_1,\dots,i_k}\in C^\infty(V)\otimes
        C^\infty(V)$ this does not have to be the case for
        $h^k(\omega)$. For instance, take
        $\omega_{i_1,\dots,i_k}(v,w)= v \sin(w)$. The same behaviour
        will occur for the explicit chain maps in
        Section~\ref{sec:chain-maps-bar-Koszul}. This is the main
        reason why one has to work with the topological completions.
    \end{enumerate}
\end{remark}

The Equations~\eqref{eq:homotopy-Koszul-complex} show that
$((K,d^K),\epsilon)$ is indeed a topologically free and \emph{finite}
resolution of $\mathcal{A}$,
\begin{equation}
    \label{eq:diagram-Koszul-resolution}
     \xymatrix{{0} 
       & \mathcal{A} \ar[l] 
       & K_0 \ar[l]_\epsilon 
       & K_1 \ar[l]_{d^K_1} 
       & \ldots \ar[l]_{d^K_2} 
       & K_k \ar[l]_{d^K_k} 
       & \ldots \ar[l]_{d^K_{k+1}} 
       & K_n \ar[l]_{d^K_{n}}
       & \ar[l] 0
      }.
\end{equation}

In an adapted notion, confer \cite[Sect.~6.12]{jacobson:1989a} for the
original one, this shows that $\mathcal{A}$ has \emph{finite
  homological dimension}.

\section{Explicit chain maps between the bar and Koszul resolution} 
\label{sec:chain-maps-bar-Koszul}

The general results concerning projective resolutions, in particular
Theorem~\ref{theorem:projective-complex-resolution}, give an
indication that there exist chain maps between the topological bar and
Koszul resolution and that they are always homotopic. For our purpose,
however, the pure existence is not enough. Since we are interested in
differential Hochschild cochains which obviously is an additional
property besides being continuous we do not try to reformulate the
general theory in the topological setting but consider explicit chain
maps and construct corresponding homotopy maps whose images under the
considered functor will turn out to respect differential cochains as
well. The following explicit formulas and constructions already appear
to some extent in \cite{bordemann.et.al:2005a:pre}. Here we provide a
self-contained presentation including all the relevant proofs.

For all $k\ge 0$ we consider the maps $F_k: K_k\longrightarrow X_k$
from \cite[Sect.~III.2$\alpha$]{connes:1994a} defined by
\begin{equation}
    \label{eq:chain-maps-Koszul-bar}
    (F_k \omega)(v, q_1, \ldots, q_k, w) 
    = \left( 
        \omega(v, w)
    \right)(q_1 - v, \ldots, q_k - v),
\end{equation}
and the maps $G_k: X_k \longrightarrow K_k$ from
\cite{bordemann.et.al:2005a:pre} defined by
\begin{equation}
    \label{eq:chain-maps-bar-Koszul}   
    \begin{split}
        (G_k \chi)(v, w)
        &= \sum_{i_1, \ldots, i_k = 1}^n 
        e^{i_1} \wedge \ldots \wedge e^{i_k} 
        \int\limits_0^1 \D t_1 \int\limits_0^{t_1}\D t_2 \cdots
        \int\limits_0^{t_{k-1}} \D t_k \\
        &\quad
        \frac{\partial^k\chi}
        {\partial q_1^{i_1} \cdots \partial q_k^{i_k}}
        (v, t_1 v + (1 - t_1) w, \ldots, t_k v + (1 - t_k) w, w).
    \end{split}
\end{equation}
In particular, $F_0 = \id_{\mathcal{A}^e} = G_0$. Note that $G_k$ is
well-defined since we assume $V$ to be convex. Clearly, these maps are
$\mathcal{A}^e$-module homomorphisms and it is a straightforward
computation which shows that $F_k$ and $G_k$ are chain maps, confer
Appendix~\ref{cha:technical-proofs} for the details. This means that
for all $k \ge 0$
\begin{equation}
    \label{eq:chain-maps}
    d^X_{k+1} \circ F_{k+1} =F_k\circ d^K_{k+1}
    \quad
    \textrm{and}
    \quad 
    d^K_{k+1} \circ G_{k+1} = G_k \circ d^X_{k+1}.
\end{equation}
Thus we have the commutative diagram of $\mathcal{A}^e$-module
morphisms
\begin{equation}
    \label{eq:diagram-bar-Koszul}
    \begin{array}{c}
        \xymatrix{{0} 
          & \mathcal{A} \ar[l] \ar @{=} [d] 
          & X_0 \ar @{=} [d]_{\id_{\mathcal{A}^e}} \ar[l]_\epsilon 
          & X_1 \ar[l]_{d^X_1} \ar @<+.5ex>[d]^{G_1} 
          & \ldots \ar[l]_{d^X_2} 
          & X_k \ar[l]_{d^X_k} \ar @<+.5ex>[d]^{G_k} 
          & X_{k+1} \ar[l]_{d^X_{k+1}} \ar @<+.5ex>[d]^{G_{k+1}} 
          & \ldots \ar[l]_{d^X_{k+2}} \\
          {0} 
          & \mathcal{A} \ar[l] 
          & K_0 \ar[l]^\epsilon 
          & K_1 \ar[l]^{d^K_1} \ar @<+.5ex>[u]^{F_1} 
          & \ldots \ar[l]^{d^K_2} 
          & K_k \ar[l]^{d^K_k} \ar @<+.5ex>[u]^{F_k} 
          & K_{k+1} \ar[l]^{d^K_{k+1}} \ar @<+.5ex>[u]^{F_{k+1}} 
          & \ldots.  \ar[l]^{d^K_{k+2}} 
        }   
    \end{array}
\end{equation}
In Appendix~\ref{cha:technical-proofs} it is further shown that
\begin{equation}
    \label{eq:chain-maps-id}
    G_k\circ F_k=\id_{K_k} \quad \textrm{for all } k\ge 0,
\end{equation}
and thus it is clear that
\begin{equation}
    \label{eq:chain-maps-projection}
    \Theta_k= F_k\circ G_k: X_k\longrightarrow X_k 
\end{equation}
for all $k\ge 0$ is a projection $\Theta_k\circ \Theta_k=\Theta_k$ on
the bar complex with $\Theta_k\circ d^X_{k+1}= d^X_{k+1} \circ
\Theta_{k+1}$. Clearly, $\Theta_0=\id_{X_0}$. The explicit form of
$\Theta_k$ is 
\begin{equation}
    \label{eq:Theta-explizit}
    \begin{split}
    &(\Theta_k \chi)(v,q_1,\ldots,q_k,w)= \sum_{i_1,\ldots,i_k=1}^n\sum_{\sigma\in S_k}
    (-1)^\sigma (q_1-v)^{i_{\sigma(1)}}\ldots
    (q_k-v)^{i_{\sigma(k)}}\\ &\qquad
    \int\limits_0^1 \D t_1 \int\limits_0^{t_1} \D t_2 \cdots
    \!\!\int\limits_0^{t_{k-1}} \D t_k\: \frac{\partial^k \chi}{\partial
      q_1^{i_1}\cdots q_k^{i_k}}(v,t_1v+(1-t_1)w,\ldots,t_kv+(1-t_k)w,w),
    \end{split}
\end{equation}
where $S_k$ denotes the set of all permutations $\sigma$ of
$(1,\ldots,k)$ and $(-1)^\sigma$ is the signum of $\sigma$. This is
obvious due to the definitions and $(e^{i_1}\wedge\dots\wedge e^{i_k})
(q_1,\dots,q_k)= \sum_{\sigma\in S_k} (-1)^\sigma
q^{i_{\sigma(1)}}\dots q^{i_{\sigma(k)}}$.

If it is possible to show that the images of the chain maps $F_k$ and
$G_k$ under the functor $\Hom_{\mathcal{A}^e}(\cdot, \mathcal{M})$,
this means their pullbacks, induce isomorphisms between the
corresponding cohomologies, the topological Koszul resolution is
identified to have the same cohomologies as the Hochschild complexes
of interest. The fact that $G_k\circ F_k=\id$ and the contravariance
of the functor already show that the induced maps $F_k^*$ are
surjective and the $G_k^*$ are injective. Thus the maps on cohomology
level are indeed isomorphisms if the chain maps $\Theta_k=F_k\circ
G_k$ are homotopic to the identity $\id_{X_K}$ by use of homotopy maps
$s_k: X_k\longrightarrow X_{k+1}$. These considerations contain two
crucial points which have to be taken into account.

\begin{remark}[Requirements of the homotopy]
    \label{remark:conditions-homotopy}
    \begin{enumerate}
    \item In contrast to the homotopy maps $h^X_k$ and $h^K_k$ used to
        show the exactness of the topological bar and Koszul
        resolution, the aspired homotopy maps $s_k:X_k\longrightarrow
        X_{k+1}$ have to be \emph{$\mathcal{A}^e$-linear} since the
        homotopy equation on the bar complex has to be carried over to
        the complex of interest by the functor
        $\Hom_{\mathcal{A}^e}(\cdot, \mathcal{M})$.
    \item In the considered framework of differential Hochschild
        complexes one has to guarantee that all involved maps with
        target $\Hom_{\mathcal{A}^e}(X_k, \mathcal{M})$, this means
        $(s_k)^*$, $(G_k)^*$, and $(\Theta_k)^*$, have values in the
        subspaces $\Hom_{\mathcal{A}^e}^{\mathrm{diff}}(X_k,
        \mathcal{M})$ of differential cochains. In contrast to the
        first point this does \emph{not} follow from abstract
        arguments and really makes it necessary to consider the
        explicit maps.
    \end{enumerate}
\end{remark}

\begin{remark}
    \label{remark:homotopy-for-K-modules}
    Regarding the resolutions as complexes of $\mathbb{K}$-modules
    instead of $\mathcal{A}^e$-modules it is very easy to find an
    explicit homotopy between $\Theta_k$ and $\id_{X_k}$. Application
    of $\Theta_k$ to the last equation of
    \eqref{eq:homotopy-bar-complex} leads to
    \begin{equation}
        \label{eq:bar-homotopy-projection}
        \Theta_k \circ h^X_{k-1} \circ d^X_k+ d^X_{k+1} \circ
        \Theta_{k+1} \circ h^X_k=
        \Theta_k. 
    \end{equation}
    Subtraction of \eqref{eq:bar-homotopy-projection} from
    \eqref{eq:homotopy-bar-complex} then yields the homotopy
    \begin{equation}
        \label{eq:wrong-homotopy}
        \id_{X_k}-\Theta_k= t_{k-1} \circ d^X_k + d^X_{k+1}\circ t_k
    \end{equation}
    with the homotopy maps $t_k= (\id_{X_{k+1}}- \Theta_{k+1}) \circ
    h^X_k$ which are not $\mathcal{A}^e$-linear.
\end{remark}

\section{Predifferential maps}
\label{sec:pre-differential-maps}

After Remark~\ref{remark:conditions-homotopy} has made clear the
requirements of the desired homotopy, we now have a closer look at a
sufficiently large class of $\mathbb{K}$-endomorphisms of the
topological bar complex $X_\bullet$ which respect the basic demand
mentioned in the second point of
Remark~\ref{remark:conditions-homotopy}. Within this class which
contains all examples considered so far, we then present an
appropriate procedure to modify arbitrary maps in order to get
$\mathcal{A}^e$-linear ones. This crucial step will then lead to an
explicit construction of the desired homotopy maps. Before defining
the mentioned class of maps, we first introduce and investigate a
certain kind of matrices which will play the role of coordinate
transformations of the convex set $V\subseteq \mathbb{R}^n$ and its
cartesian products. So, consider the set $\mathcal{B}_{sr}$ of
matrices enjoying the following convexity property
\begin{equation}
    \label{eq:convex-matrices}
    \mathcal{B}_{sr} = 
    \bigg\{
    B \in \mathrm{Mat}_{(s+2)\times(r+2)}(\mathbb{R}) 
    \; \bigg| \;
    \forall i \; \textrm{and} \; v_0, \ldots v_{r+1} \in V:
    \sum_{j=0}^{r+1} B_{ij} = 1
    \; \textrm{and} \; 
    \sum_{j=0}^{r+1} B_{ij} v_j \in V
    \bigg\}.
\end{equation}
Each such matrix $B \in \mathcal{B}_{sr}$, written as $B=(B_{ij})$
with $i=0,\dots,s+1$ and $j=0,\dots, r+1$ acts as a map $B: V^{r+2}
\longrightarrow V^{s+2}$ by
\begin{equation}
    \label{eq:matrix-as-map}
    B(v_0,\dots,v_{r+1})= \left(\sum_{j=0}^{r+1}B_{0j} v_j, \dots,
    \sum_{j=0}^{r+1}B_{(s+1)j}v_j\right), 
\end{equation}
which is the restriction of $B \otimes \id_{\mathbb{R}^n}:
\mathbb{R}^{r+2} \otimes \mathbb{R}^n \longrightarrow \mathbb{R}^{s+2}
\otimes \mathbb{R}^n$ to $V^{r+2}$. The convexity property of $B$
obviously ensures that $B(v_0, \ldots, v_{r+1})$ is indeed in
$V^{s+2}$. 

\begin{remark}[Properties and examples of $\mathcal{B}_{sr}$]
    \label{remark:convex-matrices}
    \begin{enumerate}
    \item Due to the defining properties the matrices
        $\mathcal{B}_{sr}\subseteq \mathrm{Mat}_{(s+2)\times
          (r+2)}(\mathbb{R})$ build a convex subset. In addition,
        these matrices are closed under the usual matrix
        multiplication, this means $\mathcal{B}_{sr}\cdot
        \mathcal{B}_{rk}\subseteq \mathcal{B}_{sk}$.
    \item Typical nontrivial examples of elements $B\in
        \mathcal{B}_{rs}$ are for instance given by all matrices $B$
        with rows of the form $(0, \dots, \lambda,\dots,
        (1-\lambda),\dots 0)$ for $\lambda\in [0,1]$. These examples
        also contain a sort of \emph{permutation matrices} which lead
        to maps of the form $B: (v_0,\dots, v_{r+1})\longmapsto
        (v_{\sigma(0)}, \dots, v_{\sigma(s+1)})$ with $\sigma\in
        S_{r+1}$. In particular, it is clear that the units
        $\mathbb{1}\in \mathrm{Mat}_{(s+2)\times (s+2)}(\mathbb{R})$
        are elements in $\mathcal{B}_{ss}$.
    \end{enumerate}
\end{remark}

Using this kind of matrices we now define the above mentioned class of
maps.

\begin{definition}[Predifferential maps on the topological bar
    complex] 
    \label{definition:class-of-maps}
    An endomorphism $A: X_s \longrightarrow X_r$ of the topological
    bar complex as above is said to be \emph{predifferential of
      multiorder $M=(m_0, \dots , m_{s+1})\in \mathbb{N}_0^{s+2}$} if
    it has the form
    \begin{equation}
        \label{eq:class-of-maps-explicitly}
        \begin{split}
            &(A\chi)(v, q_1, \ldots, q_r, w)
            =
            \sum_{|\alpha_0| \le m_0, \ldots, |\alpha_{s+1}| \le
              m_{s+1}}
            \hspace{-1cm}
            f^{\alpha_0 \ldots \alpha_{s+1}} (v, q_1, \ldots, q_r, w) 
            \\
            &\qquad
            \int\limits_0^1 \D t_1 \cdots \int\limits_0^1 \D t_k
            \: 
            \zeta(t_1, \ldots, t_k)
            \frac{\partial^{|\alpha_0|+ \cdots + |\alpha_{s+1}|} \chi} 
            {\partial v^{\alpha_0} \partial q_1^{\alpha_1} \cdots
              \partial q_s^{\alpha_s} \partial w^{\alpha_{s+1}}} 
            (B(t_1, \ldots, t_k) (v, q_1, \ldots, q_r, w))
        \end{split}
    \end{equation}
    with $f^{\alpha_0, \ldots, \alpha_{s+1}} \in X_r$, $\zeta: [0,
    1]^k \longrightarrow \mathbb{R}$ integrable, and a continuous map
    $B: [0, 1]^k \longrightarrow \mathcal{B}_{sr}$. The set of all
    these maps is denoted by $\tilde{\mathcal{S}}_{rs}^{k, M}$. The
    $\mathbb{R}$-vector space spanned by these elements then is
    denoted by $\mathcal{S}_{rs}^{k, M}$.
\end{definition}

\begin{remark}
    \begin{enumerate}
    \item Note that $\tilde{\mathcal{S}}_{rs}^{k, M}$ is not a vector
        space itself since the sum of two maps has not to be of the
        form \eqref{eq:class-of-maps-explicitly}. The vector space
        $\mathcal{S}_{rs}^{k, M}$ in fact consists of all
        \emph{finite} linear combinations of elements in
        $\tilde{\mathcal{S}}_{rs}^{k, M}$.
    \item Obviously, $\mathcal{S}_{rs}=
        \Union_{k=0}^{\infty}\Union_{M\in \mathbb{N}_0^{s+2}}
        \mathcal{S}_{rs}^{k, M}\subseteq \Hom_{\mathbb{R}}(X_s, X_r)$
        is a filtered subspace with
        \begin{equation}
            \label{eq:filtration-class-of-maps}
            \mathcal{S}_{rs}^{k, M} \subseteq \mathcal{S}_{rs}^{k', M'}
            \textrm { for } k \le k'  \textrm{ and  } M \le M'.
        \end{equation}
    \end{enumerate}
\end{remark}

\begin{remark}[Notation]
    \label{remark:notation-class-of-maps}
    In the following considerations and proofs it will be convenient
    to use a shorter notation for the explicit form
    \eqref{eq:class-of-maps-explicitly} of an element $A\in
    \tilde{\mathcal{S}}^{k,M}_{rs}$. With $\underline{x}=(x_0,\dots,
    x_{s+1})=(v,q_1,\dots, q_s,w)$, $\underline{t}=(t_1,\dots,t_k)$
    and $\int\limits_0^1 \D t_1 \cdots \int\limits_0^1 \D t_k \:
    \zeta(t_1, \ldots, t_k)= \int \D \mu(\underline{t})$ one simply
    writes
    \begin{equation}
        \label{eq:convenient-notation-class-of-maps}
        (A\chi)(v,\underline{q}, w)
        = \sum_{|\alpha_0| \le m_0, \ldots, |\alpha_{s+1}| \le
          m_{s+1}}
        \hspace{-1cm}
        f^{\alpha_0 \ldots \alpha_{s+1}} (\underline{x})
        \int \D\mu(\underline{t}) 
        \frac{\partial^{|\alpha_0|+ \cdots + |\alpha_{s+1}|} \chi}  
        {\partial x^{\alpha_0} \cdots \partial x^{\alpha_{s+1}}}
        (B(\underline{t}) (\underline{x})).
    \end{equation}
    When investigating $\mathcal{A}^e$-linearity, it will be useful to
    set $\underline{q}=(q_1,\dots,q_r)$ and to highlight the relevant
    arguments by $\underline{x}=(v,\underline{q},w)$.
\end{remark}

Denoting the set of matrix-valued maps used in Definition
\ref{definition:class-of-maps} by $\mathcal{B}^k_{sr}: \{ B: [0,1]^k
\longrightarrow \mathcal{B}_{sr}\: |\: B \textrm{ is continuous}\}$,
the matrix multiplication mentioned in
Remark~\ref{remark:convex-matrices} naturally induces a multiplication
$\mathcal{B}^{k_1}_{sr} \times \mathcal{B}^{k_2}_{rt}\longrightarrow
\mathcal{B}^{k_1+k_2}_{st}$ for all $k_1,k_2\in \mathbb{N}_0$ and
$s,r,t\in \mathbb{N}$ by
\begin{equation}
    \label{eq:composition-matrix-valued-function}
    (B_1 B_2) (t_1, \dots, t_{k_1+k_2})= B_1(t_1,\dots, t_{k_1})
    B_2(t_{k_1+1},\dots, t_{k_1+k_2}) \quad \textrm{for all }
    B_1\in\mathcal{B}^{k_1}_{sr}, B_2\in \mathcal{B}^{k_2}_{rt}.
\end{equation}
Using this we find the following important property of the maps
$\mathcal{S}^{k,M}_{rs}$.

\begin{lemma}[Composition of predifferential maps]
    \label{lemma:composition-class-of-maps}
    The usual composition $\circ$ of homomorphisms restricts to a map
    \begin{equation}
        \circ: \mathcal{S}_{tr}^{k_2, N} \times
        \mathcal{S}_{rs}^{k_1,M} \longrightarrow
        \mathcal{S}_{ts}^{k_1+k_2, P}, 
    \end{equation}
    with $P = (p_0, \ldots, p_{s+1})$ and $p_i = m_i + |N|$.  In
    particular, for $A_1 \in \tilde{\mathcal{S}}_{rs}^{k_1, M}$ with
    the matrix-valued function $B_1\in \mathcal{B}_{sr}$ and $A_2 \in
    \tilde{\mathcal{S}}_{tr}^{k_2, N}$ with $B_2\in\mathcal{B}_{rt}$
    one has $A_2 \circ A_1 \in \tilde{\mathcal{S}}^{k_1 + k_2,
      P}_{ts}$ with the matrix-valued function $B= B_1 B_2\in
    \mathcal{B}_{st}$ and where $P=(p_0,\dots , p_{s+1})$ is given by
    \begin{equation}
        \label{eq:estimate-order-composition}
        p_i= m_i+ \sum_{j=0}^{r+1}T(B_1(\underline{t}))_{ij} n_j 
        \le m_i+ |N|.
    \end{equation}
    There $T: \mathrm{Mat}_{(s+1)\times
      (r+1)}(\mathbb{R})\longrightarrow \mathrm{Mat}_{(s+1)\times
      (r+1)}(\{0,1\})$ is given by
    \begin{equation}
        \label{eq:matrix-zero-one}
        T(B)_{ij}= \left\{
              \begin{array}{c}
                  1\\0
              \end{array}
              \textrm{ if }
              \begin{array}{cc}
                  B_{ij}\neq 0\\
                  B_{ij}=0
              \end{array}
              \right..
    \end{equation}
\end{lemma}

\begin{proof}
    The lemma is proved if one verifies the statement for the elements
    $A_1$ and $A_2$ as mentioned. With the obvious indexes one simply
    computes 
    \begin{eqnarray*}
        \lefteqn{ ((A_2\circ A_1)\chi) (\underline{x}) }\\
        &=& 
        \sum_{|\beta_0| \le n_0, \ldots, |\beta_{r+1}| \le
          n_{r+1}}
        \hspace{-1cm}
        f_2^{\beta_0 \ldots \beta_{r+1}} (\underline{x})
        \int \D\mu(\underline{t}_2) 
        \frac{\partial^{|\beta_0|+ \dots + |\beta_{r+1}|} (A_1\chi)}  
        {\partial y^{\beta_0} \cdots \partial y^{\beta_{r+1}}}
        (B_2(\underline{t}_2) (\underline{x}))\\
        &=&
        \sum_{
          \begin{subarray}{}
              |\alpha_0| \le m_0, \ldots, |\alpha_{s+1}| \le
              m_{s+1}\\
              |\beta_0| \le n_0, \ldots, |\beta_{r+1}| \le
              n_{r+1}
          \end{subarray}
        }
        \hspace{-1cm}
        f_2^{\beta_0 \dots \beta_{r+1}} (\underline{x})
        \int \D\mu(\underline{t}_2) \int \D \mu(\underline{t}_1) \\
        && \hspace{4cm}
        \frac{\partial^{|\beta_0|+ \cdots + |\beta_{r+1}|}}  
        {\partial y^{\beta_0}  \cdots \partial y^{\beta_{r+1}}}
        \bigg|_{B_2(\underline{t}_2)(\underline{x})} 
        \left(
            f_1^{\alpha_0 \dots \alpha_{s+1}}(\underline{y}) 
            \frac{\partial^{|\alpha_0|+ \dots + |\alpha_{s+1}|} \chi}  
            {\partial x^{\alpha_0} \dots \partial x^{\alpha_{s+1}}}
            (B_1(\underline{t}_1) (\underline{y}))
        \right)\\
        &&
        = \sum_{|\gamma_0| \le p_0, \ldots, |\gamma_{s+1}| \le
          p_{s+1}}
        \hspace{-1cm}
        f^{\gamma_0 \ldots \gamma_{s+1}} (\underline{x})
        \int \D\mu(t_1,\dots,t_{k_1+k_2}) 
        \frac{\partial^{|\gamma_0|+ \dots + |\gamma_{s+1}|} \chi}  
        {\partial x^{\gamma_0} \cdots \partial x^{\gamma_{s+1}}}
        (B_1B_2(t_1,\dots, t_{k_1+k_2}) (\underline{x})),
    \end{eqnarray*}
    where a redefinition of the $t$-variables and the Leibniz rule
    yield the last form. There, the new coefficient functions
    $f^{\gamma_0 \ldots \gamma_{s+1}}$ are given by the product of 
    functions $f_2^{\beta_0 \ldots \beta_{r+1}}$ and derivatives
    of $f_1^{\alpha_0 \ldots \alpha_{s+1}}$. A counting of
    differentiations proves
    \eqref{eq:estimate-order-composition}. Then, the remaining
    assertions are obvious.
\end{proof}

The important and eponymous property of the predifferential maps $A\in
\mathcal{S}^{k,M}_{rs}$ is that their image $A^*$ under the functor
$\Hom_{\mathcal{A}^e}(\cdot ,\mathcal{M})$ respects the differential
cochains.

\begin{proposition}[Predifferential maps of the bar complex]
    \label{proposition:pre-differential-differential}
    Let $\mathcal{M}$ be an $(\mathcal{A},\mathcal{A})$-bimodule of
    order $l$. Then, every $A\in \mathcal{S}^{k,M}_{rs}$ as in
    Definition~\ref{definition:class-of-maps} which is
    $\mathcal{A}^e$-linear defines a map
    \begin{equation}
        \label{eq:pull-back-pre-differential-map}
        A^*:
        \Hom_{\mathcal{A}^e}^{\mathrm{diff},L}(X_r, \mathcal{M})
        \longrightarrow
        \Hom_{\mathcal{A}^e}^{\mathrm{diff},\tilde{L}}(X_s,
        \mathcal{M})
    \end{equation}
    via pullback respecting the differential cochains with $\tilde{L}
    = (\tilde{l}_1, \ldots \tilde{l_s})$ and
    \begin{equation}
        \label{eq:pull-back-predifferential-orders}
        \tilde{l_i} = m_i + |L| + l.
    \end{equation}
    For $A\in \tilde{\mathcal{S}}^{k,M}_{rs}$ with matrix $B\in
    \mathcal{B}_{sr}$ one has
    \begin{equation}
        \label{eq:pull-back-predifferential-orders-exactly}
        \tilde{l_i} = m_i + \sum_{j=1}^r T(B(\underline{t}))_{ij} l_j
        + T(B(\underline{t}))_{i(r+1)} l.
    \end{equation}
\end{proposition}

\begin{proof}
    By definition, it is sufficient to show that $(A^* \psi)(1\otimes
    a_1\otimes \dots \otimes a_s\otimes 1)= \psi(A(1\otimes
    a_1\otimes \dots \otimes a_s\otimes 1))$ is of the form
    \eqref{eq:differential-map}. First, one notes that
    \begin{eqnarray*}
        \lefteqn{A(1\otimes
          a_1\otimes \dots \otimes a_s\otimes 1)(v,\underline{q},w) }\\
        &=& 
        \sum_{|\alpha_1| \le m_0, \ldots, |\alpha_s| \le
          m_s}
        \hspace{-1cm}
        f^{0\alpha_1 \ldots \alpha_s 0} (v,\underline{q},w)
        \int \D\mu(\underline{t}) 
        \frac{\partial^{|\alpha_1|+ \dots + |\alpha_s|} 
          (1\otimes a_1\otimes \dots \otimes a_s\otimes 1)}  
        {\partial q_1^{\alpha_1} \cdots \partial q_s^{\alpha_s}}
        (B(\underline{t}) (v,\underline{q},w)).\\
    \end{eqnarray*}
    Application of $\psi$ in the form
    \eqref{eq:Hom-modules-elements-differential} basically means to
    build the partial derivatives and to consider the pullback with
    respect to the diagonal $\Delta_r: V\longrightarrow V^{r+2}$ in
    order to get functions in $\mathcal{A}$. The derivatives of the
    functions $f^{0\alpha_1 \ldots \alpha_s 0}$ thus simply yield
    values of functions in $\mathcal{A}$ which in the aspired form
    \eqref{eq:differential-map} contribute to the elements in
    $\mathcal{M}$. Because of the particular form of the elements
    $(1\otimes a_1\otimes \dots \otimes a_s\otimes 1)$ and the fact
    that $B(\underline{t})\circ \Delta_r =\Delta_s$, the procedure for
    the terms in the integrals yield functions of the form
    \begin{equation*}
        \frac{\partial^{|\gamma_1|+ \dots |\gamma_s|}(1\otimes
          a_1\otimes \dots 
          \otimes a_s\otimes 1)} {\partial q_1^{\gamma_1} \dots \partial
          q_s^{\gamma_s}} \circ \Delta_s 
        = \frac{\partial^{|\gamma_1|} a_1} {\partial q_1^{\gamma_1}}
        \dots \frac{\partial^{|\gamma_s|} a_s}{\partial
          q_s^{\gamma_s}}, 
    \end{equation*}
    which do not depend on the integration parameter anymore. Together
    with the inner derivatives, given by components of
    $B(\underline{t})$, the integrations can be carried out. The
    resulting constants again contribute to the elements in
    $\mathcal{M}$ and one has found the desired form
    \eqref{eq:differential-map}.  A careful counting of orders of
    differentiation finally proves the proposition.
\end{proof}

The $\mathcal{A}^e$-linearity of the map $A$ is absolutely necessary
in order to guarantee that the considered functor can be applied to
it. Out of that reason we first investigate a general procedure to
obtain $\mathcal{A}^e$-linear maps out of arbitrary ones. For this
purpose consider the \emph{extended spaces} $\tilde{X}_k=
C^\infty(V\times V\times V^k\times V\times V)\cong X_{k+2}$, $k\in
\mathbb{N}_0$, and the maps
\begin{equation}
    \label{eq:helping-restriction-map}
    p_k^*: X_k \longrightarrow \tilde{X}_k \quad \textrm{with} \quad
    (p_k^*\chi) (v',v,\underline{q},w,w') = \chi(v',\underline{q},w')
\end{equation}
and
\begin{equation}
    \label{eq:helping-extension-map}
    i_k^*: \tilde{X}_k \longrightarrow X_k \quad \textrm{with} \quad
    (i_k^* \chi) (v,\underline{q},w)= \chi(v,v,\underline{q},w,w).
\end{equation}
Obviously, the maps $i_k^*$ and $p_k^*$ are $\mathcal{A}^e$-linear and
satisfy $i_k^* \circ p_k^*=\id_{X_k}$.
For all $v',w'\in V$ and $\chi\in \tilde{X}_s$ there exists an element
$\chi_{v',w'}\in X_s$ defined by 
\begin{equation}
    \label{eq:derived-bar-complex-element}
    \chi_{v',w'}(v,\underline{q},w)= \chi(v',v,\underline{q},w,w').
\end{equation}
Then for all maps $A:X_s\longrightarrow X_r$ there
exists a corresponding map $\tilde{A}: \tilde{X}_s\longrightarrow
\tilde{X}_r$ defined by
\begin{equation}
    \label{eq:helping-derived-Ae-linear-map}
    (\tilde{A}\chi)_{v',w'}= A \chi_{v',w'}.
\end{equation}
Using this, one further defines the map $\cc{A}: X_s\longrightarrow
X_r$ by
\begin{equation}
    \label{eq:derived-Ae-linear-map}
    \cc{A}= i_r^*\circ \tilde{A} \circ p_s^*.
\end{equation}
Then one has the following result.
\begin{lemma}
    Let $A:X_s\longrightarrow X_r$ be $\mathbb{R}$-linear. Then
    $\cc{A}$ is $\mathcal{A}^e$-linear. Further, one has $\cc{A+
      B}=\cc{A}+ \cc{B}$, $\cc{\id_{X}}=\id_{X}$ and $\cc{\lambda A}=
    \lambda \cc{A}$ for $\lambda\in \mathbb{R}$.
\end{lemma}

\begin{proof}
    The properties follow directly from the $\mathcal{A}^e$-linearity
    of $\tilde{A}$ and the fact that $A\longmapsto \tilde{A}$ is an
    algebra morphism.
\end{proof}

In particular, one additionally notes that $\widetilde{A\circ B} =
\tilde{A}\circ \tilde{B}$. This is not true for $\cc{A\circ B}$ but
instead one has the following obvious results.

\begin{lemma}
    \label{lemma:properties-product-derived-maps}
    Let $A: X_r\longrightarrow X_t$ and $B: X_s\longrightarrow X_r$ be
    $\mathbb{R}$-linear.
    \begin{enumerate}
    \item If $A\circ i_r^*= i_t^*\circ \tilde{A}$, then $\cc{A\circ
          B}= A\circ \cc{B}$.
    \item If $p_r^* \circ B= \tilde{B}\circ p_s^*$, then $\cc{A\circ
          B}= \cc{A} \circ B$.
    \end{enumerate}
\end{lemma}

After this general considerations we now come back to the
predifferential maps. From the general form
\eqref{eq:class-of-maps-explicitly} one directly finds a necessary and
sufficient condition for the $\mathcal{A}^e$-linearity, namely that
there are no partial derivatives with respect to the first and last
argument and that those are reproduced under the matrix
$B(\underline{t})$.

\begin{lemma}[The subspace of $\mathcal{A}^e$-linear maps]
    \label{lemma:subspace-Ae-linear-maps}
    An element $A\in \tilde{\mathcal{S}}_{rs}^{k, M}$ as in
    \eqref{eq:class-of-maps-explicitly} with matrix-valued function
    $B\in \mathcal{B}^k_{sr}$ is $\mathcal{A}^e$-linear if and only if
    the two following conditions are satisfied:
    \begin{enumerate}
    \item The multiindex $M$ has the form $M=(0,m_1,\dots, m_s,0)$.
    \item $B$ takes values in the matrices with first row
        $(1\phantom{,}0 \phantom{,}\dots \phantom{,} 0)$ and last one
        $(0\phantom{,}\dots\phantom{,} 0\phantom{,}1)$. This means
        that $B(\underline{t})(v,\underline{q},w)=
        (v,B(\underline{t})_{0,s+1} (v,\underline{q},w), w)$ where
        $B(\underline{t})_{0,s+1} \in \mathcal{B}_{(s-2) r}$ denotes
        the matrix obtained from $B(\underline{t})$ by deleting the
        first and the last row.
    \end{enumerate}
\end{lemma}

Using this observation one can immediately define a projection onto
the subspace of $\mathcal{A}^e$-linear predifferential maps which
turns out to be the restriction of the general prescription
$A\longmapsto \cc{A}$.

\begin{proposition}[Projection onto $\mathcal{A}^e$-linear maps]
    \label{proposition:projection-Ae-linear-pre-differential}
    Let $A\in \tilde{\mathcal{S}}_{rs}^{k, M}$ be a predifferential
    map as in \eqref{eq:convenient-notation-class-of-maps}. Then the
    derived $\mathcal{A}^e$-linear map is again an element $\cc{A}\in
    \tilde{\mathcal{S}}_{rs}^{k, M}$ and given by
    \begin{equation}
        \label{eq:derived-map}
        (\cc{A}\chi)(v,\underline{q}, w)
        =
        \sum_{|\alpha_1| \le m_1, \ldots, |\alpha_{s}| \le m_{s}}
        \hspace{-1cm}
        f^{0 \alpha_1 \ldots \alpha_{s} 0} (v, \underline{q},w)
        \int \D \mu(\underline{t})
        \frac{\partial^{|\alpha_1|+ \cdots + |\alpha_{s}|} \chi}
        {\partial q_1^{\alpha_1} \cdots \partial q_s^{\alpha_s}}
        (v,(B(\underline{t})_{0,s+1}) (v,\underline{q}, w), w).
    \end{equation}
    Further, the map $A\longmapsto \cc{A}$ on all predifferential maps
    $\mathcal{S}_{rs}$ is a projection onto the $\mathcal{A}^e$-linear
    maps, this means
    \begin{equation}
        \label{eq:projection-property-Ae-linear-maps}
        \cc{\cc{A}} = \cc{A}.
    \end{equation}
    Further, for $A_1\in \mathcal{S}_{rs}$ and $A_2\in
    \mathcal{S}_{tr}$ one has
    \begin{equation}
        \label{eq:special-property-projection}
        \cc{A_1} \circ \cc{A_2} = \cc{\cc{A_1} \circ A_2}.
    \end{equation}
\end{proposition}

\begin{proof}
    With $A$ as in \eqref{eq:convenient-notation-class-of-maps} one
    gets
    \begin{equation}
        \label{eq:helping-map-pre-differential}
        (\tilde{A}\chi)(v',\underline{x}, w')
        =\hspace{-1cm} 
        \sum_{|\alpha_0| \le m_0, \ldots, |\alpha_{s+1}| \le
          m_{s+1}}
        \hspace{-1cm}
        f^{\alpha_0 \ldots \alpha_{s+1}} (\underline{x})
        \int \D\mu(\underline{t}) 
        \frac{\partial^{|\alpha_0|+ \cdots + |\alpha_{s+1}|} \chi}   
        {\partial x^{\alpha_0} \cdots \partial x^{\alpha_{s+1}}} 
        (v',B(\underline{t}) (v,\underline{q},w),w'). 
    \end{equation}
    With the definitions of the maps $i_r^*$ and $p_s^*$ this yields
    \eqref{eq:derived-map} since
    \begin{equation*}
        \frac{\partial^{|\alpha_0|+ \cdots +
            |\alpha_{s+1}|}(p_s^*\chi)}    
        {\partial v^{\alpha_0} \partial q^{\alpha_1}\cdots \partial
          q^{\alpha_s} \partial w^{\alpha_{s+1}}} (v',B(\underline{t})
        (v,\underline{q},w),w') 
        = 
        \frac{\partial^{|\alpha_1|+ \cdots +
            |\alpha_s|}\chi}    
        {\partial q^{\alpha_1}\cdots \partial
          q^{\alpha_s}} (v',B(\underline{t})_{0,s+1}
        (v,\underline{q},w),w'). 
    \end{equation*}
    The projection property
    \eqref{eq:projection-property-Ae-linear-maps} is obvious with the
    explicit form \eqref{eq:derived-map}.

    Further, the explicit forms of $\tilde{A}$ and $\cc{A}$ show that
    $\cc{A}$ satisfies $\cc{A}\circ i_s^*= i_r^*\circ \widetilde{\cc{A}}$ since
    \begin{equation*}
        \frac{\partial^{|\alpha_1|+ \cdots +
            |\alpha_s|}(i_s^*\chi)}  
        {\partial q^{\alpha_1}\cdots \partial
          q^{\alpha_s}} (v,B(\underline{t})_{0,s+1}
        (v,\underline{q},w),w) 
        = 
        \frac{\partial^{|\alpha_1|+ \cdots +
            |\alpha_s|} \chi}    
        {\partial q^{\alpha_1}\cdots \partial
          q^{\alpha_s}} (v,v,B(\underline{t})_{0,s+1}
        (v,\underline{q},w),w,w). 
    \end{equation*}
    Then \eqref{eq:special-property-projection} follows with the first
    part of Lemma \ref{lemma:properties-product-derived-maps}.
\end{proof}

\begin{lemma}
    Let $A\in \tilde{\mathcal{S}}^{k,M}_{rs}$ be as in
    \eqref{eq:convenient-notation-class-of-maps} with coefficient
    functions satisfying $f^{\alpha_0\dots
      \alpha_{s+1}}(v,\underline{q},w)= g^{\alpha_0\dots
      \alpha_{s+1}}(\underline{q})$ and matrix-valued function
    $B\in\mathcal{B}^k_{sr}$ satisfying
    $B(\underline{t})(v,\underline{q},w)=(v,C(\underline{t})(\underline{q}),
    w)$ with $C\in \mathcal{B}^k_{(s-2)(r-2)}$. Then, $\cc{A}$
    satisfies $p_r^* \circ \cc{A}= \widetilde{\cc{A}} \circ p_s^*$ and
    thus $\cc{B\circ \cc{A}}=\cc{B}\circ \cc{A}$.
\end{lemma}

\begin{proof}
    The assertions are verified with analogous arguments to those used
    in the proof of Proposition
    \ref{proposition:projection-Ae-linear-pre-differential}.
\end{proof}

 \begin{remark}[Examples]
    \label{remark:examples-predifferential}
    Note that all previously defined maps are predifferential in the
    sense of Definition \ref{definition:class-of-maps}. Indeed, for
    all $k$ one has the following.
    \begin{enumerate}
    \item $d^X_k \in \mathcal{S}_{k-1, k}^{0, (0, \ldots, 0)}$
        and $\cc{d^X_k} = d^X_k$.
    \item $\Theta_k \in \mathcal{S}_{kk}^{k, (0, 1, \ldots, 1, 0)}$
        and $\cc{\Theta^k} = \Theta^k$.
    \item $h^X_k \in \mathcal{S}_{k+1, k}^{0, (0, \ldots,
          0)}$. However, $h^X_k$ is not $\mathcal{A}^e$-linear and
        hence $\cc{h^X_k} \neq h^X_k$.
    \end{enumerate}
\end{remark}

\section{Homotopy and isomorphic differential cohomologies}
\label{sec:homotopy}

With the above preparations we can now find the desired homotopy.

\begin{remark}[First attempt]
    Having the projection onto the $\mathcal{A}^e$-linear maps, one
    could conjecture that the projections $\cc{t_k}$ of the homotopy
    map $t$ found in Remark~\ref{remark:homotopy-for-K-modules} yield
    the desired $\mathcal{A}^e$-linear maps. This is not the
    case. With the properties of the involved functions and
    \eqref{eq:special-property-projection} for $d^X_{k+1}$ application
    of the projection to \eqref{eq:wrong-homotopy} yields
    \begin{equation*}
        \id_{X_k}-\Theta_k= \cc{t_{k-1} \circ d^X_k}+ d^X_{k+1} \circ
        \cc{t_k}, 
    \end{equation*}
    which is not the desired equation since $\cc{t_{k-1} \circ
      d^X_k}\neq \cc{t_{k-1}} \circ d^X_k$.
\end{remark}

The desired $\mathcal{A}^e$-linear homotopy map can be defined
recursively.

\begin{proposition}[The $\mathcal{A}^e$-linear homotopy]
    Let the $\mathcal{A}^e$-linear maps $s_k: X_k\longrightarrow
    X_{k+1}$ for all $k\in \mathbb{N}_0$ be defined recursively by
    \begin{equation}
        \label{eq:recursive-definition-homotopy-map}
        s_k =
        \cc{
          h^X_k \circ 
          \left( \id_{X_k} - \Theta_k - s_{k-1} \circ d^X_k \right)
        }
        \quad
        \textrm{and}
        \quad
        s_0 = 0.
    \end{equation}
    This yields maps $s_k \in \mathcal{S}_{k+1, k}^{k, M(k)}$ with
    \begin{equation}
        \label{eq:homotopy-map-pre-differential}
        M(k) = (0, m_1(k), \dots, m_k(k), 0) \textrm{ and }
        m_i(k)=\binom{k-1}{i-1}\le (k-1)! \textrm{ for } i=1,\dots,k,
    \end{equation}
    which yield a homotopy for $\id_{X_k}$ and $\Theta_k$, this means
    \begin{equation}
        \label{eq:id-Theta-homotopic}
        \id_{X_k}-\Theta_k= s_{k-1}\circ d^X_k + d^X_{k+1}
        \circ s_k \quad \textrm{for all } k\ge 1.
    \end{equation}
\end{proposition}

\begin{proof}
    The assertion \eqref{eq:homotopy-map-pre-differential} is a
    consequence of Lemma \ref{lemma:composition-class-of-maps} and
    follows by induction. For the orders of differentiation one
    directly has by definition that $m_0(k)=0=m_{k+1}(k)$ for all
    $k\in \mathbb{N}$. Further, $M(0)=(0,0)$ and $M(1)=(0,1,0)=
    M(\Theta_1)$. Then, it is clear that for $k\ge 2$ the orders are
    controlled by the term $s_{k-1} \circ d^X_k$ since
    $M(\Theta_k)=(0,1,\dots,1,0)$. In order to apply
    \eqref{eq:homotopy-map-pre-differential} one can combine the
    occurring matrices in $d_k^X$ to an \emph{effective} one which is
    given by the entries $B_{ij}=\delta_{ij}+\delta_{i,j+1}$ for
    $i=0,\dots,k+1$ and $j=0,\dots,k$. Then,
    \eqref{eq:homotopy-map-pre-differential} and induction show that
    for $i=1,\dots,k$ one gets $m_i(k)=m_i(k-1)+
    m_{i-1}(k-1)=\binom{k-1}{i-1}$. In the last step one treats $i=1$
    and $i=k$ separately. The other cases follow with
    $\binom{k-2}{i-1}+ \binom{k-2}{i-2}=\binom{k-1}{i-1}$.

    The homotopy \eqref{eq:id-Theta-homotopic} is shown by induction
    with use of \eqref{eq:homotopy-bar-complex}.  For $k=1$ one has
    \begin{align*}
        d^X_2 \circ s_1 &= \cc{d^X_2 \circ h^X_1 \circ
          (\id_{X_1}-\Theta_1)}=
        \cc{(\id_{X_1}-h^X_0 \circ d^X_1) \circ (\id_{X_1}-\Theta_1)}\\
        &= \id_{X_1}-\Theta_1 -\cc{h^X_0 \circ d^X_1} + \cc{h^X_0
          \circ \Theta_0 \circ d^X_1}= \id_{X_1}-\Theta_1.
    \end{align*}
    Then, induction yields analogously
    \begin{align*}
        d^X_{k+1} \circ s_k &= \cc{(\id_{X_k}-h^X_{k-1} \circ d^X_k)
          \circ (\id_{X_k}-
          \Theta_k -s_{k-1} \circ d^X_k)}\\
        &= \id_{X_k} - \Theta_k - \cc{s_{k-1} \circ d^X_k} -
        \cc{h^X_{k-1}\circ (d^X_k- d^X_k\circ \Theta_k-
          (\id_{X_{k-1}}-
          \Theta_{k-1} -s_{k-2}\circ d^X_{k-1})\circ d^X_k)}\\
        &= \id_{X_k} - \Theta_k - s_{k-1} \circ d^X_k.
    \end{align*}
\end{proof}

With the above structures and results we can now formulate the
following important proposition.
\begin{proposition}
    \label{proposition:pullbacks-diff}
    Let $k\in \mathbb{N}_0$. The pullbacks
    \begin{equation}
        \label{eq:differential-G-star}
        G_k^*: \Hom_{\mathcal{A}^e} (K_k,\mathcal{M})
        \longrightarrow \Hom_{\mathcal{A}^e}^{\mathrm{diff},L}
        (X_k,\mathcal{M})
    \end{equation}
    only take values in the differential cochains of multiorder $L =
    (l+1, \ldots, l+1)\in \mathbb{N}^k$. With the same multiindex $L$
    one has
    \begin{equation}
        \label{eq:differential-Theta-star}
        \Theta_k^*:
        \Hom_{\mathcal{A}^e}^{\mathrm{diff}}(X_k, \mathcal{M})
        \longrightarrow
        \Hom_{\mathcal{A}^e}^{\mathrm{diff},L}(X_k, \mathcal{M}).
    \end{equation}
    Finally, for all $L \in \mathbb{N}_0^{k+1}$ one has
    \begin{equation}
        \label{eq:skStar}
        s_k^*: 
        \Hom_{\mathcal{A}^e}^{\mathrm{diff},L}(X_{k+1}, \mathcal{M})
        \longrightarrow
        \Hom_{\mathcal{A}^e}^{\mathrm{diff},\tilde{L}}(X_k, \mathcal{M}),
    \end{equation}
    where $\tilde{L} = (\tilde{l_1}, \ldots, \tilde{l_k})$ is given by
    $\tilde{l_i} = (k-1)! + |L|$.
\end{proposition}

\begin{proof}
    The statement for $G_k$ is a straightforward counting of
    degrees. For $X\in \Hom_{\mathcal{A}^e}(K_k,\mathcal{M})$ one
    simply computes $(G_k^* X)(1\otimes a_1 \otimes \dots \otimes a_k
    \otimes 1)$. This gives an expression of the form
    $\hat{a}_{i_1\dots i_k}X(e^{i_1}\wedge \dots \wedge e^{i_k})$ with
    elements $\hat{a}_{i_1\dots i_k}\in \mathcal{A}^e$ given by
    \begin{equation*}
        \hat{a}_{i_1\dots i_k}(v,w)= 
        \int\limits_0^1\D t_1 \dots \int\limits_0^{t_{k-1}} \D t_k \:
        \frac{\partial a_1}{\partial v^{i_1}} (t_1 v+ (1-t_1)w)
        \cdots \frac{\partial a_k} {\partial v^{i_k}} (t_k v+
        (1-t_k)w).
    \end{equation*}
    Using \eqref{eq:extended-structure-and-left-module-structure}, the
    product can be written as $\sum_{|\beta|\le l}\left( \Delta_0^*
        \frac{\partial^{|\beta|} \hat{a}_{i_1\dots i_k}} {\partial
          w^\beta}\right) \cdot X^\beta (e^{i_1}\wedge \dots \wedge
    e^{i_k})$ which then has the form $\sum_{|\alpha_i|\le l+1}
    \left(\frac{\partial a_1}{\partial v^{\alpha_1}} \dots
        \frac{\partial a_k}{\partial v^{\alpha_k}}\right)\cdot
    \phi^{\alpha_1\dots \alpha_k}$.

    For $\Theta_k$ the assertion follows either from the properties of
    $G_k$ or from
    Proposition~\ref{proposition:pre-differential-differential}. The
    statement concerning $s_k$ follows with the same argumentation
    where one uses the fact that all the matrices $B$ occurring in
    $s_k$ satisfy $B_{i,k+1}=0$ for all $i=1,\dots,k$.
\end{proof}

This has an immediate consequence.

\begin{theorem}[The differential Hochschild cohomologies]
    \label{theorem:isomorphisms-differential-Hochschild-cohomology}
    With the corresponding pullbacks $F_k^*$ and $G_k^*$ the
    topological bar and Koszul complexes induce the commutative
    diagram
    \begin{equation}
        \label{eq:diagram-Hom-bar-Koszul}
        \xymatrix{
          \ldots \ar[rr]^(.35){(d^X_k)^*} &&
          \Hom_{\mathcal{A}^e}^{\mathrm{diff}} (X_k,\mathcal{M}) \ar
          @<-.5ex> [d]_{F_k^*}
          \ar [rr]^{(d^X_{k+1})^*}&&
          \Hom_{\mathcal{A}^e}^{\mathrm{diff}} (X_{k+1},\mathcal{M}) \ar 
          @<-.5ex> [d]_{F_{k+1}^*}
          \ar[rr]^(.65){(d^X_{k+2})^*}&&
          \ldots\\
          \ldots \ar[rr]^(.35){(d^K_k)^*} &&
          \Hom_{\mathcal{A}^e} (K_k,\mathcal{M}) \ar @<-.5ex>
          [u]_{G_k^*} 
          \ar [rr]^{(d^K_{k+1})^*}&&
          \Hom_{\mathcal{A}^e} (K_{k+1},\mathcal{M}) \ar @<-.5ex>
          [u]_{G_{k+1}^*} 
          \ar[rr]^(.65){d^K_{k+2}}&&
          \ldots.
        }
    \end{equation}
    Further, one has a homotopy
    \begin{equation}
        \label{eq:homotopy-pullbacks}
        \id_{\Hom^{\mathrm{diff}}_{\mathcal{A}^e}(X_k,\mathcal{M})}-\Theta^*_k=
        (d^X_k)^* \circ s_{k-1}^* + (s_k)^* \circ (d^X_{k+1})^*
        \quad \textrm{for all } k\ge 1.
    \end{equation}
    Together with the isomorphism
    \eqref{eq:bijection-Hom-bar-Hochschild-diff} this yields the
    isomorphisms 
    \begin{equation}
        \label{eq:equivalence-Hochschild-Hom-Koszul-cohomology}
        \mathrm{HH}^\bullet_{\mathrm{diff}}(\mathcal{A},\mathcal{M})\cong
        \mathrm{H}(\Hom_{\mathcal{A}^e}^{\mathrm{diff}}
        (X_\bullet,\mathcal{M})) \cong
        \mathrm{H}(\Hom_{\mathcal{A}^e} (K_\bullet,\mathcal{M}))
    \end{equation}
    for the cohomology of the differential Hochschild complex.
\end{theorem}

Note that every isomorphism in
\eqref{eq:equivalence-Hochschild-Hom-Koszul-cohomology} is induced by
explicitly given maps on the level of cochains. With respect to the
application of these results in the next chapter we need the following
obvious generalization.

\begin{theorem}
    \label{theorem:isomorphism-Hochschild-bar}
    Let $\mathcal{M}^\bullet = \bigcup_{l=0}^\infty \mathcal{M}^l$ be
    a filtered $\mathcal{A}$-module, this means $\mathcal{M}^l\subset
    \mathcal{M}^{l+1}$ and $\mathcal{A} \cdot \mathcal{M}^l\subset
    \mathcal{M}^l$ for all $l\in \mathbb{N}$, such that every
    $\mathcal{M}^l$ is an $(\mathcal{A},\mathcal{A})$-bimodule of
    order $l$. Moreover, the topologies have to respect the filtration
    which means that for all $l \in \mathbb{N}$ the topology of
    $\mathcal{M}^l$ is given by the induced one from
    $\mathcal{M}^{l+1}$. Then we have:
    \begin{enumerate}
    \item The unions
        \begin{equation}
            \label{eq:extended-complexes}
            \left(\Union_{l=0}^\infty \mathrm{HC}^\bullet_{\mathrm{diff}}
                (\mathcal{A},\mathcal{M}^l),\delta \right),
            \left(\Union_{l=0}^\infty
                \Hom_{\mathcal{A}^e}^{\mathrm{diff}} (X_\bullet,
                \mathcal{M}^l), (d^X)^* \right), \: \textrm{and} \:
            \left(\Union_{l=0}^\infty
                \Hom_{\mathcal{A}^e}(K_\bullet, \mathcal{M}^l),(d^K)^* \right)
        \end{equation}
        are subcomplexes of
        $\HC^\bullet(\mathcal{A},\mathcal{M}^\bullet)$,
        $\Hom_{\mathcal{A}^e}^{\mathrm{cont}} (X_\bullet,
        \mathcal{M}^\bullet)$, and $\Hom_{\mathcal{A}^e}(K_\bullet,
        \mathcal{M}^\bullet)$, respectively.
    \item The isomorphisms
        \eqref{eq:bijection-Hom-bar-Hochschild-diff} for each $l$
        induce an isomorphism of complexes
        \begin{equation}
            \label{eq:isomorphism-complex}
            \Xi: 
            \left(
                \Union_{l=0}^\infty
                \Hom_{\mathcal{A}^e}^{\mathrm{diff}}
                (X_\bullet, \mathcal{M}^l), 
                \delta_X
            \right)  
            \longrightarrow
            \left(
                \Union_{l=0}^\infty 
                \mathrm{HC}^\bullet_{\mathrm{diff}} 
                (\mathcal{A}, \mathcal{M}^l),
                \delta 
            \right). 
        \end{equation}
    \item The pullbacks $G_k^*$, $F_k^*$, $\Theta_k^*$, and $s_k^*$
        naturally extend to the complexes
        \eqref{eq:extended-complexes}. Thus, we have induced
        isomorphisms for the corresponding cohomologies.
    \end{enumerate}
\end{theorem}

\chapter{Deformation theory of right modules on principal fibre
  bundles}
\label{cha:deformation-on-bundles}

After the general considerations concerning deformation theory of
algebras and modules as well as the detailed discussions on Hochschild
cohomologies we now come back to our initial issues. As we have seen
in the introductory Chapter~\ref{cha:gauge-theory} there is a massive
physical interest in the question how the geometry of principal fibre
bundles can be reformulated when starting with a star product on the
base manifold. From the algebraic point of view such a new formulation
is nothing but the adaption or the replacement of all algebraic
structures related to the algebra of functions on the base such that
the initial algebraic properties are preserved. The natural algebraic
joint to the star product algebra was recognized to be given by all
corresponding right modules. Within the framework of deformation
quantization built up so far this question can now be investigated.

As already seen in
Section~\ref{subsec:algebraic-properties-gauge-theory} the crucial
right modules playing an important role in classical gauge theories
are the horizontal differential forms on the considered principal
fibre bundles. These forms are of course nothing but sections of the
cotangent bundle and tensor products thereof. With respect to the
aspired generality we thus investigate the right module structures of
the sections of vector bundles over the principal bundle and, even
more general, over surjective submersions.

In the first sections of this chapter we substantiate the deformation
problem of the mentioned right modules and show that it can be treated
in the ($G$-invariant) differential setting developed in the
Chapters~\ref{cha:deformation-algebras-modules},
\ref{cha:differential-G-invariant-structures}, and \ref{cha:sheaves}.
It will become evident that this framework really describes the
deformations one is interested in. After the observation of the
additional specific properties of the relevant Hochschild complexes we
attack the central problem and compute the corresponding cohomology
groups. With the results and techniques provided in the previous
chapters we will see that it is always possible to reduce the
computation to a simpler local problem which is solvable. In the end
it turns out that all relevant cohomologies and thus all the
obstructions they encode vanish. Due to the general results of
Chapter~\ref{cha:deformation-algebras-modules} this has the immediate
consequence that the aspired deformations, especially the ones playing
a role in classical gauge theories, always exist and are uniquely
defined up to equivalence.

\section{Right modules on surjective submersions and principal fibre
  bundles}
\label{sec:general-geometric-situation}

In order to define the notion of star products on the algebra of
functions on a manifold $M$, this is necessarily equipped with a
Poisson structure as explained in
Section~\ref{sec:noncommutative-space-times}. In the following let
there be given such a Poisson manifold $M$ together with a
corresponding star product $\star$ for the formal power series
$C^\infty(M)[[\lambda]]$ of functions. In the works
\cite{kontsevich:1997:pre,kontsevich:1997a} of Kontsevich it is not
only shown that a star product always exists but also that it can be
chosen to be differential. This shall be the case in our
situation. The property means that the given deformation
$\mu=\sum_{r=0}^\infty \lambda^r \mu_r$ of the pointwise
multiplication $\mu_0$ consists of bidifferential operators $\mu_r\in
\Diffop^\bullet_{C^\infty(M)}(C^\infty(M),C^\infty(M);C^\infty(M))$. The
notion of differential operators is the one introduced in
Section~\ref{sec:algebraic-multi-diffop}. Thus one is in the general
situation of Definition~\ref{definition:algebra-type-diff} with
$\mathcal{C}=\mathcal{A}=C^\infty(M)$ where $\mu_0$ is of the
differential type of order $(0,0)$.  With respect to this algebra
structure of the differential type we now investigate the
corresponding right module structures occurring in the framework of
surjective submersions and principal fibre bundles.

\subsection{Deformation theory on surjective submersions}
\label{subsec:geometric-situation-surjective-submersion}

Let $\pp: P\longrightarrow M$ be a surjective submersion with total
space $P$ of dimension $n+k$ and the given Poisson manifold
$(M,\{\cdot,\cdot\})$ of dimension $n$ as basis. Further, let $\pe: E
\longrightarrow P$ be a vector bundle over the total space $P$ with
typical fibre $W$ of dimension $r$.

Then, the space of smooth sections $\Gamma^\infty(P,E)$ of the vector
bundle $\pe: E\longrightarrow P$ is a right $C^\infty(M)$-module with
the pointwise module structure
\begin{equation}
    \label{eq:right-module-sections}
    \rho_0(s,a)= \pp^*a\cdot s 
\end{equation}
for all $s\in \Gamma^\infty(P,E)$ and $a\in C^\infty(M)$ where one
makes use of the pointwise $C^\infty(P)$-module structure and the
pullback $\pp^*: C^\infty(M)\longrightarrow C^\infty(P)$ which is an
algebra morphism. Thus it is obvious that one is in the situation of
Definition~\ref{definition:module-type-diff} and that $\rho_0$ is of
the differential type with orders $L_\rho=l_\rho=0$. In particular,
$\Gamma^\infty(P,E)$ is a $(C^\infty(P),C^\infty(M))$-bimodule. This
and the commutativity of the pointwise multiplication make it possible
to equip $\mathcal{D}^l=
\Diffop_{C^\infty(P)}(\Gamma^\infty(P,E);\Gamma^\infty(P,E))$ with the
$(C^\infty(M), C^\infty(M))$-bimodule structure
\eqref{eq:new-bimodule-of-endomorphism}, this means
\begin{equation}
    \label{eq:bimodule-diffop}
    (a\cdot D\cdot b)(s) = \pp^*a \cdot D(\pp^*b \cdot s) 
\end{equation}
for all $a,b\in C^\infty(M)$, $s \in \Gamma^\infty(P,E)$, and $D\in
\mathcal{D}$. As already stated in
Section~\ref{subsec:Hochschild-complex-module} this convention does
not affect any considerations but is convenient and simplifies many
considerations, confer
Lemma~\ref{lemma:special-situation-free-module}.

Using these structures we are of course interested in the deformations
$\rho=\sum_{r=0}^\infty \lambda^r \rho_r$ of
\eqref{eq:right-module-sections} of the induced differential type.
This means that the $\rho_r$ are 1-cochains
\begin{equation}
    \label{eq:differential-Hochschild-sections}
    \rho_r\in \Union_{
      \begin{subarray}{c}
          L\in \mathbb{N}_0\\
          l\in \mathbb{N}_0
      \end{subarray}
    }
    \Diffop^L_{C^\infty(M)} (C^\infty(M);
    \Diffop^l_{C^\infty(P)} (\Gamma^\infty(P,E));\Gamma^\infty(P,E))
\end{equation}
in the crucial complex
\begin{equation}
    \label{eq:Hochschild-complex-sursub}
    (\HCdiff^\bullet(C^\infty(M),\mathcal{D}),\delta),
\end{equation}
where $\delta$ is induced by \eqref{eq:bimodule-diffop}. Due to the
observation made in
Corollary~\ref{corollary:Differential-Hochschild-module}, the
Hochschild differential respects the left $C^\infty(P)$-module
structure of all vector spaces $\HCdiff^k(C^\infty(M),\mathcal{D})$,
this means
\begin{equation}
    \label{eq:Hochschild-differential-functions}
    \delta (f\cdot \phi)= f\cdot \delta \phi 
\end{equation}
for all $f\in C^\infty(P)$, $\phi\in \HCdiff^\bullet(C^\infty(M),
\mathcal{D})$. Due to the results of
Chapter~\ref{cha:deformation-algebras-modules}, especially of the
Sections~\ref{sec:deformation-modules} and
\ref{sec:commutant-module-structures}, the first and the second
cohomology groups of this differential Hochschild complex represent
the obstruction for an order by order construction of the desired
deformations and the equivalence transformations between them.

The sections $\Gamma^\infty(P,E)$ can of course be restricted to open
subsets $\tilde{U}\subseteq P$ and it is clear by considering the
restriction $\pp|_{\tilde{U}}$ that the sections
$\Gamma^\infty(\tilde{U},\pe^{-1}(\tilde{U}))$ have analogous module
structures with respect to $C^\infty(\tilde{U})$ and
$C^\infty(\pp(\tilde{U}))$. The common restrictions of maps are
obviously compatible with the considered pointwise structures. Thus,
the surjective submersion together with the sheaves
\begin{equation}
    \label{eq:algebraic-structures-vector-bundle-sursub}
    \mathcal{E}= \Gamma^\infty_{PE} \quad 
    \textrm{and} \quad 
    \mathcal{A}= C^\infty_M
\end{equation}
of sections of $E$ over $P$ and functions over $M$ are a special
example of structures as in
Remark~\ref{remark:structures-sheaves-differential} and one is further
able to apply the results of
Proposition~\ref{proposition:sheaves-and-cohomology} which will be
very important for the computation of the cohomology.

\begin{remark}
    \label{remark:special-case-sections-are-functions}
    Note that for the above geometric situation of a surjective
    submersion $\pp: P\longrightarrow M$ one could already consider
    $C^\infty(P)$ as a right $C^\infty(M)$-module with
    \begin{equation}
        \label{eq:right-module-functions}
        \rho_0(f,a)= \pp^*a\cdot f 
    \end{equation}
    for all $f\in C^\infty(P)$ and $a\in C^\infty(M)$. This is of
    course a special case of \eqref{eq:right-module-sections} for the
    trivial vector bundle $E= P\times \mathbb{C}$ since
    $\Gamma^\infty(P,P\times \mathbb{C})\cong C^\infty(P)$. The
    corresponding complex then is
    $\HCdiff^\bullet(C^\infty(M),\mathcal{D})$ with $\mathcal{D}=
    \Diffop^\bullet_{C^\infty(P)}(C^\infty(P),C^\infty(P))$.
\end{remark}

\begin{remark}[The zeroth cohomology]
    As seen in Section~\ref{sec:commutant-module-structures}, the
    zeroth Hochschild cohomology is given by the commutant of the
    right module structure within the differential operators. In the
    considered geometrical context it is clear that the elements of
    the commutant are nothing but the \emph{vertical differential
      operators} $D \in
    \Diffopver(\Gamma^\infty(P,E),\Gamma^\infty(P,E))$ which are
    defined by the condition
    \begin{equation}
        \label{eq:vertical-diffops}
        D(\pp^*a \cdot s)=
        \pp^*a \cdot D(s)
    \end{equation}
    for all $a\in C^\infty(M)$ and $s\in \Gamma^\infty(P,E)$.
\end{remark}

\subsection{Deformation theory on principal fibre bundles}
\label{subsec:geometric-situation-principal-bundles}

Since every principal fibre bundle $\pp: P\longrightarrow M$ is also a
surjective submersion the above considerations still hold in this
case. However, being given the additional principal right action $\rr:
G\longrightarrow \Aut(P)$ of the structure group $G$ it is natural to
consider so-called \emph{equivariant vector bundles} over $P$. This
means that $G$ acts on the total space $E$ of the vector bundle $\pe:
E\longrightarrow P$ by vector bundle automorphisms over the principal
right action or the corresponding left action. First, we consider the
case of a right action
\begin{equation}
    \label{eq:right-action-vector-bundle}
    \rR: G\longrightarrow \Aut(\mathcal{E})
\end{equation}
with the property that $\rR_g: E\longrightarrow E$ is fibrewise
linear and
\begin{equation}
    \label{eq:right-action-over-principal-action}
    \pe\circ \rR_g = \rr_g \circ \pe
\end{equation}
for all $g\in G$.
Equation~\eqref{eq:right-action-over-principal-action} can be seen as
the $G$-invariance of $\pe$.

A right action $\rR$ on a vector bundle $E$ as above naturally induces
a corresponding left action
\begin{equation}
    \label{eq:left-action-vector-bundle}
    \lL: G\longrightarrow \Aut(E^*)
\end{equation}
on the dual vector bundle $\pe: E^*\longrightarrow P$ whose fibres
$E^*_u$ over $u\in P$ are the dual spaces of the fibres $E_u$ of
$E$. For an element $g\in G$ the action is simply given by the
transposed maps $\lL_g= \rR_g^*: E_{ug}^*\longrightarrow E^*_u$
defined by $(\lL_g v^*)(v)=v^*(\rR_g v)$. It is obvious that the so
defined maps $\lL_g: E^*\longrightarrow E^*$ again are fibrewise
linear and now satisfy
\begin{equation}
    \label{eq:left-action-over-principal-action}
    \pe\circ \lL_g = \rr_{g^{-1}} \circ \pe.
\end{equation}
Of course, such a left action induces a right action in the analogous
way. In the following it will not play any role which kind of action
on $E$, left or right, we consider. Together with the principal right
action $\rr$ both of them induce left representations of $G$ on the
space of sections of $\pp: E \longrightarrow P$. Explicitly, one has
for all $s\in \Gamma^\infty(P,E)$
\begin{eqnarray}
    \label{eq:action-on-sections-coming-from-right-action}
    g \acts s&=& \rR_{g^{-1}} \circ s \circ \rr_g \quad \textrm{or}\\ 
    \label{eq:action-on-sections-coming-from-left-action}
    g \acts s&=& \lL_g \circ s \circ \rr_g.
\end{eqnarray}

It is remarkable that the properties of the principal right action
$\rr$, this means to be \emph{free} and \emph{proper}, confer
\cite[Chap.~9]{lee:2003a} or \cite{duistermaat.kolk:2000a}, are always
inherited by the action of an equivariant vector bundle.

\begin{lemma}[Equivariant vector bundles over principal fibre bundles]  
    \label{lemma:free-proper}
    If $\pe: E\longrightarrow P$ is an equivariant vector bundle over
    a principal fibre bundle $\pp:P\longrightarrow M$ the action of
    the structure Lie group $G$ on $E$ is smooth, free and proper.
\end{lemma}

\begin{proof}
    Without loss of generality one can assume to be given a right
    action $\rR$ as in \eqref{eq:right-action-vector-bundle}.
    Smoothness is clear by definition. That $\rR$ is free, this means
    that $g\acts e \neq e$ for all $e\in E$ and all group elements
    $g\in G$ except the neutral element, is a consequence of the
    property \eqref{eq:right-action-over-principal-action} and the
    freeness of $\rr$. Properness has several equivalent
    definitions. For our purpose we use the characterization of
    \cite[Prop.~9.13]{lee:2003a}. According to the results there, the
    action $\rR$ is proper if the following condition is satisfied.
    Given a convergent sequence $\{e_n\}$ in $E$ and a sequence
    $\{g_n\}$ in $G$ such that the sequence $\{\rR_{g_n}e_n\}$
    converges, then there exists a convergent subsequence of
    $\{g_n\}$. That this really is the case can be seen very
    easily. Due to the continuity of the projection $\pe$ and the
    property \eqref{eq:right-action-over-principal-action}, the
    convergent sequences $\{e_n\}$ and $\{\rR_{g_n}e_n\}$ in $E$
    induce convergent sequences $\{ \pe (e_n)\}$ and $\{ \pe(\rR_{g_n}
    e_n)= \rr_{g_n}e_n\}$ in $P$. But then, the properness of the
    right action $\rr$ yields the assertion.
\end{proof}

Since $\rR_g$ and $\lL_g$ are fibrewise linear these definitions yield
representations and the module structures are all $G$-invariant. One
has $g\acts (fs)=(g\acts f) (g\acts s)$ where the action on $f\in
C^\infty(P)$ is given by pullback, $g\acts f= \rr_g^*f$. Regarding the
right module structure \eqref{eq:right-module-sections} as cochain and
using the induced action on maps, this implies
\begin{equation}
    \label{eq:right-module-G-invariant}
    g\acts \rho_0= \rho_0.
\end{equation}
Of course we are interested in deformations of the same type. Since
the action on $C^\infty(M)$ is trivial we are in the situation of
Remark~\ref{remark:action-only-on-module} and have to compute the
cohomology of the complex
\begin{equation}
    \label{eq:Hochschild-complex-principal-G-invariant}
    (\HCdiff^\bullet(C^\infty(M), \mathcal{D}^G),\delta).
\end{equation}
Of course, all occurring structures are $G$-invariant. In particular,
the zeroth cohomology is given by the $G$-invariant vertical
differential
operators $\Diffopver(\Gamma^\infty(P,E);\Gamma^\infty(P,E))^G$.

The properties \eqref{eq:right-action-over-principal-action} and
\eqref{eq:left-action-over-principal-action} imply that the orbits of
$\rR$ or $\lL$ are subsets of the preimages under $\pe$ of orbits of
$\rr$. This shows that the right actions are compatible with the
restrictions of maps to $\tilde{U}=\pp^{-1}(U)\subseteq P$ with open
subsets $U\subseteq M$. Altogether, the sheaves
$\mathcal{E}=\Gamma^\infty_{PE}$ and $\mathcal{A}= C^\infty_M$ now
satisfy the conditions of
Remark~\ref{remark:structures-sheaves-differential-G-inv} and we can
use the results of
Proposition~\ref{proposition:sheaves-and-cohomology-G-inv}.

\begin{remark}[Submodules of smooth sections]
    \label{remark:submodules-of-sections}
    In some relevant examples we will not be interested in
    deformations of the right $C^\infty(M)$-module
    $\Gamma^\infty(P,E)$ of all sections but rather in the
    deformations of a particular submodule
    \begin{equation}
        \label{eq:submodule-of-sections}
        \hat{\Gamma}^\infty(P,E)\subseteq \Gamma^\infty(P,E)
    \end{equation}
    consisting of sections with a further particular property. If this
    submodule gives rise to a subsheaf $\hat{\Gamma}_{PE}^\infty$ of
    substructures with all properties of $\Gamma_{PE}^\infty$ it is
    evident that all considerations are absolutely the same. In the
    following we will thus omit the extra notation.
\end{remark}

\subsection{Notation and specific properties of the relevant complexes }
\label{subsec:specific-properties-complexes}

Before we come to the computation of the cohomologies of the complexes
\eqref{eq:Hochschild-complex-sursub} and
\eqref{eq:Hochschild-complex-principal-G-invariant} we investigate
some further properties of their cochains arising in this specific
context.

\begin{remark}[Notation]
    From now on it will be convenient to use the following
    abbreviations: For any manifold $N$ and any
    $C^\infty(N)$-(bi)module $\mathcal{M}$ we set
    \begin{equation}
        \label{eq:abbreviation-Diffop-endo}
        \Diffop^\bullet (\mathcal{M})
        = \Diffop^\bullet_{C^\infty(N)} (\mathcal{M};\mathcal{M})
    \end{equation}
    and in particular
    \begin{equation}
        \label{eq:abbreviation-Diffop-functions}
        \Diffop^\bullet (N)
        = \Diffop^\bullet_{C^\infty(N)} (C^\infty(N),C^\infty(N)).
    \end{equation}
    Further, we set
    \begin{equation}
        \label{eq:abbreviation-Diffop-cochain}        
        \Diffop^L (N, \mathcal{M}) = 
        \Diffop^L_{C^\infty(N)} (\underbrace{C^\infty(N), \dots,
          C^\infty(N)}_{k\textrm{-times}}; \mathcal{M})  
    \end{equation}
    for $L = (l_1, \ldots, l_k)\in \mathbb{N}_0^k$, $k\in \mathbb{N}$,
    and
    \begin{equation}
        \label{eq:abbreviation-complex}
        \HC^\bullet(N,\mathcal{M})=\HC^\bullet (C^\infty(N),\mathcal{M}).
    \end{equation}
    For the structures of the previous section the meaning will then
    always be clear from the context.
\end{remark}

The considered cochains are given by differential operators and due to
the results of Section~\ref{sec:sheaves-Hochschild} these operators
define presheaves. Thus, the cochains can be restricted to open
subsets $\tilde{U}\subseteq P$ using the natural restriction maps
introduced in the
Sections~\ref{sec:sheaves-local-differential-operators} and
\ref{sec:sheaves-different-manifolds}. In the $G$-invariant setting
these restrictions are only allowed to open subsets $\pp^{-1}(U)$ with
$U\subseteq M$. In any case, if $U=\pp(\tilde{U})$ is the domain of a
local chart for $M$ the differential operators have the expected local
expressions. We prove the following general statement for the presheaf
of differential operators with values in differential operators.

\begin{lemma}
    \label{lemma:local-form-of-Diffops}
    Let $\tilde{U}\subseteq P$ be an open subset such that for
    $U=\pp(\tilde{U})$ there is a local chart $(U,x)$ of $M$. Further,
    let there be given $\phi \in \Diffop^L(M, \Diffop^\bullet
    (\Gamma^\infty(P,E)))$ with multiorder $L = (l_1, \ldots, l_k)\in
    \mathbb{N}_0^k$ of differentiation. Then there exist uniquely
    defined differential operators $\phi_{\tilde{U}}^{\alpha_1 \cdots
      \alpha_k}\in \Diffop^\bullet
    (\Gamma^\infty(\tilde{U},E|_{\tilde{U}}))$ with multiindices
    $\alpha_i\in \mathbb{N}_0^n$, $i=1,\ldots, k$, such that
    \begin{equation}
        \label{eq:local-form-of-Diffop}
        \phi|_{\tilde{U}} (a_1,\dots,a_k) =
        \sum_{|\alpha_1|\le l_1,\dots, 
          |\alpha_k|\le l_k}
        \left(\frac{\partial^{|\alpha_1|}a_1} {\partial x^{\alpha_1}}
            \dots 
            \frac{\partial^{|\alpha_k|}a_k}{\partial
              x^{\alpha_k}}\right) \cdot \phi_{\tilde{U}}^{\alpha_1
          \dots \alpha_k}
    \end{equation}
    for all $a_1,\ldots, a_k\in C^\infty(U)$.

    For $\phi \in \Diffop^L(M, \Diffop^\bullet
    (\Gamma^\infty(P,E))^G)$ and $\tilde{U}=\pp^{-1}(U)$ the same
    assertion is true with $G$-invariant operators
    $\phi_{\tilde{U}}^{\alpha_1 \cdots \alpha_k}\in \Diffop^\bullet
    (\Gamma^\infty(\tilde{U},E|_{\tilde{U}}))^G$.
\end{lemma}

\begin{proof}
    The proof is a straightforward generalization of the well-known
    considerations for ordinary differential operators. The assertion
    is shown by an induction over $|L|$. For $|L|=0$ one obviously has
    $\phi|_{\tilde{U}} (a_1,\dots,a_k)= (a_1 \dots a_k)\cdot
    \phi|_{\tilde{U}} (1,\dots,1)$. Due to the results of
    Chapter~\ref{cha:sheaves}, in particular the ones of
    Section~\ref{sec:sheaves-Hochschild}, the uniqueness of the local
    expression $\phi|_{\tilde{U}}$ is already clear. Thus we can
    assume that $x(U)\subseteq \mathbb{R}^n$ is convex. Now, for
    $|L|\ge 0$ consider an arbitrary point $u_0\in \tilde{U}$ and
    $\pp(u_0)=p_0\in U$. Further consider an arbitrary section $s\in
    \Gamma^\infty(\tilde{U}, E|_{\tilde{U}})$ and functions
    $a_1,\dots, a_k\in C^\infty(U)$.  The assumed convexity assures
    that for all points $p\in U$ and $t\in [0,1]$ the convex
    combination $tx(p)+ (1-t)x(p_0)$ is again in the range of the
    chart $x$ and we can use Hadamard's trick. For $l=1,\dots,k$ one
    has
    \begin{align*}
        a_l(p) &= (a_l \circ x^{-1})(x(p_0))+ \int\limits_0^1 \D t \:
        \frac{\D}{\D t}\left( (a_l\circ x^{-1})(t x(p) + (1-t)
            x(p_0))\right)\\
        &= a_l(p_0)+ \int\limits_0^1 \D t\:\frac{\partial (a_l\circ
          x^{-1})} {\partial x^i} (t x(p) + (1-t)
        x(p_0)) (x^i(p)-x^i(p_0))\\
        &= a_l(p_0) + b_i^l(p) (x^i(p)-x^i(p_0)).
    \end{align*}
    Thus one has $a_l=a_l(p_0)+ b_i^l\cdot (x^i-x^i(p_0))$. The
    functions $b_i^l\in C^\infty(U)$ given by
    $b_i^l(p)=\int\limits_0^1 \D t\:\frac{\partial (a_l\circ x^{-1})}
    {\partial x^i} (t x(p) + (1-t) x(p_0))$ satisfy
    \begin{equation}
        \label{eq:hadamard-ableitung}
        \frac{\partial^{|\alpha|} b_i^l}{\partial x^{\alpha}}(p_0) =
        \int\limits_0^1 \D t\: \frac{\partial^{|\alpha|+1} (a_i^l\circ
          x^{-1})}{\partial x^{\alpha}\partial x^i}(tx(p)+(1-t)x(p_0))
        t^{|\alpha|}|_{p=p_0}= \frac{1}{|\alpha|+1}
        \frac{\partial^{|\alpha|+1} a_l} {\partial x^{\alpha}\partial
          x^i} (p_0) 
    \end{equation}
    for all multiindices $\alpha\in \mathbb{N}_0^n$. Using the
    linearity properties one then has
    \begin{eqnarray*}
        \lefteqn{((\phi|_{\tilde{U}}(a_1,\dots,a_k))(s)) (u_0)}\\
        &=& \phi|_{\tilde{U}} (a_1(p_0)+ b_i^1\cdot
        (x^i-x^i(p_0)),a_2, \dots,a_k) (s)(u_0)\\
        &=& a_1(p_0) \phi|_{\tilde{U}} (1,a_2, \dots,a_k) (s)(u_0) +
        \left( \phi|_{\tilde{U}}\circ \mathsf{L}^{(1)}_{x^i-x^i(p_0)}
            - \mathsf{L}_{x^i-x^i(p_0)} \circ \phi|_{\tilde{U}}\right)
        (b_i^1,a_2, \dots a_k)(s)(u_0),
    \end{eqnarray*}
    where the term $\mathsf{L}_{x^i-x^i(p_0)} \circ \phi|_{\tilde{U}}$
    could be added since it vanishes in the considered expression due
    to the module structure, $((a\cdot D) (s))(u_0)= a(p_0) D(s)(u_0)$
    with $D\in \Diffop^\bullet (\Gamma^\infty(\tilde{U},
    E|_{\tilde{U}}))$ and $a\in C^\infty(U)$. After setting
    \begin{equation*}
        S^i_l= \phi|_{\tilde{U}}\circ
        \mathsf{L}^{(l)}_{x^i-x^i(p_0)} - \mathsf{L}_{x^i-x^i(p_0)} \circ
        \phi|_{\tilde{U}} \in \Diffop^{L-\mathsf{e}_l}(C^\infty(U),
        \Diffop^\bullet(C^\infty(\tilde{U}))),
    \end{equation*}
    an iteration yields
    \begin{eqnarray*}
        ((\phi|_{\tilde{U}}(a_1,\dots,a_k))(s))(u_0)
        &=& a_1(p_0)\dots a_k(p_0)
        \phi|_{\tilde{U}}(1,\dots,1)(s)(u_0)\\
        && + \sum_{l=1}^k a_1(p_0) \dots a_{l-1}(p_0)
        S^i_l(1,\dots,1,b_i^l,a_{l+1}, \dots a_k)(s)(u_0).
    \end{eqnarray*}
    Using the induction hypothesis for the $S^i_l$ and the property
    \eqref{eq:hadamard-ableitung}, the form
    \eqref{eq:local-form-of-Diffop} follows since $s$ and $u_0$ were
    arbitrary. The assertion concerning the $G$-invariance is obvious.
\end{proof}

\begin{remark}
    \label{remark:local-form-fixed-order}
    It is easy to see that Lemma~\ref{lemma:local-form-of-Diffops}
    even holds for Hochschild cochains, this means for $\phi\in
    \Diffop^L(M, \Diffop^l (\Gamma^\infty(P,E)))$ with a fixed $l\in
    \mathbb{N}_0$. Then, the assertion is true with
    $\phi_{\tilde{U}}^{\alpha_1 \cdots \alpha_k}\in \Diffop^l
    (\Gamma^\infty(\tilde{U},E|_{\tilde{U}}))$. Of course the same is
    true for $G$-invariant cochains.
\end{remark}

The local expressions in Lemma~\ref{lemma:local-form-of-Diffops} have
the interesting consequence that the considered differential operators
$\phi$ already have values in the differential operators of some fixed
order and thus are Hochschild cochains. In order to prove this we need
the following simple lemma.

\begin{lemma}
    \label{lemma:smallest-possible-order}
    Let $M$ be a manifold and let $D\in \Diffop^l(M)$ where $l\in
    \mathbb{N}_0$ is the smallest possible order of differentiation to
    choose, this means $D\notin \Diffop^{l-1}(M)$. Further let
    $\{U_i\}_{i\in I}$ be an open cover of $M$. Then the corresponding
    restrictions of the operator $D$ are differential operators
    $D|_{U_i}\in \Diffop^{l_i}(U_i)$ with $l_i\le l$ for all $i\in I$
    and there exists at least one $i_0\in I$ with $D|_{U_{i_0}}\in
    \Diffop^l(U_{i_0})$ where $l$ is the smallest possible order.
\end{lemma}

\begin{proof}
    The assertion is clear with the sheaf property of differential
    operators of some fixed order of differentiation. If there existed
    an $\tilde{l}\in \mathbb{N}_0$ with $l_i\le \tilde{l} <l$ for all
    $i\in I$ this would imply that $D\in \Diffop^{\tilde{l}}(M)$ which
    is a contradiction to the choice of $l$.
\end{proof}

\begin{lemma}
    For all $k\in \mathbb{N}$ and $L\in \mathbb{N}_0^k$ one has
    \begin{equation}
        \label{eq:diffops-values-finite-order}
        \Diffop^L (M, \Diffop^\bullet (\Gamma^\infty(P,E)))
        = \Union_{l\in \mathbb{N}_0} \Diffop^L(M,
        \Diffop^l (\Gamma^\infty(P,E))).
    \end{equation}
    The same is true for operators with values in $G$-invariant maps.
\end{lemma}

\begin{proof}
    The nontrivial inclusion to show is
    \begin{equation*}
        \Diffop^L (M, \Diffop^\bullet (\Gamma^\infty(P,E)))
        \subseteq \Union_{l\in \mathbb{N}_0} \Diffop^L(M,
        \Diffop^l (\Gamma^\infty(P,E))).
    \end{equation*}
    Let $\phi\in \Diffop^L (M, \Diffop^\bullet (\Gamma^\infty(P,E)))$
    with $L\in \mathbb{N}_0^k$ and $k\in \mathbb{N}_0$. Using the
    abbreviation $\mathcal{D}^l(P) = \Diffop^l (\Gamma^\infty(P,E))$
    one has to show the following assertion:

    There exists an $l\in \mathbb{N}_0$ such that for all $a_1,\dots,
    a_k\in C^\infty(M)$ the smallest possible
    $\tilde{l}(a_1,\dots,a_k)\in \mathbb{N}_0$ with $\phi(a_1,\dots,
    a_k)\in \mathcal{D}^{\tilde{l}(a_1,\dots,a_k)}(P)$ is bounded by
    $l$, this means $\tilde{l}(a_1,\dots,a_k)<l$.

    Now assume that this is not the case.
    % Then for all $l\in \mathbb{N}_0$ there exist $a_1^{(l)},\dots,
    % a_k^{(l)}\in C^\infty(M)$ such that the smallest possible
    % $\tilde{l}(a_1,\dots,a_k)\in \mathbb{N}_0$ with $\phi(a_1,\dots,
    % a_k)\in \mathcal{D}^{\tilde{l}(a_1,\dots,a_k)}$ is greater than
    % $l$, $\tilde{l}(a_1,\dots,a_k)\in \mathbb{N}_0> l$.
    Then there exists a strictly monotonically increasing sequence
    $\{l_n\}_{n\in \mathbb{N}}$ with $l_n\in \mathbb{N}_0$ such that
    for all $n\in\mathbb{N}$ there exist $a_1^{(n)},\dots,
    a_k^{(n)}\in C^\infty(M)$ with $\phi(a_1^{(n)},\dots,
    a_k^{(n)})\in \mathcal{D}^{l_n}(P)\setminus
    \mathcal{D}^{l_n-1}(P)$. Due to
    Lemma~\ref{lemma:smallest-possible-order} there exist charts
    $(U_n,x_n)$ of $M$ with $\phi(a_1^{(n)},\dots,
    a_k^{(n)})|_{\pp^{-1}(U_n)} \in \mathcal{D}^{l_n}(\pp^{-1}(U_n))
    \setminus \mathcal{D}^{l_n-1}(\pp^{-1}(U_n))$. Since
    $\Diffop^\bullet(\Gamma^\infty(P,E))\subset
    \Loc(\Gamma^\infty(P,E))$ is a subpresheaf one has
    $\phi(a_1,\dots, a_k)|_{\pp^{-1}(U_n)}=
    \phi|_{\pp^{-1}(U_n)}(a_1|_{U_n},\dots, a_k|_{U_n})$ and
    \eqref{eq:local-form-of-Diffop} then shows that the smallest
    possible degree of $\phi(a_1,\dots, a_k)|_{\pp^{-1}(U_n)}$ for all
    $a_1,\dots, a_k$ can not be greater than $l_n$. Thus one can
    assume that the open subsets for different minimal orders are
    disjoint, $U_n\cap U_m= \emptyset$ for $n\neq m$. Once again with
    Lemma~\ref{lemma:smallest-possible-order} we choose for all $n\in
    \mathbb{N}$ open subsets $V_n\subseteq U_n$ with
    $V_n^{\cl}\subseteq U_n$ and $\phi(a_1^{(n)},\dots,
    a_k^{(n)})|_{\pp^{-1}(V_n)} \in \mathcal{D}^{l_n} (\pp^{-1}(V_n))$
    where $l_n$ is minimal. For all $n\in \mathbb{N}$ we now choose
    $\chi_n\in C^\infty(M)$ with $\chi_n|_{V_n}=1$ and $\supp (\chi_n)
    \subseteq U_n$ and consider the well-defined functions $b_s=
    \sum_{n=1}^\infty \chi_n a_s^{(n)}\in C^\infty(M)$ with
    $b_s|_{V_n}=a_s^{(n)}|_{V_n}$ for all $s=1,\dots,k$. We then get
    \begin{equation*}
        \phi(b_1,\dots, b_k)|_{\pp^{-1}(V_n)}= \phi|_{\pp^{-1}(V_n)}
        (a_1^{(n)}|_{V_n}, \dots, a_k^{(n)}|_{V_n})=
        \phi(a_1^{(n)},\dots, a_k^{(n)})|_{\pp^{-1}(V_n)} \in
        \mathcal{D}^{l_n} (\pp^{-1}(V_n)),
    \end{equation*}
    where $l_n$ is the smallest possible order. Thus, the smallest
    possible order of $\phi(b_1,\dots, b_k)$ is at least $l_n$. Since
    $l_n \rightarrow \infty$ for $n \rightarrow \infty$, it is clear
    that $\phi(b_1,\dots,b_k)\in \Loc(\Gamma^\infty(P,E))$ has no
    well-defined degree of differentiation, $\phi(b_1,\dots,b_k)
    \notin \Diffop^\bullet(\Gamma^\infty(P,E))$. This is a
    contradiction to the assumption on $\phi$. The assertion for
    $G$-invariant maps is a direct consequence of
    \eqref{eq:diffops-values-finite-order} itself.
\end{proof}

This result shows that in the present situation the characterization
of the differential type is easy and one has not to take care of the
subtleties occurring in the general situation.

\begin{remark}[The natural notion of differential right module
    structures]
    \label{remark:notion-differential-right-module}
    The local formulas obtained in the considered situation show that
    our definitions really discuss differential right module
    structures in the natural and desired sense. If
    $\tilde{U}\subseteq P$ in Lemma~\ref{lemma:local-form-of-Diffops}
    is the domain of a local chart $(\tilde{U},z)$ on $P$ such that
    one in addition has a local $C^\infty(\tilde{U})$-module basis $\{
    s_i\}_{i=1,\dots,r}$ for $\Gamma^\infty(\tilde{U},E|_{\tilde{U}})$
    one can further consider the local expression for the operators
    $\Diffop^l(\Gamma^\infty(P,E))$, confer
    \cite[Thm.~A.5.2]{waldmann:2007a}. Regardless of a possibly
    existing $G$-invariance all $1$-cochains $\rho\in
    \Diffop^L(M,\Diffop^l(\Gamma^\infty(P,E)))$ then have the
    following local form for all $a\in C^\infty(\pp(\tilde{U}))$ and
    $s=f^is_i\in \Gamma^\infty(\tilde{U},E|_{\tilde{U}})$ with $f^i\in
    C^\infty(\tilde{U})$,
    \begin{equation}
        \label{eq:local-form-right-module}
        \rho_{\tilde{U}}(s,a)= \rho|_{\tilde{U}}(a)(s) =
        \sum_{|\alpha|\le L} 
        \left(
            \frac{\partial^{|\alpha|}a}{\partial
              x^\alpha} \cdot \phi^\alpha_{\tilde{U}}
        \right)(s)%\nonumber \\
        = \sum_{
          \begin{subarray}{c}
              |\alpha|\le L\\
              |\beta|\le l
          \end{subarray}
        }
        \left(\pp^*\frac{\partial^{|\alpha|}a} {\partial x^\alpha}
        \right) 
        \cdot 
        \frac{\partial^{|\beta|} f^i}{\partial z^\beta} \cdot
        s^{\alpha\beta}_i. 
    \end{equation}
    There, $s^{\alpha\beta}_i\in \Gamma^\infty(\tilde{U},
    E|_{\tilde{U}})$ are uniquely defined sections. Note the different
    kinds of multiindices
    $\alpha\in \mathbb{N}_0^n$ and $\beta\in \mathbb{N}_0^{n+k}$.

    Naively, one would have defined the differential right module
    structures to be operators in
    \begin{equation}
        \label{eq:naive-differential-right-modules}
        \Diffop^\bullet_{C^\infty(M)}(\Gamma^\infty(P,E),C^\infty(M);
        \Gamma^\infty(P,E)). 
    \end{equation}
    However, this does not lead to the desired local expression
    \eqref{eq:local-form-right-module}. In the principal bundle case
    the map
    \begin{equation}
        \label{eq:G-action-Diffop}
        D_g: (s,a)\longmapsto \pp^*a \cdot (g\acts s) 
    \end{equation}
    which is induced by the right action
    \eqref{eq:action-on-sections-coming-from-right-action} or
    \eqref{eq:action-on-sections-coming-from-left-action} defines such
    a differential operator $D_g\in
    \Diffop^{(0,0)}_{C^\infty(M)}(\Gamma^\infty(P,E),C^\infty(M);
    \Gamma^\infty(P,E))$ for all $g\in G$. But of course such right
    module structures containing a nonlocal shift by $g$ should not be
    called differential. These considerations justify and affirm the
    chosen approach to deformation theory of right modules of a
    particular type as introduced in
    Section~\ref{sec:algebras-modules-type}.
\end{remark}

\section{The basic steps of the computation of the cohomology}
\label{sec:basic-steps}

The strategy to compute the differential Hochschild cohomologies of
the complexes \eqref{eq:Hochschild-complex-sursub} and
\eqref{eq:Hochschild-complex-principal-G-invariant} basically consists
of four steps. Although the basic ideas are the same for both
complexes we have to discuss the two cases separately. The differences
between surjective submersions and principal fibre bundles caused by
the right action make it necessary to pursue slightly different paths.

\subsection{The basic steps for surjective submersions}
\label{subsec:basic-steps-surjective-submersion}

For surjective submersion the four steps are the following.

\begin{enumerate}
\item First, one makes use of the fact that the Hochschild complexes
    can be restricted to open subsets of an appropriate covering of
    $P$ and that one is able to consider the local situation.
\item It is a general fact that the sections of a vector bundle become
    a free module over the algebra of smooth functions on the base
    manifold if one considers the local situation with respect to a
    vector bundle chart. Using
    Lemma~\ref{lemma:special-situation-free-module}, we thus can
    consider the deformation problem for the right module of functions
    described in
    Remark~\ref{remark:special-case-sections-are-functions}.
\item Due to the fact that a surjective submersion allows a specific
    choice of charts, we can consider a surjective submersion $\pr_1:
    V\times G\longrightarrow V$ with a manifold $G$ and an open convex
    subset $V\subseteq \mathbb{R}^n$. In this case, it will be
    possible to compute the Hochschild cohomologies using the
    techniques developed in Chapter~\ref{cha:bar-koszul}.
\item In the last step we observe that the coboundaries which are used
    to compute the local cohomologies are images of cochains with
    globally bounded orders of differentiation. Thus, it will be
    possible to apply
    Proposition~\ref{proposition:sheaves-and-cohomology} and to
    compute the initial Hochschild cohomology.
\end{enumerate}

The reduction of the problem to the local situation described in the
third point makes use of the particular geometric situation. By the
constant rank theorem, confer \cite[Satz~5.4]{broecker.jaenich:1990a},
for any point $u \in P$ there exist open subsets $\tilde{U} \subseteq
P$ with $u \in \tilde{U}$ and $U \subseteq M$ with $p = \pp (u) \in U$
together with diffeomorphisms $x: U \longrightarrow V \subseteq
\mathbb{R}^n$ and $\tilde{x}: \tilde{U} \longrightarrow V \times G
\subseteq \mathbb{R}^{n+k}$, such that $\pp (\tilde{U}) = U$ and $x
\circ \pp \circ \tilde{x}^{-1} = \pr_1$. Clearly, in this case
$G\subseteq\mathbb{R}^k$ is an open subset. It is possible to choose
the open subsets in a way that $V \subseteq \mathbb{R}^n$ is
convex. Furthermore, one can achieve that in the considered situation
the vector bundle $\pe:E\longrightarrow P$ has a local trivialization
over $\tilde{U}$. This means that there exists a vector bundle chart
$\psi: E|_{\tilde{U}} \longrightarrow \tilde{U} \times W$ with $\pr_1
\circ \psi= \pe$. Such adapted local charts $(\tilde{U}, \tilde{x})$
of $P$ and $(U, x)$ of $M$ build adapted atlases and in particular
induce corresponding open coverings of $P$. These adapted structures
which will be used throughout the chapter induce the commutative
diagram
\begin{equation}
    \label{eq:diagram-surjective-submersion}
    \xymatrix{E \ar[d]_{\pe}&& E|_{\tilde{U}} \ar
      @{_{(}->}[ll]_{\iota} 
      \ar[d]_{\pe}  \ar[r]^{\psi} & **[r] \tilde{U}\times W
      \ar[dl]_{\pr_1}& \\
      P \ar[d]_{\pp} && \tilde{U} \ar @{_{(}->}[ll]_{\iota}
      \ar[d]_{\pp} \ar[rr]^{\tilde{x}}&& **[r] V\times G\subset
      \mathbb{R}^{n+k} 
      \ar[d]_{\pr_1}\\
      M && U \ar @{_{(}->}[ll]_{\iota} \ar[rr]^x &&
      **[r] V\subset \mathbb{R}^n,
    }
\end{equation}
where $\iota$ in every case denotes the embedding map. All spaces in
the first row again are vector bundles over the corresponding spaces
in the second row. Considering only the second and third row every
column of the diagram is the diagram of a surjective submersion with a
canonical right module structure for the corresponding algebras of
smooth functions as in \eqref{eq:right-module-functions}.

Using these adapted charts the reduction of the initial deformation
problem to the local one is realized as follows.

\begin{lemma}
    \label{lemma:localization}
    The restriction maps of sections and functions, the local vector
    bundle chart $\psi$ and the local charts $\tilde{x}$ and $x$ give
    rise to chain maps and chain isomorphisms
    \begin{equation}
        \label{eq:localization}
        \xymatrix{
          \HCdiff^\bullet (M,\Diffop
          (\Gamma^\infty(P,E)))\ar[d] & \\ 
          \HCdiff^\bullet (U,\Diffop
          (\Gamma^\infty(\tilde{U},E|_{\tilde{U}})))  
          \ar[r]^-{\cong} & 
          \Mat_{r\times r} \left(
              \HCdiff^\bullet (U,\Diffop(\tilde{U}))\right) \ar[dl] \\
          \HCdiff^\bullet (U,\Diffop (\tilde{U}))
          \ar[r]^-{\cong} &
          \HCdiff^\bullet
          (V,\Diffop (V\times G))
          .
        }
    \end{equation}
\end{lemma}

\begin{proof}
    As already mentioned one is in the situation of
    Section~\ref{subsec:sheaves-diff}. This in particular means that
    $\HCdiff^\bullet (M,\Diffop (\Gamma^\infty(P,E)))$ can be seen as
    a presheaf of Hochschild complexes over $P$. Thus the
    restriction map is a morphism of complexes.

    Choosing a basis $\{w_1,\dots, w_r\}$ of the typical fibre $W$ the
    homeomorphism $\psi: E|_{\tilde{U}}\longrightarrow \tilde{U}\times
    W$ gives rise to a corresponding $C^\infty(\tilde{U})$-module
    basis $\{s_1, \dots, s_r \}$ of
    $\Gamma^\infty(\tilde{U},E|_{\tilde{U}})$, thus becoming a
    finitely generated free module. With
    \eqref{eq:right-module-sections} and \eqref{eq:bimodule-diffop}
    one thus can apply Lemma~\ref{lemma:special-situation-free-module}
    and thus has justified the isomorphism with the complex of
    matrices. In Lemma~\ref{lemma:special-situation-free-module} it is
    further explained that the complex of matrices in fact is a matrix
    of one single complex. Thus the restriction to one
    of those copies is a surjective chain morphism.

    Finally, one considers the pullbacks of the inverse chart maps. It
    is simple to verify that $(\tilde{x}^{-1})^*: C^\infty(\tilde{U})
    \longrightarrow C^\infty(V\times G)$ is an algebra isomorphism
    and, in addition, a module isomorphism along the algebra
    isomorphism $(x^{-1})^*: C^\infty(U)\longrightarrow
    C^\infty(V)$. Then, the last isomorphism of complexes is a
    consequence of the functorial behaviour of differential Hochschild
    complexes as already stated in
    Remark~\ref{remark:functoriality-Hochschild-diff}.
\end{proof}

\subsection{The basic steps for principal fibre bundles}
\label{subsec:basic-steps-principal-bundles}

Considering equivariant vector bundles over principal fibre bundles we
have to modify the steps used for surjective submersions taking care
of the occurring actions of the structure group $G$.

\begin{enumerate}
\item In the $G$-invariant case one works with presheaves over the
    base manifold $M$. Thus the Hochschild complex
    \eqref{eq:Hochschild-complex-principal-G-invariant} can only be
    restricted to open subsets $U\subseteq M$ of an appropriate
    covering of $M$ and one considers the local situation with respect
    to such subsets $U$.
\item For the second step one has to work with $G$-invariant local
    module bases of the sections. Then, the additional assertion of
    Lemma~\ref{lemma:special-situation-free-module} assures that the
    initial problem for the sections can be traced back to that of the
    functions even in the $G$-invariant setting.
\item Since every principal fibre bundle has corresponding local
    bundle charts one can consider the trivial bundle $\pr_1: V \times
    G\longrightarrow V$ with the structure Lie group $G$ and an open
    convex subset $V\subseteq \mathbb{R}^n$. In this case the explicit
    computation of the cohomology will turn out to respect the
    additional $G$-invariance automatically.
\item Finally, it will be possible to apply
    Proposition~\ref{proposition:sheaves-and-cohomology-G-inv} and to
    compute the initial $G$-invariant cohomology.
\end{enumerate}

In detail, the reduction of the problem now makes use of the following
geometric situation. Given an arbitrary principal fibre bundle
$\pp:P\longrightarrow M$ it is clear by the very definition that for
any $p\in M$ there exists an open subset $U\subseteq M$ with $p\in U$
such that the bundle over $U$ is trivial, confer
Definition~\ref{definition:principal-bundle}. This means that there
exists a principal bundle chart $\varphi: P \supseteq \pp^{-1} (U)
\longrightarrow U\times G$ with the defining properties $\pr_1 \circ
\varphi =\pp$ and $\varphi \circ \rr_g = (\id_U \times r_g) \circ
\varphi$. The map $r_g: G\longrightarrow G$ is the right
multiplication with $g \in G$, this means $r_g(h) = hg$. Of course it
is possible to achieve that $U$ is the domain of a chart $x: U
\longrightarrow V \subseteq \mathbb{R}^n$ with the additional
properties that $U$ is contractible and $V\subseteq \mathbb{R}^n$ is
convex. In the following we will always work with open coverings of
$M$ consisting of such open subsets $U\subseteq M$.  Denoting the unit
element of the group by $e\in G$ one has that $U\times \{e\}$ and
$\varphi^{-1}(U\times \{e\})$ are contractible. Then it is clear that
the bundle $E$ is trivial over $\varphi^{-1}(U\times \{e\})$ such that
the sections have a module basis $\{\cc{s}_i\}_{i=1,\dots, r}$, confer
\cite[Thm.~3.4.35]{marsden.ratiu:2001}. This gives rise to a
$G$-invariant module basis $\{s_i\}_{i=1,\dots, r}$ of the sections
$\Gamma^\infty(\pp^{-1}(U),E|_{\pp^{-1}(U)})$ over $\pp^{-1}(U)= \varphi^{-1}(U\times
G)$ by setting
\begin{equation}
    \label{eq:invariant-local-module-basis}
    s_i(\varphi^{-1}(p,h))= \rR_h \cc{s}_i (\varphi^{-1}(p,e)) 
    \quad \textrm{or} \quad
    s_i(\varphi^{-1}(p,h))= \lL_{h^{-1}} \cc{s}_i (\varphi^{-1}(p,e)).  
\end{equation}
The invariance
\begin{equation}
    \label{eq:module-basis-G-inv}
    g\acts s_i=s_i
\end{equation}
for all $g\in G$ is a simple computation. Altogether, the above
structures and properties are summarized in the following commutative
diagram.

\begin{equation}
    \label{eq:diagram-principal-bundle}
    \xymatrix{
      % &E \ar'[d][dd]_{\pe} \ar@<-0.7mm>[dl]_{\rR_g} && E|_{\pp^{-1}(U)} \ar
      % @{_{(}->}[ll]_{\iota} 
      % \ar@<-2.ex>[dd]_{\pe}  \ar@<-0.7mm>[dl]_{\rR_g} &&&& \\
      % E \ar[dd]_{\pe} \ar@<-0.7mm>[ru]_{\lL_g} && E|_{\pp^{-1}(U)} \ar
      % @{_{(}->}[ll]_{\iota} \ar[dd]_(.7){\pe}
      % \ar@<-0.7mm>[ru]_{\lL_g} &&&&&\\
      &E \ar'[d][dd]_{\pe} \ar[dl]_{\rR_g}^{\lL_{g^{-1}}} &&
      E|_{\pp^{-1}(U)} \ar 
      @{_{(}->}[ll]_{\iota} 
      \ar@<-2.ex>[dd]_{\pe}  \ar[dl]_{\rR_g}^{\lL_{g^{-1}}} &&&& \\
      E \ar[dd]_{\pe} && E|_{\pp^{-1}(U)} \ar
      @{_{(}->}[ll]_{\iota} \ar[dd]_(.7){\pe}
      &&&&&\\ 
      & P \ar@{-->}[ddl]^(.7){\pp} \ar[dl]_{\rr_g} &&**[l] \pp^{-1} (U)
      \ar@{_{(}-->}'[l][ll]_{\iota} 
      \ar@{-->}[ddl]^(.7){\pp} \ar[rr]^{\varphi}
      \ar[dl]_{\rr_g} \ar@/_7mm/[uu]_(.7){s_i} &&  U\times 
      G  
      \ar@{-->}[ddl]^(.7){\pr_1} \ar[dl]_{\id_U \times r_g}
      \ar[rr]^{x\times \id_G}&& 
      V\times G 
      \ar@{-->}[ddl]^(.7){\pr_1} \ar[dl]_{\id_V \times r_g}\\
      P \ar[d]_{\pp} && \pp^{-1}(U) \ar @{_{(}->}[ll]_{\iota}
      \ar[d]_{\pp} \ar[rr]^{\varphi}
      \ar@/_7mm/[uu]_(.7){s_i} && U\times G 
      \ar[d]_{\pr_1} \ar[rr]^{x\times \id_G} &&  V\times G
      \ar[d]_{\pr_1} &\\
      M && U \ar @{_{(}->}[ll]_{\iota} \ar@{=}[rr] && U \ar[rr]^x &&
      **[r] V\subset \mathbb{R}^n.&
    }
\end{equation}
%\begin{equation}
%    \label{eq:diagram-principal-bundle-short}
%    \xymatrix{
%      E \ar@(l,u)[]^{\rr_g} \ar[dd]_{\pe} && E|_{\pp^{-1}(U)}
%      \ar@(l,u)[]^{\rr_g} \ar
%      @{_{(}->}[ll]_{\iota} \ar[dd]_{\pe}
%      &&&&&&\\
%      &&&&&&&&\\
%      P \ar@(l,u)[]^{\rr_g} \ar[d]_{\pp} && \pp^{-1}(U)
%      \ar@(l,u)[]^{\rr_g} \ar @{_{(}->}[ll]_{\iota} 
%      \ar[d]_{\pp} \ar[rr]^{\phi} \ar@/_7mm/[uu]_{\mathbf{s_i}} && U\times
%      G \ar@(ur,ul)[]_{\id_U\times r_g}
%      \ar[d]_{\pr_1} \ar[rr]^{x\times \id_G} &&  V\times G
%      \ar@(ur,ul)[]_{\id_V\times r_g} 
%      \ar[d]_{\pr_1} && \\
%      M && U \ar @{_{(}->}[ll]_{\iota} \ar@{=}[rr] && U \ar[rr]^x &&
%      **[r] V\subset \mathbb{R}^n,&&
%    }
%\end{equation}
All four triangles occurring in the diagram are of course diagrams of
principal fibre bundles. 
% Of course the $G$-invariant module basis defines a $G$-invariant
% vector bundle chart $\psi:E|_{\pp^{-1}(U)}\longrightarrow
% \pp^^{-1}(U)\times \mathbb{R}^k$ by setting $\psi(f^is_i(u))=(u,
% f^i\mathsf{e}_i)$.
With slightly different structures and morphisms the reduction of the
$G$-invariant Hochschild complex has more or less the same appearance
as before.

\begin{lemma}
    \label{lemma:localization-G-inv}
    The restriction maps of sections and functions, the $G$-invariant
    module basis, the principal bundle chart $\varphi$, and the local
    chart $x$ give rise to chain maps and chain isomorphisms
    \begin{equation}
        \label{eq:localization-G-inv}
        \xymatrix{
          \HCdiff^\bullet (M,\Diffop
          (\Gamma^\infty(P,E))^G)\ar[d] & \\ 
          \HCdiff^\bullet (U,\Diffop
          (\Gamma^\infty(\pp^{-1}(U),E|_{\pp^{-1}(U)}))^G)  
          \ar[r]^-{\cong} & 
          \Mat_{r\times r} \left(
              \HCdiff^\bullet (U,\Diffop(\pp^{-1}(U))^G)\right)
          \ar[dl] \\ 
          \HCdiff^\bullet (U,\Diffop (\pp^{-1}(U))^G)
          \ar[r]^-{\cong} &
          \HCdiff^\bullet (V,\Diffop (V\times G)^G). 
        }
    \end{equation}
\end{lemma}

\begin{proof}
    For the chain morphism given by the restriction map we now use the
    results of Section~\ref{subsec:sheaves-diff-G-inv}. With the
    $G$-invariant module basis and
    Lemma~\ref{lemma:special-situation-free-module} one gets the first
    isomorphism and the second chain morphism as in
    Lemma~\ref{lemma:localization}. For the last isomorphism one needs
    the $G$-invariant algebra morphism $(\varphi^{-1}\circ
    (x^{-1}\times \id_G))^*: C^\infty(\pp^{-1}(U))\longrightarrow
    C^\infty(V\times G)$ which is a module isomorphism along
    $(x^{-1})^*: C^\infty(U)\longrightarrow C^\infty(V)$ and
    Remark~\ref{remark:functoriality-Hochschild-diff}.
\end{proof}

\section{The local Hochschild cohomology $\HHdiff^\bullet (V,
  \Diffop(V\times G))$}
\label{sec:local-cohomology-sursub}

According to the considerations of
Section~\ref{subsec:basic-steps-surjective-submersion} we now have to
compute the Hochschild cohomologies $\HHdiff^\bullet (V,
\Diffop(V\times G))$ for the trivial surjective submersion $\pr_1:
V\times G\longrightarrow V$ with a convex subset $V\subseteq
\mathbb{R}^n$ and a smooth manifold $G$. It should be emphasized
already at the beginning that in the following we do not use the fact
that in the situation of
Section~\ref{subsec:basic-steps-surjective-submersion} $G$ is a subset
of $\mathbb{R}^k$. Instead, $G$ can be an arbitrary manifold. Thus it
will be very easy to see that the situation of
Section~\ref{subsec:basic-steps-principal-bundles} where $G$ is a Lie
group and where one has to compute the $G$-invariant Hochschild
cohomology $\HHdiff^\bullet(V,\Diffop(V\times G)^G)$ for the trivial
principal fibre bundle $\pr_1: V\times G \longrightarrow V$, can be
treated in the very same way. It will be possible to implement the
additional $G$-invariance without any restrictions.

We will often take advantage of the fact that an atlas of local charts
$(\hat{U},y)$ of $G$ induces an atlas of adapted local charts
\begin{equation}
    \label{eq:adapted-chart-V-G}
    (\tilde{U}=V\times
    \hat{U},\tilde{x}=\id_V\times y) 
\end{equation}
of $V\times G$.

The computation of the cohomology shall be performed using the
arguments and techniques of Chapter~\ref{cha:bar-koszul}, in
particular Section~\ref{sec:homotopy}. In order to do this, one first
has to show the consistency of the two notions of differential
cochains. Lemma~\ref{lemma:local-form-of-Diffops} already ensures that
the local expressions of cochains have the form as in Definition
\ref{definition:differential-maps}. In the considered case this is
sufficient for the global situation.

\begin{lemma}
    \label{lemma:form-of-cochain-on-V-G}
    Let $\phi\in \Diffop^L(V,\Diffop^l(V\times G))$ with $L\in
    \mathbb{N}_0^k$, $k,l\in \mathbb{N}_0$. Then there exist unique
    $\phi^{\alpha_1\dots \alpha_k}\in \Diffop^l(V\times G)$ such that
    \begin{equation}
        \label{eq:form-cochain-V-G}
        \phi (a_1,\dots,a_k) =
        \sum_{|\alpha_1|\le l_1,\dots, 
          |\alpha_k|\le l_k}
        \left(\frac{\partial^{|\alpha_1|}a_1} {\partial x^{\alpha_1}}
            \dots 
            \frac{\partial^{|\alpha_k|}a_k}{\partial
              x^{\alpha_k}}\right) \cdot \phi^{\alpha_1
          \dots \alpha_k}.
    \end{equation} 
    for all $a_1,\dots, a_k\in C^\infty(V)$.
    % If $\phi$ is $G$-invariant the same is true for the unique
    % $\phi^{\alpha_1\dots \alpha_k}$.
\end{lemma}

\begin{proof}
    Using the adapted charts \eqref{eq:adapted-chart-V-G} for the
    local expressions in \eqref{eq:local-form-of-Diffop} the stated
    uniqueness of the occurring $\phi_{\tilde{U}}^{\alpha_1\dots
      \alpha_k}\in \Diffop(\tilde{U})$ shows that there exist unique
    $\phi^{\alpha_1\dots \alpha_k}\in \Diffop^l(V\times G)$ with
    $\phi^{\alpha_1\dots \alpha_k}|_{\tilde{U}}=
    \phi_{\tilde{U}}^{\alpha_1\dots \alpha_k}$ such that
    \eqref{eq:form-cochain-V-G} holds.
    % The assertion concerning
    % $G$-invariance is a direct consequence of the given uniqueness.
\end{proof}

Besides this more or less obvious consistency of definitions one has
to prove that $\Diffop^{\bullet} (V\times G)$ satisfies all
requirements of the topological bimodule $\mathcal{M}$ in
Theorem~\ref{theorem:isomorphism-Hochschild-bar}.

For this purpose we have to consider an appropriate topology of the
differential operators $\Diffop^l(P)$ on a smooth manifold $P$ of
dimension $d$. The Fr\'{e}chet topology of the differential operators
of order $l\in \mathbb{N}_0$ which is induced from the one of smooth
functions on $P$ by using the local expressions will be sufficient for
this purpose. Consider all compact subsets $K\subseteq P$ such that
$K\subseteq U$ for a chart $(U,z)$ of $P$. Then, every differential
operator $D \in \Diffop^l(P)$ has the local form $D|_U=
\sum_{|\alpha|\le l} D_U^{\alpha} \frac{\partial^{|\alpha|}} {\partial
  z^{\alpha}}$ with coefficient functions $D_U^{\alpha}\in
C^\infty(U)$ and multiindices $\alpha\in \mathbb{N}_0^d$. The common
seminorms $p_{K,r}$ with $r\in \mathbb{N}_0$ inducing the Fr\'{e}chet
topology of $C^\infty(P)$, confer \eqref{eq:seminorms-functions}, are
then used to define corresponding seminorms $p^{(l)}_{U,K,r}$ on the
spaces $\Diffop(P)^l$ via
\begin{equation}
    \label{eq:seminorms-diffop}
    p^{(l)}_{U,K,r}(D)= \max_{|\alpha|\le l}(p_{K,r}(D_U^\alpha))=    
    \max_{
      \begin{subarray}{c}
          p\in K\\
          |\beta|\le r\\
          |\alpha| \le l
      \end{subarray}
    } \left| \left(\frac{\partial^{|\beta|}}{\partial z^\beta}
            D_U^\alpha \right) (p)
    \right| \quad \textrm{for all } D\in 
    \Diffop^l(P). 
\end{equation}
It is easy to see that the maps $p^{(l)}_{U,K,r}$ really are
seminorms. Thus, they induce a locally convex topology on
$\mathcal{D}^l$ with the following well-known properties.

\begin{lemma}
    \label{lemma:Diffop-topology}
    Let $P$ be a manifold and for each $l\in \mathbb{N}_0$ let
    $\Diffop^l(P)$ be equipped
    with the topology induced by the seminorms of
    \eqref{eq:seminorms-diffop}. Then the following statements hold:
    \begin{enumerate}
    \item The locally convex spaces $\Diffop^l(P)$ are Fr\'{e}chet
        spaces, this means they are Hausdorff, metrizable and
        complete.
    \item The embedding map $\Diffop(P)^l \hookrightarrow
        \Diffop(P)^k$ for all $l\le k$ is a continuous embedding with
        closed image. In particular, the topology of $\Diffop(P)^l$ is
        the one induced from $\Diffop(P)^k$.
    \item The composition of differential operators is a
        continuous map $\Diffop(P)^l\times
        \Diffop(P)^{l'}\longrightarrow \Diffop(P)^{l+l'}$ for all
        $l,l'\in \mathbb{N}_0$.
    \item With respect to the natural Fr\'{e}chet topology of
        $C^\infty(P)$ the spaces $\Diffop(P)^l$ are topological
        $(C^\infty(P),C^\infty(P))$-bimodules with the left module
        structure $a\cdot D=\mathsf{L}_a\circ D$ and the right module
        structure $D\cdot a= D\circ \mathsf{L}_a$ for $D\in
        \Diffop(P)$ and $a\in C^\infty(P)$.
    \item If $P$ is a principal fibre bundle the $G$-invariant
        differential operators are closed subspaces
        $\Diffop^l(P)^G\subseteq \Diffop^l(P)$ for all $l\in
        \mathbb{N}_0$.
    \end{enumerate}
\end{lemma}

\begin{proof}
    It is obvious that the locally convex spaces $\Diffop^l(P)$ are
    topological vector spaces. For the Hausdorff property one only has
    to show that for all $0\neq D\in \Diffop^l(P)$ there exist
    seminorms with $p_{U,K,r}^{(l)}(D)\neq 0$. This is clear with
    $K=\{p\}$ for appropriate $p\in U\subseteq P$ and $r=0$.

    By definition, a manifold is in particular a topological space
    which is second countable, this means the topology has a countable
    basis. Thus one can choose a countable covering $\{U_n\}_{n\in
      \mathbb{N}}$ of $P$ with chart domains and, in addition, a
    countable covering of each $U_n$ with compact subsets $K\subseteq
    U_n\subseteq P$, confer \cite[App.~A.1]{waldmann:2007a}. Since
    the locally convex topology is already induced by the
    corresponding countable set of seminorms $\{p_n\}_{n\in
      \mathbb{N}}$ the topology is metrizable, for instance with the
    metric $d(D,D')=\max_n \frac{1}{2} \frac{p_n(D-D')}{1+p_n(D-D')}$
    for all $D,D'\in \Diffop^l(P)$. The countable set of seminorms
    inducing the topology further implies that the topology is first
    countable, this means that for all $D\in \Diffop^l(P)$ there
    exists a countable set of neighbourhoods such that any
    neighbourhood contains one of those. For example, the open balls
    of rational radius with respect to the seminorms $p_n$ build such
    a set. Thus, $\Diffop^l(P)$ is complete if and only if any Cauchy
    sequence converges. Due to \eqref{eq:seminorms-diffop} a Cauchy
    sequence $\{D_n\}$ of operators induces Cauchy sequences
    $\{D_{U,n}^\alpha\}$ of coefficient functions with respect to the
    Fr\'{e}chet topology of $C^\infty(U)$ for the relevant chart
    $(U,z)$. Then, these sequences converge to smooth functions
    $D_U^\alpha\in C^\infty(U)$. Considering different charts $(U,z)$
    and $(U',z')$ with $U\cap U'\neq \emptyset$ the local coefficient
    functions are related by $D_{U,n}^\alpha|_{U\cap U'}=
    \sum_{|\beta|\le l} D_{U',n}^\beta|_{U\cap U'} \cdot
    f^\alpha_\beta$ where $f^{\alpha}_{\beta}\in C^\infty(U\cap U')$
    contains the Jacobi matrices of the coordinate change. Since the
    functions are a topological algebra this behaviour is the same for
    the limits $D_U^\alpha$ and $D_{U'}^\alpha$. Thus, they determine
    a well-defined differential operator $D\in \Diffop^l(P)$ given by
    the local expressions. By construction, the sequence $\{D_n\}$
    converges to $D$. Thus, $\Diffop^l(P)$ is a Fr\'{e}chet space.

    The next three assertions are easy consequences of
    \eqref{eq:seminorms-diffop}. If $D\in \Diffop(P)^l$ and $l\le k$
    one has $p^{(k)}_{U,K,r}(D)= p^{(l)}_{U,K,r}(D)$. This and the
    first point imply the second assertion. The third and fourth part
    are a simple consequence of the Leibniz rule. With the obvious
    notation the third part follows with $p^{(l+l')}_{U,K,r}(D\circ
    D')\le p^{(l)}_{U,K,r}(D) p^{(l')}_{U,K,r+l}(D')$ and for the
    fourth part one has $p^{(l)}_{U,K,r}(f \cdot D)\le p_{K,r}(f)
    p^{(l)}_{U,K,r}(D)$ and $p^{(l)}_{U,K,r}(D\cdot f)=
    p^{(l)}_{U,K,r}(D) p_{K,r+l}(f)$.

    For the last statement let $D_n\in \Diffop^l(P)^G$ be a sequence
    of differential operators converging to $D\in \Diffop^l(P)$ with
    respect to the considered topology. The fact that for all compact
    subsets $K\subseteq U$ of a chart domain $U\subseteq P$ one has
    $p_{K,r}(D_n(f))\le p_{U,K,r}^{(l)}(D_n) p_{K,r+l}(f)$, implies
    that for all $f\in C^\infty(P)$ the sequence $D_n(f)$ converges to
    $D(f)$ with respect to the usual Fr\'{e}chet topology of
    $C^\infty(P)$. Since we only work with smooth group actions the
    right action $\rr_g: P\longrightarrow P$ is a diffeomorphism for
    all $g\in G$. Thus, $\rr_g^*: C^\infty(P) \longrightarrow
    C^\infty(P)$ is a continuous map and the equation
    $\rr_g^*(D_n(f))= D_n(\rr_g^* f)$ implies $\rr_g^*(D(f))=
    D(\rr_g^* f)$ for all $f\in C^\infty(P)$. This shows $g\acts D=D$
    which proves the assertion.
\end{proof}

The next lemma makes use of the canonical global coordinates
$x^i:\mathbb{R}^n \supseteq V \longrightarrow \mathbb{R}$,
$i=1,\dots,n$, which extend to maps $x^i\in C^\infty(V\times G)$.

\begin{lemma}
    \label{lemma:form-Diffop-V-G}
    Every $D\in \Diffop^l(V\times G)$ with $l\in \mathbb{N}_0$ is of
    the form
    \begin{equation}
        \label{eq:form-Diffop-V-G}
        D = \sum_{r=0}^l \frac{1}{r!} D^{i_1 \dots i_r}
        \frac{\partial^r}{\partial x^{i_1}\dots \partial x^{i_r}} 
    \end{equation}
    with uniquely defined vertical operators $D^{i_1\dots i_r}\in
    \Diffopverg{l-r}(V\times G)$ which are symmetric in the indices
    $i_1,\dots, i_r\in \{1,\dots, n\}$. 
    % Clearly, for $G$-invariant
    % $D$ the same is true for the $D^{i_1\dots i_r}$.
\end{lemma}

\begin{proof}
    With the adapted local charts \eqref{eq:adapted-chart-V-G} the
    considered $D \in \Diffop^l(V\times G)$ has the local expression 
    \begin{eqnarray}
        \label{eq:local-Diffop-V-G}
        D|_{\tilde{U}}
        &=& \sum_{s=0}^l\sum_{r=0}^s \frac{1} {r! (s-r)!}
        D_{\hat{U}}^{i_1\dots i_r j_1 \dots j_{s-r}}
        \frac{\partial^r} {\partial x^{i_1} \dots \partial x^{i_r}} 
        \frac{\partial^{s-r}} {\partial y^{j_1} \dots \partial
          y^{j_{s-r}}} \nonumber \\
        &=& \sum_{r=0}^l\sum_{s=r}^l \frac{1} {r! (s-r)!}
        D_{\hat{U}}^{i_1\dots i_r j_1 \dots j_{s-r}}
        \frac{\partial^{s-r}} {\partial y^{j_1} \dots \partial
          y^{j_{s-r}}} 
        \frac{\partial^r} {\partial x^{i_1} \dots \partial x^{i_r}},
    \end{eqnarray}
    with unique functions $D_{\hat{U}}^{i_1\dots i_r j_1 \dots
      j_{s-r}}\in C^\infty(V\times \hat{U})$ which are symmetric in
    each group of indices $i_1,\dots,i_r\in \{1,\dots, n\}$ and
    $j_1,\dots,j_{s-r}\in \{1,\dots, k\}$.  The behaviour under a
    change of the charts $(\hat{U},y)$ of $G$ shows that the locally
    given vertical differential operators
    \begin{equation}
        \label{eq:local-vertical-diffop}
        D_{\hat{U}}^{i_1\dots i_r}:= \sum_{s=r}^l \frac{1} {(s-r)!}
        D_{\hat{U}}^{i_1\dots i_r j_1 \dots j_{s-r}}
        \frac{\partial^{s-r}} {\partial y^{j_1} \dots \partial
          y^{j_{s-r}}}\in \Diffopverg{l-r}(V\times \hat{U})
    \end{equation}
    uniquely define global ones $D^{i_1\dots i_r}\in
    \Diffopverg{l-r}(V\times G)$ with $D^{i_1\dots i_r}|_{V\times
      \hat{U}}= D^{i_1\dots i_r}_{\hat{U}}$, thus yielding
    \eqref{eq:form-Diffop-V-G}.
\end{proof}

The two lemmas above now imply the applicability of
Theorem~\ref{theorem:isomorphism-Hochschild-bar}.

\begin{lemma}
    \label{lemma:conditions-satisfied}
    The filtered space $\Diffop^\bullet (V \times G)$ satisfies the
    conditions of
    Theorem~\ref{theorem:isomorphism-Hochschild-bar}. Hence one has
    \begin{equation}
        \label{eq:isomorphism-Hochschild-bar-surj-sub}
        \mathrm{HC}^\bullet_{\mathrm{diff}} (V,
        \Diffop(V\times G)) 
        \cong 
        \Union_{l=0}^\infty \Hom_{\mathcal{A}^e}^{\mathrm{diff}}
        (X_\bullet, \Diffop^l(V\times G))
    \end{equation}
    and
    \begin{equation}
        \label{eq:equivalence-Hochschild-Koszul-cohomology}
        \mathrm{HH}^\bullet_{\mathrm{diff}}
        (V,\Diffop(V\times G)) 
        \cong  
        \mathrm{H} \left(
            \Union_{l=0}^\infty \Hom_{\mathcal{A}^e}
            (K_\bullet,\Diffop^l(V\times G))
        \right).
    \end{equation}
    % If $G$ is a Lie group $\Diffop^\bullet(V\times G)$ can be
    % replaced
    % by $\Diffop^\bullet(V\times G)^G$ to get analogous assertions
    % with $G$-invariant operators.
\end{lemma}

\begin{proof}
    Since $\pr_1^*: C^\infty(V)\longrightarrow C^\infty(V\times G)$ is
    a continuous map with respect to the corresponding Fr\'{e}chet
    topologies, Lemma~\ref{lemma:Diffop-topology} immediately implies
    for $P=V\times G$ that the spaces $\Diffop^l (V \times G)$ are
    topological $(C^\infty(V),C^\infty(V))$-bimodules with the
    structures corresponding to \eqref{eq:bimodule-diffop}.
    Lemma~\ref{lemma:form-Diffop-V-G} and the Leibniz rule show that
    the right $C^\infty(V)$-module structures satisfy condition
    \eqref{eq:right-module-structure-differential}, so all
    $\Diffop^l(V\times G)$ are $(C^\infty(V),C^\infty(V))$-bimodules
    of order $l\in \mathbb{N}_0$. Then, the stated isomorphisms
    follow.
\end{proof}

So far, the problem of computing the Hochschild cohomology has only
been paraphrased and reformulated for other complexes. But now, it is
possible to compute the cohomology of $\bigcup_{l=0}^\infty
\Hom_{\mathcal{A}^e} (K_\bullet, \Diffop^l (V \times G))$
explicitly.

It is a crucial point to recognize that the differential operators on
$C^\infty(V\times G)$ and thus the relevant complex
$\Union_{l=0}^\infty\Hom_{\mathcal{A}^e}(K_\bullet,\Diffop^l(V\times
G))$ have a particular grading. With Lemma~\ref{lemma:form-Diffop-V-G}
it follows immediately that any differential operator $D\in
\Diffop^l(V\times G)$ can be written as a direct sum $D= \sum_{r=0}^l
D_r$ where the $D_r$ are given by
\begin{equation}
    \label{eq:Diffop-r}
    D_r= \frac{1}{r!} D^{i_1\dots i_r} \frac{\partial^r}{\partial
      x^{i_1}\dots \partial x^{i_r}}= 
    \sum_{
      \begin{subarray}{c}
          \alpha\in \mathbb{N}_0^n\\
          |\alpha|=r
      \end{subarray}
      }
      D_r^\alpha \frac{\partial^r}{\partial x^\alpha}
      = 
      \sum_{
      \begin{subarray}{c}
          \alpha\in \mathbb{N}_0^n\\
          |\alpha|=r
      \end{subarray}
    }
    D_r^\alpha \frac{\partial^r}{\partial (x^1)^{\alpha_1} \dots
      (x^n)^{\alpha_n}} 
\end{equation}
with vertical operators $D^{i_1\dots i_r}$ and $D_r^\alpha\in
\Diffopverg{l-r}(V\times G)$ satisfying $D_r^\alpha(\pr_1^* a \cdot
f)= \pr_1^*a\cdot D_r^\alpha(f)$ for all $a\in C^\infty(V)$ and $f\in
C^\infty(V\times G)$. The $D_r\in \Diffop^l(V\times G)$ are obviously
specified by the property that the order of differentiation with
respect to the \emph{horizontal} direction $V$ is exactly given by
$r$. The subset of all such differential operators in
$\Diffop^l(V\times G)$ of the form \eqref{eq:Diffop-r} shall be
denoted by $\Diffop^l_r(V\times G)$. Then for all $l\in \mathbb{N}_0$
one has the decomposition
\begin{equation}
    \label{eq:decomposition-diffops-on-product}
    \Diffop^l(V\times G) = \bigoplus_{r=0}^l
    \Diffop^l_r(V\times G),
\end{equation}
showing that $\Diffop^l(V\times G)$ is a graded space. The
decomposition \eqref{eq:decomposition-diffops-on-product} naturally
transfers to
\begin{equation}
    \label{eq:decomposition-Hom-Koszul-complex}
    \Hom_{\mathcal{A}^e} (K_\bullet,\Diffop^l(V\times G))=
    \bigoplus_{r=0}^l \Hom_{\mathcal{A}^e}
    (K_\bullet,\Diffop^l_r(V\times G)) \quad \textrm{for all
    } l \in \mathbb{N}.
\end{equation}

\begin{definition}[The degree of horizontal differentiation]
    By linear extension of
    \begin{equation}
        \label{eq:degree-Hom}
        \deg \psi =r \psi \quad \textrm{for all } \psi \in
        \Hom_{\mathcal{A}^e}(K_\bullet,\Diffop^l_r(V\times G))
    \end{equation}
    one obtains a map 
    \begin{equation}
        \label{eq:degree-map-complex}
        \deg: \bigcup_{l=0}^\infty
        \Hom_{\mathcal{A}^e} (K_\bullet, \Diffop^l (V \times G))
        \longrightarrow \bigcup_{l=0}^\infty \Hom_{\mathcal{A}^e}
        (K_\bullet, \Diffop^l (V\times G)),
    \end{equation}
    which, by the obvious reasons, is called the \emph{degree of
      horizontal differentiation}.
\end{definition}

With $\Hom_{\mathcal{A}^e}(K_0,\Diffop^l_r(V\times G))\cong
\Diffop^l_r(V\times G)$ one obviously has a corresponding map $\deg:
\Diffop^\bullet (V\times G)\longrightarrow \Diffop^\bullet(V\times G)$
defined by
\begin{equation}
    \label{eq:degree}
    \deg D_r = r D_r \quad \textrm{for all }
    D_r\in \Diffop^l_r(V\times G),
\end{equation}
such that $(\deg \psi)(\omega)= \deg (\psi(\omega))$ for all $\psi \in
\Hom_{\mathcal{A}^e} (K_k, \Diffop^l (V\times G))$ and $\omega \in
K_k$ with $l, k \in \mathbb{N}_0$.

The $(C^\infty(V), C^\infty(V))$-bimodule structure of
$\Diffop^\bullet(V\times G)$ as in \eqref{eq:bimodule-diffop} induces
a corresponding $C^\infty(V\times V)$-module structure according to
Proposition \ref{proposition:Continuous-Hochschild}. For the latter
structure we can formulate the following simple but useful lemma.

\begin{lemma}
    \label{lemma:coordinate-functions}
    Let $x^i: V \longrightarrow \mathbb{R}$ be the coordinate
    functions with respect to the canonical basis $\{e_i\}_{i = 1,
      \ldots, n}$ of $\mathbb{R}^n$ and let $\xi^i = x^i \otimes 1 - 1
    \otimes x^i \in C^\infty(V \times V)$ as in
    \eqref{eq:functions-xi}. Then, for all $D \in \Diffop^\bullet
    (V\times G)$ one has
    \begin{equation}
        \label{eq:commutator-Ae-differentiation}
        (\xi^i \cdot D) \circ \frac{\partial} {\partial x^j}
        -\xi^i \cdot \left(D \circ \frac{\partial} {\partial
              x^j}\right) 
        = \delta^i_jD,
        \quad \textrm{for all } i,j=1,\dots, n
    \end{equation}
    and
    \begin{equation}
        \label{eq:degD-short}
        \sum_{j=1}^n 
        \left( 
            \left(
                (-\xi^j) \cdot D
            \right)
            \circ \frac{\partial}{\partial x^j}
        \right)
        =
        \deg D.
    \end{equation}
    Moreover, one finds that for all $r\le l\in \mathbb{N}_0$ and
    $i=1,\dots,n$ the left multiplication $l_{\xi^i}$ with $\xi^i$ is
    a map
    \begin{equation}
        \label{eq:Ae-multiplication-degrees}
        l_{\xi^i}: \Diffop^l_r(V\times G) \longrightarrow
        \Diffop^{l-1}_{r-1}(V\times G). 
    \end{equation}
\end{lemma}

\begin{proof}
    With the bimodule structure \eqref{eq:bimodule-diffop} it is
    obvious that $x^i\cdot (D\circ \frac{\partial} {\partial x^j})
    =(x^i\cdot D)\circ \frac{\partial} {\partial x^j} $ and a simple
    computation shows $(D\circ \frac{\partial}{\partial x^j}) \cdot
    x^i=(D\cdot x^i) \circ \frac{\partial}{\partial x^j}
    +\delta^i_jD.$ Then, \eqref{eq:commutator-Ae-differentiation}
    follows by \eqref{eq:regularity-condition}. Due to the linearity
    of both sides it is sufficient to show \eqref{eq:degD-short} for
    $D\in \Diffop^l_r(V\times G)$ as in \eqref{eq:Diffop-r}. Then, one
    first computes for all $f\in C^\infty(V\times G)$
    \begin{eqnarray*}
        \lefteqn{ \left(\sum_{j=1}^n (D\cdot x^j)\circ
              \frac{\partial}{\partial x^j}\right)(f)}\\  
        &=&
        \sum_{|\alpha|=r} D_r^\alpha \frac{\partial^r}{\partial
          (x^1)^{\alpha_1} \dots \partial (x^n)^{\alpha_n}}
        \left(\sum_{j=1}^n \pr_1^* x^j \cdot \frac{\partial f}{\partial
              x^j}\right)\\  
        &=&
        \sum_{|\alpha|=r} D_r^\alpha
        \left(\sum_{j=1}^n \pr_1^* x^j \cdot\frac{\partial^r}{\partial 
              (x^1)^{\alpha_1} \dots \partial (x^n)^{\alpha_n}} 
            \frac{\partial f}{\partial x^j}  
            + \sum_{j=1}^n \alpha_j \frac{\partial^r f} {\partial 
              (x^1)^{\alpha_1} \dots \partial (x^n)^{\alpha_n}}
        \right)\\
        &=& 
        \sum_{j=1}^n \pr_1^* x^j \cdot D \left( \frac{\partial f}
            {\partial x^j} \right)
        + rD(f) 
        = 
        \left( \sum_{j=1}^n (x^j\cdot D)\circ \frac{\partial}{\partial
              x^j} + \deg D \right) (f),
    \end{eqnarray*}
    where one has used the Leibniz rule  
    \begin{equation}
        \label{eq:Leibniz-particular}
        \frac{\partial^{\alpha_j}} {\partial (x^j)^{\alpha_j}} 
        \left(\pr_1^* x^j \cdot \frac{\partial f}{\partial x^j} \right)
        = \sum_{t=0}^{\alpha_j}\binom{\alpha_j}{t} \frac{\partial^t
          \pr_1^*x^j}{\partial (x^j)^t} \cdot
        \frac{\partial^{\alpha_j-t+1} f} 
        {\partial (x^j)^{\alpha_j-t+1}}= \pr_1^* x^j \cdot
        \frac{\partial^{\alpha_j}} {\partial (x^j)^{\alpha_j}} 
        \frac{\partial f}{\partial x^j} + \alpha_j
        \frac{\partial^{\alpha_j} f} {\partial (x^j)^{\alpha_j}}.
    \end{equation}
    This yields
    \begin{equation}
        \label{eq:degD}
        \sum_{j=1}^n 
        \left(
            (D \cdot x^j) \circ \frac{\partial}{\partial x^j} 
            - (x^j \cdot D) \circ \frac{\partial}{\partial x^j}
        \right) 
        =
        \deg D,
    \end{equation}
    and \eqref{eq:degD-short} follows again with
    \eqref{eq:regularity-condition}. The property
    \eqref{eq:Ae-multiplication-degrees} is a simple consequence of
    the Leibniz rule and \eqref{eq:regularity-condition}.
\end{proof}

After the above preparing considerations we now define the crucial
maps which will lead to a homotopy equation.

\begin{definition}
    \label{definition:homotopy-maps-Koszul}
    We define the map
    \begin{equation}
        \label{eq:homotopy-map}
        \delta_K^{-1}:
        \Union_{l=0}^\infty \Hom_{\mathcal{A}^e} 
        (K_\bullet, \Diffop^l (C^\infty(V \times G)))
        \longrightarrow 
        \Union_{l=0}^\infty \Hom_{\mathcal{A}^e}
        (K_{\bullet-1}, \Diffop^l (C^\infty(V \times G)))  
    \end{equation}
    by the linear extension of the maps
    \begin{equation}
        \label{eq:homotopy-map-explicitly}
        (\delta_K^{-1})^k \psi
        =
        \begin{cases}
            \frac{1} {k+r} (\delta_K^*)^k\psi & \textrm{for} \; k\ge
            1\\
            0 & \textrm{for} \; k=0
        \end{cases}
    \end{equation}    
    for $\psi \in \Hom_{\mathcal{A}^e}
    (K_k,\Diffop^l_r(C^\infty(V\times G)))$ where the maps
    \begin{equation}
        \label{eq:deltastarDef}
        (\delta_K^*)^k: 
        \Union_{l=0}^\infty \Hom_{\mathcal{A}^e} 
        (K_k, \Diffop^l (C^\infty(V \times G))) 
        \longrightarrow 
        \Union_{l=0}^\infty \Hom_{\mathcal{A}^e}
        (K_{k-1}, \Diffop^l (C^\infty(V \times G)))
    \end{equation}
    are defined by %$(\delta_K^*)^0=0$ and
    \begin{equation}
        \label{eq:homotopy-map-Hom-Koszul}
        ((\delta_K^*)^k\psi) 
        (e^{i_1} \wedge \cdots \wedge e^{i_{k-1}})
        = 
        - \sum_{j=1}^n
        \psi (e^j \wedge e^{i_1} \wedge \cdots \wedge e^{i_{k-1}})
        \circ \frac{\partial} {\partial x^j}
    \end{equation}
    for all $i_s \in \{1,\ldots, n\}$ with the dual basis
    $\{e^i\}_{i=1\dots, n}$ of the canonical basis $\{e_i\}_{i=1\dots,
      n}$ of $\mathbb{R}^n$.
\end{definition}

With \eqref{eq:boundary-operator-Koszul-complex-short} for the case of
the canonical basis of $\mathbb{R}^n$ the coboundary operators
\begin{equation}
    \label{eq:coboundary-Koszul}
    \delta_K^k = (d^K_{k+1})^*
\end{equation}
of the considered complex are given by
\begin{equation}
    \label{eq:differential-Hom-Koszul}
    (\delta_K^k\psi)
    (e^{i_1} \wedge \cdots \wedge e^{i_{k+1}})
    =
    \sum_{j=1}^n \xi^j \cdot\psi \left(
        \ins_{\mathrm{a}}(e_j) 
        (e^{i_1}\wedge \cdots \wedge e^{i_{k+1}}) 
    \right) 
\end{equation}
for $\psi\in \Union_{l=0}^\infty \Hom_{\mathcal{A}^e}(K_k, \Diffop^l
(C^\infty(V \times G)))$ and $i_s \in \{1,\ldots, n\}$ with $\xi^j$ as
in \eqref{eq:functions-xi}.

\begin{remark}
    \label{remark:Koszul-homotopy-diff}
    Note that by definition for all $r\le l\in \mathbb{N}_0$ and $k\ge
    1$ one has
    \begin{equation}
        \label{eq:homotopy-degrees}
        \delta_K^{-1}:
        \Hom_{\mathcal{A}^e}(K_k ,\Diffop^l_r(V\times G)) 
        \longrightarrow
        \Hom_{\mathcal{A}^e}(K_{k-1},\Diffop^{l+1}_{r+1}(V\times
        G)).
    \end{equation}
    and, even for $k\ge 0$, 
    \begin{equation}
        \label{eq:differential-degrees}
        \delta_K:
        \Hom_{\mathcal{A}^e}(K_k,\Diffop^l_r(V\times G)) 
        \longrightarrow
        \Hom_{\mathcal{A}^e}(K_{k+1},\Diffop^{l-1}_{r-1} (V\times
        G))
    \end{equation}
    due to \eqref{eq:Ae-multiplication-degrees}. It is further
    remarkable that in contrast to $\delta_K$ the maps $\delta_K^*$
    and $\delta_K^{-1}$ are \emph{not} $\mathcal{A}^e$-linear.
\end{remark}

\begin{proposition}[The crucial homotopy]
    \label{proposition:computation-of-cohomology}
    The $\mathbb{K}$-linear map $\delta_K^{-1}$ yields an explicit
    homotopy for the identity map of the considered complex $\left(
        \bigcup_{l=0}^\infty \Hom_{\mathcal{A}^e}
        (K_\bullet,\Diffop^l(C^\infty(V\times G))), \delta_K
    \right)$. This means that for $k\ge 1$
    \begin{equation}
        \label{eq:homotopy}
        \delta_K^{k-1}\circ (\delta_K^{-1})^k +
        (\delta_K^{-1})^{k+1}\circ \delta_K^k  =
        \id.
    \end{equation}
    Consequently, the corresponding cohomologies are trivial,
    this means 
    \begin{equation}
        \label{eq:Koszul-cohomology-vanishes}
        \mathrm{H}^k \left(
            \Union_{l=0}^\infty
            \Hom_{\mathcal{A}^e}(K, \Diffop^l (V \times
            G)) 
        \right)
        = \{0\} \quad \textrm{for } k \ge 1. 
    \end{equation}
\end{proposition}

\begin{proof}
    With the formulas of Lemma~\ref{lemma:coordinate-functions} one
    computes for $k\ge 1$
    \begin{eqnarray*}
        \lefteqn{
          \left(
              (\delta_K^{k-1}\circ (\delta_K^*)^k +
              (\delta_K^*)^{k+1}\circ \delta_K^k)\psi
          \right)
          (e^{i_1}\wedge \dots \wedge e^{i_k})} \\ 
        &=& 
        - \xi^i\cdot 
        \left(
            \psi (e^j \wedge \insa(e_i) (e^{i_1}\wedge \dots \wedge 
            e^{i_k})) \circ \frac{\partial}{\partial x^j}
        \right) 
        - \left( 
            \xi^i \cdot \psi (\insa(e_i) (e^j\wedge e^{i_1}\wedge \dots 
            \wedge e^{i_k}))
        \right) 
        \circ \frac{\partial} {\partial x^j}\\
        &=& 
        - \xi^i\cdot 
        \left(
            \psi (e^j \wedge \insa(e_i) (e^{i_1}\wedge \dots \wedge
            e^{i_k})) \circ \frac{\partial}{\partial x^j}
        \right) 
        + \left( 
            (-\xi^j) \cdot \psi (e^{i_1}\wedge
            \dots \wedge e^{i_k})
        \right) 
        \circ \frac{\partial} {\partial x^j}\\
        && + \left( 
            \xi^i \cdot \psi (e^j \wedge \insa(e_i) (e^{i_1}\wedge
            \dots \wedge e^{i_k}))
        \right) 
        \circ \frac{\partial} {\partial x^j}\\
        &=& 
        \deg \psi(e^{i_1}\wedge \dots \wedge e^{i_k}) + \psi(e^j\wedge
        \insa(e_j)(e^{i_1}\wedge \dots \wedge e^{i_k}))\\
        &=& 
        ((\deg + k\id )\psi)(e^{i_1}\wedge \dots \wedge e^{i_k}).
    \end{eqnarray*}
    %For $k=0$ a similar computation can be made using that by
    %definition the coboundary operators vanish, $\delta_K^{-s}=0$ for
    %$s\in \mathbb{N}$.
    This shows that
    \begin{equation}
        \label{eq:homotopy-deg}
        \delta_K^{k-1}\circ (\delta_K^*)^k + (\delta_K^*)^{k+1}\circ
        \delta_K^k  = \deg + k \id
    \end{equation}
    for all $k\in \mathbb{N}$. Due to
    Remark~\ref{remark:Koszul-homotopy-diff} and the definition of
    $\delta^{-1}$ this yields \eqref{eq:homotopy}. The last statement
    \eqref{eq:Koszul-cohomology-vanishes} is an immediate consequence.
\end{proof}

Proposition~\ref{proposition:computation-of-cohomology} now already
states that the local Hochschild cohomologies are trivial. According
to the explanations of
Proposition~\ref{proposition:sheaves-and-cohomology} this result
extends to the global situation if all cocycles in a determined
Hochschild module are images of cochains with globally bounded orders
of differentiation. In the present case this is true since there
exists an explicit homotopy map $\delta^{-1}$ yielding the desired
property.

\begin{theorem}[The local Hochschild cohomology $\HHdiff^\bullet
    (V,\Diffop(V\times G))$]
    \label{theorem:local-Hochschild-cohomologies-surjective-submersions}
    Let $V \subset \mathbb{R}^n$ be an open and convex subset and let
    $G$ be a manifold. Then one has for all $k \in \mathbb{N}$:
    \begin{enumerate}
    \item The $k$-th differential Hochschild cohomology with values in
        $\Diffop(V \times G)$ is trivial
        \begin{equation}
            \label{eq:cohomology-groups-local}
            \HHdiff^k (V,\Diffop(V\times G))= \{0\}. 
        \end{equation}
    \item There exists an explicit homotopy
        \begin{equation}
            \label{eq:explicit-homotopy}
            \delta^{k-1} \circ (\delta^{-1})^k+ (\delta^{-1})^{k+1}
            \circ \delta^k=\id,
        \end{equation}
        where the maps $(\delta^{-1})^k: \HCdiff^k (V,
        \Diffop(V \times G)) \longrightarrow \HCdiff^{k-1}
        (V, \Diffop (V \times G))$ are given by
        \begin{equation}
            \label{eq:explicit-homotopy-map}
            (\delta^{-1})^k
            = \Xi^{k-1} \circ 
            \left(
                (G_{k-1})^* \circ (\delta_K^{-1})^{k} \circ (F_k)^* 
                + (s_{k-1})^*
            \right) \circ (\Xi^k)^{-1}.
        \end{equation}
    \item In particular, for any multiindex $L = (l_1, \ldots,
        l_k) \in \mathbb{N}_0^k$ one has
        \begin{equation}
            \label{eq:explicit-homotopy-map-diff}
            (\delta^{-1})^k:
            \Diffop^L (V, \Diffop^l (V \times G))
            \longrightarrow \Diffop^{\tilde{L}} (V,
            \Diffop^{l+1} (V\times G)), 
        \end{equation} 
        where the new multiindex $\tilde{L} = (\tilde{l}_1, \ldots,
        \tilde{l}_{k-1}) \in \mathbb{N}_0^{k-1}$ is given by
        \begin{equation}
            \label{eq:new-multi-index}
            \tilde{l}_i = 
            \max\{(k-2)! + |L|, l+2\}\le (k-2)! + |L|+ l+2
            \quad \textrm{for all } i = 1, \ldots, k-1.
        \end{equation}
    \end{enumerate}
\end{theorem}

\begin{proof}
    The proof of Equation~\eqref{eq:explicit-homotopy} is a simple
    computation which makes use of \eqref{eq:homotopy},
    \eqref{eq:id-Theta-homotopic} and the properties of the involved
    functions. Then, Equation~\eqref{eq:cohomology-groups-local} is
    trivial. The third assertion \eqref{eq:explicit-homotopy-map-diff}
    follows with Proposition~\ref{proposition:pullbacks-diff} and
    Remark~\ref{remark:Koszul-homotopy-diff} by counting the orders of
    differentiation.
\end{proof}

As stated in the beginning of the section all previous considerations
can be reformulated in the $G$-invariant setting. 

\begin{remark}[The $G$-invariant local cohomology $\HHdiff^\bullet(V,
    \Diffop(V\times G)^G)$]
    \label{remark:local-situation-G-inv}
    If $G$ is a Lie group all statements of this section are still true
    when replacing $\Diffop(V\times G)$ by $\Diffop(V\times G)^G$. In
    detail, this is the case because of the following reasons.
    \begin{enumerate}
    \item The assertions of the
        Lemmas~\ref{lemma:form-of-cochain-on-V-G} and
        \ref{lemma:form-Diffop-V-G} still hold for the $G$-invariant
        operators. This is a direct consequence of the
        Equations~\eqref{eq:form-cochain-V-G} and
        \eqref{eq:form-Diffop-V-G} due to the stated uniqueness.
    \item As seen in the last part of
        Lemma~\ref{lemma:Diffop-topology}, $\Diffop^l(V\times G)^G$ is
        a closed subspace of $\Diffop^l(V\times G)$ for all $l\in
        \mathbb{N}_0$. With $\pr_1^*: C^\infty(V)\longrightarrow
        C^\infty(V\times G)^G$ this implies that
        $\Diffop^\bullet(V\times G)^G$ satisfies the conditions of
        Theorem~\ref{theorem:isomorphism-Hochschild-bar}. Thus,
        Lemma~\ref{lemma:conditions-satisfied} can be reformulated.
    \item The differential operators $\frac{\partial}{\partial x^j}\in
        \Diffop^1(V\times G)^G$ for $j=1,\dots,n$ which occur in the
        defining Equation~\eqref{eq:homotopy-map-Hom-Koszul} for the
        map $\delta_K^*$ are $G$-invariant. Further, the decomposition
        \eqref{eq:decomposition-diffops-on-product} and the map $\deg$
        are $G$-invariant. Thus, the derived map $\delta_K^{-1}$
        satisfies
        \begin{equation}
            \label{eq:homotopy-map-G-inv}
            \delta_K^{-1}: 
            \Union_{l=0}^\infty \Hom_{\mathcal{A}^e}
            (K_\bullet, \Diffop^l(V\times G)^G) \longrightarrow 
            \Union_{l=0}^\infty \Hom_{\mathcal{A}^e}
            (K_{\bullet-1}, \Diffop^l(V\times G)^G).
        \end{equation}
        Moreover, this implies that the maps $(\delta^{-1})^k$ in
        Theorem~\ref{theorem:local-Hochschild-cohomologies-surjective-submersions}
        have restrictions
          \begin{equation}
            \label{eq:explicit-homotopy-map-diff-G-inv}
            (\delta^{-1})^k:
            \Diffop^L (V, \Diffop^l (V \times G)^G)
            \longrightarrow \Diffop^{\tilde{L}} (V,
            \Diffop^{l+1} (V\times G)^G).
        \end{equation}
        Thus,
        Theorem~\ref{theorem:local-Hochschild-cohomologies-surjective-submersions}
        is also true in the $G$-invariant case and one has
        \begin{equation}
            \label{eq:cohomology-groups-local-G-inv}
            \mathrm{HH}^k_{\mathrm{diff}}
            (V,\Diffop(V\times G)^G)= \{0\} \quad \textrm{for all }
            k\ge 1. 
        \end{equation}
    \end{enumerate}
\end{remark}

\section{Existence and uniqueness of the relevant deformations}
\label{sec:global-cohomologies-deformations}

In the last step of our computation we use the obtained local results
and apply the Propositions~\ref{proposition:sheaves-and-cohomology}
and \ref{proposition:sheaves-and-cohomology-G-inv}. This is only
possible due to the found explicit homotopy map $\delta^{-1}$ which
guarantees that all coboundaries are images of cochains with
\emph{uniformly bounded degrees of differentiation}. Altogether, we
can reverse the steps of the Lemmas~\ref{lemma:localization} and
\ref{lemma:localization-G-inv} and find the following central
theorems.

\begin{theorem}[Hochschild cohomology for surjective submersions]
    \label{theorem:Hochschild-cohomologies-surjective-submersions}
    Let $\pe: E\longrightarrow P$ be a vector bundle over a surjective
    submersion $\pp: P\longrightarrow M$. Then,
    \begin{equation}
        \label{eq:cohomology-groups-surjective-submersion}
        \mathrm{HH}^k_{\mathrm{diff}}
        (M,\Diffop^\bullet(\Gamma^\infty(P,E))) 
        = \left\{
        \begin{array}{c}
            \Diffopverg{\bullet}(\Gamma^\infty(P,E)) \\ 
            \{0\} 
        \end{array}
        \right.
        \textrm{for}
        \begin{array}{c}
            k=0\\
            k\ge 1.
        \end{array}
    \end{equation}
    More specifically, every $\phi \in \Diffop^L(M,
    \Diffop^l(\Gamma^\infty(P,E)))$ with $L = (l_1, \ldots, l_k) \in
    \mathbb{N}_0^k$, $k \ge 1$, and $\delta \phi = 0$ is of the form
    \begin{equation}
        \label{eq:cocycle-exact-surjective-submersion}
        \phi = \delta \Theta
    \end{equation}
    with $\Theta \in \Diffop^{\tilde{L}} (M,
    \Diffop^{l+1}(\Gamma^\infty(P,E)))$ and $\tilde{L}$ as in
    \eqref{eq:new-multi-index}.
\end{theorem}

\begin{proof}
    The proof is now an easy consequence of
    Proposition~\ref{proposition:sheaves-and-cohomology} and
    Theorem~\ref{theorem:local-Hochschild-cohomologies-surjective-submersions}
    using a partition of unity which is subordinate to the atlas of
    adapted charts $\{\tilde{U}_i\}_{i\in I}$ as used in
    \eqref{eq:diagram-surjective-submersion}.
\end{proof}

Analogously, one gets the following result for principal fibre
bundles.

\begin{theorem}[Hochschild cohomology for principal fibre bundles]
    \label{theorem:Hochschild-cohomologies-principal-fibre-bundle}
    Let $\pp: P\longrightarrow M$ be a principal fibre bundle with
    structure Lie group $G$ and let $\pe: E\longrightarrow P$ be an
    equivariant vector bundle. Then one has
    \begin{equation}
        \label{eq:cohomology-groups-principal-fibre-bundle}
        \mathrm{HH}^k_{\mathrm{diff}}
        (M,\Diffop^\bullet(\Gamma^\infty(P,E))^G) 
        = \left\{
        \begin{array}{c}
            \Diffopverg{\bullet}(\Gamma^\infty(P,E))^G \\ 
            \{0\} 
        \end{array}
        \right.
        \textrm{for}
        \begin{array}{c}
            k=0\\
            k\ge 1.
        \end{array}
    \end{equation}
    More specifically, every $\phi \in \Diffop^L(M,
    \Diffop^l(\Gamma^\infty(P,E))^G)$ with $L = (l_1, \ldots, l_k) \in
    \mathbb{N}_0^k$, $k \ge 1$, and $\delta \phi = 0$ is of the form
    \begin{equation}
        \label{eq:cocycle-exact-principal-bundle}
        \phi = \delta \Theta
    \end{equation}
    with $\Theta \in \Diffop^{\tilde{L}} (M,
    \Diffop^{l+1}(\Gamma^\infty(P,E))^G)$ and $\tilde{L}$ as in
    \eqref{eq:new-multi-index}.
\end{theorem}

Due to these results and those of the
Sections~\ref{sec:deformation-modules} and
\ref{sec:commutant-module-structures} we can now find the existence
and uniqueness up to equivalence of the aspired deformations of right
module structures. Of course we always consider differential and
$G$-invariant differential deformations in the sense of the
Definitions~\ref{definition:deformation-algebra},
\ref{definition:deformation-right-module},
\ref{definition:algebra-type-diff}, \ref{definition:module-type-diff}
and \ref{definition:G-invariant-type}.

\begin{theorem}[Deformations on surjective submersions]
    \label{theorem:deformations-surjective-submersion}
    Let $\pe: E\longrightarrow P$ be a vector bundle over a surjective
    submersion $\pp: P\longrightarrow M$ with a Poisson manifold $M$
    as basis. Further, let $\star$ be a star product on $M$.  Then
    there always exists a differential deformation of the right module
    structure \eqref{eq:right-module-sections} of the sections
    $\Gamma^\infty(P,E)$ which is unique up to equivalence.
\end{theorem}

\begin{theorem}[Deformations on principal fibre bundles]
    \label{theorem:deformations-principal-bundles}
    Let $\pe: E\longrightarrow P$ be an equivariant vector bundle over
    a principal fibre bundle $\pp: P\longrightarrow M$ with a Poisson
    manifold $M$ as basis. Further, let $\star$ be a star product on
    $M$. Then there always exists a $G$-invariant differential
    deformation of the right module structure
    \eqref{eq:right-module-sections} of the sections
    $\Gamma^\infty(P,E)$ which is unique up to equivalence.
\end{theorem}

\begin{remark}[Notation]
    In analogy to $a\star b= \mu(a,b)= \sum_{r=0}^\infty \lambda^r
    \mu_r(a,b)$ we simply write
    \begin{equation}
        \label{eq:notation-right-module-bullet}
        s\bullet a=\rho(s,a)= \sum_{r=0}^\infty \lambda^r \rho_r(s,a)
    \end{equation}
    for the deformed right module structure. All structures of the new
    right module structure are simply denoted by
    \begin{equation}
        \label{eq:notation-right-module-all}
        (\Gamma^\infty(P,E)[[\lambda]], \bullet)_{(C^\infty(M)[[\lambda]],\star)}.
    \end{equation}
\end{remark}

From the considerations of this chapter, in particular
Section~\ref{sec:basic-steps}, it becomes clear that the special case
of functions $C^\infty(P)$ on the total space $P$ is crucial for all
deformations as above. This is the reason why these deformations are
called deformation quantizations of the bundles themselves. Besides
having clarified the occurring notions we have found the proofs for the
Theorems~\ref{theorem:existence-uniqueness-sursub} and
\ref{theorem:existence-uniqueness-principal}.

\begin{corollary}[Deformation quantization of surjective submersions]
    \label{corollary:DQ-sursub}
    Every surjective submersion over a Poisson manifold with a given
    differential star product admits a deformation quantization which
    is unique up to equivalence.
\end{corollary}

\begin{corollary}[Deformation quantization of principal fibre bundles]
    \label{corollary:DQ-principal}
    Every principal fibre bundle over a Poisson manifold with a given
    differential star product admits a deformation quantization which
    is unique up to equivalence.
\end{corollary}

\begin{example}[Tangent and cotangent bundles]
    For a surjective submersion $\pp:P\longrightarrow M$ one can
    consider the tangent and cotangent bundles $TP\longrightarrow P$
    and $T^*P\longrightarrow P$ as well as arbitrary tensor products
    thereof. Then the corresponding tensor fields have a deformed
    right module structure with respect to a star product $\star$ on
    $M$.

    If $P$ is a principal fibre bundle one has a right action $\rR$ of
    the structure group $G$ on the tangent bundle $TP$ given by the
    tangent map
    \begin{equation}
        \label{eq:right-action-tangent-map}
        \rR_g= T\rr_g: T_uP\longrightarrow T_{ug}P
    \end{equation}
    and the induced left action $\lL_g=(T\rr_g)^*$ on the cotangent
    bundle. Then it is obvious that the induced actions
    \eqref{eq:action-on-sections-coming-from-right-action} and
    \eqref{eq:action-on-sections-coming-from-left-action} on the
    sections are nothing but the usual pullbacks of vector fields and
    differential forms. Thus, the actions are algebra automorphisms
    with respect to the tensor product.
\end{example}

With respect to the applications in classical gauge theories we find
the following important theorem where all results are summarized.

\begin{theorem}[Deformation quantization of horizontal forms]
    Let $\pp:P\longrightarrow M$ be a principal fibre bundle with
    structure Lie group $G$ and right action $\rr$. Further, let the
    base be a Poisson manifold $M$ with a corresponding differential
    star product $\star$.

    Then there exists a right $(C^\infty(M)[[\lambda]],\star)$-module
    structure $\bullet$ of the horizontal forms
    \begin{equation}
        \label{eq:horizontal-forms-formal-series}
        \Gamma^\infty_{\hor}(P,\Dach{\bullet} T^*P)[[\lambda]]
    \end{equation}
    with the following properties:
    \begin{enumerate}
    \item The module structure $\bullet$ is a differential deformation
        of the pointwise module, this means
        \begin{equation}
            \label{eq:deformation-forms}
            \alpha \bullet a = \pp^*a \cdot \alpha + \sum_{r=1}^\infty
            \lambda^r \rho_r(\alpha,a) 
        \end{equation}
        for $\alpha\in \Gamma^\infty_{\hor}(P,\Dach{\bullet}
        T^*P)[[\lambda]]$ and $a\in C^\infty(M)[[\lambda]]$ with
        $\rho_r(\alpha,a)=\rho_r(a)(\alpha)$ for differential
        operators $\rho_r\in \Union_{L,l\in \mathbb{N}_0}\Diffop^L
        (M,\Diffop^l(\Gamma^\infty_{\hor}(P,\Dach{\bullet}
        T^*P)))$.
    \item The module structure is $G$-invariant, this means that for
        all $g\in G$
        \begin{equation}
            \label{eq:deformation-forms-G-inv}
            \rr_g^*(\alpha\bullet a)= (\rr_g^* \alpha)\bullet a
        \end{equation}
    \item The deformation $\bullet$ is unique up to equivalence. 
    \end{enumerate}
\end{theorem}

\begin{proof}
    One has to show that the horizontal forms
    $\Gamma^\infty_{\hor}(P,\Dach{t} T^*P)$ for all $t\in
    \mathbb{N}_0$ give rise to a subsheaf of substructures as in
    Remark~\ref{remark:submodules-of-sections}. The defining property
    of such forms $\alpha$ to vanish if one argument is vertical, this
    means $\insa(V)\alpha=0$ if $V\in \Gamma^\infty(VP)$, is not
    affected by the right module structure
    \eqref{eq:right-module-sections}. The sheaf structure is clear and
    the pullback $\rr_g^*$ of forms is in fact an action on the
    horizontal forms since $T\pp\circ T\rr_g=T\pp$ implies that
    $T\rr_g: VP\longrightarrow VP$. Note that in this case for a local
    chart $(U,x)$ of $M$ the tensor products of pullbacks $\pp^*\D
    x^i$ for $i=1,\dots,n=\dim M$ build a
    $C^\infty(\pp^{-1}(U))$-module basis of
    $\Gamma_{\hor}^\infty(\pp^{-1}(U),\Dach{t} T^*\pp^{-1}(U))$ which
    is obviously $G$-invariant. Application of
    Theorem~\ref{theorem:Hochschild-cohomologies-principal-fibre-bundle}
    shows that
    \begin{equation}
        \label{eq:cohomology-horizontal-forms-trivial}
        \HHdiff^k(M,\Diffop(\Gamma_{\hor}^\infty(P,\Dach{t}
        T^*P))^G)=\{0\} 
    \end{equation}
    for all $k\ge 1$ and the assertion is again a consequence of the
    general results of Chapter~\ref{cha:deformation-algebras-modules}.
\end{proof}

\begin{remark}
    With the given definitions it is obvious that in the present
    situations Proposition~\ref{proposition:unit-acts-as-unit} can be
    applied. Thus, for all deformed right module structures with
    respect to a star product $\star$ the unit function
    $1\in C^\infty(M)$ acts by the identity. This means that
    \begin{equation}
        \label{eq:unit-function-acts-as-identity}
        s\bullet 1= s \quad \textrm{and} \quad f\bullet 1=f
    \end{equation}
    for all $s\in\Gamma^\infty(P,E)[[\lambda]]$ and $f\in
    C^\infty(P)[[\lambda]]$ as above.
\end{remark}

\chapter{The commutants of the deformations}
\label{cha:symbol-calculus}

For the further investigations of the deformations discussed in
Chapter~\ref{cha:deformation-on-bundles}, in particular the
computation of the commutants within the differential operators, we
will make use of the well-known symbol calculus for differential
operators. The basic idea of symbol calculus comes from the general
transition from a filtered vector space $V=\Union_{l=0}^\infty V_l$ to
the corresponding graded vector space $W=\bigoplus_{l=0}^\infty W_l$
by setting $V_0=W_0$ and $W_l=V_l/V_{l-1}$ for $l\in \mathbb{N}$. For
the filtered spaces of differential operators this procedure yields
the spaces of the so-called \emph{symbols}. In many important
geometric examples these symbols are isomorphic to other geometric
structures, mostly tensor fields. Then, many crucial properties of a
differential operator in fact can be expressed by according properties
of the corresponding geometric structures. For many concrete problems
and investigations this is a very useful approach.

In the first section of this short chapter we repeat the basic concept
for differential operators $\Diffop^\bullet(M)$ on the smooth
functions on a manifold $M$. The presentation and discussion of the
relevant structures, in particular the symmetrized covariant
derivative, result in the observation that every differential operator
can be identified with a unique series of symmetric multivector
fields. Basically, these considerations are found in
\cite[Sect.~5.4.1]{waldmann:2007a} or
\cite[Anhang~B]{neumaier:1998a}. Further information on symbol
calculus can be found in \cite{pflaum:1998b}. For our purposes we
formulate the generalized assertions for the differential operators
$\Diffop^\bullet(M, C^\infty(P))$ of the functions on the base
manifold $M$ and the total space $P$ of a surjective submersion
$\pp:P\longrightarrow M$ and for the differential operators
$\Diffop^\bullet(\Gamma^\infty(M,E))$ of the sections of a vector
bundle $\pe: E\longrightarrow M$.
Section~\ref{sec:covariant-derivatives} is dedicated to the
investigation of the crucial symmetrized covariant derivatives. It is
shown that in the situations of interest there exist adapted covariant
derivatives such that the isomorphisms established by symbol calculus
preserve further information.  As a first application of these
observations it is shown in Section~\ref{sec:respecting-fibration}
that there always exist deformation quantizations that respect a
simple algebraic property of the pullback $\pp^*:
C^\infty(M)\longrightarrow C^\infty(P)$.  In combination with the
results of Section~\ref{sec:commutant-module-structures} the developed
symbol calculus finally allows a detailed investigation of the
commutants of the deformed structures which will be performed in
Section~\ref{sec:commutants}.

\section{Symbol calculus for differential operators}
\label{sec:symbol-functions}

The starting point of the considerations is the already mentioned fact
that with respect to a chart $(U,x)$ of a manifold $M$ any
differential operator $D\in\Diffop^l(M)$ of degree $l\in \mathbb{N}_0$
on the functions of $M$ has the local form
\begin{equation}
    \label{eq:local-form-diffop-functions}
    D|_U = \sum_{r=0}^l \frac{1}{r!} D_{U}^{i_1 \dots i_r}
    \frac{\partial^r} {\partial x^{i_1} \dots \partial x^{i_r}},
\end{equation}
with uniquely defined functions $D_{U}^{i_1 \dots i_r}\in C^\infty(U)$
for $r=0,\dots, l$ which are symmetric in the indices $i_1,\dots,
i_r$. Conversely, if a $\mathbb{K}$-linear map $D:
C^\infty(M)\longrightarrow C^\infty(M)$ has the local form
\eqref{eq:local-form-diffop-functions} for an atlas of $M$ it is a
differential operator $D\in \Diffop^l(M)$. Due to the transformation
behaviour under a change of charts it follows that these functions for
$r=l$ define a global symmetric tensor field $\sigma_l(D)\in
\Gamma^\infty(M,S^kTM)$, the so-called \emph{principal symbol} of
$D$. This is defined by its restrictions
\begin{equation}
    \label{eq:leading-symbol}
    \sigma_l(D)|_U=\frac{1}{l!} D_{U}^{i_1 \dots i_l}
    \frac{\partial}{\partial x^{i_1}} \vee \dots \vee
    \frac{\partial}{\partial x^{i_l}}
\end{equation}
to the chart domains of a corresponding atlas of $M$ where $\vee$ is
the symmetric tensor product which for every $\mathbb{K}$-vector space
$V$ is defined by $v_1\vee \dots \vee v_l= \sum_{\tau\in
  S_l}v_{\tau(1)}\otimes \dots \otimes v_{\tau(l)}$. Of course, it is
$\sigma_l(D)=0$ if and only if $D\in \Diffop^{l-1}(M)$.

For the further considerations one has to choose a covariant
derivative $\nabla^M= \nabla^{TM}$ on $TM$, this means a
$\mathbb{K}$-linear map $\nabla^M: \Gamma^\infty(M,TM)\times
\Gamma^\infty(M,TM)\longrightarrow \Gamma^\infty(M,TM)$ with the
properties $\nabla_{aX}Y= a \nabla^M_XY$ and $\nabla_X(aY)=
a\nabla^M_XY+ X(a) Y$ for all $a\in C^\infty(M)$ and $X,Y\in
\Gamma^\infty(M,TM)$. Given any symmetric $l$-form $\gamma\in
\Gamma^\infty(M,S^lT^*M)$ the symmetrized covariant derivative
$\mathsf{D}_M\gamma\in \Gamma^\infty(M,S^{l+1}T^*M)$ induced by
$\nabla^M$ is defined by
\begin{equation}
    \label{eq:symmetrized-covariant-derivative-global-form}
    (\mathsf{D}_M\gamma)(X_1,\dots, X_{l+1})=
    \sum_{s=1}^{l+1}(\nabla^M_{X_s}\gamma)
    (X_1,\dots,\stackrel{s}{\wedge}, \dots, X_{l+1}),
\end{equation}
where the dual covariant derivative and its extension to the symmetric
tensor product are again denoted by $\nabla^M$. The notation
$\stackrel{s}{\wedge}$ indicates that the argument $X_s$ is
omitted. For $l=0$, this means for $a\in C^\infty(M)$, one has
$\mathsf{D}_Ma=\D a$. Due to the defining properties of a covariant
derivative it is clear that $\nabla^M$ is a differential operator
$\nabla^M\in \Diffop^{(0,1)}(\Gamma^\infty(M,TM), \Gamma^\infty(M,TM);
\Gamma^\infty(M,TM))$. Likewise, it is easy to see that the
symmetrized covariant derivative is a differential operator
$\mathsf{D}_M\in \Diffop^1(\Gamma^\infty(M,S^lT^*M);
\Gamma^\infty(M,S^{l+1}T^*M))$ for all $l\in \mathbb{N}_0$ and thus
can be restricted to open subsets $U\subseteq M$. The local expression
in a chart $(U,x)$ of $M$ turns out to be
\begin{equation}
    \label{eq:symmetrized-covariant-derivative-explicit-local}
    \mathsf{D}_M|_U= \D x^i \vee \nabla_{\frac{\partial}{\partial x^i}}.
\end{equation}
Thus, $\mathsf{D}_M$ is a derivation of degree one of the symmetric
algebra $(\bigoplus_{l=0}^\infty \Gamma^\infty(M,S^lT^*M),\vee)$.
Note that the local expression
\eqref{eq:symmetrized-covariant-derivative-explicit-local} may also
serve as a definition for $\mathsf{D}_M$.

Defining $\mathsf{D}_M^{(l)}= \frac{1}{l!}\mathsf{D}_M^l$ for $l\in
\mathbb{N}_0$ a simple induction shows that for $a\in C^\infty(M)$
\begin{equation}
    \label{eq:multiple-symmetrized-covariant-derivative}
    (\mathsf{D}_M^{(l)}a)|_U= \frac{1}{l!} \left( \frac{\partial^l a|_U}  
        {\partial x^{i_1}\dots \partial x^{i_l}} + 
        \Gamma_{i_1\dots i_l}(a|_U)\right) \D x^{i_1}\vee \dots \vee \D
    x^{i_l},
\end{equation}
where $\Gamma_{i_1\dots i_l}(a|_U)$ depends linearly on $a|_U$ and its
partial derivatives at most up to order $l-1$. Using the derivation
property of $\mathsf{D}_M$ another simple induction yields
\begin{equation}
    \label{eq:symmetrized-covariant-derivative-derivation}
    \mathsf{D}_M^{(l)}(a\cdot b)= \sum_{r=0}^l
    (\mathsf{D}_M^{(r)}a ) \vee
    (\mathsf{D}_M^{(l-r)}b )
\end{equation}
for all $a,b\in C^\infty(M)$. Since a tensor field $T_l\in
\Gamma^\infty(M,S^lTM)$ can locally be written as
\begin{equation}
    \label{eq:symmetric-multi-vector-field-local}
    T_l|_U= \frac{1}{l!} T^{i_1\dots i_l}\frac{\partial}{\partial x^{i_1}}
    \vee \dots \vee \frac{\partial}{\partial x^{i_l}}
\end{equation}
with unique functions $T^{i_1\dots i_l}\in C^\infty(U)$ one can thus
define a differential operator $D_{T_l}\in\Diffop^l(M)$ by
\begin{equation}
    \label{eq:diffop-tensor-field}
    D_{T_l} a= \frac{1}{l!} \left\langle T_l,
        \mathsf{D}_M^{(l)}a\right\rangle, 
\end{equation}
where $\langle\cdot, \cdot \rangle$ denotes the \emph{natural pairing}
of symmetric multivector fields with symmetric differential forms
defined by the local expressions, $\langle T_l,\gamma\rangle|_U=
T_l|_U(\D x^{i_1},\dots, \D x^{i_l})
\gamma|_U(\frac{\partial}{\partial x^{i_1}},\dots,
\frac{\partial}{\partial
  x^{i_l}})$. Equation~\eqref{eq:diffop-tensor-field} in fact defines
a differential operator since one locally has
\begin{eqnarray*}
    \label{eq:local-diffop-tensor-field}
    (D_{T_l} a)|_U &=& 
    \frac{1}{l!} 
    \left\langle \frac{1}{l!}
        T^{i_1\dots i_l}\frac{\partial}{\partial x^{i_1}} 
        \vee \dots \vee \frac{\partial}{\partial x^{i_l}},
        \frac{1}{l!} 
        \left( \frac{\partial^l a|_U}{\partial x^{i_1}\dots \partial
              x^{i_l}} +  \Gamma_{i_1\dots i_l}(a|_U)\right) \D
        x^{i_1}\vee \dots \vee  \D x^{i_l}\right\rangle\\
    &=& \frac{1}{l!} T^{i_1\dots i_l}
    \left( \frac{\partial^l
          a|_U}{\partial x^{i_1}\dots \partial x^{i_l}} + 
        \Gamma_{i_1\dots i_l}(a|_U)\right),
\end{eqnarray*}
confer Lemma~\ref{lemma:hom-local-diffops}. Due to the properties of
$\Gamma_{i_1\dots i_j}(a|_U)$ one finds that $\sigma_l(D_{T_l})=T_l$.
Thus, for all $D\in \Diffop^l(M)$ one has $D-D_{\sigma_l(D)}\in
\Diffop^{l-1}(M)$ and an induction finally yields the following
well-known proposition.

\begin{proposition}[Symbol calculus for $\Diffop(M)$]
    \label{proposition:symbol-calculus}
    Let $\nabla^M$ be a covariant derivative on a manifold $M$. Then,
    every differential operator $D\in \Diffop^l(M)$ of order $l$ can
    be identified with a unique series $T_0,\dots, T_l$ of symmetric
    multivector fields $T_j\in \Gamma^\infty(M,S^jTM)$, $j=0,\dots,l$,
    yielding
    \begin{equation}
        D=\sum_{j=0}^l D_{T_j},
    \end{equation}
    where $T_0= D(1)$ and $T_l=\sigma_l(D)$ are independent of the
    choice of $\nabla^M$.  The corresponding \emph{symbol map}
    $\sigma$, defined by
    \begin{equation}
        \label{eq:symbol-map-explicit}
        \sigma(D)=T_0+\dots +T_l\in  \bigoplus_{s=0}^l
        \Gamma^\infty(M,S^s TM)
    \end{equation}
    then is a vector space isomorphism
    \begin{equation}
        \label{eq:symbol-map}
        \sigma: \Diffop^\bullet(M)\longrightarrow
        \Gamma^\infty(M,S^\bullet TM)= \bigoplus_{l=0}^\infty
        \Gamma^\infty(M,S^l TM).
    \end{equation}
\end{proposition}

If $\pp:P\longrightarrow M$ is a surjective submersion one can not
only consider the above explained symbol calculus on the base and on
the total space separately but also a combined version. The following
proposition shows that a covariant derivative $\nabla^M$ on $TM$ and
an additional choice of a connection $TP=VP\oplus HP$ as in
Appendix~\ref{cha:bundle-geometry} give rise to a symbol calculus for
the differential operators $\Diffop^\bullet(M, C^\infty(P))$ where
$C^\infty(P)$ is equipped with the obvious $C^\infty(M)$-module
structure which is induced by the pullback $\pp^*$.

\begin{proposition}[Symbol calculus for $\Diffop^\bullet(M,
    C^\infty(P))$]
    \label{proposition:symbol-calculus-M-P}
    Let $\pp:P\longrightarrow M$ be a surjective submersion and let
    $\nabla^M$ be a covariant derivative on $TM$. Then, depending on
    $\nabla^M$ and a choice of a connection $TP=VP\oplus HP$ every
    differential operator $D\in \Diffop^l(M,C^\infty(P))$ of order
    $l\in \mathbb{N}_0$ can be identified with a unique series
    $T_0,\dots, T_l$ of horizontal symmetric multivector fields
    $T_j\in \Gamma^\infty(P,S^jHP)$, $j=0,\dots,l$ via $ D=
    \sum_{j=0}^l D_{T_j}$ where
    \begin{equation}
        \label{eq:symbol-calculus-M-P}
        D_{T_j}(a)= \frac{1}{j!} 
        \langle T_j, \pp^* \mathsf{D}_M^{(j)}a \rangle.
    \end{equation}
\end{proposition}

\begin{proof}
    The proof is a straightforward generalization of the
    considerations above. As a special case of the results of
    Lemma~\ref{lemma:local-form-of-Diffops} and
    Remark~\ref{remark:local-form-fixed-order} one finds that every
    differential operator $D\in \Diffop^l(M,C^\infty(P))$
    has the local form
    \begin{equation}
        \label{eq:local-form-diffop-M-P}
        D|_{\pp^{-1}(U)}(a)= \sum_{r=0}^l 
        D_{\pp^{-1}(U)}^{i_1\dots i_r}\cdot
        \pp^*\left(
            \frac{\partial^r a} {\partial x^{i_1} \dots \partial
              x^{i_r}}  
        \right)
    \end{equation}
    for all $a\in C^\infty(U)$ where $(U,x)$ is a local chart of $M$
    and where $D_{\pp^{-1}(U)}^{i_1\dots i_r}\in
    C^\infty(\pp^{-1}(U))$. For $r=l$ these functions define a global
    horizontal tensor field $\sigma_l(D)\in \Gamma^\infty(P,S^lHP)$
    via
    \begin{equation}
        \label{eq:leading-symbol-M-P}
        \sigma_l(D)|_{\pp^{-1}(U)} = \frac{1}{l!} D_{\pp^{-1}(U)}^{i_1\dots i_l} 
        \left(\frac{\partial}{\partial x^{i_1}}\right)^{\hlift}
        \vee \dots \vee  
        \left(\frac{\partial}{\partial x^{i_l}}\right)^{\hlift}.
    \end{equation}
    Locally, the horizontal lifts
    $\left\{\left(\frac{\partial}{\partial x^i}
        \right)^{\hlift}\right\}_{i=1,\dots,n}$ are a module basis of
    $\Gamma^\infty(P,HP)$ and the dual basis is given by $\{\pp^* \D
    x^i\}_{i=1,\dots, n}$. This and
    \eqref{eq:multiple-symmetrized-covariant-derivative} imply that
    every tensor field $T_j\in \Gamma^\infty(P,S^jHP)$ gives rise to a
    differential operator $D_{T_j}\in \Diffop^j(M,C^\infty(P))$ by the
    natural pairing \eqref{eq:symbol-calculus-M-P}.  Then,
    $\sigma_j(D_{T_j})=T_j$ leads to $D-D_{\sigma_l(D)}\in
    \Diffop^{l-1}(M,C^\infty(P))$ and a simple induction yields the
    assertion.
\end{proof}

Besides the generalization made in
Proposition~\ref{proposition:symbol-calculus-M-P} we will need another
one for the differential operators
$\Diffop^\bullet(\Gamma^\infty(M,E))$ of the sections of a vector
bundle $\pp: E\longrightarrow M$. The adaptions are the following. One
now always considers charts $(U,x)$ of $M$ such that $E|_U$ is trivial
and $\Gamma^\infty(U, E|_U)$ has a local module basis
$\{e_i\}_{i=1,\dots, k}\subseteq \Gamma^\infty(E|_U)$ with its dual
counterpart $\{e^i\}_{i=1,\dots, k}\subseteq
\Gamma^\infty(U,E^*|_U)$. Then, every differential operator $D\in
\Diffop^l(\Gamma^\infty(M,E))$ has the local form
\begin{equation}
    \label{eq:diffop-section-local}
    D|_U(s)= \sum_{r=0}^l \frac{1}{r!} \frac{\partial^r s^i}{\partial
      x^{i_1}\dots \partial x^{i_r}} D^{i_1\dots i_r}_{Ui}
\end{equation}
for all local sections $s=s^i e_i\in \Gamma^\infty(U,E|_U)$ with
coefficients $s^i\in C^\infty(U)$ and where $D^{i_1\dots i_r}_{Ui}\in
\Gamma^\infty(U,E|_U)$. This induces a principal symbol
$\sigma_l(D)\in \Gamma^\infty(M,S^lTM\otimes E^*\otimes E)$ by
\begin{equation}
    \label{eq:principal-symbol-sections}
    \sigma_l(D)|_U= \frac{1}{l!} \frac{\partial}{\partial x^{i_1}} 
    \vee \dots \vee \frac{\partial}{\partial x^{i_l}}\otimes
    e^i\otimes D^{i_1\dots i_l}_{Ui}.
\end{equation}
Note that $E^*\otimes E$ is isomorphic to the endomorphism bundle
$\End(E)$. With $\nabla^M$ as before and an additional covariant
derivative $\nabla^E: \Gamma^\infty(M,TM)\times
\Gamma^\infty(M,E)\longrightarrow \Gamma^\infty(M,E)$ on $E$ one gets
a new one on the tensor product $TM\otimes E$, again denoted by
$\nabla^E$. The symmetrized version $\mathsf{D}_E$ thereof is defined
in analogy to \eqref{eq:symmetrized-covariant-derivative-global-form}
and \eqref{eq:symmetrized-covariant-derivative-explicit-local}. But
now it is a $(\Gamma^\infty(M,S^\bullet T^*M), \vee)$-module
derivation of $\Gamma^\infty(M,S^\bullet T^*M \otimes E)$ satisfying
\begin{eqnarray}
    \label{eq:covariant-derivative-module-derivation}
    \mathsf{D}_E(\alpha \vee \beta \otimes s)&=& 
    (\mathsf{D}_M \alpha) \vee (\beta \otimes s) + \alpha \vee
    \mathsf{D}_E(\beta\otimes s) \nonumber \\
    &=& (\mathsf{D}_M \alpha) \vee (\beta \otimes s) + \alpha \vee
    (\mathsf{D}_M \beta) \otimes s + \alpha \vee \beta \otimes
    \mathsf{D}_E s
\end{eqnarray}
for all $\alpha,\beta \in \Gamma^\infty(M,S^\bullet T^*M)$ and $s\in
\Gamma^\infty(M,E)$. The analogue of
\eqref{eq:multiple-symmetrized-covariant-derivative} now reads
\begin{equation}
    \label{eq:symm-covariant-derivative-sections}
    (\mathsf{D}^{(l)}s)|_U = \frac{1}{l!} \left( \frac{\partial^l
          s^i}{\partial x^{i_1} \dots \partial x^{i_l}} +
        \Gamma^i_{i_1\dots i_l} (s|_U) \right) \D x^{i_1}\vee \dots \vee \D
    x^{i_l} \otimes e_i. 
\end{equation}
With the local form
\begin{equation}
    \label{eq:tensor-field-symbol-section}
    T_l|_U= \frac{1}{l!} T^{i i_1\dots i_l}_j \frac{\partial}{\partial
      x^{i_1}} \vee \dots \vee \frac{\partial}{\partial x^{i_l}}
    \otimes e^j \otimes e_i
\end{equation}
of any tensor field $T_l\in \Gamma^\infty(M,S^lTM \otimes E^* \otimes
E)$ and analogous considerations to the previous cases one finds the
following extended version of symbol calculus.

\begin{proposition}[Symbol calculus for $\Diffop(\Gamma^\infty(M,E))$]
    \label{proposition:symbol-calculus-sections}
    Let $\pp: E \longrightarrow M$ be a vector bundle and let
    $\nabla^M$ and $\nabla^E$ be covariant derivatives on $M$ and $E$,
    respectively. Then, every differential operator $D\in
    \Diffop^l(\Gamma^\infty(M,E))$ of order $l$ can be identified with
    a unique series $T_0, \dots T_l$ of tensor fields $T_j\in
    \Gamma^\infty(M,S^jTM \otimes E^*\otimes E)\cong
    \Gamma^\infty(M,S^jTM \otimes \End(E))$, $j=0, \dots, l$, such
    that $D= \sum_{j=0}^l D_{T_j}$, and where $T_l= \sigma_l(D)$ is
    independent of the choices of $\nabla^M$ and $\nabla^E$. Thus
    one has a vector space isomorphism
    \begin{equation}
        \label{eq:symbol-map-sections}
        \sigma: \Diffop^\bullet(\Gamma^\infty(M,E)) \longrightarrow
        \Gamma^\infty(M,S^\bullet TM \otimes \End(E)).
    \end{equation}
\end{proposition}

\section{Covariant derivatives on surjective submersions and principal
  fibre bundles}
\label{sec:covariant-derivatives}

As seen above, the choice of covariant derivatives is crucial for the
symbol calculus of differential operators. In the following we prove
the existence of covariant derivatives with particular properties and
discuss some resulting consequences which will be used in the later
applications. The basic idea of the following lemma comes from the
related assertions in \cite[Prop.~4.3]{bordemann:2004a:pre} and
\cite[Satz~2.3.20]{weiss:2006a}.

\begin{lemma}[Adapted covariant derivative]
    \label{lemma:covariant-derivative-sursub}
    Let $\pp:P\longrightarrow M$ be a surjective submersion and
    $TP=VP\oplus HP$ be a connection. Further, let $\nabla^M$ be a
    torsion-free covariant derivative on $TM$. Then there always
    exists a covariant derivative $\nabla^P$ on $TP$ with the
    following properties.
    \begin{enumerate}
    \item $\nabla^P$ is torsion-free, this means
        \begin{equation}
            \label{eq:torsion-free}
            \nabla^P_Z W-\nabla^P_W Z= [Z,W]
        \end{equation}
        for all $Z,W\in \Gamma^\infty(P,TP)$.
    \item $\nabla^P$ respects the vertical bundle, this means 
        \begin{equation}
            \label{eq:respect-vertical-bundle}
            \nabla^P_Z V\in \Gamma^\infty(P,VP) \quad \textrm{for all }
            V \in \Gamma^\infty(P,VP), Z\in \Gamma^\infty(P,TP).
        \end{equation}
    \item $\nabla^P$ and $\nabla^M$ satisfy
        \begin{equation}
            \label{eq:related-covariant-derivatives}
            T\mathsf{p}\circ \nabla^P_{X^{\hlift}}Y^{\hlift} =
            \left(\nabla^M_X Y \right)\circ \mathsf{p}
        \end{equation}
        for all $X,Y \in \Gamma^\infty(M,TM)$.
    \end{enumerate}
    If $P$ is a principal fibre bundle with structure group $G$ and
    $TP=VP\oplus HP$ is a principal connection, this means $T\rr_g
    H_uP= H_{ug}P$ for all $u\in P$ and $g\in G$, it is possible to
    achieve the following additional property.
    \begin{enumerate}
        \setcounter{enumi}{3}
    \item $\nabla^P$ is $G$-invariant, this means
        \begin{equation}
            \label{eq:covariant-derivative-G-invariant}
            \rr_g^* \nabla^P_Z W = \nabla^P_{\rr_g^* Z} \rr_g^*W
        \end{equation}
        for all $Z,W\in \Gamma^\infty(P,TP)$ and $g\in G$.
    \end{enumerate}
\end{lemma}

\begin{proof}
    The vertical bundle $VP$ is integrable, this means $[V,W]\in
    \Gamma^\infty(P,VP)$ for all $V,W\in \Gamma^\infty(P,VP)$. Due to
    this fact it is possible to choose a torsion-free covariant
    derivative $\nabla^{VP}: \Gamma^\infty(P,VP)\times
    \Gamma^\infty(P,VP)\longrightarrow \Gamma^\infty(P,VP)$. This can
    now be extended to a covariant derivative $\nabla^P$ on $TP$. Due
    to the Leibniz rule it is possible to define $\nabla^P$ by the
    values on vertical and horizontal vector fields. Of course one
    defines
    \begin{equation}
        \label{eq:covariant-VP}
        \nabla^P_VW=\nabla^{VP}_V W
    \end{equation}
    for $V,W\in \Gamma^\infty(P,VP)$. Using the connection in the form
    of the corresponding projection $\mathcal{P}: TP\longrightarrow
    VP$ onto the vertical space one further sets
    \begin{equation}
        \label{eq:covariant-derivative-VP-HP-general}
        \nabla^P_VH= (\id-\mathcal{P})[V,H] \quad \textrm{and}\quad
        \nabla^P_HV= \mathcal{P} [H,V] 
    \end{equation}
    where $H\in \Gamma^\infty(P,HP)$ and $V\in
    \Gamma^\infty(P,VP)$. In addition, this definition can be made
    locally by using horizontal lifts instead of arbitrary horizontal
    vectors since there exist local module bases of
    $\Gamma^\infty(P,HP)$ consisting of horizontal lifts. The
    well-known facts concerning related vector fields yielding
    \eqref{eq:consequences-related-fields} then show that
    \begin{equation}
        \label{eq:covariant-derivative-VP-HP}
        \nabla^P_VW=\nabla^{VP}_V W, \quad \nabla^P_V X^{\hlift}=0,
        \quad \textrm{and } \nabla^P_{X^{\hlift}}V=[X^{\hlift},V].
    \end{equation}
    Note that these equations also yield consistent definitions since
    $(aX)^{\hlift}= (\pp^*a) X^{\hlift}$ and $V(\pp^*a)=0$ for
    all $a\in C^\infty(M)$. For the last case of horizontal lifts we
    now use $\nabla^M$ and set
    \begin{equation}
        \label{eq:covariant-derivative-HP}
        \nabla^P_{X^{\hlift}} Y^{\hlift}= (\nabla^M_X
        Y)^{\hlift} + \frac{1}{2} 
        \left( 
            [X^{\hlift}, Y^{\hlift}]- [X,Y]^{\hlift}
        \right).
    \end{equation}
    for all $X,Y\in \Gamma^\infty(M,TM)$.  A simple calculation shows
    that the so defined $\nabla^P$ is torsion-free and it is obvious
    that $\nabla^P$ respects the vertical bundle. With
    \eqref{eq:consequences-related-fields} one finally verifies
    equation \eqref{eq:related-covariant-derivatives}.  For principal
    fibre bundles one simply chooses a covariant derivative
    $\nabla^{VP}$ that additionally is $G$-invariant. Since $\rr_g^*
    \Gamma^\infty(P,VP)=\Gamma^\infty(P,VP)$ and with
    $\rr_g^*X^{\hlift}=X^{\hlift}$ the above defined covariant
    derivative then is easily verified to be $G$-invariant. Confer
    \cite[Satz~2.3.20]{weiss:2006a} for an explicit construction of
    all this.
\end{proof}

\begin{remark}
    \begin{enumerate}
    \item Since $\nabla^P$ as in
        Lemma~\ref{lemma:covariant-derivative-sursub} is torsion-free
        and respects the vertical bundle the covariant derivative of a
        horizontal lift in vertical direction is vertical, this means
        \begin{equation}
            \label{eq:cov-derivative-horizontal-lift-vertical}
            \nabla^P_{V}X^{\hlift}\in \Gamma^\infty(P,VP) \quad
            \textrm{for all } V\in \Gamma^\infty(P,VP), X\in
            \Gamma^\infty(M,TM).
        \end{equation}
        This is obvious with $\nabla^P_{V} X^{\hlift}= \nabla^P_{X^{\hlift}}V+
        [V,X^{\hlift}]$ and \eqref{eq:consequences-related-fields}.
    \item Note that a $G$-invariant covariant derivative $\nabla^P$ on
        $TP$ already defines a covariant derivative $\nabla^M$ by
        \eqref{eq:related-covariant-derivatives}. If $\nabla^P$ is
        additionally torsion-free the same is true for $\nabla^M$.
    \end{enumerate}
\end{remark}

\begin{lemma}
    \label{lemma:covariant-derivative-pull-back}
    Let $\nabla^M$ and $\nabla^P$ be torsion-free covariant
    derivatives as in Lemma~\ref{lemma:covariant-derivative-sursub}.
    Then, for all $l \in \mathbb{N}$ the corresponding symmetrized
    covariant derivatives $\mathsf{D}_P^{(l)}$ and
    $\mathsf{D}_M^{(l)}$ are related by
    \begin{equation}
        \label{eq:symmetric-covariant-derivatives-pullback}
        \mathsf{D}_P^{(l)} \circ \mathsf{p}^*=
        \mathsf{p}^*\circ \mathsf{D}_M^{(l)}.
    \end{equation}
\end{lemma}

\begin{proof}
    One has to show that $\mathsf{D}_P (\pp^*\gamma)=
    \pp^*(\mathsf{D}_M\gamma)$ for all forms $\gamma\in
    \Gamma^\infty(M,S^lT^*M)$. This can be done by evaluation on
    arbitrary points $u\in P$ and vectors $Z_1(u),\dots, Z_{l+1}(u)\in
    T_uP$ which of course can be seen as values of vector fields
    $Z_1,\dots, Z_{l+1}\in \Gamma^\infty(P,TP)$. Due to the defining
    equation \eqref{eq:symmetrized-covariant-derivative-global-form}
    one has to show the equality of
    \begin{eqnarray}
        \label{eq:covariant-derivative-pullback}
        \left(\mathsf{D}_P (\pp^*\gamma)\right)_u (Z_1(u),\dots,
        Z_{l+1}(u)) 
        &=&
        \sum_{s=1}^{l+1} 
        \Bigg( 
        Z_s(u)\left((\pp^*\gamma)(Z_1, \dots,
            \stackrel{s}{\wedge},\dots, Z_{l+1})\right) \\
        && 
        - \sum_{i\neq s}
        (\pp^*\gamma) (Z_1,\dots, \nabla^P_{Z_s} Z_i,
        \dots, \stackrel{s}{\wedge}, \dots,Z_{l+1})(u)
        \Bigg) \nonumber
    \end{eqnarray}
    and
    \begin{equation}
        \label{eq:pullback-covariant-derivative}
        \left(\pp^*(\mathsf{D}_M \gamma)\right)_u (Z_1(u),\dots, Z_{l+1}(u))=
        (\mathsf{D}_M \gamma)_{\pp(u)} (T_u\pp Z_1(u),\dots, T_u\pp Z_{l+1}(u)).
    \end{equation}
    This can be done by a case differentiation since it is clear with
    the connection $TP=VP\oplus HP$ that the relevant vectors can be
    seen as the direct sum $Z_i(u)=V_i(u)+ X_i^{\hlift}(u)$ of values
    of vertical vector fields $V_i\in \Gamma^\infty(P,VP)$ and
    horizontal lifts $X_i^{\hlift}\in \Gamma^\infty(P,HP)$ of vector
    fields $X_i\in \Gamma^\infty(M,TM)$. Since $\nabla^P$ respects the
    vertical bundle and with
    \eqref{eq:cov-derivative-horizontal-lift-vertical} it follows that
    both \eqref{eq:covariant-derivative-pullback} and
    \eqref{eq:pullback-covariant-derivative} are zero if at least one
    vector field is vertical and all others are horizontal
    lifts. Finally, let all vector fields be horizontal lifts. Then,
    both expressions turn out to be
    \begin{equation*}
        \sum_{s=1}^{l+1}
        \left(
            X_s(\gamma(X_1\dots, \stackrel{s}{\wedge},\dots,
            X_{l+1}))- \sum_{i\neq s} \gamma(X_1,\dots,
            \nabla^M_{X_s}X_i,\dots, \stackrel{s}{\wedge},\dots,
            X_{l+1})
        \right)(\pp(u)),
    \end{equation*}
    since $X_s^{\hlift}(\pp^*a)=\pp^*(X_s(a))$ for all $a\in
    C^\infty(M)$.
\end{proof}

For a principal fibre bundle the $G$-invariance of $\nabla^P$ induces
the $G$-invariance of the corresponding symmetrized covariant
derivative.

\begin{lemma}
    \label{lemma:symm-cov-der-G-invariant}
    Let $\pp:P\longrightarrow M$ be a principal fibre bundle and let
    $\nabla^P$ be a $G$-invariant covariant derivative on $TP$. Then,
    for all $l\in \mathbb{N}_0$ the corresponding symmetrized
    covariant derivative satisfies
    \begin{equation}
        \label{eq:symmetrized-covariant-derivative-G-invariance}
        \mathsf{D}_P^{(l)} \circ \rr_g^*= \rr_g^* \circ \mathsf{D}_P^{(l)}
    \end{equation}
    for all $g\in G$.
\end{lemma}

\begin{proof}
    With the compatibility of the pullbacks $\rr_g^*$ and the
    $G$-invariance of $\nabla^P$ one computes with
    \eqref{eq:symmetrized-covariant-derivative-global-form} for all
    $\gamma\in\Gamma^\infty(P,S^lT^*P)$ and vector fields $Z_1,\dots,
    Z_{l+1}\in \Gamma^\infty(P,TP)$
    \begin{eqnarray*}
        \lefteqn{ \left(\mathsf{D}_P (\rr^*_g \gamma)\right) (\rr^*_g
          Z_1,\dots, \rr^*_g Z_{l+1})}\\
        &=& \sum_{s=1}^{l+1}
        \left(
            (\rr^*_g  Z_s) (\rr^*_g \gamma(
            \rr^*_g Z_1,\dots,\stackrel{s}{\wedge},\dots,
            \rr^*_g Z_{l+1})) - \sum_{i\neq s} \rr^*_g \gamma ( \rr^*_g 
            Z_1, \dots, \nabla^P_{\rr^*_g  Z_s}  \rr^*_g Z_i, \dots,
            \stackrel{s}{\wedge}, \dots,  \rr^*_g Z_{l+1})
        \right)\\
        &=&  \rr^*_g \sum_{s=1}^{l+1}
        \left(
            Z_s (\gamma(Z_1,\dots,\stackrel{s}{\wedge},\dots,
            Z_{l+1})) - \sum_{i\neq s} \gamma (Z_1, \dots,
            \nabla^P_{Z_s} Z_i, \dots, 
            \stackrel{s}{\wedge}, \dots,Z_{l+1})
        \right)\\
        &=& (\rr^*_g  (\mathsf{D}_P \gamma)) (\rr^*_g Z_1,\dots,
        \rr^*_g Z_{l+1}).
    \end{eqnarray*}
    This shows
    \eqref{eq:symmetrized-covariant-derivative-G-invariance}.
\end{proof}

The above result concerning $G$-invariance has a useful generalization
with respect to the covariant derivatives occurring in
Proposition~\ref{proposition:symbol-calculus-sections}. For this
purpose we need the following lemma.

\begin{lemma}[$G$-invariant covariant derivatives for equivariant
    vector bundles]
    Every equivariant vector bundle $\pe:E\longrightarrow P$ over a
    principal fibre bundle $\pp:P\longrightarrow M$ with structure
    group $G$ admits a $G$-invariant covariant derivative $\nabla^E:
    \Gamma^\infty(P,TP)\times \Gamma^\infty(P,E)\longrightarrow
    \Gamma^\infty(P,E)$, this means with
    \begin{equation}
        \label{eq:invariant-cov-derivative-bundle}
        g\acts \nabla^E_Zs= \nabla^E_{\rr_g^*Z} (g\acts s)
    \end{equation}
    for all $Z\in \Gamma^\infty(P,TP)$, $s\in \Gamma^\infty(P,E)$, and
    $g\in G$.
\end{lemma}

\begin{proof}
    The assertion is nothing but a special case of
    \cite[Cor.~B.38]{guillemin.ginzburg.karshon:2002a} stating that
    every equivariant vector bundle with a proper group action on the
    base manifold admits an invariant connection. In our case the
    principal right action clearly satisfies the required condition.
\end{proof}

\begin{lemma}
    \label{lemma:symm-cov-der-G-invariant-vector-bundle}
    Let $\pe: E\longrightarrow P$ be an equivariant vector bundle over
    the principal fibre bundle $\pp:P\longrightarrow M$ as in
    Section~\ref{subsec:geometric-situation-principal-bundles}. Further,
    let $\nabla^E$ and $\nabla^P$ be $G$-invariant covariant
    derivatives on $E$ and $TP$, respectively. Then, the induced
    covariant derivative $\nabla^E$ on $TP\otimes E$ is $G$-invariant
    with respect to the representation on the sections $\alpha\otimes
    s\in \Gamma^\infty(P,S^\bullet T^*P\otimes E)$ given by
    \begin{equation}
        \label{eq:action-tensor-product}
        g \acts (\alpha\otimes s)=
        (\rr^*_g \alpha) \otimes (g\acts s). 
    \end{equation}
    Moreover, the same is true for the corresponding symmetrized
    covariant derivative, this means one has
    \begin{equation}
        \label{eq:symmetrized-covariant-derivative-G-inv-bundle}
        \mathsf{D}_E^{(l)} (g \acts s)= g\acts (\mathsf{D}_E^{(l)} s)
    \end{equation}
    for all $l\in \mathbb{N}_0$, $s \in \Gamma^\infty(E)$, and $g\in G$.
\end{lemma}

\begin{proof}
    The proof is an analogous computation to the one in the proof of
    Lemma~\ref{lemma:symm-cov-der-G-invariant} making use of the
    involved definitions. In particular, one needs the fact that for
    $\gamma \in \Gamma^\infty(P,S^l T^*P\otimes E)$ and $Z,Z_1,\dots,
    Z_l\in \Gamma^\infty(P,TP)$ one has
    \begin{equation}
        \label{eq:covariant-derivative-pairing-sections}
        (\nabla^E_Z \gamma) (Z_1,\dots, Z_l)=
        \nabla^E_Z(\gamma (Z_1,\dots ,Z_l))- \sum_{i=1}^l
        \gamma(Z_1,\dots, \nabla^P_Z Z_i, \dots, Z_l)
    \end{equation}
    and
    \begin{equation}
        \label{eq:compatibility-action-pairing}
        g \acts (\gamma(Z_1,\dots, Z_l))= (g\acts \gamma)(\rr^*_g
        Z_1,\dots, \rr^*_g Z_l).
    \end{equation}
    % Explicitly, for all $\gamma\in
    % \Gamma^\infty(S^lT^*P \otimes E)$ and $Z_1,\dots, Z_{l+1}\in
    % \Gamma^\infty(TP)$ the slightly modified computation now reads
    % \begin{eqnarray*}
    %     \lefteqn{\left(\mathsf{D}_E (g \acts \gamma)\right)
    %       (\rr^*_g Z_1,\dots, \rr^*_g Z_{l+1})}\\
    %     &=& \sum_{s=1}^{l+1} \left( \nabla^E_{\rr^*_g Z_s} ((g \acts
    %         \gamma) (\rr^*_g Z_1,\dots,\stackrel{s}{\wedge},\dots,
    %         \rr^*_g Z_{l+1})) - \sum_{i\neq s} (g \acts \gamma)
    %         (\rr^*_g Z_1, \dots, \nabla^P_{\rr^*_g Z_s} \rr^*_g Z_i,
    %         \dots, \stackrel{s}{\wedge}, \dots, \rr^*_g Z_{l+1})
    %     \right)\\
    %     &=& g \acts \left( \nabla^E_{Z_s}
    %         (\gamma(Z_1,\dots,\stackrel{s}{\wedge},\dots, Z_{l+1}))
    %         - \sum_{i\neq s} \gamma (Z_1, \dots, \nabla^P_{Z_s} Z_i,
    %         \dots, \stackrel{s}{\wedge}, \dots,Z_{l+1})
    %     \right)\\
    %     &=& \left(g \acts (\mathsf{D}_P \gamma)\right) (\rr^*_g
    %     Z_1,\dots, \rr^*_g Z_{l+1}).
    % \end{eqnarray*}
\end{proof}

The above results have the following important consequence.

\begin{proposition}
    \label{proposition:symbol-calculus-verticality-invariance}
    Let $\pp: P\longrightarrow M$ be a surjective submersion and let
    $\nabla^M$ and $\nabla^P$ be torsion-free covariant derivatives on
    $TM$, and $TP$ respectively, as in
    Lemma~\ref{lemma:covariant-derivative-sursub}. Further, let
    $\nabla^E$ be an arbitrary covariant derivative on the vector
    bundle $\pe: E\longrightarrow P$. Then, the corresponding symbol
    map
    \begin{equation}
        \label{eq:symbol-map-sections-submersion}
        \sigma: \Diffop^\bullet(\Gamma^\infty(P,E)) \longrightarrow
        \Gamma^\infty(P,S^\bullet TP \otimes \End(E)) 
    \end{equation}
    restricts to a vector space isomorphism
    \begin{equation}
        \label{eq:isomorphism-vertical-operators}
        \sigma: \Diffopverg{\bullet}(\Gamma^\infty(P,E)) \longrightarrow
        \Gamma^\infty(P,S^\bullet VP \otimes \End(E))
    \end{equation}
    from the vertical differential operators to the vertical symmetric
    multivector fields with values in the endomorphism bundle
    $\End(E)$.

    Moreover, if $\pe: E\longrightarrow P$ is an equivariant vector
    bundle over the principal fibre bundle $\pp: P\longrightarrow M$
    and if the covariant derivatives $\nabla^P$ and $\nabla^E$ are
    $G$-invariant, the symbol map $\sigma$ is also $G$-invariant. This
    means that
    \begin{equation}
        \label{eq:symbol-map-G-invariance}
        g\acts \sigma(D)= \sigma (g\acts D)
    \end{equation}
    with the naturally induced representations of $G$ on the
    differential operators $D\in \Diffop^\bullet(\Gamma^\infty(P,E))$
    and the sections $\Gamma^\infty(P,S^\bullet TP\otimes \End(E))$.
\end{proposition}

\begin{proof}
    In order to prove the isomorphism
    \eqref{eq:isomorphism-vertical-operators} consider $D=
    \sum_{j=0}^l D_{T_j}\in \Diffop^l(\Gamma^\infty(P,E))$ with
    $T_j\in \Gamma^\infty(P,S^jTP \otimes \End(E))$. First assume that
    $T_j\in \Gamma^\infty(S^jVP\otimes \End(E))$ are vertical for all
    $j=1,\dots, l$. By definition this yields for all $s\in
    \Gamma^\infty(P,E)$ and $a\in C^\infty(M)$ that $j! D_{T_j}( \pp^*
    a \cdot s)= \langle T_j, \mathsf{D}_E^{(j)}(\pp^*a \cdot s)
    \rangle = \langle T_j, \sum_{r=0}^j \pp^* \mathsf{D}_M^{(r)}\pp^*a
    \vee \mathsf{D}_E^{(j-r)} s \rangle= \pp^*a \cdot \langle T_j,
    \mathsf{D}_E^{(j)}s \rangle$. This shows that $D$ is vertical. The
    other implication is slightly more technical. If $D(\pp^*a\cdot
    s)= \pp^*a \cdot D(s)$ this implies
    \begin{eqnarray*}
        0= \sum_{j=1}^l \frac{1}{j!} \langle T_j, \sum_{r=0}^{j-1}
        \pp^* \mathsf{D}_M^{(j-r)} a \vee \mathsf{D}_E^{(r)} s\rangle
        &=&
        \sum_{j=1}^l \sum_{r=0}^{j-1} \frac{1}{j!} \binom{j}{r}
        \langle \langle T_j, \pp^* \mathsf{D}_M^{(j-r)}a \rangle,
        \mathsf{D}^{(r)}_E s\rangle \\
        &=& 
        \sum_{r=0}^{l-1} \frac{1}{r!} \langle \langle \sum_{j=r+1}^l
        \frac{1}{(j-r)!} T_j, \pp^* \mathsf{D}_M^{(j-r)} a\rangle,
        \mathsf{D}_E^{(r)} s\rangle
    \end{eqnarray*}
    for all $s\in \Gamma^\infty(P,E)$ and $a\in C^\infty(M)$. Due to the
    symbol calculus this shows that
    \begin{equation*}
        0= \sum_{j=r+1}^l \frac{1}{(j-r)!} \langle T_j, \pp^*
        \mathsf{D}_M^{(j-r)} a\rangle \quad \textrm{for all} \quad
        r=0, \dots, l-1.
    \end{equation*}
    For $r=l-1$ one gets $\langle T_l, \pp^* \mathsf{D}_M a \rangle$
    which locally for a chart $(U,x)$ of $M$ implies that $\inss(\pp^*
    \D x^i) T_l=0$. Then, the successive treatment of the cases $r=l-2$
    to $r=0$ yields that $T_1,\dots, T_l$ can not have any horizontal
    contribution showing that $T_0+ \dots + T_l\in
    \bigoplus_{r=0}^l\Gamma^\infty(P,S^rVP \otimes \End(E))$. This
    proves the isomorphism
    \eqref{eq:isomorphism-vertical-operators}.

    With the additional assumptions,
    Lemma~\ref{lemma:symm-cov-der-G-invariant-vector-bundle}, and the
    properties of the naturally involved representations of $G$ a
    simple computation for $T_j \in \Gamma^\infty(P,S^j TP
    \otimes \End(E))$ gives $(g\acts D_{T_j})(s)= g\acts
    \frac{1}{j!}\langle T_j, \mathsf{D}_E^{(j)}(g^{-1}\acts s) \rangle
    = \frac{1}{j!}\langle g\acts T_j, \mathsf{D}_E^{(j)} s
    \rangle$. There one has used the induced representation of $G$ on
    the sections $\phi\in \Gamma^\infty(P,\End(E))$. If the action on
    $E$ is a right action $\rR_g$ one has $(g\acts \phi)= \rR_{g^{-1}}
    \circ \phi_{\circ \rr_g} \circ \rR_g$ with the obvious meaning. In
    any case, the resulting equation
    \begin{equation}
        \label{eq:symbol-map-equivariant}
        g\acts D_{T_j}= D_{g\acts T_j}
    \end{equation}
    for all $j\in \mathbb{N}_0$ is nothing but the stated
    $G$-invariance of the symbol map $\sigma$.
\end{proof}

\section{Deformations that preserve the fibration}
\label{sec:respecting-fibration}

In Section~\ref{sec:global-cohomologies-deformations} the question
concerning existence and uniqueness up to equivalence of differential
and $G$-invariant deformations has been answered completely. In this
section we concentrate on a further property of deformation
quantizations of surjective submersions and principal fibre bundles,
this means the deformed right module structures $\bullet$ of the
functions $C^\infty(P)[[\lambda]]$ with respect to a star product
algebra $(C^\infty(M)[[\lambda]],\star)$.

\begin{definition}[Fibration preserving deformation]
    Let $\bullet$ be a deformation quantization of a surjective
    submersion or a principal fibre bundle $\pp:P\longrightarrow
    M$. The structure $\bullet$ is said to \emph{preserve the
      fibration} if
    \begin{equation}
        \label{eq:respecting-fibration}
        (\pp^*a)\bullet b= \pp^*(a\star b)
    \end{equation}
    for all $a,b\in C^\infty(M)[[\lambda]]$.
\end{definition}

Such deformations are relevant to consider since they respect a
classically given algebraic property. In other words, the property
states that the pullback $\pp^*: C^\infty(M)[[\lambda]]\longrightarrow
C^\infty(P)[[\lambda]]$ in the deformed situation still is a module
morphism along the identity when regarding the algebra
$(C^\infty(M)[[\lambda]], \star)$ as a module over itself.

In the subsequent considerations we show that there always exist
deformations that preserve the fibration. Before doing this we note
that there exists an equivalent characterization of the property
\eqref{eq:respecting-fibration}.

\begin{lemma}
    A deformation quantization $\bullet$ preserves the fibration if
    and only if it satisfies the equation
    \begin{equation}
        \label{eq:equivalent-to-respecting-fibration}
        1\bullet a= \pp^*a
    \end{equation}
    for all $a\in C^\infty(M)$.
\end{lemma}

\begin{proof}
    The assertion is obvious due to the given right module property
    and the property $1\star a=a$ of any star product.
\end{proof}

Given a deformation $\bullet$, the aim will be to find an appropriate
equivalence transformation such that the transformed deformation
preserves the fibration. In order to find such a series of
differential operators we make use of the symbol calculus as explained
in the Propositions~\ref{proposition:symbol-calculus} and
\ref{proposition:symbol-calculus-M-P}.

\begin{lemma}
    \label{lemma:function-induced-diffop}
    Let $D\in \Diffop^L(M,\Diffop^l(P))$ be a differential operator
    with $L\in \mathbb{N}_0^k$, $k\in \mathbb{N}$, $l\in\mathbb{N}_0$.
    Then, every function $f\in C^\infty(P)$ gives rise to a
    differential operator $D_f\in \Diffop^L(M, C^\infty(P))$ which is
    defined by
    \begin{equation}
        \label{eq:function-induced-diffop}
        D_f(a_1,\dots, a_k)= D(a_1,\dots,a_k)(f)
    \end{equation}
    for all $a_1,\dots,a_k\in C^\infty(M)$.
    % Moreover, if $D\in \Diffop^L(M,\Diffop^l(P)^G)$ takes values in
    % the $G$-invariant differential operators of a principal fibre
    % bundle, the differential operator $D_f$ induced by a
    % $G$-invariant function $f=\pp^*a$ with $a\in C^\infty(M)$
    % satisfies
    % \begin{equation}
    %     \label{eq:function-induced-diffop-invariant}
    %     \rr_g^*\circ D_f= D_f.
    % \end{equation}
\end{lemma}

\begin{proof}
    With the present algebraic definition of differential operators
    the proof is a straightforward induction over $r=|L|$.
    % The assertion concerning $G$-invariance is clear, too.
\end{proof}

With the preceding lemmas it is now possible to show that for every
deformation quantization $\bullet$ one can explicitly construct a
particular equivalence transformation which relates the both sides of
Equation~\eqref{eq:equivalent-to-respecting-fibration}.

\begin{lemma}
    \label{lemma:special-diffop}
    Let $\bullet$ be a deformation quantization of a surjective
    submersion $\mathsf{p}: P \longrightarrow M$. Then there exists a
    formal series $T = \id + \sum_{r=1}^\infty \lambda^r
    T_r$ of differential operators $T_r \in
    \Diffop(P)$ such that
    \begin{equation}
        \label{eq:Diffop-deformation-quantization}
        T(\mathsf{p}^*a)
        = 1 \bullet a  
    \end{equation}
    for all $a\in C^\infty(M)$.

    If $\bullet$ is a deformation quantization of a principal fibre
    bundle it can be achieved that $T$ is a series of $G$-invariant
    differential operators.
\end{lemma}

\begin{proof}
    The deformed module structure $\bullet$ is given by a formal
    series $\rho =\sum_{r=0}^\infty \lambda^r\rho_r$ of differential
    operators $\rho_r\in \Diffop^{l_r}(M,\Diffop^{m_r}(P))$ with
    $l_r,m_r\in
    \mathbb{N}_0$. Lemma~\ref{lemma:function-induced-diffop} then
    shows that $D_r(a)=\rho_r(a)(1)$ defines a differential operator
    $D_r\in \Diffop^{l_r} (M,C^\infty(P))$ for all $r\in
    \mathbb{N}_0$. After the choice of a connection $HP= VP\oplus HP$
    and a torsion-free covariant derivative $\nabla^M$ on $TM$
    Proposition~\ref{proposition:symbol-calculus-M-P} shows that
    \begin{equation}
        \label{eq:diff-M-P}
        D_r(a) = \sum_{s=0}^{l_r} \frac{1}{s!}
        \left\langle 
            T^r_s, \mathsf{p}^* \mathsf{D}_M^{(s)} a
        \right\rangle 
        \quad \textrm{with unique} \quad T^r_s \in
        \Gamma^\infty(P,S^sHP). 
    \end{equation}
    Using a torsion-free covariant derivative $\nabla^P$ on $TP$ with
    the properties of Lemma~\ref{lemma:covariant-derivative-sursub} we
    define $T_r \in \Diffop^{l_r}(P)$ by
    \begin{equation}
        \label{eq:diff-P}
        T_r(f)
        =  \sum_{s=0}^{l_r} \frac{1}{s!}
        \left\langle T^r_s, \mathsf{D}_P^{(s)} f \right\rangle. 
    \end{equation}
    By the definition of a deformation quantization one has
    $D_0=\pp^*\in \Diffop^0(M,C^\infty(P))$ and with
    $\mathsf{D}_P^{(0)}=\id$ this implies that $T_0=\id$. Due to
    property \eqref{eq:symmetric-covariant-derivatives-pullback} one
    finds $T_r(\mathsf{p}^* a)= D_r(a)$ which shows that the so
    constructed series $T=\sum_{r=0}^\infty\lambda^r T_r$ has the
    required property.

    In the principal fibre bundle case the proof is slightly
    different. The given $G$-invariance of the $\rho_r$ implies that
    $\rr_g^*(D_r(a))=D_r(a)$. Thus, one has $D_r=\pp^*\circ D^M_r$
    with differential operators $D^M_r\in \Diffop^{l_r}(M)$ for which
    we apply the ordinary symbol calculus with respect to $\nabla^M$
    as above. Choosing a principal connection and using the canonical
    extension of the horizontal lift to multivector fields one gets
    \begin{equation}
        \label{eq:diff-M-P-principal-bundle}
        D_r(a) = \pp^*\sum_{s=0}^{l_r} \frac{1}{s!}
        \left\langle 
            T^r_s, \mathsf{D}_M^{(s)} a
        \right\rangle = \sum_{s=0}^{l_r} \frac{1}{s!}
        \left\langle 
            \left(T^r_s\right)^{\hlift}, \mathsf{D}_P^{(s)} a
        \right\rangle
        \quad \textrm{with unique} \quad T^r_s \in
        \Gamma^\infty(M,S^sTM). 
    \end{equation}
    Then one defines the $T_r\in \Diffop^{l_r}(P)$ in analogy to
    \eqref{eq:diff-P}. If the used $\nabla^P$ is $G$-invariant the
    compatibility of the natural pairing with the pullback $\rr_g^*$,
    the $G$-invariance of the horizontal lifts, and equation
    \eqref{eq:symmetrized-covariant-derivative-G-invariance} show that
    the $T_r$ are $G$-invariant.
\end{proof}

With this preceding lemma we finally come to the aspired result.

\begin{proposition}[Fibration preserving deformations]
    \label{proposition:fibration-sursub}
    Every surjective submersion and every principal fibre bundle
    admits a deformation quantization which preserves the fibration.
\end{proposition}

\begin{proof}
    Let $\bullet$ be an arbitrary deformation quantization with
    respect to a given star product $\star$. Since the map $T$ in
    Lemma~\ref{lemma:special-diffop} has all properties of an
    equivalence transformation, $f\bullett a= T^{-1}\left(Tf \bullet a
    \right)$ defines a new deformation quantization $\bullett$. With
    \eqref{eq:Diffop-deformation-quantization} we then get
    $(\mathsf{p}^*a) \bullett b =T^{-1}(1\bullet (a\star b))=
    \mathsf{p}^*( a\star b)$ for all $a,b \in C^\infty(M)$.

    The additional $G$-invariance of the map $T$ for principal fibre
    bundles guarantees the $G$-invariance of $\bullett$ if $\bullet$
    has this property.
\end{proof}

\begin{remark}
    Note that the above proof and Lemma~\ref{lemma:special-diffop} not
    only yield the existence of deformations which preserve the
    fibration. Moreover, the above considerations in fact provide an
    explicit procedure to find such a deformation when starting with
    an arbitrary one by constructing the corresponding equivalence
    transformation.
\end{remark}

\section{The commutants of the deformed right modules}
\label{sec:commutants}

As seen in Chapter~\ref{cha:deformation-algebras-modules}, the fact
that the crucial Hochschild cohomologies are trivial not only implies
that deformations exist and are unique up to equivalence, but also
that one has a lot of information concerning the corresponding
commutants. In the subsequent considerations we apply the general
results of Section~\ref{sec:commutant-module-structures} to the
particular situation of vector bundles $\pe: E\longrightarrow P$ over
surjective submersions and principal fibre bundles.

First of all we concentrate on surjective submersions
$\pp:P\longrightarrow M$ and deformations $\bullet$ as in
Theorem~\ref{theorem:deformations-surjective-submersion}. Due to the
vanishing first cohomology group
\begin{equation}
    \label{eq:vanishing-first-cohomology-group}
    \HHdiff^1(M,\Diffop(\Gamma^\infty(P,E)))=\{0\},
\end{equation}
the results of Proposition~\ref{proposition:description-commutant}
ensure that every choice of a deformation quantization $\bullet$ and a
decomposition
\begin{equation}
    \label{eq:decomposition-diffop}
    \Diffop(\Gamma^\infty(P,E))=\Diffopver(\Gamma^\infty(P,E))\oplus
    \overline{\Diffopver(\Gamma^\infty(P,E))}
\end{equation}
of all differential operators into the classical commutant and a
complementary subspace lead to an isomorphism
\begin{equation}
    \label{eq:iso-commutant-surj-sub}
    \rho': 
    \Diffopver(\Gamma^\infty(P,E)) [[\lambda]]
    \longrightarrow 
    \{ D\in \Diffop(\Gamma^\infty(P,E)) [[\lambda]] \: | \: D(s\bullet
    a) = D(s)\bullet a\} 
\end{equation}
onto the commutant of the deformed module structure with $\rho'=\id +
\sum_{r=1}^\infty \lambda^r \rho'_r$. According to
Corollary~\ref{corollary:deformation-commutant} this further induces
an associative deformation $\left(\Diffopver (\Gamma^\infty(P,E))
    [[\lambda]], \mu' \right)$ of the algebra
$(\Diffopver(\Gamma^\infty(P,E)),\circ)$ of vertical differential
operators with the usual composition of maps. Using the definitions
\begin{equation}
    \label{eq:bulletp-starp}
    A \bulletp s = \rho'(A)s
    \quad 
    \textrm{and}
    \quad A \starp B= \mu'(A,B)
\end{equation}
for all $A, B \in \Diffopver(\Gamma^\infty(P,E)) [[\lambda]]$ and $s
\in \Gamma^\infty(P,E)[[\lambda]]$, the sections
$C^\infty(P)[[\lambda]]$ inherit a $(\starp, \star)$-bimodule
structure which is shortly denoted by
\begin{equation}
    \label{eq:bimodule-surjective-submersion}
    {}_{(\Diffopver (\Gamma^\infty(P,E))[[\lambda]],\starp)} 
    (\bulletp,\Gamma^\infty(P,E)[[\lambda]],
    \bullet)_{(C^\infty(M)[[\lambda]],\star)}. 
\end{equation}

In the special case of functions $C^\infty(P)$, this means when
considering deformation quantizations of surjective submersion one
gets the following result, see
\cite[Prop.~4.18]{bordemann.neumaier.waldmann.weiss:2007a:pre}

\begin{proposition}[Mutual commutants for the module $C^\infty(P)$ of
    functions]
    \label{proposition:mutual-commutants}
    The commutant within the differential operators of
    $\Diffopver(P)[[\lambda]]$ acting via $\bulletp$ on
    $C^\infty(P)[[\lambda]]$ is given by $C^\infty(M)[[\lambda]]$
    acting via $\bullet$. Thus the two algebras in
    \eqref{eq:bimodule-surjective-submersion} are mutual commutants.
\end{proposition}

\begin{proof}
    Let $A = \sum_{r=0}^\infty \lambda^r A_r\in
    \Diffop^\bullet(P)[[\lambda]]$ satisfy $A(D \bulletp f) = D
    \bulletp (Af)$ for all $D \in \Diffopver(P)[[\lambda]]$ and $f \in
    C^\infty(P)[[\lambda]]$. It follows that in zeroth order $A_0$
    commutes with all undeformed left multiplications with functions
    $D \in C^\infty(P)$. Since $C^\infty(P)$ is a unital algebra this
    is only possible if $A_0 \in C^\infty(P)$. Further we know that
    $A_0$ commutes with all $D \in \Diffopver(P)$ in zeroth order, so
    $A_0$ is constant in fibre directions. This means nothing but $A_0
    = \mathsf{p}^*a_0$ for some $a_0 \in C^\infty(M)$.  Since $A(D
    \bulletp f) - (D \bulletp f) \bullet a_0 = D \bulletp (Af - f
    \bullet a_0)$ and since $Af - f \bullet a_0$ has vanishing zeroth
    order, a simple induction shows that $Af=f\bullet a$ for a series
    $a=\sum_{r=0}^\infty\lambda^r a_r\in C^\infty(M)[[\lambda]]$.
\end{proof}

For $G$-invariant deformations we can apply the results of
Section~\ref{sec:G-invariant-types}. In order to achieve a
$G$-invariant bimodule structure
\eqref{eq:bimodule-surjective-submersion} it is necessary to find a
$G$-invariant decomposition \eqref{eq:decomposition-diffop}.  As we
will see now,
Proposition~\ref{proposition:symbol-calculus-verticality-invariance}
guarantees that this crucial condition can be satisfied. Taking a
$G$-invariant covariant derivative $\nabla^E$ on the vector bundle
$\pe: E\longrightarrow P$ one obtains a $G$-invariant symbol map
$\sigma$ as in \eqref{eq:symbol-map-sections-submersion}. The choice
of a principal connection, this means a $G$-invariant decomposition
$TP= VP \oplus HP$ as in Appendix~\ref{cha:bundle-geometry}, naturally
induces a $G$-invariant decomposition
\begin{equation}
    \label{eq:invariant-decomposition-vector-fields}
    \Gamma^\infty(P,S^\bullet TP \otimes \End(E))=
    \Gamma^\infty(P,S^\bullet VP \otimes \End(E)) \oplus 
    \cc{\Gamma^\infty(P,S^\bullet VP \otimes \End(E))}.
\end{equation}
The complementary subspace $\cc{\Gamma^\infty(P,S^\bullet VP
  \otimes \End(E))}$ consists of all tensor fields having a
nontrivial horizontal component. The invariance of the composition is
clear with the properties $T\rr_g HP= HP$ and $T\rr_g VP=VP$ and the
given action on the vector fields via pullback with the principal
right action. With the symbol map $\sigma$ from
Proposition~\ref{proposition:symbol-calculus-verticality-invariance}
one can define
\begin{equation}
    \label{eq:complement-by-symbol-calculus}
    \cc{\Diffopver(\Gamma^\infty(P,E))}=
    \sigma^{-1}(\cc{\Gamma^\infty(P,S^\bullet VP \otimes \End(E))}).
\end{equation}
Due to the isomorphism \eqref{eq:isomorphism-vertical-operators} and
the $G$-invariance of $\sigma$ this yields a $G$-invariant
decomposition \eqref{eq:decomposition-diffop}, this means
\begin{equation}
    \label{eq:invariant-complement}
    G \acts \cc{\Diffopver(\Gamma^\infty(P,E))} \subseteq
    \cc{\Diffopver(\Gamma^\infty(P,E))}. 
\end{equation}
Note that the classical commutant $\Diffopver(\Gamma^\infty(P,E))$ is
$G$-invariant by the compatibility of the differential type with any
$G$-action, confer the observations of the
Remarks~\ref{remark:representations-and-cohomology} and
\ref{remark:G-invariant-differential-type}. Since the cohomology
groups $\HHdiff^1(M,\Diffop(\Gamma^\infty(P,E)))=\{0\}$ and
$\HHdiff^1(M,\Diffop(\Gamma^\infty(P,E))^G)=\{0\}$ are trivial one can
apply Proposition~\ref{proposition:commutant-G-invariance}, confer
\cite[Thm.~5.8]{bordemann.neumaier.waldmann.weiss:2007a:pre}.

\begin{theorem}
    \label{theorem:bimodule-equivariant-bundle}
    Let $\pe: E\longrightarrow P$ be an equivariant vector bundle over
    the principal fibre bundle $\pp:P \longrightarrow M$ and let
    $\bullet$ be a $G$-invariant deformation of the module structure
    of the sections $\Gamma^\infty(P,E)$ with respect to a star
    product $\star$ as in
    Theorem~\ref{theorem:deformations-principal-bundles}. Then there
    exists a bimodule structure
    \begin{equation}
        \label{eq:bimodule-equivariant-bundle}
        {}_{(\Diffopver (\Gamma^\infty(P,E))[[\lambda]],\starp)} 
        (\bulletp,\Gamma^\infty(P,E)[[\lambda]],
        \bullet)_{(C^\infty(M)[[\lambda]],\star)}
    \end{equation}
    as in \eqref{eq:bimodule-surjective-submersion} with the further
    property that $\starp$ and $\bulletp$ are $G$-invariant, this
    means
    \begin{equation}
        \label{eq:deformed-commutant-algebra-G-invariant}
        g\acts (A \starp B)= (g\acts A)\starp (g\acts B)
    \end{equation}
    and
    \begin{equation}
        \label{eq:deformed-commutant-left-module-G-invariant}
        g\acts (A \bulletp s)= (g\acts A)\bulletp (g\acts s)
    \end{equation}
    for all $A,B\in \Diffopver(\Gamma^\infty(P,E))[[\lambda]]$, $s\in
    \Gamma^\infty(P,E)[[\lambda]]$ and $g\in G$.
    Moreover, the deformed bimodule structure
    \eqref{eq:bimodule-equivariant-bundle} is unique up to
    $G$-invariant equivalence. 
\end{theorem}

\begin{remark}[Infinitesimal gauge transformations]
    Interpreting the deformed module structures in the context of
    classical gauge theories the consideration of the commutant within
    the differential operators gives rise to an extension of the
    algebra of infinitesimal gauge transformations. For two such
    $V_1,V_2\in \gau(P)= \Gamma^\infty(P,VP)^G\subseteq
    \Diffopverg{\mathrm{1}}(P)$ the deformed action on $f\in
    C^\infty(P)[[\lambda]]$ gives
    \begin{equation}
        \label{eq:deformed-action-trafo}
        V_1 \bulletp (V_2\bulletp f)- V_2\bulletp (V_1\bulletp f)=
        (V_1\starp V_2- V_2\starp V_1) \bulletp f= \Lie_{[V_1,V_2]} f
        + O(\lambda). 
    \end{equation}
    In general, the higher order contributions are nontrivial. If the
    former algebraic relation \eqref{eq:inf-gauge-trafo-commutator}
    shall still be valid in the deformed situation the Lie algebra of
    infinitesimal gauge transformations $\gau(P)$ has to be extended
    to the deformed algebra $(\Diffopver(P)[[\lambda]],\starp)$ of
    vertical differential operators of which $\gau(P)$ is no longer a
    Lie subalgebra.
\end{remark}

\chapter{Deformation quantization of associated vector bundles}
\label{cha:associated-vector-bundles}

The importance of principal fibre bundles in differential geometry
comes not least from the fact that every vector bundle can be seen as
an associated vector bundle with respect to some principal fibre
bundle. In this chapter we investigate how the deformations discussed
in Chapter~\ref{cha:deformation-on-bundles} behave under this process
of association and how they induce corresponding deformations on the
level of vector bundles. As a main result, it will turn out that every
deformation quantization of a principal fibre bundle naturally induces
a corresponding deformation quantization of any associated vector
bundle in the well-known sense that is recalled and discussed in the
first section. This observation makes use of the isomorphism between
the invariant vector-valued functions on the principal fibre bundle
and the sections of the associated bundle. Even more general, we can
show that any invariant deformed module structure of the horizontal
forms of the principal fibre bundle induces a corresponding
deformation of the forms on the base manifold with values in the
associated bundle.

\section{Deformation quantization of vector bundles}
\label{sec:quantization-vector-bundles}

It is a well-known fact that the sections $\Gamma^\infty(M,E)$ of a
vector bundle $\pe:E\longrightarrow M$ have the structure of a
finitely generated projective module over the smooth functions
$C^\infty(M)$ of the base manifold $M$. Moreover, the vector bundle is
already determined by the sections, \cite{swan:1962a}, and the theorem
of Serre and Swan \cite[Thm.~2.10]{gracia.bondia:2000a} states that
the category of vector bundles over a manifold $M$ is equivalent to
the category of finitely generated projective modules over
$C^\infty(M)$. This one-to-one correspondence naturally implies what
shall be understood by a deformation quantization of a vector bundle
$\pe: E\longrightarrow M$ with respect to some star product $\star$ on
the Poisson manifold $M$. It is nothing but a deformed module
structure
\begin{equation}
    \label{eq:quantization-vector-bundle}
    (\Gamma^\infty(M,E)[[\lambda]],\bullet)_{(C^\infty(M)[[\lambda]],\star)}
\end{equation}
in the sense of
Definition~\ref{definition:deformation-right-module}. The explicit
definition and the corresponding well-known results can be found in
\cite{bursztyn.waldmann:2000b} and \cite{waldmann:2005a}. In our case
where $\star$ is a differential star product one of course demands the
same for the deformation $\bullet$. In the form $s \bullet a=s\cdot a+
\sum_{r=1}^\infty \rho_r(s,a)$, this means that the $\rho_r$ are
bidifferential operators
\begin{equation}
    \label{eq:differential-for-vector-bundle}
    \rho_r\in \Diffop^{(L_r,l_r)}(C^\infty(M),\Gamma^\infty(M,E);
    \Gamma^\infty(M,E)) \cong \Diffop^{L_r}(M,
    \Diffop^{l_r} (\Gamma^\infty(M,E))) 
\end{equation}
for some $L_r,l_r\in \mathbb{N}_0$ and all $r\ge 1$. The stated
isomorphism in \eqref{eq:differential-for-vector-bundle} is verified
by a simple induction and shows that in the considered case the
meaning of a differential deformation is always clear, confer
Remark~\ref{remark:notion-differential-right-module}. The classical
commutant of the considered right module structure
$\Gamma^\infty(M,E)$ is given by the sections of the endomorphism
bundle $\Gamma^\infty(M, \End(E))\cong \Diffop^0(\Gamma^\infty(M,E))$.
Thus, one has the bimodule structure
\begin{equation}
    \label{eq:undeformed-bimodule-vector-bundle}
    {}_{(\Gamma^\infty(M,\End(E)),\circ)}
    \Gamma^\infty(M,E)_{(C^\infty(M),\cdot)}. 
\end{equation}
With the results of Section~\ref{sec:projective-modules} and
Proposition~\ref{proposition:projective-differential-modules} we can
apply all considerations of
Chapter~\ref{cha:deformation-algebras-modules} and find a simple proof
for the main assertion of the following theorem, see
\cite{bursztyn.waldmann:2000b}, \cite[Thm.~1]{waldmann:2005a}.

\begin{theorem}[Deformation quantization of vector bundles] 
    \label{theorem:quantization-vector-bundles}
    Every vector bundle $\pe: E\longrightarrow M$ over a Poisson
    manifold with a star product $\star$ has a corresponding
    deformation quantization which is unique up to equivalence.
    Moreover, there exists a deformed bimodule structure
    \begin{equation}
        \label{eq:deformed-bimodule-vector-bundle}
        {}_{(\Gamma^\infty(M,\End(E))[[\lambda]], \star'_E)}
        (\bullet'_E,\Gamma^\infty(M,E)[[\lambda]],\bullet) 
        _{(C^\infty(M)[[\lambda]],\star)},
    \end{equation}
    and for a fixed star product $\star$ the deformed bimodule
    structures $\bullet'_E, \bullet$ and the algebra structure
    $\star'_E$ are unique up to equivalence.

    In addition, the two deformed algebras
    $(\Gamma^\infty(M,\End(E))[[\lambda]],\star'_E)$ and
    $(C^\infty(M)[[\lambda]],\star)$ are mutual commutants. In other
    words, the commutant
    $\End_{C^\infty(M)[[\lambda]]}(\Gamma^\infty(M,E)[[\lambda]],\bullet)$
    of the right module structure is isomorphic to
    $(\Gamma^\infty(M,\End(E))[[\lambda]]$ as
    $\mathbb{C}[[\lambda]]$-module and as algebra via $\star'_E$.
\end{theorem}

Note that only the last assertion does not follow from the general
considerations of Chapter~\ref{cha:deformation-algebras-modules}.

\section{Deformation quantization of associated vector bundles}
\label{sec:deformation-associated-bundle}

Now we investigate the above situation for associated vector bundles
$E= P\times_G V$ as explained in
Appendix~\ref{cha:bundle-geometry}. The following considerations
basically follow \cite{bordemann.neumaier.waldmann.weiss:2007a:pre}
but the more conceptual proceeding shows that some results can be
obtained in a slightly simpler way.

In addition to the results of
Proposition~\ref{proposition:associated-vector-bundles} we make the
following simple observation.

\begin{lemma}[The endomorphism bundle of an associated vector bundle]
       \label{lemma:endomorphisms-associated-bundle}
       Let $\pp:P\longrightarrow M$ be a principal fibre bundle,
       $\rep: G\longrightarrow \Aut(V)$ be a representation of the
       structure group $G$ on a finite dimensional vector space $V$
       and let $P\times_G V$ denote the associated vector bundle. Then
       the induced representation of $G$ on $\End(V)$ yields a vector
       bundle isomorphism
        \begin{equation}
            \label{eq:isomorphism-associated-endomorphism-bundle}
            P\times_G \End(V) \cong \End(P\times_G V)
        \end{equation}
        over the identity $\id_M$ for the endomorphism bundle of the
        associated bundle, given by
        \begin{equation}
            \label{eq:isomorphism-associated-endomorphism-bundle-explicit}
            [u,L]([u,v])= [u,L(v)]
        \end{equation}
        for equivalence classes $[u,L] \in P\times_G \End(V)$ and
        $[u,v]\in P\times_G V$. Thus, one has the isomorphism
        \begin{equation}
            \label{eq:isomorphism-sections-associated-endomorphism-bundle}
            (C^\infty(P)\otimes
            \End(V))^G \cong \Gamma^\infty(M,\End(P\times_G V)).
        \end{equation}
        With the isomorphism $(C^\infty(P)\otimes V)^G\cong
        \Gamma^\infty(M,P\times_G V)$ the application of a section of
        the endomorphism bundle to a section of the bundle reads
        \begin{equation}
            \label{eq:section-endomorphism-bundle-on-section}
            (f \otimes L)(h\otimes v)= fh \otimes L(v)
        \end{equation}
        for all $f\otimes L\in (C^\infty(P)\otimes
        \End(V))^G$ and $h\otimes v \in (C^\infty(P)\otimes V)^G$.
    \end{lemma}

\begin{proof}
    The isomorphism
    \eqref{eq:isomorphism-associated-endomorphism-bundle} is obvious
    with the given identification
    \eqref{eq:isomorphism-associated-endomorphism-bundle-explicit}. Then,
    the remaining statements follow easily using the corresponding
    isomorphisms of invariant vector-valued functions on $P$ and
    sections on the associated bundle, confer
    \eqref{eq:isomorphism-sections-associated-explicit}.
\end{proof}

In the following it will be shown that a deformation quantization of
the principal fibre bundle $P$ as considered in
Chapter~\ref{cha:deformation-on-bundles} leads to the bimodule
structure~\eqref{eq:deformed-bimodule-vector-bundle} for any
associated vector bundle $P\times_G V$, confer
\cite[Lemma~6.1]{bordemann.neumaier.waldmann.weiss:2007a:pre}.

\begin{lemma}[Deformed associated bimodules]
    \label{lemma:extended-bimodule}
    The $G$-invariant bimodule structure 
    \begin{equation}
        \label{eq:bimodule-principal-again}
        {}_{( \Diffopver (P)
          [[\lambda]],\starp)} (\bulletp,C^\infty(P)
        [[\lambda]],\bullet)_{(C^\infty(M)[[\lambda]],\star)} 
    \end{equation}
    as in Theorem~\ref{theorem:bimodule-equivariant-bundle} yields a
    bimodule structure
    \begin{equation}
        \label{eq:extended-bimodule-invariance}
        {}_{((\Diffopver (P)\otimes \End(V))^G [[\lambda]],\starp)}
        (\bulletp, (C^\infty(P)\otimes V)^G[[\lambda]],
        \bullet )_ 
        {(C^\infty(M) [[\lambda]],\star)}.
    \end{equation}
\end{lemma}

\begin{proof}
    First of all it is clear that the structures $\starp,\bulletp$,
    and $\bullet$ naturally induce corresponding operations with
    respect to the tensor products $(\Diffopver
    (P)\otimes \End(V))[[\lambda]]$ and $(C^\infty(P)\otimes
    V)[[\lambda]]$, yielding a bimodule structure
    \begin{equation}
        \label{eq:extended-bimodule}
        {}_{((\Diffopver (P)\otimes \End(V)) [[\lambda]],\starp)}
        (\bulletp, (C^\infty(P)\otimes V)[[\lambda]],
        \bullet )_ 
        {(C^\infty(M) [[\lambda]],\star)}.
    \end{equation}
    For factorising elements $A\otimes L, B\otimes K\in \Diffopver
    (P)\otimes \End(V)$, $f\otimes v\in C^\infty(P)\otimes V$ and
    $a\in C^\infty(M)$ the definitions are $(A\otimes L) \starp
    (B\otimes K) = (A\starp B) \otimes (L\circ K)$, $(A\otimes L)
    \bulletp (f \otimes v)= (A \bulletp f) \otimes L(v)$, and $(f
    \otimes v)\bullet a= (f\bullet a) \otimes v$ with some abuse of
    notation.
\end{proof}

Considering only the right module structure, the isomorphism
$(C^\infty(P)\otimes V)^G\cong \Gamma^\infty(M, P\times_G V)$
immediately implies the following theorem, see
\cite[Thm.~6.2]{bordemann.neumaier.waldmann.weiss:2007a:pre}.

\begin{theorem}[Deformation quantization of associated vector bundles]
    \label{theorem:quantization-associated-vector-bundle}
    By \eqref{eq:extended-bimodule-invariance} every deformation
    quantization $\bullet$ of a principal fibre bundle induces a
    corresponding deformation quantization of any associated vector
    bundle.
\end{theorem}

In order to find the relation between the left module structures in
\eqref{eq:extended-bimodule-invariance} and
\eqref{eq:deformed-bimodule-vector-bundle} we observe that every
$A\otimes L\in (\Diffopver(P)\otimes \End(V))^G$ can be seen as an
element $A\otimes L \in \Gamma^\infty(M,\End(P\times_G V))$ by setting 
\begin{equation}
    \label{eq:diffop-as-endo}
    (A \otimes L)(h\otimes v)= A(h)\otimes L(v)
\end{equation}
for all $h\otimes v\in (C^\infty(P)\otimes V)^G$. With
\eqref{eq:isomorphism-sections-associated-endomorphism-bundle} it
follows that this identification yields a surjective map
\begin{equation}
    \label{eq:surjection-endomorphisms}
    (\Diffopver(P) \otimes \End(V))^G \longrightarrow
    \Gamma^\infty(M,\End(P\times_G V))
\end{equation}
and that the following diagram commutes,
\begin{equation}
    \label{eq:diagram-diffops-endos}
    \xymatrix @R=+1.7pc @C=-4pc {
      & (\Diffop_{\mathrm{ver}}(P)\otimes \End(V))^G \ar@{->>}[dr]& \\
      (C^\infty(P)\otimes \End(V))^G \ar@{<->}[rr]^-{\cong}
      \ar@{^{(}->}[ur] & & \Gamma^\infty(M,\End(P\times_G V)).
    }
\end{equation}
That every element $A\otimes L\in
(\Diffop_{\mathrm{ver}}(P)\otimes \End(V))^G$ by
\eqref{eq:diffop-as-endo} acts in the same way as an element $f\otimes
L'\in C^\infty(P)\otimes \End(V)$ by
\eqref{eq:section-endomorphism-bundle-on-section} is plausible due to
the following consideration. The arguments $h\otimes v\in
(C^\infty(P)\otimes V)^G$ satisfy $h(u.g)v=h(u)\rep_{g^{-1}} v$ for
all $u\in P$ and $g\in G$. Thus, every vertical differentiation of $h$
is nothing but a linear transformation of $v$. Explicitly, one has
$\ddto h(u.\exp(-t\xi))v= h(u) \rep'_\xi v$ for each $\xi\in \lieg$ in
the Lie algebra and the induced Lie algebra representation
$\rep':\lieg\longrightarrow \End(V)$.

The bimodule \eqref{eq:extended-bimodule-invariance} shows that for
every $D\in (\Diffopver(P)\otimes \End(V))^G[[\lambda]]$ there exists
a unique element $\phi(D)\in
\End_{C^\infty(M)[[\lambda]]}(\Gamma^\infty(M,P\times _G
V)[[\lambda]], \bullet)$ in the commutant of the right module
structure with $D \bulletp s = \phi(D)s$ for all $s\in
\Gamma^\infty(M,P\times _G V)[[\lambda]]$. The map $ \phi:
(\Diffopver(P) \otimes \End(V))^G[[\lambda]] \longrightarrow
\End_{C^\infty(M)[[\lambda]]}(\Gamma^\infty(M,P\times _G V) [[\lambda]],
\bullet)$ is surjective since this is already true for the restriction
to $ (C^\infty(P) \otimes \End(V))^G[[\lambda]]$. This assertion is a
slightly stronger version of
\cite[Lemma~6.3]{bordemann.neumaier.waldmann.weiss:2007a:pre}.

\begin{lemma} 
    \label{lemma:surjection-on-commutant}
    The map
    \begin{equation}
        \label{eq:surjection-commutant}
        \phi: (C^\infty(P) \otimes \End(V))^G[[\lambda]]
        \longrightarrow 
        \End_{C^\infty(M)[[\lambda]]}(\Gamma^\infty(M,P\times
        _G V) [[\lambda]], \bullet)
    \end{equation}
    is surjective.
\end{lemma}

\begin{proof}
    Let $K=\sum_{r=0}^\infty \lambda^r K_r \in
    \End_{C^\infty(M)[[\lambda]]} (\Gamma^\infty(M,P\times _G V)
    [[\lambda]], \bullet)$ be an element in the commutant. By
    definition it is clear that $K_0\in \Gamma^\infty(M,\End(P\times_G
    V))\subseteq (\Diffopver(P) \otimes \End(V))^G$. With $\phi(K_0)=
    \sum_{r=0}^\infty\lambda^r \phi(K_0)_r$ it follows by the
    definition of $\phi$ that $K_0=\phi(K_0)_0$. Thus, the element $K-
    \phi(K_0)= \sum_{r=1}^\infty \lambda^r K_r^{(1)}$ of the commutant
    begins in order $\lambda$. Due to the
    $\mathbb{C}[[\lambda]]$-linearity of $\phi$, iteration proves the
    lemma.
\end{proof}

As cited in Theorem~\ref{theorem:quantization-vector-bundles} there is
an isomorphism $\psi: \End_{C^\infty(M)[[\lambda]]}
(\Gamma^\infty(M,P\times_G V) [[\lambda]], \bullet) \longrightarrow
\Gamma^\infty(M,\End(P\times_G V))[[\lambda]]$ such that for each
$D\in (C^\infty(P)\otimes \End(V))^G[[\lambda]]$ one finally has
\begin{equation}
    \label{eq:actions-commutants}
    D\bulletp s= \phi(D)s= \psi(\phi(D)) \mathbin{\bulletp_E} s
\end{equation}
for all $s\in \Gamma^\infty(M,P\times_GV)[[\lambda]]$. The left module
properties in \eqref{eq:deformed-bimodule-vector-bundle} and
\eqref{eq:extended-bimodule-invariance} then imply the following
result for the map $\gamma=\psi\circ \phi$, confer
\cite[Thm.~6.4]{bordemann.neumaier.waldmann.weiss:2007a:pre}.

\begin{theorem}[The associated deformed commutant] 
    \label{theorem:associated-commutant}
    Let $\bullet$ denote a deformation quantization of a principal
    fibre bundle as well as the induced deformation quantization of an
    associated vector bundle $E=P\times_G V$. Then, for all structures
    $\star'_E$ and $\bullet'_E$ as in
    \eqref{eq:deformed-bimodule-vector-bundle} there exists a
    surjective algebra homomorphism
    \begin{equation}
        \label{eq:surjection-endos}
        \gamma: ((\Diffopver(P) \otimes \End(V))^G[[\lambda]],\starp) 
        \longrightarrow  
        (\Gamma^\infty(M,\End(E)) [[\lambda]], \star'_E)
    \end{equation}
    such that
    \begin{equation}
        \label{eq:property-surjection-endos}
        D \bulletp s = \phi(D)\bullet'_E s
    \end{equation}
    for all $D\in (\Diffopver(P) \otimes
    \End(V))^G[[\lambda]]$ and $s\in \Gamma^\infty(M,E)[[\lambda]]$.
\end{theorem}

\begin{remark}
    A further investigation that can be found in
    \cite{bordemann.neumaier.waldmann.weiss:2007a:pre} shows that the
    extension to the vertical differential operators is necessary. In
    general, $(C^\infty(P) \otimes \End(V))^G$ is \emph{not} deformed
    into a \emph{subalgebra} of $((\Diffopver(P)
    \otimes \End(V))^G[[\lambda]], \star')$.
\end{remark}

In analogy to the above discussions for functions one gets the following
statement for differential forms when using the isomorphism
\eqref{eq:isomorphisms-forms-associated}. 

\begin{proposition}[Deformation quantization of associated forms]
    \label{proposition:associated-forms}
    Every $G$-invariant deformed right module structure
    \begin{equation}
        \label{eq:deformed-bimodule-forms}
        % {}_{(\Diffopver(\Gamma^\infty_\hor
        %   (\Gamma^\infty(P,\Dach{k}T^*P)))[[\lambda]], \starp)}
        (
        % \bulletp,
        \Gamma^\infty_\hor(P,\Dach{k} T^*P)[[\lambda]], \bullet)_{(C^\infty(M)[[\lambda]],\star)}
    \end{equation}
    of the horizontal forms of a principal fibre bundle induces a
    deformation
    \begin{equation}
        \label{eq:deformed-bimodule-associated-forms}
        (\Gamma^\infty(M, \Dach{k}T^*M\otimes
        (P\times_G V)), \bullet)_{(C^\infty(M)[[\lambda]],\star)}
    \end{equation}
    for the differential forms on the base manifold with values in
    any associated vector bundle.
\end{proposition}

\begin{proof}
    The proof is a straightforward generalization of
    Lemma~\ref{lemma:extended-bimodule} where one makes use of the
    isomorphism $(\Gamma^\infty_\hor(P,\Dach{k}T^*P) \otimes V)^G
    \cong \Gamma^\infty(M, \Dach{k}T^*M\otimes (P\times_G V))$.
\end{proof}

% Local Variables: 
% TeX-master: "Dissertation"
% End: 

%% file: anhang.tex
\chapter{Geometry of principal fibre bundles}
\label{cha:bundle-geometry}

This appendix gives a short overview over some concepts of the
differential geometry of Lie groups and principal fibre bundles used
in the main text. In order to concentrate on the relevant facts it is
assumed that the reader is already familiar with the fundamental
concepts in differential geometry.  The presented notions can be found
in the references
\cite{duistermaat.kolk:2000a,michor:2008,kobayashi.nomizu:1963a,waldmann:2007a}.
With respect to the aspired physical applications in gauge theories
all necessary aspects of differential geometry can also be found in
\cite{bleecker:2005,naber:1997,naber:2000}.

\section{Lie groups and group actions}
\label{sec:group-actions}

In this first section we collect some basic facts on Lie groups and
their actions and representations. All the details omitted in this
overview can be found in \cite{duistermaat.kolk:2000a} or in any other
book on differential geometry treating Lie groups. By definition, a
Lie group $G$ is a manifold such that the group multiplication and
taking the inverse are smooth maps. The Lie algebra of the group is
given by the tangent space $\lieg= T_eG$ at the neutral element $e\in
G$ equipped with the following Lie bracket. Let $l_gh=gh$ for $g,h\in
G$ denote the left multiplication. Since for any $\xi\in \lieg$ there
exists a unique \emph{left-invariant vector field} $X^\xi=l_g^*X^\xi$
with $X^\xi(e)=\xi$, explicitly given by $X^\xi(g)=T_el_g(\xi)$, the
Lie bracket of vector fields gives rise to the definition $[\xi,\eta]=
[X^\xi,X^\eta](e) \in \lieg$ for $\xi,\eta\in \lieg$. Any $X^\xi$ has
a complete flow $\Fl^\xi_t$ with $\Fl^\xi_t\circ l_g= l_g\circ
\Fl^\xi_t$. The \emph{exponential map} $\exp: \lieg \longrightarrow G$
is defined by $\exp(\xi)=\Fl^\xi_1(e)$. Then the map
$\mathbb{R}\longrightarrow G$ with $t\longmapsto \exp(t\xi) =
\Fl^{t\xi}_1(e) = \Fl^\xi_t(e)$ is a \emph{one-parameter group}, thus
having the important properties $\exp((t+s)\xi)=\exp(t\xi)\exp(s\xi)$
and $\exp(0)=e$. Further, one has $\ddto \exp(t\xi)=\xi$.

If $M$ is a manifold a \emph{left action of $G$ on $M$} is given by a
smooth map $\phi: G\times M \longrightarrow M$, often denoted by
$\phi(g,p)=\phi_g(p)=\phi^p(g)$ with smooth maps $\phi_g:
M\longrightarrow M$ and $\phi^p: G\longrightarrow M$, such that
$\phi_e=\id_M$ and $\phi_g\circ \phi_h= \phi_{gh}$. For a \emph{right
  action} one clearly has $\phi_g\circ \phi_h= \phi_{hg}$. An action
is called \emph{free} if $\phi_g(p)=p$ implies that $g=e$. Further, an
action is called \emph{proper} if the map $G\times M\longrightarrow
M\times M$ with $(g,p)\longmapsto (\phi_g(p),p)$ is proper, this means
if any compact subset has a compact preimage. It should be noted that
properness has many different but equivalent formulations, confer
\cite[Chap.~9]{lee:2003a}. In order to describe the infinitesimal
version of an action one considers the \emph{fundamental vector
  fields} $\xi_M\in \Gamma^\infty(M,TM)$ with respect to elements
$\xi\in \lieg$ which are defined by
\begin{equation}
    \label{eq:fundamental-vector-field}
    \xi_M(p)=\ddto \phi_{\exp(t\xi)}p = T_e\phi^p \xi
\end{equation}
for all $p\in M$. The map $\phi': \lieg \longrightarrow
\Gamma^\infty(M,TM)$ with $\xi\longmapsto \xi_M$ then satisfies
$[\xi_M,\eta_M]= -[\xi,\eta]_M$ for left actions and $[\xi_M,\eta_M]=
[\xi,\eta]_M$ for right actions.

A \emph{representation} of a Lie group $G$ on a finite dimensional
vector space $V$ is a left action $\rep: G\times V \longrightarrow V$
such that each $\rep_g: V\longrightarrow V$ is a linear map.
Obviously, a representation can be seen as a group homomorphism
\begin{equation}
    \label{eq:Lie-group-representation}
    \rep: G\longrightarrow \Aut(V)
\end{equation}
from $G$ to the group of automorphisms $\Aut(V)$ of $V$ with the
composition of maps as group structure. With the identification $T_v
V=V$ the fundamental vector fields can be seen as linear maps
$\xi_V:V\longrightarrow V$. Then, the infinitesimal version $\rep':
\xi \longmapsto \xi_V$ of the representation is a \emph{Lie algebra
  representation}, this means a \emph{Lie algebra homomorphism}
\begin{equation}
    \label{eq:Lie-algebra-representation}
    \rep': \lieg \longrightarrow \End(V),
\end{equation}
where the Lie bracket $[\cdot,\cdot]_{\End(V)}$ for the endomorphisms
is the usual commutator. The crucial point is that $\rep'([\xi,\eta])=
[\xi,\eta]_V= -[\xi_V,\eta_V]_{\Gamma^\infty(V,TV)}= [\xi_V,
\eta_V]_{\End(V)}$. With $\End(V)=T_{\id_V}\Aut(V)$ the induced
infinitesimal representation of \eqref{eq:Lie-group-representation}
can be seen as tangent map $\rep'= T_e \rep$.

Important examples for Lie group actions and representations are the
following. Any Lie group $G$ acts on itself from the left by the
\emph{conjugation}
\begin{equation}
    \label{eq:conjugation}
    \Conj: G\times G \longrightarrow G, \quad \quad \Conj_g h=
    ghg^{-1}. 
\end{equation}
This induces the \emph{adjoint representation $\Ad$ of $G$ on its Lie
  algebra $\lieg$}
\begin{equation}
    \label{eq:adjoint-group-representation}
    \Ad: G\longrightarrow \End(\lieg), \quad \quad \Ad_g=T_e\Conj_g,
\end{equation}
and the \emph{adjoint representation $\ad$ of the Lie algebra $\lieg$
  on itself} as infinitesimal version which turns out to be
\begin{equation}
    \label{eq:adjoint-lie-algebra-representation}
    \ad= \Ad': \lieg \longrightarrow \End(\lieg), \quad \quad \ad_\xi
    \eta= [\xi,\eta]. 
\end{equation}
The adjoint representation $\Ad$ acts by Lie algebra automorphisms
which explicitly means that $\Ad_g[\xi,\eta]= [\Ad_g \xi,\Ad_g \eta]$.

\section{Principal fibre bundles and surjective submersions}
\label{sec:principal-bundles}

The fundamental geometric structure discussed in this work is the one
of a principal fibre bundle. This is nothing but a particular
\emph{$G$-bundle}, confer \cite{michor:2008}, for a Lie group $G$. The
explicit definition is the following.

\begin{definition}[Principal fibre bundle]
    \label{definition:principal-bundle}
    A \emph{principal fibre bundle} $(\pp,P,M,G)$ consists of a smooth
    mapping $\pp:P\longrightarrow M$ called \emph{projection} between
    two manifolds where $P$ is called \emph{total space} and $M$ is
    called \emph{base space}, together with a Lie group $G$ referred
    to as \emph{structure group} such that the following assertions
    hold.
    \begin{enumerate}
    \item For each $p\in M$ there exists an open neighborhood
        $U\subseteq M$ and a diffeomorphism $\varphi: P|_U=
        \pp^{-1}(U)\longrightarrow U\times G$ such that the diagram
        \begin{equation}
            \label{eq:diagram-principal-bundle-chart}
            \xymatrix{
              P|_U \ar[r]^{\varphi} \ar[d]_{\pp} & U\times G
              \ar[dl]^{\pr_1} \\ 
              U&
            }
        \end{equation}
        commutes where $\pr_1$ is the projection onto the first
        component.\\
        Such a pair $(U,\varphi)$ is referred to as a \emph{principal
          bundle chart}.
    \item There exists a family $\{U_i,\varphi_i\}_{i\in I}$ of
        principal bundle charts with some index set $I$ such
        that $\{U_i\}_{i\in I}$ is an open cover of $M$ and for any
        point $p\in M$ with $p\in U_{ij}= U_i\cap U_j\neq \emptyset$
        for some $i,j\in I$ and $g\in G$ one has
        \begin{equation}
            \label{eq:transition-principal-bundle}
            (\varphi_i\circ \varphi_j^{-1}) (p,g)= (p,
            \varphi_{ij}(p)g), 
        \end{equation}
        where the \emph{transition function} $\varphi_{ij}:
        U_{ij}\longrightarrow G$ is a smooth map and
        $\varphi_{ij}(p)g$ denotes the group multiplication in $G$.
        Such a family $\{U_i,\varphi_i\}_{i\in I}$ is referred to as
        \emph{principal bundle atlas}.
    \end{enumerate}
    If the structures are clear one simply refers to the total space
    $P$ as corresponding principal fibre bundle.
\end{definition}

The definition has the following immediate consequences.

\begin{proposition}[Cocycle conditions and right action]
    \label{proposition:cocycle-right-action}
    \begin{enumerate}
    \item The transition functions of a principal fibre bundle satisfy
        the so-called \emph{cocycle conditions}
        \begin{equation}
            \label{eq:cocycle-conditions}
            \varphi_{ij}(p).\varphi_{jk}(p)= \varphi_{ik}(p) \quad
            \textrm{and} \quad \varphi_{ii}(p)=e.
        \end{equation}
    \item Every principal fibre bundle is equipped with a smooth right
        action $\rr: P\times G \longrightarrow P$ of $G$ on the fibres
        of $P$ which is free and proper. This right action is
        well-defined by the local definition in a principal bundle
        chart $(U,\varphi)$ and reads
        \begin{equation}
            \label{eq:definition-principal-action}
            \rr(\varphi^{-1}(p,h),g)= \varphi^{-1}(p,hg)
        \end{equation}
        for all $p\in U$ and $h,g\in G$. Sometimes, we will use the
        abbreviated notation $\rr(u,g)=u.g$.
    \end{enumerate}
\end{proposition}

The principal right action is the fundamental structure of any
principal fibre bundle. With the following well-known notion of a
submersion one can in fact give an alternative definition of a
principal bundle based on the right action.

\begin{definition}[Submersion]
    \label{definition:submersion}
    A map $\pp: P\longrightarrow M$ is said to be a \emph{submersion
      at $u\in P$} if the tangent map $T_u\pp: T_uP\longrightarrow
    T_{\pp(u)}M$ at the point $u$ is surjective. The map $\pp$ is
    called a \emph{submersion} if it is a submersion at each point
    $u\in P$.
\end{definition}

Given a free and proper right action $\rr: P\times G \longrightarrow
P$ of a Lie group on a manifold $P$ the quotient $M=P/G$ is again a
smooth manifold and the corresponding projection $\pp:P\longrightarrow
M$ is a surjective submersion, confer
\cite[Thm.~1.11.4]{duistermaat.kolk:2000a}. According to
\cite[Lemma~18.3]{michor:2008}, $P$ is a principal fibre bundle with
structure Lie group $G$ and right action $\rr$. Thus one can identify
principal fibre bundles with free and proper actions on smooth
manifolds. Since it is the case for any surjective submersion any
principal bundle admits local sections, this means smooth maps
$\sigma: U\longrightarrow P$ with $\pp\circ \sigma= \id_U$ for
sufficiently small and open subsets $U\subseteq M$.

\section{Connections}
\label{sec:connections}

Before we present the notion of a connection on a fibre bundle
together with some of its relevant issues which all can be found in
\cite[Chap.~IV]{michor:2008}, we clarify the general notion of
\emph{vector-valued differential forms}. For a vector bundle $\pe:
E\longrightarrow M$ the $E$-valued $k$-forms, $k=0,\dots, \dim E$, are
given by
\begin{equation}
    \label{eq:vector-valued-forms}
    \Gamma^\infty(M,\Dach{k} T^*M\otimes E),
\end{equation}
where one makes use of the usual tensor product of vector bundles.
For a vector space $V$, the $V$-valued $k$-forms on $M$ are defined by
\begin{equation}
    \label{eq:vector-space-valued-forms}
    \Gamma^\infty(M,\Dach{k} T^*M\otimes (M\times V))=
    \Gamma^\infty(M,\Dach{k} T^*M)\otimes V,
\end{equation}
where $M\times V$ is the trivial vector bundle.

In the most general case a \emph{connection} on a manifold $P$ is
given by a vector-valued one-form $\mathcal{P}\in
\Gamma^\infty(P,T^*P\otimes TP)\cong \Gamma^\infty(P,\End(TP))$ which
is a fibre projection $\mathcal{P}\circ \mathcal{P}=\mathcal{P}$ when
seen as a map $\mathcal{P}: TP \longrightarrow TP$. The image $\image
\mathcal{P}\subseteq TP$ is called \emph{vertical space} and $\ker
\mathcal{P}\subseteq TP$ is referred to as \emph{horizontal space}.

\subsection{Connections on surjective submersions}
\label{subsec:connection-submersion}

The tangent bundle of the total space $P$ of a surjective submersion
$\pp:P\longrightarrow M$ is canonically equipped with a subbundle
$VP$, the so-called \emph{vertical bundle} consisting of all
\emph{vertical vectors} $V\in TP$ with $T\pp V=0$. In other words,
\begin{equation}
    \label{eq:vertical-bundle}
    VP= \ker T\pp \subseteq TP
\end{equation}
is the kernel of the tangent map $T\pp$ of the projection $\pp$.

In this case, a \emph{connection} on $P$ is a projection $\mathcal{P}$
onto the vertical bundle $\image \mathcal{P}= VP$. The kernel defines
the so-called \emph{horizontal bundle} $HP= \ker \mathcal{P}$ which is
a complementary subbundle $HP$ to $VP$ such that
\begin{equation}
    \label{eq:connection}
    TP = VP \oplus HP.
\end{equation}
It is remarkable that the connection is already determined by the
choice of a horizontal bundle $HP$. The \emph{curvature} $\mathcal{R}$
of the connection $\mathcal{P}$ is the vector-valued two-form
$\mathcal{R}\in \Gamma^\infty(P,\Dach{2}T^*P\otimes VP)$ defined by
\begin{equation}
    \label{eq:curvature}
    \mathcal{R}(Z,W)= \mathcal{P}[(\id_{TP}-\mathcal{P})Z,
    (\id_{TP}-\mathcal{P})W].
\end{equation}

An important structure in the context of connections of surjective
submersions is the horizontal lift of vector fields.

\begin{definition}[Horizontal lift]
    Let $\pp:P\longrightarrow M$ be a surjective submersion with a
    connection given by a horizontal bundle $HP$. Then the
    \emph{horizontal lift} of a vector field $X\in
    \Gamma^\infty(M,TM)$ is the uniquely defined horizontal vector
    field $X^\hlift\in \Gamma^\infty(P,HP)$ which is $\pp$-related to
    $X$, this means
    \begin{equation}
        \label{eq:horizontal-lift}
        T\pp \circ X^{\hlift}= X \circ \pp.
    \end{equation}
\end{definition}

% Given a smooth map $\phi: M\longrightarrow M'$ between two manifolds
% and vector fields $X,Y\in\Gamma^\infty(TM)$ and $X',Y'\in
% \Gamma^\infty(TM')$ with $T\phi\circ X=X'\circ \phi$ and $T\phi\circ
% Y=Y'\circ \phi$ it is a well-known fact concerning related vector
% fields that $T\phi\circ [X,Y]= [X',Y']\circ \phi$. For the above
% situation this has some useful consequences, which are summarized
% together with other useful properties of horizontal lifts and vertical
% vector fields.
It is helpful for many applications in the main text to summarize some
of the obvious properties of horizontal lifts and vertical vector
fields.

\begin{lemma}[Horizontal lifts and vertical vector fields]
    \label{lemma:horizontal-lifts-vertical-vector-fields}
    Let $\pp:P\longrightarrow M$ be a surjective submersion with a
    connection $TP=VP\oplus HP$. Then for all $X,Y\in
    \Gamma^\infty(M,TM)$, $V,W\in \Gamma^\infty(P,VP)$, $a\in
    C^\infty(M)$ and $\alpha\in \Gamma^\infty(M,T^*M)$ the following
    equations hold.
    \begin{gather}
        \label{eq:consequences-related-fields}
        T\pp \circ [X^{\hlift},Y^{\hlift}]= [X,Y]\circ \pp, \quad
        T\pp\circ [V,X^{\hlift}]=0, \quad
        \quad  T\pp\circ [V,W]=0\\
        \label{eq:horizontal-lift-module-homo}
        (aX)^\hlift = (\pp^*a) X^\hlift,\\
        \label{eq:pullback-forms-lifts-vertical}
        (\pp^*\alpha)(X^\hlift)= \pp^*(\alpha(X)), \quad \quad
        (\pp^*\alpha)(V)=0,\\
        \label{eq:pullback-functions-lifts-vertical}
        X^\hlift (\pp^*a)= \pp^*(X(a)), \quad \quad V(\pp^*a)=0.
    \end{gather}
\end{lemma}

The first equation of \eqref{eq:consequences-related-fields} yields
\begin{equation}
    \label{eq:horizontal-lift-and-curvature}
    [X^\hlift, Y^\hlift] - [X,Y]^\hlift= \mathcal{R}(X^\hlift, Y^\hlift),
\end{equation}
which shows that the curvature of the connection is the obstruction
against the horizontal lift $\hlift:
\Gamma^\infty(M,TM)\longrightarrow \Gamma^\infty(P,HP)$ being a Lie
algebra homomorphism. However,
Equation~\eqref{eq:horizontal-lift-module-homo} points out that it is
always a $C^\infty(M)$-module homomorphism.

Let $\Gamma^\infty_\hor(P,\Dach{\bullet} T^*P)$ denote the
\emph{horizontal forms} which by definition vanish under the insertion
of any vertical vector. Then, for any connection $HP$ with
\emph{horizontal projection} $\Gamma^\infty(P,TP)\ni Z\longmapsto
(\id_{TP}-\mathcal{P})Z= Z^H\in \Gamma^\infty(P,HP)$ one can define
the \emph{covariant exterior derivative}
\begin{equation}
    \label{eq:cov-ext-derivative}
    \D_{HP}: \Gamma^\infty(P,\Dach{\bullet} T^*P) \longrightarrow
    \Gamma^\infty_\hor(P,\Dach{\bullet+1} T^*P)
\end{equation}
by $(\D_{HP}\alpha)(Z_1,\dots, Z_{k+1})= (\D \alpha)(Z_1^H, \dots,
Z_{k+1}^H)$ for all $\alpha\in \Gamma^\infty(P,\Dach{k} T^*P)$,
$k=1,\dots,\dim P$, and $Z_1,\dots, Z_{k+1}\in \Gamma^\infty(P,TP)$.

\subsection{Principal connections}
\label{subsec:connection-principal}

If $\pp:P\longrightarrow M$ is a principal fibre bundle with structure
Lie group $G$ a connection is said to be a \emph{principal connection}
if it is $G$-invariant. For the projection $\mathcal{P}$ this means
$T\rr_g \circ \mathcal{P}= \mathcal{P}\circ T\rr_g$ and the horizontal
spaces then have the property $T\rr_g H_uP= H_{u.g}$ for all $g\in G$
and $u\in P$. In this case it easily follows that the horizontal lifts
are $G$-invariant, this means
\begin{equation}
    \label{eq:invariance-horizontal-lift}
    \rr_g^*X^\hlift= X^\hlift
\end{equation}
for all $g\in G$.

The exponential map $\exp: \lieg \longrightarrow G$ has a natural
generalization $\exp: C^\infty(P,\lieg)^G\longrightarrow
C^\infty(P,G)^G$ from the $G$-invariant $\lieg$-valued functions $\Xi
\in C^\infty(P,\lieg)^G$ with $\Xi \circ \rr_g=\Ad_{g^{-1}} \circ \Xi$
to the $G$-valued functions $F \in C^\infty(P,G)^G$ with $F\circ
\rr_g= \Ad_{g^{-1}} \circ F$ defined by $\exp(\Xi)(u)= \exp(\Xi(u))$
for all $u\in P$. In analogy to the fundamental vector fields
$\xi_P\in \Gamma^\infty(P,VP)$ for the principal right action which
clearly are vertical, one defines vector fields $\Xi_P\in
\Gamma^\infty(P,VP)$ for functions $\Xi\in C^\infty(P,\lieg)$ by
$\Xi_P(u)=(\Xi(u)_P)(u)$ which have the flow
$\Fl^{\xi_P}_t=\rr_{\exp{t\Xi}}$. The behaviour under the principal
action is given by $\rr_g^*\Xi_P= (\Ad_g\circ \Xi \circ
\rr_g)_P$ and $\rr_g^*\xi_P= (\Ad_g\xi)_P$.

The fundamental vector fields $\xi_P$ of the principal right action
give rise to a vector bundle isomorphism $P\times \lieg \cong VP$ by
$(u,\xi)\longmapsto \xi_P(u)=T_e\rr^u \xi$ showing that the vertical
bundle is always trivial. As a consequence, a principal connection
$\mathcal{P}$ induces a $\lieg$-valued one-form $\omega \in
\Gamma^\infty(P,T^*P)\otimes \lieg$ by $\omega(Z)= (T_e\rr^u)^{-1}
\mathcal{P}Z$ with $Z\in TP$. Such a form $\omega$ has particular
properties and already determines the connection by $HP=\ker \omega$
and $\mathcal{P}(Z)=(\omega(Z))_P$ for $Z\in
\Gamma^\infty(P,TP)$. This gives rise to the following definition.

\begin{definition}[Connection one-form]
    \label{definition:connection-one-form}
    A $\lieg$-valued one-form $\omega\in \Gamma^\infty(P,T^*P)\otimes
    \lieg$ is called \emph{connection one-form} if it has the
    following properties.
    \begin{enumerate}
    \item $\omega$ reproduces the generators of the fundamental vector
        fields, 
        \begin{equation}
            \label{eq:connection-form-reproducing-generators}
            \omega(\xi_P)= \xi \quad \textrm{for all } \xi\in \lieg.  
        \end{equation}
    \item $\omega: TP\longrightarrow \lieg$ is $G$-invariant, this
        means
        \begin{equation}
            \label{eq:connection-one-form-invariant}
            \omega \circ T \rr_g = \Ad_{g^{-1}} \circ  \omega. 
        \end{equation}
    \end{enumerate}
    The space of all connection one-forms is denoted by $\mathcal{C}$.
\end{definition}

It can be easily verified that the space $\mathcal{C}$ of connection
one-forms is an affine vector space over
$\Gamma^\infty_\hor(P,T^*P)\otimes \lieg$. The covariant exterior
derivative with respect to $\omega\in \mathcal{C}$ is denoted by
$\D_\omega$. If $\rep: V\longrightarrow \Aut(V)$ is a representation
of the structure Lie group $G$ on a finite dimensional vector space
$V$ it yields a map 
\begin{equation}
    \label{eq:covariant-exterior-derivative}
    \D_\omega: (\Gamma^\infty(P,\Dach{k} T^*P)\otimes V)^G
    \longrightarrow 
    (\Gamma^\infty_\hor(P,\Dach{k+1} T^*P)\otimes V)^G,
\end{equation}
where $(\Gamma^\infty(P,\Dach{k} T^*P)\otimes V)^G$ denotes the
$G$-invariant $V$-valued forms $\alpha$ with respect to the induced
left action $\rr_g^*\otimes \rep_g$, this means $\rr_g^* \alpha=
\rep_{g^{-1}} \circ \alpha$. Note that the connection one-form
$\omega$ itself is a particular invariant $\lieg$-valued form in the
above sense. For $\lieg$-valued forms one can define a Lie bracket
$[\cdot,\cdot]_\wedge$ by bilinear extension of
\begin{equation}
    \label{eq:bracket-g-valued-forms}
    [\alpha \otimes \xi, \beta \otimes \eta]_\wedge= \alpha\wedge
    \beta \otimes [\xi,\eta]
\end{equation}
for all $\alpha,\beta\in \Gamma^\infty(P,\Dach{\bullet}T^*P)$ and
$\xi,\eta\in \lieg$ where $\wedge$ is the anti-symmetric tensor
product of differential forms.

In analogy to the connection one-form, the curvature $\mathcal{R}$
gives rise to a $\lieg$-valued \emph{curvature two-form $\Omega\in
  \Omega^2(P,\lieg)^G_\hor$} which is defined by $\Omega(Z,W)=-(T_e
\rr^u)^{-1}\mathcal{R}(Z,W)$ and satisfies $\mathcal{R}(Z,W)=
-(\Omega(Z,W))_P$. It turns out that the curvature two-form for a
connection one-form $\omega$ is given by the \emph{structure equation}
\begin{equation}
    \label{eq:connection-curvature-forms}
    \Omega= \D_\omega \omega= \D \omega +
    \frac{1}{2}[\omega,\omega]_\wedge. 
\end{equation}
The definition of the covariant exterior derivative and $\D \circ
\D=0$ immediately imply the \emph{Bianchi identity}
\begin{equation}
    \label{eq:Bianchi-identity}
    \D_\omega \Omega= \D \Omega+ [\omega,\Omega]_\wedge= 0.
\end{equation}

\section{Associated vector bundles}
\label{sec:associated-vector-bundles}

Principal fibre bundles have the important property that any
representation of the structure group on a vector space $V$ gives rise
to a corresponding vector bundle with typical fibre $V$. The following
proposition contains some of the most important facts used in this
work. The proofs of the well-known assertions can all be found in
\cite{michor:2008}.

\begin{proposition}[Associated vector bundles]
    \label{proposition:associated-vector-bundles}
    Let $\pp: P\longrightarrow M$ be a principal fibre bundle with
    right action $\rr$ and $\rep: G\longrightarrow \Aut(V)$ be a left
    representation of the structure group $G$ on a finite dimensional
    vector space $V$. Moreover, let $\rR: (P\times V)\times
    G\longrightarrow P\times V$ be the right action with
    $\rR((u,v),g)= (\rr_gu, \rep_{g^{-1}}v)$. Then the following
    assertions hold.
    \begin{enumerate}
    \item The quotient $P\times_G V = (P\times V)/ G$ given by
        the set of all orbits of $\rR$ has a unique structure of a
        smooth manifold such that the projection $P\times V
        \longrightarrow (P\times V)/ G$ is a surjective
        submersion. The equivalence class of an element $(u,v)$ is
        denoted by $[u,v]$.
    \item The map $\pe: P\times_G V \longrightarrow M$ which is
        well-defined by the prescription $\pe([u,v])= \pp(u)$
        defines a vector bundle with typical fibre $V$ which is called
        the \emph{associated vector bundle}.
    \end{enumerate}
    Altogether, one has the following commutative diagram
    \begin{equation}
        \label{eq:diagram-associated-vector-bundle}
        \xymatrix{
          P\times V \ar[r] \ar[d]_{\pr_1} & P\times_G V \ar[d]^{\pe}\\
          P \ar[r]^\pp & M.
        }
    \end{equation}
    \begin{enumerate}
        \setcounter{enumi}{2}
    \item With the choice of a principal connection $\omega$ there is
        a vector space isomorphism
        \begin{equation}
            \label{eq:isomorphisms-forms-associated}
            s: (\Gamma^\infty_\hor(P,\Dach{k}T^*P) \otimes V)^G
            \longrightarrow \Gamma^\infty(M, \Dach{k}T^*M\otimes
            (P\times_G V)),
        \end{equation}
        which is well-defined by the prescription 
        \begin{equation}
            \label{eq:isomorphisms-forms-associated-explicit}
            (s(\alpha))(X_1,\dots, X_k)(p)= [u, \alpha_u(X_1^\hlift(u),
            \dots, X_k^\hlift(u))]
        \end{equation}
        for $\alpha\in (\Gamma^\infty_\hor(P,\Dach{k}T^*P) \otimes
        V)^G$, $X_1,\dots, X_k\in \Gamma^\infty(M,TM)$, $p\in M$, and
        some $u\in P$ with $\pp(u)=p$. As a special case this
        yields the isomorphism
        \begin{equation}
            \label{eq:isomorphism-sections-associated}
            s: (C^\infty(P)\otimes V)^G \longrightarrow
            \Gamma^\infty(M, P\times_G V)
        \end{equation}
        of $C^\infty(M)$-modules between the $G$-invariant $V$-valued
        functions on $P$ and the smooth sections of the associated
        bundle. Explicitly, the formulas are
        \begin{equation}
            \label{eq:isomorphism-sections-associated-explicit}
            s_f(p)= (s(f))(p)= [u,f(u)]
        \end{equation}
        and 
        \begin{equation}
            \label{eq:iso-sections-module-iso}
            a\cdot s_f= s_{\pp^*a \cdot f}
        \end{equation}
        for $f \in (C^\infty(P)\otimes V)^G$ and $a\in C^\infty(M)$.
    \item With the isomorphism
        \eqref{eq:isomorphism-sections-associated} for the sections of
        an associated bundle $E=P\times_G V$ any principal connection
        $\omega$ induces a covariant derivative
        \begin{equation}
            \label{eq:covariant-derivative}
            \nabla^E: \Gamma^\infty(M,TM)\times \Gamma^\infty(M,E)
            \longrightarrow \Gamma^\infty(M,E)
        \end{equation}
        by
        \begin{equation}
            \label{eq:covariant-derivative-explicit}
            \nabla^E_Xs_f= s_{X^\hlift f}
        \end{equation}
        for all $X\in\Gamma^\infty(M,TM)$ and $f\in C^\infty(P,V)^G$.
        The curvature $\mathcal{R}\in
        \Gamma^\infty(P,\Dach{2}T^*P\otimes VP)^G$ of the connection
        induces the curvature $\mathcal{R}^E\in \Gamma^\infty(M,
        \Dach{2}T^*M \otimes \End(E))$ of the covariant derivative and
        one has
        \begin{eqnarray}
            \label{eq:curvature-covariant-derivative}
            \mathcal{R}^E(X,Y)s_f &=&
            s(\mathcal{R}(X^\hlift,Y^\hlift)f)\nonumber \\
            &=& s((\omega([X^\hlift, Y^\hlift]))_Pf) \\
            &=& \left (\nabla^E_X \nabla^E_Y -\nabla^E_Y \nabla^E_X -
                \nabla^E_{[X,Y]} \right) s_f.\nonumber
        \end{eqnarray}
    \end{enumerate}
\end{proposition}

\chapter{Basic concepts of homological algebra}
\label{cha:homological-algebra}

In the following there are presented some basic notions of homological
algebra which are used in the main text. For further details, in
particular for the notion of categories and functors, confer
\cite[Chap.~1,3, and 6]{jacobson:1989a}.

\section{Complexes, cohomology and homotopy}
\label{sec:complexes}

\begin{definition}[Complex]
    \label{definition:chain-complex}
    Let $\mathsf{R}$ be a ring. An \emph{$\mathsf{R}$-complex} $(C,d)$
    is given by two $\mathbb{Z}$-indexed sets $C=\{C_i\}_{i\in
      \mathbb{Z}}$ of $\mathsf{R}$-modules and $d=\{d_i\}_{i\in
      \mathbb{Z}}$ of $\mathsf{R}$-homomorphisms $d_i:
    C_i\longrightarrow C_{i-1}$ such that
    \begin{equation}
        \label{eq:differential-chain-complex}
        d_{i-1} \circ d_i=0
    \end{equation}
    for all $i\in \mathbb{Z}$.

    If $(C,d)$ and $(C',d')$ are complexes a \emph{(chain)
      homomorphism} of $C$ into $C'$ is an indexed set
    $\alpha=\{\alpha_i\}_{i\in \mathbb{Z}}$ of homomorphisms
    $\alpha_i: C_i\longrightarrow C'_i$ such that the diagram
    \begin{equation}
        \label{eq:chain-homomorphism}
        \xymatrix{C_{i-1} \ar[d]_{\alpha_{i-1}} & C_i
          \ar[d]^{\alpha_i} \ar[l]_{d_i}\\
        C'_{i-1} & C'_i  \ar[l]^{d'_i}}
    \end{equation}
    commutes for all $i\in \mathbb{Z}$. Briefly, this is denoted as
    \begin{equation}
        \label{eq:chain-homomorphism-short}
        \alpha \circ d= d' \circ \alpha.    
    \end{equation}
\end{definition}

This definition naturally induces the category
$\mathsf{R}$-$\mathbf{comp}$ of $\mathsf{R}$-complexes which is an
abelian category, confer \cite[Def.~6.1 and 6.2]{jacobson:1989a},
since the homomorphisms between two
complexes have the structure of an abelian group.

Given a complex $(C,d)$ the elements of the submodules $Z_i= \ker d_i$
given by the kernel are called \emph{$i$-cycles} and the elements in
the image $B_i= \image d_{i-1}$ which is a submodule of $Z_i$ because
of \eqref{eq:differential-chain-complex} are called
\emph{$i$-boundaries}. In addition, $d$ is often called \emph{boundary
  operator}.

The resulting $\mathsf{R}$-module $H_i= H_i(C)= Z_i/B_i$ is called
\emph{$i$-th homology module} or \emph{$i$-th homology group} of the
complex $(C,d)$. It follows that a map $\alpha_i: C_i\longrightarrow
C'_i$ of a chain homomorphism $\alpha$ satisfies
$\alpha_i(Z_i)\subseteq Z'_i$ and $\alpha_i (B_i)\subseteq B'_i$ so
there is an induced map $H_i(\alpha): H_i(C)\longrightarrow
H_i(C')$. Altogether one gets additive functors, the so-called
\emph{$i$-th homology functors} $H_i$ from the category
$\mathsf{R}$-$\mathbf{comp}$ of $\mathsf{R}$-chain complexes to the
category $\mathsf{R}$-$\mathbf{mod}$ of
$\mathsf{R}$-modules.

A complex with $C_i=0$ for $i<0$ is called \emph{positive} or
\emph{chain complex} and the elements in a non-vanishing $C_i$ are
referred to as \emph{$i$-chains}. If $C_i=0$ for $i>0$ the complex is
called \emph{negative} or \emph{cochain complex}. In this case one
denotes $C_{-i}$ as $C^i$ whose elements are \emph{$i$-cochains} and
$d_{-i}$ as $d^i$, the \emph{coboundary operators} or
\emph{differentials}. Analogously, one defines the \emph{$i$-cocycles}
and \emph{$i$-coboundaries} as the elements of $Z^i=\ker d^i$ and
$B^i= \image d^{i-1}$, respectively, and finally the \emph{$i$-th
  cohomology} group $H^i=Z^i/B^i$.

\begin{definition}[Homotopy]
    \label{definition:homotopy}
    Let $\alpha$ and $\beta$ be chain homomorphisms of a complex
    $(C,d)$ into a complex $(C',d')$. Then $\alpha$ is said to be
    \emph{homotopic} to $\beta$ if there exists an indexed set
    $s=\{s_i\}_{i\in \mathbb{Z}}$ of module morphisms $s_i:
    C_i\longrightarrow C'_{i+1}$ such that
    \begin{equation}
        \label{eq:homotopy-formula}
        \alpha_i-\beta_i=d'_{i+1}\circ s_i + s_{i-1} \circ d_i.
    \end{equation}
\end{definition}

With some abuse of notation the set of maps $s$ is often referred to
as \emph{(chain) homotopy map}. It is easy to verify that homotopy,
indicated by $\alpha \sim \beta$ is an equivalence relation and that
there even exists a multiplication. If $\alpha,\beta:
(C,d)\longrightarrow (C',d')$ are homotopic with homotopy map $s$ and
$\gamma,\delta: (C',d')\longrightarrow (C'',d'')$ are homotopic via
$t$, then one easily proves the homotopy $\gamma\circ \alpha\sim
\delta\circ \beta$ with homotopy maps $u_i= \gamma_{i+1}\circ s_i +
t_i \circ \beta_i$.  The relevance of homotopy lies in the fact that
homotopic chain morphisms $\alpha\sim \beta$ yield the same maps on
the homology groups, this means $H_i(\alpha)= H_i(\beta)$.

\section{Projective resolutions and derived functors}
\label{sec:resolutions-derived-functors}

Projective resolutions are a crucial structure for the investigation
of (co)homologies in the purely algebraic framework. Before we come to
the definition it is necessary to recall the notion of a projective
module, confer \cite[Chap.~1, $\mathsection$2A]{lam:1999a}.

\begin{definition}[Projective module] 
    An $\mathsf{R}$-module $P$ is called \emph{projective} if for
    every homomorphism $f:P\longrightarrow N$ and every epimorphism
    $p:M\longrightarrow N$ there exists a homomorphism
    $h:P\longrightarrow N$ such that $f= p\circ h$.
\end{definition}

The defining property can be visualized in the diagram
\begin{equation}
    \label{eq:projective-module}
    \xymatrix{& P \ar@{-->}[dl]_h \ar[d]^f & \\
      M \ar[r]_p & N \ar[r] & 0.
    }
\end{equation}

\begin{remark}[Equivalent definitions]
    \label{remark:projective-module}
    There exist the following further equivalent characterizations of
    a projective (right) $\mathsf{R}$-module $P$, see \cite[Chap.~1,
    $\mathsection$2A]{lam:1999a}.
    \begin{enumerate}
    \item There exists another $\mathsf{R}$-module $Q$ such that the
        direct sum $P\oplus Q$ is a free $\mathsf{R}$-module, this
        means $P\oplus Q \cong \mathsf{R}^{(I)}=\bigoplus_{i\in I}
        \mathsf{R}_i$ with $\mathsf{R}_i=\mathsf{R}$ for each $i$ in
        some index set $I$.
    \item There exists a \emph{dual basis}, this means there exist
        families of elements $\{e_i\}_{i\in I}\subseteq P$ and of
        $\mathsf{R}$-linear functionals $\{\epsilon^i\}_{i \in
          I}\subseteq \Hom_{\mathsf{R}}(P,\mathsf{R})$ for some index
        set $I$ such that for any $p\in P$ one has $\epsilon^i(p)=0$
        for all except finitely many $i\in I$ and
        \begin{equation}
            \label{eq:dual-basis-equation}
            p=\sum_{i \in I} e_i \epsilon^i(p).
        \end{equation}
    \end{enumerate}
\end{remark}

Now it is possible to define the important notion of a projective
resolution.

\begin{definition}[Complex over and resolution of a module]
    \label{definition:resolution}
    Let $M$ be an $\mathsf{R}$-module. A \emph{complex} over $M$ is a
    positive complex $C=(C,d)$ together with a homomorphism $\epsilon:
    C_0\longrightarrow M$, called an \emph{augmentation}, such that
    $\epsilon \circ d_1=0$. The complex $(C,\epsilon)$ over $M$ is
    called a \emph{resolution} of $M$ if the sequence
    \begin{equation}
        \label{eq:sequence-resolution}
        \xymatrix{ 0 & M \ar[l] & C_0 \ar[l]_{\epsilon} & C_1
          \ar[l]_{d_1} & \dots \ar[l] & C_{n-1} \ar[l]  &
          C_n \ar[l]_{d_n} & \dots \ar[l]
          }
    \end{equation}
    is exact, this means if $\image d_n= \ker d_{n-1}$ for all $n\in
    \mathbb{N}$, $\image d_1= \ker \epsilon$ and $\epsilon$ is
    surjective. Further, the complex $(C,\epsilon)$ over $M$ is called
    \emph{projective}, if every $C_i$ is projective.
\end{definition}

Projective resolutions and even free resolutions with free $C_i$
always exist, confer \cite[Sect.~6.5]{jacobson:1989a}. For a resolution
the exactness of \eqref{eq:sequence-resolution} implies that the
homologies $H_i(C)=0$ are trivial for $i>0$ and that $H_0(C)= C_0/
\ker \epsilon \cong M$. Concerning the above structures one now has
the following important result.

\begin{theorem}
    \label{theorem:projective-complex-resolution}
    Let $(C,\epsilon)$ be a projective complex over the module $M$ and
    let $(C', \epsilon')$ be a resolution of the module $M'$. Further
    let $\mu: M\longrightarrow M'$ be a homomorphism. Then there
    exists a chain homomorphism $\alpha$ of the complex $C$ into $C'$
    such that $\mu \circ \epsilon = \epsilon'\circ
    \alpha_0$. Moreover, any two such homomorphisms $\alpha$ and
    $\beta$ are homotopic.
\end{theorem}

The proof consists of a multiple application of the above stated
property of projective modules. The assertion of
Theorem~\ref{theorem:projective-complex-resolution} belongs to the
diagram

\begin{equation}
    \label{eq:chain-morphism-resolution}
    \xymatrix{
      0& M \ar[l] \ar[d]_\mu &C_0 \ar[l]_\epsilon
      \ar[d]_{\alpha_0}&\dots \ar[l] &C_{n-1} \ar[l] \ar[d]_{\alpha_{n-1}}&C_n
      \ar[l]_{d_n} \ar[d]_{\alpha_n}&\dots \ar[l] \\
      0& M' \ar[l]&C'_0 \ar[l]_{\epsilon'}&\dots \ar[l]&C'_{n-1} \ar[l]&C'_n
      \ar[l]_{d'_n}& \dots. \ar[l]
      }
%     \xymatrix{ \dots \ar[r]& C_n \ar[r]^{d_n} \ar[d]_{\alpha_n} &
%       C_{n-1} \ar[r] \ar[d]_{\alpha_{n-1}} & \dots \ar[r]& C_0
%       \ar[r]^\epsilon \ar[d]_{\alpha_0} & M \ar[r] \ar[d]_\mu &0\\
%       \dots \ar[r]& C'_n \ar[r]^{d'_n} & C'_{n-1} \ar[r] &
%       \dots \ar[r]& C'_0 \ar[r]^{\epsilon'} & M' \ar[r] &0.
%     }
\end{equation}

An important application of the above structures is the definition of
derived functors of an additive functor $F$ from the category
$\mathsf{R}$-$\mathbf{mod}$ to the category $\mathbf{Ab}$ of abelian
groups. Note that a functor between categories where the sets of
morphisms all are abelian groups, is called \emph{additive} if the
corresponding map between the morphisms is additive, confer
\cite[Sect.~3.1]{jacobson:1989a}.

For an $\mathsf{R}$-module $M$ let
\begin{equation}
    \label{eq:projective-resolution}
    \xymatrix{0 & M \ar[l] & C_0 \ar[l]_{\epsilon} & C_1 \ar[l]_{d_1}
      & \dots \ar[l]_{d_2}}     
\end{equation}
be a projective resolution of $M$. Application of a covariant additive
functor $F$ then yields a sequence of homomorphisms of abelian groups
\begin{equation}
    \label{eq:functor-sequence}
    \xymatrix{0 & F(M) \ar[l] & F(C_0) \ar[l]_{F(\epsilon)} & F(C_1)
      \ar[l]_{F(d_1)} & \dots \ar[l]_{F(d_2)}}
\end{equation}
since $F(d_i)\circ F(d_{i-1})= F(d_i \circ d_{i-1})=0$. This sequence
is not necessarily exact, so there can occur nontrivial homology
groups $H_n(F(C))$. Now consider a projective resolution
$(C',\epsilon')$ of a different module $M'$ with a corresponding
homomorphism $\mu: M\longrightarrow M'$. Due to Theorem
\ref{theorem:projective-complex-resolution} there exists a chain
homomorphism $\alpha: C\longrightarrow C'$ with $\mu \circ \epsilon =
\epsilon' \circ \alpha_0$ and two such homomorphism $\alpha,\beta$ are
always homotopic. The functor $F$ then yields chain homomorphisms
$F(\alpha)$ with $F(\mu) \circ F(\epsilon)= F(\epsilon')\circ
F(\alpha_0)$ and since it is additive the homotopy $\alpha\sim \beta$
translates to $F(\alpha)\sim F(\beta)$.  In particular, if $M=M'$ and
$\mu=\id$ it follows that $H_n(F(\alpha)): H_n(F(C))\longrightarrow
H_n(F(C'))$ is an isomorphism.

Altogether it follows that for the given covariant functor $F$ and
every $n\in \mathbb{N}_0$ one has found another additive functor
$L_nF$, the so called \emph{$n$-th left derived functor} from
$\mathsf{R}$-$\mathbf{mod}$ to $\mathbf{Ab}$, given by
\begin{eqnarray}
    \label{eq:left-derived-functor-objects}
    L_nF(M)&=& H_n(F(C)),\\
    \label{eq:left-derived-functor-morphisms}
    L_nF(\mu)&=& H_n(F(\alpha))
\end{eqnarray}
for every $\mathsf{R}$-module $M$ and every morphism $\mu:
M\longrightarrow M'$. In particular, one has $L_0F(M)= F(C_0)/\image
F(d_1)$. Due to the above results the right hand sides of
\eqref{eq:left-derived-functor-objects} and
\eqref{eq:left-derived-functor-morphisms} do not depend on the
resolution of $M$ and the chain homomorphism $\alpha$ \emph{over}
$\mu$. If the functor $F$ is contravariant its application to
\eqref{eq:projective-resolution} yields a cochain complex
\begin{equation}
    \label{eq:functor-sequence-contravariant}
    \xymatrix{0 \ar[r] & F(M) \ar[r]^{F(\epsilon)} & F(C_0)
      \ar[r]^{F(d_1)} & F(C_1) 
      \ar[r]^{F(d_2)} & \dots . }
\end{equation}
With the same arguments as above and using the cohomology functors
$H^n$ one gets the \emph{$n$-th right derived functor} $R^nF$. There
one has in particular $R^0F(M)= \ker F(d_1)$.

The following example of derived functors is essential for the
discussion of Hochschild cohomologies as explained in
Section~\ref{sec:hochschild-Ext}. Let $\mathsf{R}$ be a ring and $N$
be an $\mathsf{R}$-module. Then consider the contravariant additive
functor
\begin{equation}
    \label{eq:hom-functor}
    \mathrm{hom}(\cdot, N): \textrm{$\mathsf{R}$-$\mathbf{mod}$}
    \longrightarrow \mathbf{Ab},
\end{equation}
given by
\begin{equation}
    \label{eq:hom-functor-objects}
    \mathrm{hom}(\cdot, N)(M) = \Hom_{\mathsf{R}}(M,N)
\end{equation}
for all $\mathsf{R}$-modules $M$ and
\begin{equation}
    \label{eq:hom-functor-morphisms}
    \mathrm{hom}(\cdot, N)(\mu: M\longrightarrow M') = \left(\mu^*:
    \begin{array}{rcl}
        \Hom_{\mathsf{R}}(M',N) &\longrightarrow& \Hom_{\mathsf{R}}(M,N) \\
        \nu &\longmapsto& \mu^*\nu=\nu \circ \mu
    \end{array}
    \right)
\end{equation}
for all homomorphisms $\mu:M\longrightarrow M'$.
The $n$-th right derived functor of $\mathrm{hom}(\cdot, N)$ is
denoted by
\begin{equation}
    \label{eq:Ext}
    \Ext^n_{\mathsf{R}}(\cdot,N)= R^n \mathrm{hom}(\cdot, N),
\end{equation}
and its value for a module $M$ with an arbitrary projective resolution
$(C,\epsilon)$ is thus given by
\begin{equation}
    \label{eq:Ext-on-module}
    \Ext^n_{\mathsf{R}}(M,N)= H^n(\Hom_{\mathsf{R}}(C,N)).
\end{equation}
If the context is clear, this means for fixed modules $M$ and $N$, we
refer to the latter as the \emph{$n$-th $\Ext$ group}. Since
$\mathrm{hom}(\cdot, N)$ is \emph{left exact} the exactness of
$C_1\rightarrow C_0 \xrightarrow{\epsilon} M \rightarrow 0$ implies
that of $0 \rightarrow \Hom_{\mathsf{R}}(M,N) \xrightarrow{\epsilon^*}
\Hom_{\mathsf{R}}(C_0,N) \rightarrow \Hom_{\mathsf{R}}(C_1,N)$ and it
follows that
\begin{equation}
    \label{eq:Ext-zero}
    \Ext^0(M,N)\cong \Hom_{\mathsf{R}}(M,N).
\end{equation}

% Concluding this short overview we quote the following theorem which
% provides further equivalent characterizations of projective modules.

% \begin{theorem}[Projective modules and $\Ext$ groups] 
%     The following conditions on a module $M$ are equivalent.
%     \begin{enumerate}
%     \item $M$ is projective.
%     \item $\Ext^n(M,N)=0$ for all $n\ge 1$ and all modules $N$.
%     \item $\Ext^1(M,N)=0$ for all $N$.
%     \end{enumerate}
% \end{theorem}

\chapter{Explicit chain maps of the bar and Koszul complex}
\label{cha:technical-proofs}

This appendix provides the missing computations concerning the maps
$F_k: K_k\longrightarrow X_k$ and $G_k: X_k \longrightarrow K_k$
defined by the Equations~\eqref{eq:chain-maps-Koszul-bar} and
\eqref{eq:chain-maps-bar-Koszul} in Section~\ref{sec:Koszul}. First,
it will be shown that they are indeed chain maps between the bar and
Koszul complex. With the definitions
\eqref{eq:boundary-operator-bar-complex} and
\eqref{eq:boundary-operator-Koszul-complex} for the differentials
$d^X_k$ and $d^K_k$ and the notation as in Section~\ref{sec:Koszul}
one computes for $k\ge 1$
\begin{eqnarray*}
    \lefteqn{ (d^X_k F_k \omega)(v, q_1, \dots , q_{k-1}, w)}\\
    &=& 
    (F_k\omega)(v,v,q_1, \dots, q_{k-1},w) + \sum_{i=1}^{k-1}
    (-1)^i(F_k\omega)(v, q_1,  \dots, q_i, q_i, \dots, q_{k-1},w)\\
    &&+(-1)^k(F_k\omega)(v, q_1, \dots, q_{k-1}, w,w) \\
    &=& \omega(v,w)(0, q_1-v, \dots, q_{k-1}-v)+ \sum_{i=1}^{k-1} (-1)^i
    \omega(v,w)(q_1-v, \dots, q_i-v, q_i-v, \dots, q_{k-1}-v)\\
    &&+ (-1)^k \omega(v,w)(q_1-v, \dots, q_{k-1}-v, w-v)\\
    &=& (-1)^k(-1)(-1)^{k-1} \omega(v,w)(v-w, q_1-v, \dots, q_{k-1}-v)\\
    &=& (d^K_k \omega)(v,w)(q_1-v,\dots, q_{k-1}-v) \\
    &=& (F_{k-1} d^K_k \omega)(v,w) (v,q_1, \dots, q_{k-1},w).
\end{eqnarray*}
This shows the desired equation $d^X_k \circ F_k=F_{k-1} \circ d^K_k$.
For $G_k$ the computation is longer. In a first step one has
\begin{align*}
    &(d^K_k G_k \chi) (v,w)\\
    = &\insa(v-w) (e^{i_1}\wedge \dots \wedge e^{i_k}) \int\limits_0^1
    \D t_1 \dots \int\limits_0^{t_{k-1}} \D t_k \frac{\partial^k
      \chi}{\partial q_1^{i_1} \dots \partial q_k^{i_k}}
    (v, t_1 v+ (1-t_1)w, \dots,t_k v +(1-t_k)w ,w)\\
    = &\sum_{j=1}^k(-1)^{j+1} e^{i_1}\wedge \dots
    \stackrel{i_j}{\wedge}\dots \wedge e^{i_k} (v^{i_j}-w^{i_j})
    \int\limits_0^1 \D t_1 \dots \int\limits_0^{t_{k-1}} \D t_k
    \frac{\partial^k \chi}{\partial q_1^{i_1} \dots \partial
      q_k^{i_k}}
    (v, q_1(t_1), \dots , q_k(t_k),w)\\
    = &\sum_{j=1}^k(-1)^{j+1} e^{i_1}\wedge \dots
    \stackrel{i_j}{\wedge}\dots \wedge e^{i_k} \int\limits_0^1 \D t_1
    \dots \int\limits_0^{t_{k-1}} \D t_k \frac{\partial}{\partial
      t_j}\left( \frac{\partial^{k-1} \chi} {\partial q_1^{i_1} \dots
          \stackrel{j}{\wedge}\dots \partial q_k^{i_k}} (v, q_1(t_1),
        \dots , q_k(t_k),w) \right)
\end{align*}
where $\stackrel{j}{\wedge}$ means to leave out the $j$-th term of an
expression. Here and in the following we use the abbreviation
$q_r(t_i)=t_iv+(1-t_i)w$ for $r,i=1,\dots, k$ where the first index
$r$ denotes the position of the argument. For $j=2,\dots ,k-1$ we use
the fact that for $f\in C^\infty(\mathbb{R}^2)$ one has
$\int\limits_0^{t_j} \D t_{j+1} \frac{\partial f} {\partial t_j} (t_j,
t_{j+1})= \frac{\D} {\D t_j}\left( \int\limits_0^{t_j}\D t_{j+1}
    f(t_j, t_{j+1})\right)- f(t_j, t_j)$ and that subsequent
integration yields
\begin{equation}
    \label{eq:integration-general}
    \int\limits_0^{t_{j-1}} \D t_j\int\limits_0^{t_j} \D
    t_{j+1} \frac{\partial f} {\partial t_j} (t_j, t_{j+1})=
    \int\limits_0^{t_{j-1}}\D t_{j+1} \: f(t_{j-1}, t_{j+1})- 
    \int\limits_0^{t_{j-1}} \D t_j \: f(t_j, t_j).
\end{equation}
For $j=1$ the same argument is applied in the form
\begin{equation}
    \label{eq:integration-one}
    \int\limits_0^1 \D t_1\int\limits_0^{t_1} \D
    t_2 \frac{\partial f} {\partial t_1} (t_1, t_2)=
    \int\limits_0^1 \D t_2 \: f(1, t_2)- 
    \int\limits_0^1 \D t_1 \: f(t_1, t_1)
\end{equation}
noting that $q_1(1)=v$. In the case $j=k$ the integration can be
carried out directly. After setting $t_0=1$ the second summands in the
above expressions \eqref{eq:integration-general} and
\eqref{eq:integration-one} have the same structure whereby we can
combine them into one sum.
\begin{align*}
    &(d_k G_k \chi) (v,w)\\
    = & \phantom{+} e^{i_2}\wedge \dots \wedge e^{i_k} 
    \int\limits_0^1\D t_2 \dots \int\limits_0^{t_{k-1}} \D t_k 
    \frac{\partial^{k-1} \chi}{\partial q_2^{i_2}\dots \partial
      q_k^{i_k}}
    (v, v, q_2(t_2), \dots , q_k(t_k), w)\\
    &+ \sum_{j=2}^{k-1} (-1)^{j+1} 
    e^{i_1}\wedge \ldots \stackrel{i_j}{\wedge}\ldots \wedge e^{i_k} 
    \int\limits_0^1\D t_1 \dots \int\limits_0^{t_{j-2}} \D t_{j-1}
    \int\limits_0^{t_{j-1}} \D t_{j+1} \int\limits_0^{t_{j+1}} \D
    t_{j+2} \dots \int\limits_0^{t_{k-1}} \D t_k  \\
    & \quad \quad \frac{\partial^{k-1} \chi}{\partial q_1^{i_1}\dots
      \stackrel{j}{\wedge}\dots \partial
      q_k^{i_k}}
    (v,q_1(t_1), \dots, q_{j-1}(t_{j-1}), q_j(t_{j-1}), q_{j+1}(t_{j+1})
    \dots , q_k(t_k), w)\\ 
    & - \sum_{j=1}^{k-1} (-1)^{j+1} 
    e^{i_1}\wedge \ldots \stackrel{i_j}{\wedge}\ldots \wedge e^{i_k} 
    \int\limits_0^1\D t_1 \dots \int\limits_0^{t_{j-1}} \D t_j
    \int\limits_0^{t_j} \D t_{j+2} \int\limits_0^{t_{j+2}} \D
    t_{j+3} \dots \int\limits_0^{t_{k-1}} \D t_k  \\
    & \quad \quad \frac{\partial^{k-1} \chi}{\partial q_1^{i_1}\dots
      \stackrel{j}{\wedge}\dots \partial
      q_k^{i_k}}
    (v,q_1(t_1), \dots, q_j(t_j), q_{j+1}(t_j), q_{j+2}(t_{j+2}) \dots
    , q_k(t_k), w)\\ 
    & + (-1)^{k+1} 
    e^{i_1}\wedge \ldots \wedge e^{i_{k-1}} 
    \int\limits_0^1\D t_1 \dots \int\limits_0^{t_{k-2}} \D t_{k-1}\\
    & \quad \left( \frac{\partial^{k-1} \chi}{\partial q_1^{i_1}\dots
          \partial
          q_{k-1}^{i_{k-1}}}
        (v,q_1(t_1), \dots, q_{k-1}(t_{k-1}), q_k(t_{k-1}), w)\right.\\
    & \quad \left.- \frac{\partial^{k-1} \chi}{\partial q_1^{i_1}\dots
          \partial
          q_{k-1}^{i_{k-1}}}
        (v,q_1(t_1), \dots,q_{k-1}(t_{k-1}), w, w)\right).
\end{align*}
Now we rename the indices $(i_1,\dots, i_k)$ and the variables $(t_1,
\dots t_k)$ in the first three terms such that the expressions
$e^{i_1}\wedge \ldots \wedge e^{i_{k-1}}$ and the integrals
$\int\limits_0^1\D t_1 \dots \int\limits_0^{t_{k-2}} \D t_{k-1}$ are
the same for all summands and do not depend on the summation
parameter $j$. As a consequence, the derivatives as well as the arguments
of $\chi$ are changed. After a change $j'=j-1$ of the summation index
in the second term the penultimate expression can be seen as a summand
and raises the upper bound of the resulting sum. This gives
\begin{align*}
    &(d^K_k G_k \chi) (v,w)\\
    = &e^{i_1}\wedge \ldots \wedge e^{i_{k-1}} 
    \int\limits_0^1\D t_1 \dots \int\limits_0^{t_{k-2}} \D t_{k-1}\\
    &\left( 
        \frac{\partial^{k-1} \chi} {\partial q_2^{i_1}\dots \partial
          q_k^{i_{k-1}}}
        (v, v, q_2(t_1), \dots , q_k(t_{k-1}), w)\right.\\
        &+ \sum_{j=1}^{k-1} (-1)^j 
        \frac{\partial^{k-1} \chi} {\partial q_1^{i_1}\dots
          \partial q_j^{i_j} \partial q_{j+2}^{i_{j+1}} \dots \partial
          q_k^{i_{k-1}}} 
        (v,q_1(t_1), \dots, q_j(t_j), q_{j+1}(t_j), q_{j+2}(t_{j+1}),
        \dots , q_k(t_{k-1}), w)\\ 
        &+ \sum_{j=1}^{k-1} (-1)^j 
        \frac{\partial^{k-1} \chi} {\partial q_1^{i_1}\dots
          \partial q_{j-1}^{i_{j-1}} \partial q_{j+1}^{i_j}
          \dots \partial  q_k^{i_{k-1}}}
        (v,q_1(t_1), \dots, q_j(t_j), q_{j+1}(t_j), q_{j+2}(t_{j+1}),
        \dots , q_k(t_{k-1}), w)\\ 
        &\left. +(-1)^k 
        \frac{\partial^{k-1} \chi} {\partial q_1^{i_1}\dots
          \partial q_{k-1}^{i_{k-1}}} 
        (v,q_1(t_1), \dots, q_{k-1}(t_{k-1}),w, w)
    \right).\\
    &= e^{i_1}\wedge \ldots \wedge e^{i_{k-1}} 
    \int\limits_0^1\D t_1 \dots \int\limits_0^{t_{k-2}} \D t_{k-1}
    \frac{\partial^{k-1} d^X_k\chi} {\partial q_1^{i_1}\dots
      \partial q_{k-1}^{i_{k-1}}} 
    (v,q_1(t_1), \dots, q_{k-1}(t_{k-1}), w)\\
    &= (G_{k-1} d^X_k \chi)(v,w).
\end{align*}
The penultimate equation becomes clear with the definition
\eqref{eq:boundary-operator-bar-complex} of the boundary operator
$d^X_k$ of the bar complex and so we have shown $d^K_k \circ
G_k=G_{k-1} \circ d^X_k$.

The verification of $G_k \circ F_k= \id_{K_k}$ for $k\ge 0$ is a
simple computation. With \eqref{eq:chain-maps-Koszul-bar} and
\eqref{eq:element-Koszul-complex} one has
\begin{eqnarray*}
    (F_k\omega)(v,q_1,\dots, q_k,w)&=& \frac{1}{k!} \omega_{j_1\dots
      j_k}(v,w) e^{j_1}\wedge \dots \wedge e^{j_k}(q_1-v, \dots,
    q_k-v)\\
    &=& \omega_{j_1\dots
      j_k}(v,w) (q_1^{j_1}-v^{j_1}) \dots (q_k^{j_k}-v^{j_k}).
\end{eqnarray*}
With the same notation as above, this and $\int\limits_0^1 \D t_1
\dots \int\limits_0^{t_{k-1}} \D t_k=\frac{1}{k!}$ gives
\begin{eqnarray*}
    (G_k F_k \omega)(v,w) &=& e^{i_1}\wedge\dots \wedge e^{i_k}
    \int\limits_0^1 \D t_1  \dots \int\limits_0^{t_{k-1}} \D t_k
    \frac{\partial^k (F_k \omega)} {\partial q_1^{i_1} \dots \partial
      q_k^{i_k}} (v, q_1(t_1), \dots, q_k(t_k),w)\\
    &=& e^{i_1}\wedge\dots \wedge e^{i_k} \int\limits_0^1\D t_1 \dots
    \int\limits_0^{t_{k-1}} \D t_k \: \omega_{i_1\dots i_k} (v,w)\\
    &=& \frac{1}{k!} \omega_{i_1\dots i_k}(v,w)
    e^{i_1}\wedge\dots\wedge e^{i_k}= \omega(v,w)
\end{eqnarray*}
which proves the statement.

% Local Variables: 
% TeX-master: "Dissertation"
% End: 

%% file: paper.tex
\chapter*{Publication}
\label{cha:publications}

\addcontentsline{toc}{chapter}{Publication}

\begin{itemize}
\item {\sc Bordemann, M., Neumaier, N., Waldmann, S., Weiss, S.:}
    {\it Deformation quantization of surjective submersions and
      principal fibre bundles.} Preprint {\bf arXiv: 0711.2965} (2007),
    confer \cite{bordemann.neumaier.waldmann.weiss:2007a:pre}. To
    appear in: {\it Journal f{\"{u}}r die reine und angewandte
      Mathematik (Crelles Journal)}
\end{itemize}

% Local Variables: 
% TeX-master: "Dissertation"
% End: 

%% file: danke.tex
\chapter*{Acknowledgement}
\label{cha:acknowledgement}

\addcontentsline{toc}{chapter}{Acknowledgement}
\thispagestyle{empty}

It is a pleasure for me to thank all those people and institutions
that supported me in various ways during the last three years. Many
thanks go to

\begin{itemize}
\item my supervisor Stefan Waldmann for the excellent supervision of
    this PhD thesis and the very good working atmosphere during all
    the years I know him.
\item Nikolai Neumaier for the many helpful discussions and for
    reading the proof.
\item Martin Bordemann for the valuable initial impulses which
    motivated this thesis and the joint publication together with him,
    Nikolai, and Stefan.
\item Prof.~Dr.~Hartmann R{\"{o}}mer and his successor
    Prof.~Dr.~Stefan Dittmaier for the kind acceptance in their
    department.
\item Florian Becher, Svea Beiser, and Matteo Carrera for the friendly
    working atmosphere in the office. Many thanks also go to all the
    other group members, in particular Maximilian Hanusch and Dominic
    Maier for reading the proof and Klaus Zimmermann for the IT
    support.
\item Henrique Bursztyn for giving me the opportunity to work with him
    at the \emph{Instituto Nacional de Matem{\'{a}}tica Pura e
      Aplicada} in Rio de Janeiro for half a year. In this context my
    thanks also go to Alejandro Cabrera and all the other
    mathematicians I met there for the very valuable
    discussions. Furthermore, I would like to thank all people I met
    in Brazil for their hospitality, their support in any situation,
    and for the great time I had in Brazil.
\item the \emph{Studienstiftung des deutschen Volkes} for the
    financial support.
\item Andrea Wachten for her love and understanding.
\item my parents Ursula and Josef as well as my brother Josef and his
    son Mario for their unbounded support which always encouraged me to
    continue my way.
\end{itemize}

\clearpage

% Local Variables: 
% TeX-master: "Dissertation"
% End: 